\documentclass[11pt, amssymb,amsfonts]{amsart}
\usepackage{amsrefs}
\usepackage{float}
\makeatletter
\def\l@section{\@tocline{1}{10pt plus0pt}{0pt}{}{\bfseries}}

\def\@tocline#1#2#3#4#5#6#7{\relax
    \ifnum #1>-1
  \ifnum #1>\c@tocdepth 
  \else
    \par \addpenalty\@secpenalty
    \begingroup \hyphenpenalty\@M
    \@ifempty{#4}{%
      \@tempdima\csname r@tocindent\number#1\endcsname\relax
    }{%
      \@tempdima#4\relax
    }%
    \parindent\z@ \leftskip#3\relax \advance\leftskip\@tempdima\relax
    \rightskip\@pnumwidth plus4em \parfillskip-\@pnumwidth
    #5\leavevmode\hskip-\@tempdima #6\nobreak\relax
    \hfil\hbox to\@pnumwidth{\@tocpagenum{#7}}\par
    \nobreak
    \endgroup
  \fi
\fi}

\makeatother
\usepackage{algorithm}
\usepackage{algorithmic}


\usepackage{hyperref}
\usepackage{a4wide}
\usepackage[a4paper, margin=2.3cm]{geometry}
\usepackage{amsmath,amssymb,amsthm}
\usepackage{color}
\usepackage{parskip}
\usepackage{graphicx}
\usepackage{hyperref}
\usepackage{parskip}
\usepackage{setspace}
\usepackage{amsmath}
\usepackage{mathtools}
\usepackage{amssymb}
\usepackage{cases}
\usepackage{algorithm}
\usepackage{algorithmic} 
\usepackage{todonotes}
\usepackage{stmaryrd}
\usepackage{IEEEtrantools}
\usepackage[nocompress]{cite}
\usepackage{tikz}
\usetikzlibrary{patterns,snakes}
\newcommand{\norm}[1]{ \left|  #1 \right| }
\newcommand{\Norm}[1]{ \left\|  #1 \right\| }

\def\dist{\hbox{\rm dist \,\,}}

\def\Var{\hbox{\bf Var}}

\def\P{{\hbox{\bf P}}}
\def\E{{\hbox{\bf E}}}

\font \roman = cmr10 at 10 true pt

\def\rank{\hbox{\rm rank}}

\def\be#1{ \begin{equation}\label{#1} }

\def\bas{\begin{align*}}
\def\eas{\end{align*}}
\def\bi{\begin{itemize}}
\def\ei{\end{itemize}}

\def\dist{{\hbox{\roman dist}}}

\def\emph#1{{\it #1}}
\def\textbf#1{{\bf #1}}

\def\Span{\hbox{ \rm Span} \,\,\, }

\def \small {\scriptsize}

\theoremstyle{plain}
 \theoremstyle{plain}
  \newtheorem{theorem}{Theorem}
  \numberwithin{theorem}{section}

  \numberwithin{problem}{section}

  \newtheorem{lemma}{Lemma}
  
  \numberwithin{notation}{section}
  \numberwithin{lemma}{section}
  \newtheorem{corollary}{Corollary}
   \numberwithin{corollary}{section}

\theoremstyle{remark}
  \newtheorem{remark}{Remark}
  \numberwithin{remark}{section}

\theoremstyle{definition}
  \newtheorem{definition}{Definition}
\numberwithin{definition}{section}

\theoremstyle{definition}
  
  \numberwithin{summary}{section}
  
  \newtheorem{question}{Question}
  \numberwithin{question}{section}

\DeclareUnicodeCharacter{2113}{\ensuremath{\ell}}
\DeclareUnicodeCharacter{2082}{\ensuremath{\infty}}
\DeclareUnicodeCharacter{221E}{\ensuremath{\infty}}
\begin{document}
\include{psfig}

\setcounter{tocdepth}{2}


\title[New matrix perturbation bounds with relative strength]{New matrix perturbation bounds with relative strength: Perturbation of eigenspaces }
\pagenumbering{arabic}

\author{Phuc Tran, Van Vu }
\thanks{Department of Mathematics, The University of Hong Kong (HKU) \\
  \textit{Email:} \texttt{trandangphuc234@gmail.com, vuhavan001@gmail.com}.}
\date{}

\begin{abstract} 
Matrix perturbation bounds (such as Weyl and Davis--Kahan) are used abundantly in many areas of mathematics and data science. Many bounds (such as the above two) involve the spectral norm of the noise 
matrix and are sharp in worst-case analysis.   

\vskip2mm 

In order to refine these classical bounds, we introduce a new parameter, which we refer to as the relative strength. This parameter measures the strength of the action of the noise matrix on the relevant eigenvectors of the ground matrix. It has turned out that in a number of situations, we can use the relative strength as a replacement for the spectral norm (which can be seen as the absolute strength). 
 This has led to a number of notable improvements under 
certain sets of assumptions, which are frequently met in practice. A representative example is the case when the noise matrix is random.  
\vskip2mm 

For the purpose of our study, we introduce a new method of analysis, which combines the classical contour integral argument with new (combinatorial) ideas. This method is robust and of independent interest.

\vskip2mm 
In the current paper, we focus on the perturbation of eigenspaces (Davis--Kahan type results). Perturbation bounds for eigenspaces are essential in statistics and theoretical computer science, and thus deserve a special treatment. Furthermore,  this will lay the ground for the more technical treatment of general matrix functionals, which appears in a future paper.

 \vskip5mm 


\end{abstract} 
\maketitle

{\it Acknowledgment.} We would like to thank  R. Vershynin, T. Tao, and the reviewers for their valuable comments and suggestions, which greatly helped improve the presentation and quality of the paper.

\newpage
\tableofcontents
\newpage

\section{Introduction}

\subsection{Perturbation bounds with relative strength} 

Let $A$ and $E$ be symmetric real matrices of size $n$, and  $\tilde A:= A+E$. We view $A$ as the truth (or data/ground) matrix and $E$ as noise. 
Consider the spectral decomposition 
$$A =\sum_{i=1}^n \lambda_i u_i u_i^T , $$ where $\lambda_i$ are the eigenvalues and $u_i$ are the corresponding eigenvectors. We order the eigenvalues decreasingly,  $\lambda_1 \ge \lambda_2 \dots \geq \lambda_n$. We also order the singular values of $A$ in a similar manner 
$\sigma_1 \ge \sigma_2 \ge \dots \ge \sigma_n$. It is well known that 
$\{\sigma_1, \dots, \sigma_n \} = \{| \lambda_1| , \dots,|  \lambda_ n| \} $. Therefore, $\sigma_i = | \lambda_{\pi (i) } | $ for a permutation $\pi \in S_n$. We use the notation $\tilde \lambda_i, \tilde u_i,$ etc., for $\tilde A$, with the same ordering. 

\begin{remark} To simplify the presentation, we temporarily assume that the eigenvalues (singular values) are different, so the eigenvectors (singular vectors) are well-defined (up to signs).  Our argument works for matrices with multiple eigenvalues,  and the results
will be stated in this generality. Our main results extend to rectangular matrices; see Section \ref{section:recMa}. \end{remark} 

A perturbation bound measures the difference $g(\tilde A) - g(A)$ for a 
matrix functional $g$. Perturbation theory is a well-established area of mathematics \cite{Book1, KVbook, Kato1, SS1, WW1, Higham2008}, with a wide range of applications across all branches of computer science, data science, and statistics.

{\it \noindent Examples.} Important examples of matrix functionals \(g(A)\) include the leading eigenvalues and singular values, the leading eigenvectors, the eigenspace spanned by the leading \(p\) eigenvectors, and the best rank-\(p\) approximation.

Many classical perturbation bounds, such as Weyl's inequality for eigenvalues and the Davis--Kahan theorem for eigenspaces, are expressed in terms of the spectral (operator) norm of the perturbation matrix \(E\). While these bounds are sharp, they do not take into account the interaction between \(E\) and \(A\). In this sense, they are sharp from a worst-case perspective, where \(E\) may act adversarially on the (relevant) eigenvectors of \(A\); see the discussion in the next subsection. The spectral norm can be viewed as measuring the absolute strength of the perturbation. Indeed,
\[
\|E\|
=
\max_{\|u\|=\|v\|=1}|u^T Ev|.
\]

To refine these classical perturbation bounds, we introduce the notion of
{\it relative strength}, a parameter that measures the strength of the action of the
noise matrix \(E\) on the relevant eigenvectors of the ground matrix \(A\).
We show that, in a number of perturbation bounds, the spectral norm can be replaced by this relative strength, leading to significant improvements under certain sets of assumptions that frequently arise in applications. An important example is when \(E\) is random, which we will discuss in detail.

In order to achieve these new bounds, we introduce a novel method, which we call {\it the method of combinatorial contour expansion}, to bound $  g(\tilde A)- g(A),$ for a wide range of functions $g$.  This method is robust and of independent interest. In many cases, it has led us to the true nature of the perturbation in question, as the results will be essentially sharp.

  In this paper, we focus on the perturbation of eigenspaces, which is important in a number of applications and provides a simple setting to demonstrate our method; see Sections~\ref{section:applications} and \ref{section:method}, respectively. We will treat other matrix functionals in a subsequent paper \cite{TranVu2}.

 \noindent\textbf{Notation.} We use standard asymptotic notation such as $O, o, \Theta, \Omega, \tilde{O},$ etc., under the assumption that $n \rightarrow  \infty$.  For a matrix $M$, $\|M \|$ and $\|M\|_F$ denote its operator (spectral) norm and Frobenius norm, respectively. Logarithms have a natural base. 

\subsection{Perturbation of eigenspaces: Davis--Kahan bound} \label{section:classical} 
The standard perturbation bound for eigenspaces is the Davis--Kahan theorem \cite{DKoriginal, Book1}. Let
\(S\subset\{1,\ldots,n\}\), and let \(\Pi_S\) and \(\tilde{\Pi}_S\) denote the orthogonal projections onto the subspaces spanned by
\(\{u_i:i\in S\}\) and \(\{\tilde u_i:i\in S\}\), respectively.


\begin{theorem}[Davis--Kahan] \label{DKgeneral}  Let  $\delta_S :=\min \{ | \lambda_i - \lambda_j | , i \in S, j \notin S \}$ be the spectral separation. Then 
\begin{equation}    \|\tilde{\Pi}_S - \Pi_S \| \leq  \frac{\pi}{2} \frac{\Norm{E}}{\delta_{S}}.
\end{equation}  
\end{theorem} 
The original formulation uses the perturbed spectral separation \(\min_{i\in S,\,j\notin S}|\lambda_i-\widetilde{\lambda}_j|\); see Theorem~VII.3.2 in \cite{Book1}.  This version is easier to use in applications.

In practice, the most important special case is when $S= \{1, \dots, p \}$.  In this case, we obtain the following bound for the perturbation of $\Pi_p$ and $\tilde{\Pi}_p$, the orthogonal projections onto the $p$-leading eigenspaces of $A$ and $\tilde{A}$ respectively.

\begin{corollary}  \label{DKp} Let $\delta_p := \lambda_p -\lambda_{p+1}$. Then 
\begin{equation}
\Norm{\tilde{\Pi}_p - \Pi_{p}} \leq \frac{\pi}{2}  \frac{\|E\|}{\delta_p}.
\end{equation}

\end{corollary}

\begin{remark}  Since the trivial bound is 1, Corollary \ref{DKp} is only nontrivial under the gap condition 
\begin{equation} \label{gaptonoise0} \delta_S \ge \frac{\pi}{2 }   \| E \| . \end{equation} 
In fact, if one wants the bound in Theorem \ref{DKgeneral} to be  of order  $O(\epsilon)$, then one needs a strong gap assumption
\begin{equation} \label{gaptonoise1} \delta_S =\Omega \left(\frac{1}{\epsilon}   \| E \| \right) . \end{equation}

\end{remark}

To discuss the sharpness of the Davis--Kahan bound, let us consider the following example:
 $$A= \tiny{\begin{pmatrix} 
     1 & 0 & 0 & 0 &  0 & \cdots & 0\\
     0 & 1- \epsilon & 0 & 0 &  0 & \cdots & 0  \\
     0 & 0 & 1 - 2\epsilon  & 0  & 0 & \cdots & 0  \\
     0 & 0 & 0 & 1 -3 \epsilon  &0 & \cdots & 0  \\
         0 & 0 & 0 & 0 &  1/2 & \cdots & 0 \\
    \cdots \\
     0 & 0 & 0 & 0 &  0 & \cdots & 1/2
 \end{pmatrix}},  $$
 $$ E:= \tiny{
 \begin{pmatrix}
     0 & 0 & 0 & 0 & 0 & \cdots & 0 \\
     0 & 5 \epsilon & 0 & 0 & 0 & \cdots & 0  \\
     0 & 0 & 5\epsilon & 0  & 0 &  \cdots & 0  \\
     0 & 0 & 0 & 0 & 0 & \cdots & 0  \\
     0 & 0 & 0 & 0 & 0 &  \cdots & 0  \\
       \cdots \\
    0 & 0 & 0 & 0 & 0 &  \cdots & 0
 \end{pmatrix}}, \, \tilde{A} = A+E= \tiny{ \begin{pmatrix}
     1 & 0 & 0 & 0 & 0  & \cdots & 0\\
     0 & 1+ 4 \epsilon & 0 & 0 &  0 & \cdots & 0 \\
     0 & 0 & 1 + 3\epsilon & 0 & 0 & \cdots & 0 \\
     0 & 0 & 0 & 1 -3 \epsilon  &0 & \cdots & 0 \\
     0 & 0 & 0 & 0 &  1/2& \cdots & 0 \\
    \cdots \\
    0 & 0 & 0 & 0 &  0 & \cdots & 1/2
 \end{pmatrix}}, $$
for some $0<\epsilon<1/3$. Here and throughout, $e_1,\ldots,e_n$ denote the standard basis vectors; for instance, $e_1=(1,0,\ldots,0)$.

Regardless of how small $\epsilon$ is, the two leading eigenvectors change from $e_1,e_2$ (in that order) to $e_2,e_3$.  Thus, both $\Pi_1$ and $\Pi_2$ change radically. This shows that the Davis--Kahan bound is sharp up to a constant factor.

 However, we also notice that this scenario is very rare: among all $n$ dimensions, the noise matrix $E$ acts
{\it adversarially} exactly on the two relevant eigenvectors of $A$. 
Such an alignment is unlikely for a generic perturbation. For instance, if the nonzero entries of $E$ were placed at random, its effect on the leading eigenspaces would typically be much smaller.

\section{New perturbation bounds for eigenspaces } \label{section: main}

To ease the presentation, we first focus on the special case when $S=\{1, \dots, p\} $. In this case, $\Pi_S := \sum_{i \in S} u_i u_i^T = \sum_{i=1}^p u_i u_i^T$ is the orthogonal projection onto the eigenspace spanned by the first $p$ eigenvectors. The treatment of general $S$ appears in Subsection \ref{subsection:generalS}.

\subsection{A special case: $A$ has low rank $r$.} 
For the sake of presentation, we first consider the case where \(A\) has low rank \(r\ll n\). In practice, this setting is of considerable interest, as datasets arising in many applications are low-rank or approximately low-rank; see Section~\ref{section:applications}. 
To ensure that \(\Pi_S\) is well defined, we assume that the first $p$
eigenvalues $\lambda_1 \ge \lambda_2 \ge \dots \ge \lambda_p$  and the eigengap $\delta_p:= \lambda_p -\lambda_{p+1} $ are positive. If an eigenvalue \(\lambda_i\) has multiplicity \(k\), we choose the corresponding eigenvectors to be any orthonormal basis of its \(k\)-dimensional eigenspace.

We start by defining {\it skewness} parameters, which measure the interaction of $E$ with the eigenvectors of $A$. The word "skewness" has a geometric meaning.  For instance, $x$ shows how the vectors $u_i$ and $Eu_j$ are "skewed" from each other. The case where they align (as in the example given in the introduction) corresponds to the worst-case analysis, and we try to use the skewness to our advantage. 

\noindent\begin{definition} \label{xyz0} 
Define

\begin{itemize} 
 
 \item  $x:= \max _{i,j} | u_i ^T E u_j | $,   where  $ 1\le i,j \le r. $

 
\item $y := \max_{1 \le i < j \le r }\frac{| u_i^T E (I- U_{r} U_{r } ^T) Eu_j| } {\lambda_p }, $  where $U_{r}$ is the $n \times r$ matrix whose columns are the eigenvectors of the nontrivial eigenvalues of $A$.  

\item   $w:= \max \|u_i^T E\| $,  $1\le i \le r $. 

\end{itemize} \end{definition}


  We define the  {\it relative strength} of $E$ (with respect to this setting) as follows 
\begin{equation} \label{resnorm} 
 \| E \|_{A,p}:= \sqrt {r} x + \sqrt{r} y. \end{equation}

 In this definition, $x$ plays the role of the first-order term, and $y$ the second-order term, while \(\sqrt r\) serves as a normalization factor. For a trivial bound 
 on $x$ and $y$, we can use $x \le \| E\|$, and $y  \le \frac{\| E \|^2}{\lambda_p}$.
If one uses these trivial bounds, one can recover the Davis--Kahan theorem in a slightly worse form.  Our key point is that \(x\) and \(y\) capture the actual interaction between \(E\) and the eigenvectors of \(A\), and in many settings considered throughout this paper, they can lead to much better bounds. As an illustration, at the end of this section, we show that for random $E$, this improvement can yield a gain of order \(\sqrt n\).

\begin{theorem} \label{deterministicHRfirstKeigenvalues} 
Assume that  $ {(\bf C0):}  \,\,\,
\delta_p:= \lambda_p -\lambda_{p+1}  \ge  \max \{ 12 rx ,  \frac{144 r w^2}{  \lambda_p  } \} . $
Then
\begin{equation} \label{finalboundHRfirstKeigenvalues}  \| \tilde  \Pi_{p} - {\Pi}_{p} \| \le  12 \sqrt{p}  \frac{\| E\|_{A,p} }{\delta_p} + 12 \sqrt p \frac{ \|E\|}{\lambda_p} .
\end{equation} 
\end{theorem}

 Let us comment on the terms in our bound.  The first term $\frac{\| E\|_{A,p}}{\delta_p} $ (ignoring the  $12 \sqrt p$ factor) 
realizes our goal, as in this term we replace the spectral 
norm $\| E\|$  in the Davis--Kahan bound 
$\frac{\| E\|}{\delta_p} $ by the relative strength $\| E \|_{A,p}$. The  second term $ \frac{\| E\| } {\lambda_p }$   is a necessary 
adjustment. At the intuitive level, one should expect the appearance of 
this {\it noise-to-signal} term. Speaking in engineering terms, once the intensity of noise is too large compared to a signal (eigenvalue), it could destroy the eigenvector of that signal completely. The following result from random matrix theory  \cite{B-GN1} provides rigorous support. 

\begin{theorem} If $E$ is a random matrix with independent, zero-mean, normally distributed entries and $\|A\| =c \|E\|$, for a constant $c<1/2$, then with probability $1-o(1)$,
$$\| \tilde{u}_1 \tilde{u}_1^T - u_1 u_1^T\| = 1 - o(1).$$
\end{theorem}

This term $\frac{\| E\|} {\lambda_p }$ also improves the Davis--Kahan bound, as $\delta_p:= \lambda_p - \lambda_{p+1} \le \lambda_p $. In many applications, the leading eigenvalues $\lambda_1, \dots, \lambda_p$ are guaranteed to be large, while we are uncertain about the gap. 

Theorem \ref{deterministicHRfirstKeigenvalues} and 
Theorems \ref{deterministicHR} and \ref{deterministicHRgeneral} below are sharp up to the $12 \sqrt p$ factor; see the detailed discussion in  Appendix \ref{section:sharpness}. In fact, our analysis provides a way to write down the perturbation almost precisely;  see the discussion at the end of Section \ref{section:method}.

For a numerical illustration, consider the case where $E$ is a random Wigner matrix and $p=O(1)$. It is well known 
from random matrix theory that $\| E \| = \Theta (\sqrt n)$ with probability $1- o(1)$; see \cite{Ver1book,tao2012topics}. On the other hand, it is easy to prove that, with probability \(1-o(1)\), 
$x= O( \log r )$ and $y= O \big(\frac {\sqrt n \log r + r \log r }{\lambda_p} \big)$; see Appendix \ref{AppendixA}. These estimates improve upon the trivial bounds on \(x\) and \(y\) discussed before Theorem~\ref{deterministicHRfirstKeigenvalues} by a factor of \(n^{1/2-o(1)}\). It follows that, with probability $1-o(1)$,
$$ \| \tilde  \Pi_{p} - {\Pi}_{p} \| = O \left(\frac{r^{1/2} \log r }{\delta_p} +  \frac{ \|E\|} {\lambda_p } \right) = O \left(\frac{r^{1/2} \log r }{\delta_p} +  \frac{ \sqrt{n}} {\lambda_p } \right).$$ 
In comparison, Corollary~\ref{DKp}, which follows from the Davis--Kahan theorem, gives
$$ \| \tilde  \Pi_{p} - {\Pi}_{p} \|= O \left(\frac{\| E \| }{ \delta_p} \right) = O \left(\frac{\sqrt n}{\delta_p } \right),  $$
with probability $1-o(1)$.
Thus, our bound improves upon Davis--Kahan whenever
\(r\leq n^{1-\alpha}\) for some constant \(\alpha>0\) and
\(\delta_p=o(\lambda_p)\).

A fact of equal importance is that one can apply Theorem~\ref{deterministicHRfirstKeigenvalues} without the strong gap assumptions \eqref{gaptonoise0} and \eqref{gaptonoise1}. In the present setting, it suffices to have $\delta_p \gg r \cdot \max\{\log r, \frac{n}{\lambda_p}\} $. In many applications, $r=O(1)$ and $\lambda_p=\Theta(n)$ (see, e.g., \cite{UT1, CP1, JP1, J1, CT1, abbe2020entrywise}, and Section \ref{section:applications} for further discussion), and hence this assumption requires a very modest separation between  $\lambda_p$ and $\lambda_{p+1}$.

\subsection{The arbitrary rank case} \label{section:arbitraryrank} We now turn to the case where \(A\) has arbitrary rank, while still taking \(S=\{1,\dots,p\}\). We determine the relevant eigenvectors as follows. First, fix a \textit{safe distance} \(\bar\lambda\), and let \(r\) be the smallest index such that
$
\lambda_p-\lambda_{r+1}>\bar\lambda.
$
The motivation is the following:  if \(\bar\lambda\) is sufficiently large, the perturbation of the eigenvectors corresponding to eigenvalues beyond this safe distance should have little or no influence on the perturbation of the projection \(\Pi_S\). Thus, we regard only the first \(r\) eigenvectors \(u_1,\dots,u_r\) as {\it relevant}  to the problem. We now define the skewness parameters:
\noindent \begin{definition} \label{xyzP} 
 Define

\begin{itemize}
 
 \item  $x:= \max _{i,j} | u_i ^T E u_j | $,   where  $ 1\le i, j \le r . $

 
\item $y :=\max _{i,j,k}  |u_i^T E \left( \sum_{l \ge r+1}  \frac{u_l u_l^T}{\lambda_k - \lambda_l} \right) E u_j | $,  where $ 1 \le i \neq j \le r $ and $1 \le k \le p$.

\item   $w:= \max \|u_i^T E\| $,  where   $1 \le i \le r.$

\end{itemize} \end{definition}

We set the relative strength as 
$$ \| E\| _{A,p, \bar \lambda} := \sqrt r x + \sqrt r y. $$

\begin{theorem} \label{deterministicHR}
    Assume that  $ {(\bf C0):}  \,\,\,
\delta_p \ge  \max \{ 12 rx ,  \frac{144 r w^2}{  \bar{\lambda} }\} . $
Then
\begin{equation} \label{deterministicHRbound}  \| \tilde  \Pi_{p} - {\Pi}_{p} \| \le  12 \sqrt{p}  \frac{\| E\|_{A,p, \bar{\lambda}} }{\delta_p} + 12 \sqrt p \frac{ \|E\|}{\bar{\lambda}}.
\end{equation} 
\end{theorem}

Compared to the previous theorem, we need to modify the (second-order) parameter \(y\) to account for contributions from the entire spectrum. In Sections~\ref{section:randomnoise}--\ref{section:applications}, we show that this contribution is negligible in several natural settings. Both Theorem~\ref{deterministicHR} and its rectangular version, Corollary~\ref{cor: recPip}, follow as corollaries of the next theorem, which is the main result of this paper.


\subsection  {Treatment of general $S$}  \label{subsection:generalS} Now we discuss the most general case when $A$ can have  any rank and $S$ is an arbitrary set of $p$ indices. Let $\Lambda_S:=\{ \lambda_i| i \in S\}$, and set 
$$ \delta_S := \min_{i \in S, j \notin S } | \lambda_i -\lambda_j | > 0,$$
 which guarantees that $\Pi_S$ is well defined. We set a  {\it safe  distance }
$\bar \lambda \ge \delta_S$ and let $N_{\bar \lambda} (S)$ be a set of indices 
 such that  for all $j \notin N_{ \bar \lambda }, \dist (\lambda_j, \Lambda_S ) > \bar \lambda$. Let $r$ be the cardinality of $N_{\bar \lambda}(S)$, $r := | N_{\bar \lambda} (S)|$.  Notice here that we allow  $\dist (\lambda_i, \Lambda_S ) $ to be larger than $\bar \lambda$ for some $i \in N_{\bar \lambda } (S)$. This flexibility is useful in some applications.  We now define the skewness parameters accordingly:
\noindent \begin{definition} \label{xyz} 
 Define

\begin{itemize}
 
 \item  $x:= \max _{i,j} | u_i ^T E u_j | $,   where  $ i,j \in N_{\bar{\lambda}}(S). $

 
\item $y :=\max _{i,j,k}  |u_i^T E \left( \sum_{l \notin N_{\bar{\lambda}}(S)} \frac{u_l u_l^T}{\lambda_k - \lambda_l} \right) E u_j | $,  where $ i \neq j \in N_{\bar{\lambda}}(S)$ and $k \in S$.

\item   $w:= \max \|u_i^T E\| $,  where   $i \in N_{\bar{\lambda}}(S).$

\end{itemize} \end{definition}

We set the {\it relative strength} of $E$ (with respect to this general setting) as 
\begin{equation} \label{resnorm} 
 \| E \|_{A,S, \bar \lambda}:= \sqrt rx + \sqrt r y. \end{equation}

The statement of the general theorem is slightly more involved than the previous ones. We first construct a simple closed contour \(\Gamma\) that contains exactly the eigenvalues indexed by \(S\), such that the distance from any eigenvalue of \(A\) to \(\Gamma\) is at least \(\delta_S/2\), and the distance from any \(\lambda_j\), \(j\notin N_{\bar\lambda}(S)\), to \(\Gamma\) is at least \(\bar\lambda/2\).
A simple choice is to construct \(\Gamma\) from horizontal and vertical segments, with the vertical segments bisecting the intervals between the eigenvalues in \(\Lambda(S)\) and their nearest neighbors to the left and right outside \(\Lambda(S)\). Let \(\widetilde S\) denote the set of indices corresponding to the eigenvalues of \(\tilde A\) lying inside \(\Gamma\), and define
$$
\tilde{\Lambda}_{\tilde S}
:=
\{\tilde{\lambda}_i:i\in\tilde S\},
\,\,\text{and}\,\,\,
\tilde{\Pi}_{\tilde S}
:=
\sum_{i\in\tilde S}\tilde u_i\tilde u_i^T.
$$

\begin{theorem} \label{deterministicHRgeneral} 
Assume that  $ {(\bf D0)}:  \,\,\,
\delta_S \ge  \max \{ 12 rx , \frac{144 r w^2}{  \bar \lambda }\}. $
Then,
\begin{equation} \label{finalboundHRgeneral}  \| \tilde  \Pi_{\tilde S} - {\Pi}_{S} \| \le  12 \sqrt{p}  \frac{\| E\|_{A,S, \bar \lambda} }{\delta_S} + 12 \sqrt p \frac{ \|E\|}{\bar \lambda}.
\end{equation} 
\end{theorem}


As in the discussion following Theorem~\ref{deterministicHRfirstKeigenvalues}, the key point is that we replace the spectral norm \(\|E\|\) in Theorem~\ref{DKgeneral} by the relative strength \(\|E\|_{A,S,\bar\lambda}\), at the cost of the adjustment term \(\frac{\|E\|}{\bar\lambda}\).
Again, a clear illustration is when $E$ is random, a setting we study in detail in Section~\ref{section:randomnoise}. In this case, the corollaries of Theorem~\ref{deterministicHRgeneral} yield improvements over Theorem~\ref{DKgeneral} and several other existing results, and settle a question posed by O'Rourke, Vu, and Wang in~\cite{OVK13}.

 We would like to point out a delicate point: in Theorem \ref{deterministicHRgeneral}, the two sets $S$ and $\tilde S$ are not necessarily the same. Thus, the LHS in Theorem \ref{deterministicHRgeneral} is not necessarily the perturbation 
$\tilde{\Pi}_{S} - \Pi_S$. However, nothing is lost. If $S = \tilde S$, then we are done. 
If $S \neq \tilde S$, and the bound in Theorem \ref{deterministicHRgeneral} is $o(1)$, then 
$\| \tilde \Pi_{S}- \Pi_S\| =1 -o(1)$. 


It is often easy to force $S=\tilde{S}$, using Theorem \ref{deterministicHRgeneral} and a properly constructed contour. Let us illustrate this by deducing Theorem~\ref{deterministicHR} from Theorem~\ref{deterministicHRgeneral}. We can assume the bound in Theorem \ref{deterministicHRgeneral} is meaningful, namely that its RHS  is less than 1. By the pigeonhole principle, the two spaces that we project onto must have the same dimension, which means $|\tilde S |=|S|$. We construct the contour  $\Gamma$ to be a rectangular 
box going very far to the right (passing $\|A \| + \| E \|$, say). By this construction and Weyl's inequality, the eigenvalues of  $\tilde A$ inside $\Gamma$  must be the largest $k$ eigenvalues for some $k$.  But since $p = |S| =| \tilde S|$, the contour must contain exactly the first $p$ eigenvalues of $\tilde A$.

Theorem~\ref{deterministicHR} follows from Theorem~\ref{deterministicHRgeneral} by setting
$
N_{\bar\lambda}(S)=\{1,2,\dots,r\}.
$
To obtain Theorem~\ref{deterministicHRfirstKeigenvalues}, we use the flexibility in the definition of \(N_{\bar\lambda}(S)\) and take it to consist of the indices of all nonzero eigenvalues of \(A\).

Finally, let us make a remark on Assumption {\bf D0}.
In the proof, we can formally strengthen this assumption as follows. Since for any two orthogonal projections $\tilde{\Pi}_{\tilde{S}}$ and $\Pi_S$,
\(
\big\|\tilde\Pi_{\tilde{S}}-\Pi_S\big\|\leq 1, 
\)
Theorem \ref{deterministicHRgeneral} is trivial if \( 12\sqrt{p}\frac{\|E\|}{\bar \lambda}
> 1.\)
Thus, we can, for free, 
assume that  $\frac{\sqrt{p}\|E\|}{\bar \lambda}
 \leq \frac{1}{12},$  and add this assumption to {\bf D0}. After this addition and a simple normalization, Assumption {\bf D0} can be written
as 
\begin{equation} \label{D0: full version}
    \max \left\lbrace \frac{\sqrt{p}\|E\|}{\bar{\lambda}}, \frac{rx}{\delta_S}, \frac{\sqrt{r} w}{\sqrt{\bar{\lambda} \delta_S}} \right\rbrace \leq \frac{1}{12}.
\end{equation}
We will use this extended Assumption \textbf{D0} in the proofs. 

\subsection{The structure of the paper}

We divide the rest  of the paper into three parts: 
A discussion on related results and sharpness of our new results; Methodology and proofs; and Applications. These parts are independent and can be read in non-linear order.

The first part consists of Section~\ref{subsection: related results JW} and Appendix~\ref{section:sharpness}. In Section~\ref{subsection: related results JW}, we compare our results with several recent refinements of the Davis--Kahan bound. In Appendix~\ref{section:sharpness}, we discuss the necessity of the terms and assumptions in our results and construct examples showing that our bounds are sharp up to constant factors.

The second part (Methodology and proofs) consists of four sections. 

\begin{itemize}

\item In Section~\ref{section:method}, we introduce our method of combinatorial contour expansion.

\item In Section~\ref{section:proof1}, we establish three main lemmas and use them to prove Theorem~\ref{deterministicHRgeneral}.

\item In Section~\ref{section: proof of main lemmas}, we first prove a combinatorial lemma (Lemma~\ref{contour1-HR} in Subsection~\ref{subsec: combi lemma}), which provides the key tool for defining the combinatorial profiles introduced earlier. We then prove the three main lemmas through a detailed analysis of these profiles.

\item In Section~\ref{section:proof3}, we prove the main theorems stated in Section~\ref{section:randomnoise} and Section~\ref{sec: extension}. To maintain the flow of the presentation, we defer several technical parts of these proofs to the appendices.

\end{itemize}

The third part (Applications) has two sections. 

\begin{itemize}

\item In Section~\ref{section:randomnoise}, we derive several new bounds for random perturbations \(E\). Among others, these results partially settle an open problem posed in~\cite{OVK13}. Random perturbations arise naturally in many applications and also provide a particularly clear setting in which to illustrate the improvement obtained by using the relative strength.

\item In Section~\ref{section:applications}, we discuss two applications in statistics and computer science: the estimation of eigenvectors and eigenspaces of covariance matrices, and matrix completion. Due to the technical nature of these applications, we focus on the main ideas and refer readers to subsequent papers for detailed treatments.

\end{itemize}


\section{Related results } \label{subsection: related results JW}

In recent years, various refinements of the Davis--Kahan theorem have been obtained under additional assumptions
\cite{EBW1, HLMNV1, HE1, Ip1, JW1, JW2, KX1, KL1, OVK13, OVK22, Vu1, Wu1, Z1, tran2025davis, xie2024euclidean}.
In this section, we discuss some of the most recent results and compare them with the bounds obtained in this paper. To simplify the presentation and comparison, we restrict our discussion to the important case that
$
S=\{1,\dots,p\}$ and $ p=O(1).
$
We focus here on deterministic perturbations; the case of random \(E\) will be treated in Section~\ref{section:randomnoise}.

Studies from \cite{HE1, KL1} showed that if 
\begin{equation} \label{JWassumption} 2 \| E \| \le \delta_p, \end{equation} 
then 
$$\tilde \Pi_p- \Pi_p= \sum_{i \le p } \sum_{j  > p } \frac{1}{\lambda_i - \lambda_j} \left( u_i u_i^T E u_j u_j^T + u_j u_j^T E u_i u_i^T \right)  + T . $$

The first term is the linear term in the Taylor 
expansion (see Subsection \ref{subsec: contour representation}), while  $T$ represents the remainder of the series. For $T$, the best-known estimate was (see \cite{KL1} and also the discussion following (1.3) in \cite{JW1}) 
\begin{equation} \label{Tbound}  \| T \| =  O \bigg[ \bigg( \frac{\| E\| } {\delta_p } \bigg)^2 \bigg] . \end{equation}
Notice that  $\| u_i u_i^T E u_j u_j^T  \|  = | u_i^T E u_j | $. Thus, the most convenient way to bound  the first term in \eqref{JWassumption}
is to have 
$$ \|  u_i u_i^T E u_j u_j^T + u_j u_j^T E u_i u_i^T  \| \le 2 \bar x $$ where 
$$\bar x =\max_{1 \le i, j \le n} | u_i^T E u_j |. $$
With this estimate, we have 
\begin{equation} \label{KLbound}  
\| \tilde \Pi_p -\Pi_p \| = O \bigg[  \bar x \sum_{ 1 \le i \le p < j \le n } \frac{1}{\lambda_i- \lambda_j }  +   \bigg(\frac{\| E\| } {\delta_p } \bigg)^2 \bigg]. \end{equation}

Under a similar assumption, in \cite{tran2025davis}, we recently proved   
\begin{theorem} \label{theorem:BT}
  Let $r \ge p$ be the smallest integer such that $ \frac{\norm{\lambda_p}}{2} \leq \norm{\lambda_p - \lambda_{r+1}} $, and set  $x: =\max_{i,j \le r} | u_i^T E u_j | $.
 Assume furthermore that  $$4\|E\| \leq \delta_p: = \lambda_p - \lambda_{p+1}  \leq \frac{\norm{\lambda_p}}{4}. $$  Then,
\begin{equation}
\big\|\tilde{\Pi}_p - \Pi_p\big\| = O \left[ \frac{\|E\|}{\norm{\lambda_p}} \log \left(\frac{6 \sigma_1}{\delta_p} \right) + \frac{r^2 x}{\delta_p} \right].
\end{equation}

\end{theorem}
The weakness of these results, as already noted in \cite{JW1},  is the strong signal-to-gap assumption \eqref{JWassumption}, which is basically the same as \eqref{gaptonoise0}. 
 
 In a recent paper, Jirak and Wahl  \cite{JW1} worked out a different assumption and proved the following theorem
 for positive semi-definite matrices. Set 
$$x_0:= \max_{i,j \geq 1} \frac{\norm{ u_i^T E u_j}}{\sqrt{\lambda_i \lambda_j}},$$
$$\textbf{r}_p:= \sum_{1 \leq i \leq p} \frac{\lambda_i}{ \min_{j > p}|\lambda_i -\lambda_j|} + 
\sum_{j >p} \frac{\lambda_j}{\min_{1 \leq i \leq p} |\lambda_j -\lambda_i|}.$$

 \begin{theorem}[Theorem 1, \cite{JW1}]  \label{JWtheorem}
  Let $A$ be a positive semi-definite matrix. If 
  \begin{equation} \label{JWnewassumption}
      \textbf{r}_p \leq \frac{1}{8x_0},
  \end{equation}
  then, 
  \begin{equation}
      \Norm{\tilde \Pi_p- \Pi_p} = \textstyle O \bigg( x_0 \sqrt{\sum_{\substack{1 \leq i \leq p \\ j > p}}  \frac{\lambda_i \lambda_j}{(\lambda_i -\lambda_j)^2}  }  \,\,\,\bigg).
  \end{equation}
 \end{theorem}

\medskip

Let us now compare our new bound, Theorem~\ref{deterministicHR}, with the recent refinements of the Davis--Kahan theorem mentioned above.

{\it \noindent Comparison with \eqref{KLbound}.}
As mentioned above, a significant advantage of Theorem~\ref{deterministicHR} is that it does not require the gap-to-noise condition~\eqref{JWassumption}. Even under this assumption, our bound is stronger when the leading term is \(\frac{x}{\delta_p}\). Indeed, in this case, Theorem~\ref{deterministicHR} gives
$
O\left(\frac{x}{\delta_p}\right).
$
By definition, \(x\leq \bar x\) (see Definition~\ref{xyz}), while \(\frac{1}{\delta_p}\) is only one of the \(p(n-p)\) terms in
$
\sum_{\substack{1\leq i\leq p\\p<j\leq n}}
\frac{1}{\lambda_i-\lambda_j}.
$
Thus, $
\frac{x}{\delta_p}
\leq
\bar x
\sum_{\substack{1\leq i\leq p\\p<j\leq n}}
\frac{1}{\lambda_i-\lambda_j}.
$
Due to the large number of terms in this sum, the first term in~\eqref{KLbound} can be significantly larger than \(\frac{x}{\delta_p}\).

\medskip

{\it \noindent Comparison with Theorem~\ref{theorem:BT}.}
As in the previous comparison, a significant advantage of Theorem~\ref{deterministicHR} is that it does not require the strong gap-to-noise condition~\eqref{JWassumption}. The bound in Theorem~\ref{theorem:BT} may be simpler to apply, as it does not involve the quantity \(y\). On the other hand, it contains an additional logarithmic factor and an additional factor \(r^{3/2}\). The proof of Theorem~\ref{theorem:BT} was, in fact, our first attempt to analyze the contour representation of \(\tilde\Pi_{\tilde S}-\Pi_S\). See Subsection~\ref{subsec: contour representation} for further details.

\medskip

{\it \noindent Comparison with Theorem~\ref{JWtheorem}.} One can notice that the definition of $x_0$ contains the term $\frac{ | u_n^T E u_n|}{\lambda_n } $, where  $\lambda_n =\sigma_n$ is the least singular value.  Thus, if $\lambda_n$ is close to zero (the matrix is near singular), then $x_0$ is very large and \eqref{JWnewassumption} is hard to satisfy. Our results, on the other hand,  do not involve  $\sigma_n$ (or any non-important eigenvalues). Intuitively,  if $\lambda_1$ and $\lambda_n$ are far apart, we should not care about  $\lambda_n$  in estimating the perturbation of $u_1$; see the discussion at the beginning of Subsection \ref{section:arbitraryrank} about safe distance, which explains this point in more detail. 

Another advantage of Theorem~\ref{deterministicHR} is that it does not require \(A\) to be positive semi-definite. Even when 
\(A\) is restricted to be positive semi-definite, the two results are not directly comparable, as both depend on the full spectrum of \(A\), but in different ways. A natural way to compare them is therefore to consider several representative spectral distributions and examine the resulting bounds.

Without loss of generality, we assume that \(A\) is positive definite; otherwise,
Theorem~\ref{JWtheorem} is trivial, as its right-hand side is infinite.  We may then write $\lambda_i =\lambda_1 h(i) $, where $1 \ge h(i) > 0 $ is a monotone decreasing function describing the decay of the
eigenvalues of \(A\).
We will consider three natural types of decay

\begin{itemize}

\item Polynomial decay: $h(i)= i^{-c}$, for some positive constant $c$.

\item Exponential decay: $h(i) = e^{-ci} $, for some positive constant $c$. 

\item Logarithmic decay: $h(i)= \log^{-c} (i),$ for some positive constant $c$.

\end{itemize} 


Next, we assume that the noise matrix $E$ is random (which is the most important application in both papers). This assumption makes it possible to estimate the terms $|u_i ^T E u_j |$. Indeed, in the case when $E$ is Wigner, with high probability, any
term $|u_i^T E u_j|$ (and also  the quantity $x$ and $\bar x$ defined before) are $O(\log n)= \tilde O(1)$.
 This also determines $x_0$, as one can easily show that in this case $x_0= O \big(\frac{\log n}{\lambda_n} \big)$, with high probability. 
The bound in Theorem \ref{JWtheorem} simplifies to 
\begin{equation} \label{newJWbound} 
\tilde O \bigg[ \frac{1}{ \lambda_n }  \bigg(\sum_{\substack{1 \leq i \leq p \\ p < j \leq n}} \frac{\lambda_i \lambda_j } {(\lambda_i -\lambda_j )^2} \bigg)^{1/2}  \bigg] = 
\tilde O \bigg[ \frac{1}{ \lambda_1 h(n) }  \Big( \sum_{\substack{1 \leq i \leq p \\ p < j \leq n}} \frac{h(i)h(j) } {(h(i) -h(j))^2}  \Big)^{1/2}  \bigg]. 
\end{equation}
In the application of Theorem \ref{deterministicHR}, we choose $r$ so that $\bar \lambda= \frac{1}{2} \lambda_p$. In all the natural 
spectral distributions we consider, we can choose $r= c(p)$, a constant depending on $p$.
As we set $p =O(1)$, the bound in our Theorem \ref{deterministicHR} becomes 
\begin{equation} \label{newTVbound}
\tilde{O} \left( \frac{\| E \| } { \lambda_p}  + \frac{1}{\delta_p} + \frac{y}{\delta_p} \right). 
\end{equation}
Taking into account the facts that $\| E\| = \Theta (\sqrt n)$ and $\lambda_i= h(i) \lambda_1$, we can rewrite \eqref{newTVbound} as 
\begin{equation} \label{newTVbound2}
\tilde{O} \bigg[ \frac{\sqrt n  } { \lambda_1 }  + \frac{1}{\lambda_1( h(p)-h(p+1)) } + \frac{1}{\lambda_1 (h(p)-h(p+1)) }  \sum_{l>p} \frac{1}{\lambda_1 (h(p)- h(l) )} \bigg]. 
\end{equation}

\textit{Case 1: $h(i)= i^{-c}$} To estimate the  RHS of \eqref{newJWbound}, notice that 
$$ \frac{1}{ \lambda_1 h(n) }  \Big( \sum_{\substack{1 \leq i \leq p \\ p<j \leq n}} \frac{h(i)h(j) } {(h(i) -h(j))^2}  \Big)^{1/2} 
\ge \frac{1}{\lambda_1} n^c \bigg(\sum_{j >p} j^{-c} \big)^{1/2} 
= \Theta \bigg(\frac{n^c (1+ n ^{1-c} )^{1/2}}{\lambda_1}  \bigg) = \Theta \bigg(\frac{n^c+ n^{1/2+c/2}}{\lambda_1} \bigg) . $$
Therefore, the RHS of \eqref{newJWbound} is at least 
\begin{equation} \label{JWpolynom} \tilde{\Omega} \bigg(\frac{n^c+ n^{1/2+c/2}}{\lambda_1}   \bigg).\end{equation}
On the other hand, our bound \eqref{newTVbound2} is at most 
\begin{equation*}
    \begin{split}
\tilde{O} \bigg[ \frac{\sqrt n  } { \lambda_1 }  + \frac{1}{\lambda_1 } + \frac{1}{\lambda_1^2 }  \sum_{n\geq l>p} \frac{1}{ (p^{-c}- l^{-c} )} \bigg] & = \tilde{O} \bigg[ \frac{\sqrt{n}}{\lambda_1} + \frac{1}{\lambda_1} + \frac{n}{\lambda_1^2} \bigg] = \tilde{O} \bigg( \frac{\sqrt{n}}{\lambda_1} \bigg), 
    \end{split}
\end{equation*}
which improves  \eqref{JWpolynom} by a factor of $\max \{n^{c-1/2}, n^{c/2} \}.$

\textit{Case 2: $h(i) = \exp(-ci) $.} By a similar estimation, the RHS of Theorem \ref{JWtheorem} is at least 
\begin{equation} \label{JWexpo} \tilde{\Omega} \bigg( \frac{\exp (nc)}{\lambda_1} \bigg).\end{equation}
On the other hand, the bound
in Theorem \ref{deterministicHR} is at most 
$$ \tilde{O} \bigg( \frac{\sqrt{n}}{\lambda_1} + \frac{1}{\lambda_1} + \frac{1}{\lambda_1^2} \sum_{n \geq l>p} \frac{1}{e^{-pc} -e^{-lc}} \bigg) = \tilde{O} \bigg( \frac{\sqrt{n}}{\lambda_1}  \bigg),$$
which improves \eqref{JWexpo} by a  factor of $\frac{e^{nc}}{\sqrt{n}}$.

\textit{Case 3: $h(i)= \log^{-c} (i)$. } In this case, both bounds are $\tilde O \left( \frac{\sqrt n }{\lambda_1} \right)$.
(The computation is a bit more tedious than the other two cases, so we leave it as an exercise.) 

\noindent\textbf{Other related works.} 
In \cite{Ip1}, the author replaced the original eigenvalue gap $\lambda_p -\tilde{\lambda}_{p+1}$ with either the absolute separation or the relative separation with respect to the $p$-leading eigenspaces of $(A,\tilde{A})$, 
generalizing the Davis--Kahan bound for arbitrary square matrices.  
Many other works (e.g., \cite{EBW1, AV, CTP1, abbe2020entrywise, lyu2022exact})  consider the perturbation in the infinity norm or the two-to-infinity norm, which is critical for several well-known problems, such as finding hidden partitions and matrix completion. The results concerning random noise $E$ will be discussed later in Section \ref{section:randomnoise}.

To conclude this section, let us also point out that our method is robust and can be used to treat a number of matrix functionals; see Section \ref{section:method} and our representative works \cite{TranVulowrank, TranVUeigenvalue, TranVishnoiVu2025, TranVu2}.



\section{The method of combinatorial contour expansion }\label{section:method}

In this section, we describe our method, which works for a variety of matrix functionals. 
 
\subsection{Contour representation.}  \label{subsec: contour representation} 
 We start with the general idea of using a contour integral to represent matrix perturbation. 
 Let us recall the Cauchy integral formula~\cite{CAbook, agarwal2011introduction}.
 
 \begin{theorem}[Cauchy integral formula] \label{theo: Cauchy}
Let $\Gamma$ be a simple closed contour, and let $f$ be an analytic function in the whole simply connected domain $S$ containing $\Gamma$. Then for any $k \ge 1$
\begin{equation}\label{Cauchy0} 
   \frac{1}{2 \pi {\bf i}} \int_{\Gamma} \frac{f(z)}{(z-\lambda)^k}\,dz 
   = 
   \begin{cases}
      \frac{f^{(k-1)}(\lambda)}{(k-1)!}, & \lambda \text{ inside } \Gamma, \\[4pt]
      0, & \lambda \text{ outside } \Gamma.
   \end{cases}  
\end{equation}
 Here and later, {\bf i} denotes $\sqrt {-1} $, $0! =1$, and $f^{(l)}(z)$ is the $l^{th}$ derivative. 
\end{theorem}

Let  $\Gamma $ be a contour containing $\lambda_i, i \in S$ where $S$ is a subset of $\{1, \dots, n\}$, 
 and assume that all $\lambda_j, j \notin S$ are outside $\Gamma$. 
 By the Cauchy integral formula, applying for $f \equiv 1$, one obtains the  following contour formula (often referred to as the Riesz projection formula) \cite{Book1, Kato1}:
 \begin{equation} \label{contour-formula} 
  \sum_{ i \in S} u_i u_i ^T =  \frac{1} {2 \pi {\bf i }} \int_{\Gamma}    (z-A)^{ -1} dz.  \end{equation} 
Applying a similar argument to $\tilde A$, we have\footnote{If some eigenvalues of $\tilde A$ lie on the left vertical edge of $\Gamma$, we may shift this edge by a random infinitesimal amount.} 
 \begin{equation} \label{contour-formula1} 
  \frac{1}{2 \pi {\bf i}}  \int_{\Gamma}   (z-\tilde A)^{-1} dz  = \sum_{ i \in \tilde S}  \tilde u_i  \tilde u_i ^T .  \end{equation} Thus, we obtain a contour representation for the perturbation 
 \begin{equation}  \label{contourrep}
 \tilde \Pi_{ \tilde S} - \Pi_S =  \frac{1} {2 \pi {\bf i} } \int_{\Gamma}   [(z-\tilde A)^{-1}- (z- A)^{-1} ]  dz. 
  \end{equation}

Using  the resolvent formula 
\(
 M^{-1} - (M+N)^{-1}= (M+N)^{-1} N M^{-1}
\)
and the fact that $z- A = (z - \tilde A) +E$, we obtain 
\begin{equation*} (z-\tilde A)^{-1} - (z- A)^{-1} = (z-A)^{-1}  E (z-\tilde A)^{-1} .  \end{equation*} 
Applying this identity repeatedly, we obtain, at least formally,
\begin{equation} \label{TaylorEx}
(z-\tilde A)^{-1}-(z-A)^{-1}
=
\sum_{s=1}^{\infty}(z-A)^{-1}\big[E(z-A)^{-1}\big]^s.
\end{equation}
To justify this expansion, one needs to verify that the series on the RHS converges. We defer this issue for the moment and refer to Subsection~\ref{TaylorExpansion} for a rigorous justification.

Now set 
\begin{equation}
F_s := \int_{\Gamma}  (z-A)^{-1} [ E (z-A)^{-1} ]^s dz.
\end{equation}
We have 
$$ \tilde \Pi_ {\tilde S}- \Pi_S=  \frac{1} {2 \pi {\bf i} } \sum_{s=1}^{\infty} F_s.$$
Therefore, 
\begin{equation} \label{bound0} 
\big\| \tilde \Pi _{\tilde S} - \Pi _S \big\| \le \frac{1}{2\pi} \sum_{s=1}^{\infty} \Norm{F_s}.
\end{equation}

Several previous works have used the Taylor expansion idea; see, for instance, \cite{Kato1, xia2021normal, KX1, KL1, ding2024esprit, JW1}. 
The key matter is how to bound $\| F_s \|$. In \cite[Chapter 2]{Kato1}, one simply uses 
\begin{equation}\label{trivialF_s}
    \| F_s \|  \le  \|E\|^{s}  \int_{\Gamma} \|(z-A)^{-1}\|^{s+1} |dz|  \le \left( \frac{\|E\|}{\delta_S} \right)^s  \int_{\Gamma} \| (z- A)^{-1} \| |dz| \le C_{\Gamma} \left( \frac{\|E\|}{\delta_S} \right)^s, 
\end{equation}
where $C_\Gamma$ is a constant depending on $\Gamma$. Summing over $s$, we obtain the classical Davis--Kahan bound, up to a constant factor.

A refinement is to compute the first-order term
\(
\int_{\Gamma}(z-A)^{-1}E(z-A)^{-1}\,dz
\)
explicitly and bound the remaining terms by
\(\textstyle
O\big[\big(\frac{\|E\|}{\delta_S}\big)^2\big].
\)
This approach leads to \eqref{KLbound} in \cite{KX1}, as discussed in the previous section. More generally, one could compute the first \(k\) terms explicitly and bound the remainder by
\(
O\big[\big(\frac{\|E\|}{\delta_S}\big)^{k+1}\big].
\)
However, the resulting expressions become increasingly complicated as \(k\) grows, while the strong gap-to-noise condition~\eqref{gaptonoise0} remains unavoidable.

In the early stages of our study, we exploit the repetitive nature of the series to work out a bootstrapping argument to obtain Theorem \ref{theorem:BT}. However, this approach still requires the strong gap-to-noise condition; see \cite{tran2025davis, TranVishnoiVu2025} for details. 

\subsection {A new argument}
Our new approach to estimating \(\|F_s\|\) is substantially different, involving two additional expansions and a combinatorial analysis. Recall that
\(
(z-A)^{-1}
=
\sum_{i=1}^n\frac{1}{z-\lambda_i}u_i u_i^T.
\)
To control \(F_s\), we first decompose each resolvent factor \((z-A)^{-1}\) appearing in \(F_s\) into two parts,
\[
P:=\sum_{i\in I}\frac{1}{z-\lambda_i}u_i u_i^T,
\qquad
Q:=\sum_{j\in J}\frac{1}{z-\lambda_j}u_j u_j^T,
\]
where \(I\) is the set of \textit{important} indices and \(J\) is its complement. Consequently,
$$
(z-A)^{-1}\big[E(z-A)^{-1}\big]^s
=
(P+Q)\big[E(P+Q)\big]^s.
$$
We now expand this product as the sum of $2^{s+1} $ terms, each of which is of the form $ R_1 E R_2 \dots, E R_{s+1} $,  where $R_i$ is either $P$ or $Q$. This is the second expansion in our proof, beyond the Taylor expansion. 

At first glance, the situation seems worse, as the expansion produces an enormous number of terms. However, we have two key observations which will help with the analysis. First, notice that the only term which does not involve any important direction is $ QE \dots EQ $ (which involves only $Q$ operators). Using the formula of $Q$, we further expand this product into the sum of $(n- r)^{s+1} $ terms, each  of the form 
$$ \prod_{l _i \in J} \frac{1}{ z- \lambda_{l_i}  }  u_{l_1} u_{l_1} ^T E  u_{l_2} u_{l_2} ^T E \cdots E u_{l_s} u_{l_s} ^T . $$
The key point here is that all eigenvalues $\lambda_{l_i}  $ involved in this form are  {\it outside } the contour. Thus, the integral of every product $\prod_{l_i \in J}   \frac{1}{ z -\lambda_{l_i  }  }$ is zero. 
Therefore, $ \int_{\Gamma}  QE \dots EQ dz$ vanishes completely. 

For the remaining terms, we are not so lucky. They simply will not vanish. Our second key observation is that any of these terms contains many $P$'s, each of which is adjacent to an $E$. We will rely crucially on this property to exploit the interaction of the important directions 
(contained in $P$) with the noise matrix $E$. 

For a concrete example, let us consider the term $ PEP$ ($s=1$). Writing $P= \sum_{i \in I} \frac {1}{ z-\lambda_i } u_i u_i^T $ and expanding the product (this is the third expansion in our analysis), we obtain 
$$ PEP= \sum_{i_1, i_2 \in I }  \frac{1} { (z-\lambda_{i_1} ) (z -\lambda_{i_2} ) }  u_{i_1} u_{i_1} ^T  E u_{i_2} u_{i_2} ^ T =  \sum_{i_1, i_2 \in I }  \frac{1} { (z-\lambda_{i_1} ) (z -\lambda_{i_2} ) } 
 u_{i_1} u_{i_2} ^T  (u_{i_1}^T  E u_{i_2}) . $$
Since $u_{i_1}, u_{i_2} $ are unit vectors, a  trivial bound for $|u_{i_1} ^T E u_{i_2}  |$ is $\| E \| $. However, as mentioned before,  if $E$ is in a general position with respect to 
both vectors $u_{i_1} $ and $u_{i_2} $, then we expect a saving.  This observation naturally motivates the critical ``skewness'' parameter \(x\) introduced in Definition~\ref{xyz}. 

The main technical part of our proof is to handle $ (P+Q) [E (P+Q)] ^s $ for a general $s$. After applying the second and third expansions, this expression becomes an (enormous) sum of terms of the form $$h(z) u_{l_1} u_{l_1} ^T E u_{l_2} u_{l_2} ^T E \dots E u_{l_{s+1} } u_{l_{s+1} } ^T, $$ where 
$$h(z)= \frac{1}{ \prod_{i_l \in I}  (z- \lambda_{i_l}) \times \prod_{j_{l'} \in J}  (z- \lambda_{j_{l'}})}.$$
Moreover, the function \(h(z)\) varies from term to term. The main technical challenge is that the complex-analytic and linear-algebraic parts of these terms cannot be estimated separately: their contributions must be tracked simultaneously, as the losses arising from one part are compensated by gains from the other.

There are many ways to integrate $h(z)$. 
For instance, after many trials, we have found that the most natural approach, using the residue theorem,  is a {\it wrong } way.
We end up designing an algorithm to execute this integration (Lemma \ref{contour1-HR}). A key point here is that we can use this algorithm to obtain a combinatorial profile of the above term, in the form of a graph. We next use this combinatorial profile to split the whole collection of terms into many subsums, which we will be able to estimate separately, by a combination of analytic, linear algebraic, and combinatorial tools; see Section \ref{section:proof2} for a more detailed discussion. For concrete examples, see the derivations of the bounds for $\|F_1\|$ and $ \|F_2\|$ in Subsections \ref{F_1Es} and \ref{subsec: F2}, respectively.

\subsection{Sharpness of the method} \label{subsec: sharpness of method} Our argument  often leads to sharp results.  Technically, we  are going to  estimate $F_1 $ and $F_2$ almost precisely, and obtain sharp bounds, 
 $\| F _1\| \le h_1$ and $\| F_2  \| \le h_2 $.
 The quantity $h:= h_1 +h_2$ will contain all the terms in our main theorem.  We next show that the remaining $F_s$ satisfies  
 \begin{equation} \label{Fs} \| F_s \| \le 2^{3-s} h,  \end{equation}  or some estimate of this sort (with $2^{3-s}$ replaced by another fast decaying function). 
  Combining these estimates yields
 $$ \sum_{s=1}^{\infty} \| F_s\| \le 3 h.$$

 Once this plan has been executed properly, we obtain the desired bound, which is sharp, as the computation concerning  $F_1$ and $F_2$ is almost exact, and their contribution dominates the rest. 
See Section~\ref{section:sharpness} for a further discussion of the sharpness of the main results.

\subsection{Generalization}  Consider the 
spectral decomposition $A= U \Lambda U^T$ of $A$,  where $\Lambda$ is the diagonal matrix with the $i$th entry being  $\lambda_i$. Let $S$ be a subset of 
$\{1, \dots, n \}$. 
Define 
\begin{equation} \label{f(A)def} g(A):= f_S(A)=  U f_S(\Lambda ) U ^T , \end{equation} where $f_S(\Lambda) $ is the diagonal matrix with the $i$th entry being  $f(\lambda_i)$ if $i \in S$ and zero otherwise. This is one of the most common matrix functionals considered in the literature; see \cite{Book1, Higham2008}.  The case $f \equiv 1$ gives the 
orthogonal projection onto an eigenvector subspace spanned 
by the eigenvectors defined by $S$, which is the topic of this paper. 
Another important example is when $f(z)= z$ and $S$ defines the set of the $p$ largest singular values of $A$. In this case, $g(A)$ is the best rank $p$ approximation of $A$. These two operators alone are already of huge importance 
with respect to  applications; see, for instance, \cite{ GoLo, KVbook, KMO2, KMO1, CT1, CP1, JP1, J1, abbe2020entrywise, C1, JW1, JW2}.

In a forthcoming paper \cite{TranVu2}, we extend our method to treat the perturbation 
$\| f_S (\tilde A)- f _S(A) \|$. The above scheme extends naturally, and $F_s$ becomes 
$\int_{\Gamma} f(z) (z-A) ^{-1} [ E (z-A)^{-1} ]^s dz $, thanks to the Cauchy integral formula. However, the estimation becomes much more involved, as the derivatives of $f$ do not vanish.

\section{Eigenspace perturbation with respect to random noise}\label{section:randomnoise}

In this section, we consider the case where \(E\) is random, with independent but not necessarily identically distributed entries. In particular, we allow the entries to have different variances, and the resulting variance profile will play an important role in our analysis.

In Subsection~\ref{DKdiscussion}, we discuss several previous results on this topic and formulate a problem that partially motivated our study. We then present new results that provide a partial solution to this problem and yield improvements in a broad range of settings.

\subsection{Existing results for eigenspace perturbation}\label{DKdiscussion}  

The case where $A$ has a low rank (say, $r$) and $E$ is random was first studied by the second author in \cite{Vu1}. The leading intuition in \cite{Vu1}  is that with respect to $A$, $E$  should behave like an $r$ dimensional 
random matrix, which would yield quantitative improvements. 

In the ideal case where $E$ is a Wigner matrix, its norm is 
$\Theta (\sqrt n )$ (with high probability). Thus, Theorem \ref{DKp} yields a bound of order  $\frac{\sqrt n}{ \delta_p} $. In the spirit of the above intuition, we would expect to replace $\sqrt  n$ by $\sqrt r$, or (perhaps more realistically) by some quantity depending polynomially on $r$ and $\log n$. In \cite{Vu1}, the second author proved

\begin{theorem}[\cite{Vu1}] \label{Vu1} For any positive constants $\alpha_1, \alpha_2$ there is a positive constant $C$ such that the following holds. Assume that $r:= \rank A \leq n^{1 - \alpha_1}$ and $0< \lambda_1 \leq n^{\alpha_2}$. Let $E$ be a random Bernoulli matrix. Then with high probability 
\begin{equation*}
\Norm{\tilde{\Pi}_1 - \Pi_1 }^2 \leq C \bigg( \frac{\sqrt{r \log n}}{\delta_1}+ \frac{\| E\| ^2}{\delta_1 \lambda_1} + \frac{\|E \| }{\lambda_1} \bigg).
\end{equation*}

\end{theorem}

The first term in the bound is $\frac{\sqrt{r \log n}}{\delta_1}$, mirroring the above intuition. In many cases, this happens to be the leading term. 
However, we have an extra square on the LHS.  Continuing in this direction, about 10 years ago, O'Rourke, Wang, and the second author \cite{OVK13} obtained the following stronger and more general theorem.  
\begin{definition} We say the $m \times n$ random matrix $E$ is $(C_1, c_1, \gamma)$-concentrated if for any unit vectors $u \in \mathbb{R}^m, v \in \mathbb{R}^n$, and every $t >0$, 
\begin{equation*}
\P \left( \norm{u^T E v} > t \right) \leq C_1 \exp \left( -c_1 t^\gamma \right).
\end{equation*}
\end{definition}
\begin{theorem}[O'Rourke, Wang, and Vu \cite{OVK13}] \label{OVW2} 
Assume that $E$ is $(C_1,c_1,\gamma)$-concentrated for a trio of constants $C_1,c_1,\gamma > 0$ and $A$ has rank $r$. Let $1 \leq p \leq r$ be an integer. Assume that $\lambda_p > 0$. Then for any $t>0$, 
\begin{equation*}
\Norm{\tilde{\Pi}_p - \Pi_p } \leq 4 \sqrt{2 p} \bigg( \frac{\|E\|}{\lambda_p}+ \frac{t r^{1/\gamma}}{\delta_p}+ \frac{\|E\|^2}{\lambda_p \delta_p} \bigg),
\end{equation*} 
with probability at least $1 - 6 C_1 9^j \exp \left( -c_1 \frac{\delta_p^\gamma}{8^\gamma} \right) - 2C_1 9^{2r} \exp \left( -c_1r \frac{t^\gamma}{4^\gamma} \right). $

\end{theorem}
In the case when $E$ is Wigner with sub-Gaussian 
entries, we have $\gamma=2$. If we further assume that $p= O(1)$, then Theorem \ref{OVW2}  implies that with probability at least $.99$
$$ \Norm{\tilde{\Pi}_p - \Pi_p }= O \bigg( \frac{\|E\|}{\lambda_p}+ \frac{ \sqrt r}{\delta_p}+ \frac{\|E\|^2}{\lambda_p \delta_p} \bigg).$$
We notice that in this bound, besides the term $\frac{ \sqrt r } {\delta _p}$, which corresponds to our dimension reduction intuition, there are two other terms. The term 
$\frac{\| E \| }{\lambda_p} $ is the noise-to-signal ratio, and is necessary. Studies in random matrix theory
\cite{ABP1, B-GN1} showed that when this ratio is large, 
the sum $A +E$ behaves like a totally random matrix, so all information from $A$ is lost.  More mysterious is the 
$\frac{\|E\|^2}{\lambda_p \delta_p}$ term, which also appeared in Theorem \ref{Vu1} (which had a  different proof). 
Unlike the other two, there is no natural explanation for this term.  Should it be there, or is it just an artifact of the proof? 
\begin{question} \label{question1}
Is the third  term $\frac{\|E\|^2}{\lambda_p \delta_p}$ necessary? 
\end{question}

Motivated by this question, very recently O'Rourke, Wang, and the second author  \cite{OVK22} 
considered the case when $E$ is GOE. In this case, they showed that only the first two terms are sufficient. (Their original setting is for rectangular matrices; the study of the symmetric case is actually simpler, using the  same arguments.) In particular, they showed (under mild assumptions) that  with high probability,
\begin{equation} \label{OVW2023}
\Norm{\tilde{\Pi}_p - \Pi_p } = O \left( \frac{r^2 \log n} {\delta_p}+ \frac{\|E\|}{\lambda_p} \right) .
\end{equation} 
The same result holds for the general Wigner matrix, 
but the proof requires new technical lemmas. This gives the impression that the answer to Question \ref{question1} 
is negative. 

In  \cite{KX1}, Koltchinskii and  Xia 
also considered Gaussian noise and a general $A$ (which can have arbitrary rank). Among others, they showed (under some assumptions) that with probability $1-o(1),$
\begin{equation*}
\Norm{ \tilde{\Pi}_p - \Pi_p} = O \bigg[   \bigg(\frac{\|E\|}{\delta_{(p)}} \bigg)^2 + \frac{\sqrt{\log n}}{\delta_{(p)}} \bigg].
\end{equation*}

If we could ignore the term $\frac{\sqrt{\log n}}{\delta_{(p)}}$ on the RHS, then this is a quadratic improvement over the original Davis--Kahan theorem 
(Theorem \ref{DKp}). On the other hand, here we still need the gap-to-noise assumption \eqref{gaptonoise0}. 

The bounds of  Jirak and Wahl (see Section \ref{subsection: related results JW}) also apply well for the random setting; see \cite{JW1, JW2} for more details.  For many other related works,  see \cite{XZ1, Wa1, CCF1} and the references therein. 

\subsection{New results} \label{section:randomper2}


Before stating the main results \footnote{We would like to thank the reviewers for pointing out the mistakes in the previous versions of the theorems in this section.}, we list a few ensembles of random matrices with increasing generality. 

\begin{itemize}

\item Wigner matrices: This is the most well-known class of random matrices in the literature. Wigner matrices are Hermitian matrices with upper-triangular entries being independent (but not necessarily iid) random variables with mean 0 and variance 1. The diagonal entries are independent with mean 0 and variance 2.  (This class contains GOE and GUE as a subcase). 

\item Generalized Wigner matrices: Generalized Wigner matrices are Hermitian matrices with upper diagonal and diagonal entries being independent (but not necessarily iid) random variables with mean 0 and variance
$\sigma_{ij}^2, 1 \le i \le j \le n$. The second moments $\sigma_{ij}^2$ do not need to be the same, but are comparable, namely, there is a constant $C>0$ such that $\sigma_{ij}^2/ \sigma_{lk}^2 \le C$ for all
$i,j,k,l$. Finally $\sum_j \sigma_{ij}^2= \sum_j \sigma_{lj} ^2 $ for any $i$ and $l$. The last condition means the expectations of the square length of the rows are the same.

\item Regular matrices: We further extend the class of Generalized Wigner matrices by dropping the condition that the moments are comparable. Thus, a regular random matrix is a Hermitian matrix with upper diagonal and diagonal entries being independent (but not necessarily iid) random variables with mean 0 and variance
$\sigma_{ij}^2, 1 \le i \le j \le n$, where $\sum_j \sigma_{ij}^2= \sum_j \sigma_{lj} ^2 $ for any $i$ and $l$.

\end{itemize}

As already mentioned, in the case when  $E$ is random, the critical skewness parameters $x$, $y$, and hence the relative strength $\|E\|_{A,p,\bar{\lambda}}$ (see Definition \ref{xyz}) are relatively easy to estimate, and we can achieve explicit and significant savings compared to classical bounds. 

Regarding Question~\ref{question1}, the new results in this section reveal a subtle phenomenon. When $E$ is regular (a class that includes Wigner matrices), the answer to Question~\ref{question1} is indeed negative, consistent with the Gaussian case treated by O'Rourke et al.; see~\eqref{OVW2023}. On the other hand, the answer is positive in general: the error term in question is \emph{not} an artifact of the proof. It is genuinely necessary and originates from the third error term in Theorem~\ref{deterministicHR}. In the regular case, this term vanishes because of an algebraic cancellation; see the discussion leading to~\eqref{subtle}.

Let $E = (\xi_{ij})_{i,j=1}^{n}$ be an $ n \times n$ real symmetric random matrix, whose upper-triangular entries are  random variables $\xi_{ij}, 1\le i \le j \le n $, 
  which are independent, with mean zero and variance $\sigma^2_{ij}$. 
If  with probability 1, 
$$ |\xi_{ij}| \leq K,\,\,\, \forall 1\leq i \le j \leq n,$$ then we say that $E$ is $K$-bounded. If furthermore $\sigma_{ij}^2 \le \sigma^2$, for all $1 \le i \le j \le n$ also holds, we 
say that $E$ is $(K,\sigma)$-bounded. Define 
$$m_4:= \max_{1 \leq i\leq n} \frac{1}{n} \sum_{j=1}^n \E \xi_{ij}^4. $$
The parameter $m_4$ is the maximal row-wise average of the fourth moments of the entries of $E$. This quantity arises naturally in a second-moment
estimate for quadratic forms involving $E$.

\begin{remark} \label{remark: truncation} If a random variable $\xi$  has infinite support,  then we can truncate it 
at a properly chosen level $K$ where $\P (| \xi| \ge K) $ is negligible. After this, we can treat $\xi$  as a $K$-bounded variable.  For instance, $\mathcal{N}(0,1)$ is basically $100\sqrt {\log n}  $ bounded. 
It is easy to show that 
with probability at least $1 -n^{-98}$ (which is sufficiently large in most applications), a random matrix with iid $\mathcal{N}(0,1)$ entries is $100 \sqrt {\log n} $
bounded. See the introduction of \cite{TVuFourmoment} for more discussion about the truncation idea.

We recall the definition of sub-Gaussian random variable, which
appears frequently in both theory and practice. 
Similar to the Gaussian case, the truncation idea works well for this random variable.

\begin{definition} \label{def: subgau} A random variable $\xi$ is called sub-Gaussian if it has mean zero and there is a positive constant $C$ such that 
\[
\|\xi\|_{\psi_2}:= \inf \big\{ c > 0: \E \exp \big(\xi^2/c^2 \big) \leq 2 \big\} \leq C.
\]
This particularly implies that for every $t > 0$, 
    $\P (|\xi| \geq t) \leq 2 \exp(-t^2/C^2).$
\end{definition}
\end{remark}

\subsection{The low rank case} We first consider the case where $A$ has low rank and $\lambda_p > 0$. In this case, it is convenient to choose the parameter $r$ in Theorem \ref{deterministicHR} to be exactly the rank of $A$.  Given this choice of $r$, we have that $\lambda_j =0$ , for all $j \ge r+1$. We also set the safe distance $\bar{\lambda}= .99 \lambda_p$. 
Define 
\begin{equation} \label{epsilon12} \epsilon_1:=\frac{\|E\|}{\lambda_p  } \,\,\,\text{and}\,\,\, \epsilon_2=\frac{1} {\delta_p }.\end{equation}
Similar to previous sections, we now define our key assumption 
\begin{equation}
(\textbf{C2}): \,\,\, \delta_p\geq \max \bigg\{12 r t_1, \frac{144 r \|E\|^2}{\lambda_p}\bigg\}.
\end{equation}
Here $t_1$ is a parameter that will be a  part of the theorem. Finally, we set 
 $$\mu:= \max_{1 \leq k <l \leq r} \bigg\{ \bigg|\sum_{i=1}^n u_{ki} u_{li} \sum_{j=1}^n \sigma_{ji}^2 \bigg| \bigg\}.$$

Notice that we  can also write $\mu$ as 
\[
\mu:= \max_{1 \leq k <l \leq r} \big|\E (E u_k) \cdot (E u_l)  \big| = \max_{1 \leq k <l \leq r} \big|\E ( u_k^T E^T E u_l)  \big| = \max_{1 \leq k<l \leq r}|u_k^T \E E^2 u_l|.
\]
Thus, \(\mu\) measures the largest off-diagonal interaction of the second-moment matrix \(\E E^2\) with respect to the eigenvectors
\(u_1,\ldots,u_r\) of \(A\).


\begin{theorem} \label{mainTh01}
 Let $A$ be a symmetric matrix with rank $r$ and $E$ be a 
symmetric $(K,\sigma)$-bounded random matrix. Assume that $\lambda_p > 0$ and that Assumption \textbf{C2} holds with probability $1- \alpha_1$ for some choice of $t_1$. Then for any $t_2 >0$, with probability at least  
$$ 1 - \alpha_1- 2r^2 \exp\left( \frac{-t_1^2/2}{2\sigma^2 + Kt_1} \right) - \frac{5r(r-1)n m_4}{2t_2^2}, $$
we have 
$$\Norm{\tilde{\Pi}_p - \Pi_p} \leq 12 \sqrt{p} \left( \epsilon_1 + \sqrt{r} t_1 \epsilon_2 + \frac{ \sqrt{r} (t_2+\mu+rt_1^2)}{\|E\|} \epsilon_1 \epsilon_2 \right).$$
\end{theorem}


We can derive Theorem \ref{mainTh01}  from Theorem \ref{deterministicHR} by combining the inequality 
\begin{equation} \label{yboundbyx}
    y=\max_{\substack{1\leq i<j \leq r\\ l \leq p}} \bigg| \frac{u_i^T E (I -\sum_{k \leq r}u_k u_k^T)E u_j}{\lambda_l} \bigg|
\leq
\max_{1\leq i<j\leq r}
\frac{|u_i^T E^2u_j|}{\lambda_p}
+
\frac{rx^2}{\lambda_p},
\end{equation}
with standard probabilistic tools such as Chebyshev's inequality, to obtain suitable high-probability bounds on \(x\) and \(|u_i^T E^2u_j|\).  The formal proof requires some careful bookkeeping due to the number of parameters involved, and is deferred to Subsection~\ref{subsec: proof random} and
Appendix~\ref{subsec: moment easy}.


In a typical application of Theorem \ref{mainTh01}, the parameters
$\sigma, K, n, m_4 $ and $\mu$ are given or can be estimated explicitly. Once these parameters are determined, we optimize $\alpha_1$, $t_1$ and $t_2$ under the constraint that Assumption {\bf C2} holds with probability at least $1 -\alpha_1$, with an eye on our priority (having a better probability bound 
or a better perturbation bound).


Notice that when $E$ is regular,  $\sum_{j=1}^n \sigma_{ji}^2 $ does not depend on $i$. 
Since \(u_k\) and \(u_l\) are orthogonal for \(k\neq l\), it follows that \(\mu=0\). In this case, 
the bound in Theorem~\ref{mainTh01} reduces to 
\begin{equation} \label{subtle} \Norm{\tilde{\Pi}_p - \Pi_p} \leq 12 \sqrt{p} \bigg( \epsilon_1 + \sqrt{r} t_1 \epsilon_2 + \frac{ \sqrt{r} t_2 +r^{3/2} t_1^2}{\|E\|} \epsilon_1 \epsilon_2 \bigg).\end{equation}


By a routine computation (see Appendix \ref{app: 5.3 to concretebound} for details), we deduce the following corollary.  

\begin{corollary} \label{lowrank-random1}
Let  $A$ be a symmetric matrix with rank $r$ and $E$ be a 
random regular matrix with sub-Gaussian entries satisfying $\max_{1 \leq i,j \leq n}\|\xi_{ij}\|_{\psi_2} \leq C$. Assume that $\lambda_p >0$.  There are constants $C_1, C_2 > 0$, depending on $C$, such that if  
\[
\delta_p \geq C_1r \cdot \max \bigg\{\sqrt{r} \log^{3/2} n,  \frac{n}{\lambda_p} \bigg\},
\]
then, with probability $1-o(1)$,
$$\Norm{\tilde{\Pi}_p - \Pi_p} \leq C_2 \sqrt{p} \left( \frac{\sqrt{n} } {\lambda_p} + 
\frac{ \sqrt {r} \log^{3/2} n } {\delta _p } \right).$$
\end{corollary}

This corollary provides a negative answer to Question \ref{question1}.




\begin{remark}
    In many applications, $p$ is of order $\Theta(1)$ and the RHS can be written as $\tilde{O} \left( \frac{\|E\|}{\lambda_p} + \frac{\sqrt{r}}{ \delta_p} \right).$ As discussed above, these two terms are inevitable.    \end{remark}

\subsection{The full rank case}
Now we address the general case when $A$ can have any rank and spectral distribution. In this case, we choose the parameter $r$  (which is no longer the rank of $A$)  following Theorem \ref{deterministicHR}. We need to modify the definition of our parameters accordingly. First, let 
$$\epsilon_1 = \frac{\| E \| }{ \lambda_p - \lambda_{r+1} } \,\,\,\text{and}\,\,\, \epsilon_2 = \frac{1}{\delta_p}.$$
The analog of Assumption \textbf{C2} is 
$$(\textbf{D2}):\,\,\,\,\delta_p \geq \bigg\{ 12 r t_1, \frac{144 r \|E\|^2}{\lambda_p-\lambda_{r+1}} \bigg\}.$$
Next, we set 
$$s_{jj}^{(\bar p)} := (\lambda_p -\lambda_{r+1} ) \cdot \sum_{l >r} \frac{u_{lj}^2}{\lambda_{\bar{p}}-\lambda_l},\,\,\, \text{for each}\,\,\, 1 \leq j \leq n, 1 \leq \bar{p} \leq p,$$
and let 
$$ \mu_{A,E}:=  \max_{\substack{1 \leq k < l  \leq r \\ 1\leq \bar{p} \leq p}} \bigg\{ \bigg|\sum_{i=1}^n u_{ki} u_{li} \sum_{j=1}^n \sigma_{ji}^2 s_{jj}^{(\bar{p})} \bigg| \bigg\}.$$
When \(r=\operatorname{rank}(A)\) and \(\lambda_p>0\), the parameter
\(\mu_{A,E}\) can be controlled by the simpler quantity \(\mu\). More precisely,
\[
\mu_{A,E}
\leq
\mu+r\sigma^2
\leq
\mu+rt_1^2.
\]
Thus, in the low-rank setting, one may use \(\mu\) in place of
\(\mu_{A,E}\), which leads to a simpler computation.


We obtain the following generalization of Theorem \ref{mainTh01}. 

\begin{theorem} \label{mainTh-002} 
Let $A$ be a symmetric matrix and $E$ be a symmetric $(K,\sigma)-$bounded random matrix with independent upper-triangular entries. Assume that Assumption \textbf{D2} holds with probability $1- \alpha_1$ for some choice of $t_1$. Then, for any $t_2 >0$, with probability at least 
$$1 - \alpha_1- 2r^2 \exp\left( \frac{-t_1^2/2}{2\sigma^2 + Kt_1} \right) - \frac{ 32 r(r-1)p n \sigma^2 (\sigma^2+K^2)}{t_2^2}, $$ we have 
$$\Norm{\tilde{\Pi}_p - \Pi_p} \leq 12 \sqrt{p} \left( \epsilon_1 + \sqrt{r} t_1 \epsilon_2 +  \frac{ \sqrt{r} (t_2+\mu_{A,E})}{\|E\|} \epsilon_1 \epsilon_2 \right).$$
\end{theorem}

\section {Applications} \label{section:applications}
In many important applications in theoretical computer science and data science, the input matrix has the form
\(
\widetilde A=A+E,
\)
where \(A\) has low rank (or approximately low rank) and \(E\) is random, providing an ideal setting to apply our theorems. We have also observed that, in most applications, \(E\) is not standard white noise, but rather a problem-specific random matrix arising from the underlying model. 

In what follows, we consider the following two problems, both of which are among central problems in data science:

\begin{itemize}
\item Estimating eigenvectors and eigenspaces of a hidden covariance matrix using those of a sample covariance matrix. 

\item Matrix completion from missing and noisy data.

\end{itemize}

We discuss how these problems fit into our framework and highlight quantitative improvements obtained from our new results. 
Many other applications of the
same spirit will appear in future papers.


\subsection{Estimating eigenvectors and eigenspaces of a hidden covariance matrix}  \label{section:spiked} 

\subsubsection{The general problem} Let $X = [ \xi_1, \,\, \xi_2,\,...,\,\, \xi_d]^T$ be a zero-mean random vector of dimension $d$ with  covariance matrix $A = (a_{ij})_{1 \leq i, j \leq d},$ where $a_{ij} = a_{ji} = \mathbf{Cov}(\xi_i, \xi_j).$  
Let  $X_1, X_2, \dots, X_n$ be $n$ iid samples of $X$ (referred to as the observed data).
We consider the sample covariance matrix
$$\textstyle \tilde{A} := \frac{1}{n} \sum_{i=1}^{n} X_i X_i^T.$$  
Consider the spectral decomposition of $A$, $ A = \sum_{i=1}^{d} \lambda_i u_i u_i^T$, in which \( \lambda_1 \geq \lambda_2 \geq \dots \geq \lambda_d  \ge 0\) are the eigenvalues of \( A \), with corresponding orthonormal eigenvectors \( u_i, 1 \leq i \leq d \).
Let $H_p :=\Span \{u_1, \dots, u_p \} $ be the eigenspace spanned by the $p$ leading eigenvectors and $\Pi_p$ be the orthogonal projection onto $H_p$. We define $\tilde{\lambda}_i, \tilde{u}_i, \tilde{H}_p$ and $\tilde{\Pi}_p$ similarly, with respect to $\tilde{A}$.

It is common practice in statistics and data science to use spectral parameters of the sample covariance matrix $\tilde A$ to estimate those of the (unknown) covariance matrix $A$. A central task here is to determine the accuracy of these estimates.  In particular, one would like to know under which conditions these estimators are {\it consistent } (i.e., the error term tends to zero with probability $1-o(1)$) (as $d$ tends to infinity); see \cite{JL1, JW2, jung2009pca, shen2016statistics, srivastava2013covariance, RV1, adamczak2010quantitative}. In this section,  we focus on the estimation of the leading eigenvectors and eigenspaces. Technically speaking, we want to estimate the errors
$\|\tilde{u}_p \tilde{u}_p^T - u_p u_p^T\|$ and 
$\|\tilde{\Pi}_p - \Pi_p\|$, for small $p$. 

There is a large literature studying 
$\|\tilde{u}_p \tilde{u}_p^T - u_p u_p^T\| $ ($\|\tilde{\Pi}_p - \Pi_p\|$); see the books \cite{pearson1901liii, hotelling1933analysis, jolliffe2011principal, lorenz1956empirical, BYZ1, WW1} and the references therein.  The problem has been completely solved in the classical setting when $d$ is fixed. Among others, it was shown that if $d$ is fixed and $n$ tends to infinity, then the error tends to zero \cite{AndersonPCA2003, dauxois1982asymptotic}. However, in the modern setting (high-dimensional statistics),  where $d$ is large ($d \rightarrow \infty$), and  $n$ can be smaller than $d$, the situation is far from well understood. 

In the high-dimensional setting, one way to attack the problem is to use methods from random matrix theory to obtain limiting results for the error terms; see, for example,  \cite{J1, BYZ1, FW1, BKYY1, BY2}. By the nature of the approach, bounds obtained by this method are asymptotically sharp. However, this approach requires strong assumptions, in particular $n$ and $d$ must be of the same order of magnitude, and so are all the eigenvalues of $A$;  see \cite{TranvuCovariance} for details. Furthermore, the limit is always a positive constant, which means that the estimators are inconsistent;  see, for example, \cite[Theorem 1]{JL1} or \cite[Section 11 - Theorem 11.5]{BYZ1}.

Another approach uses perturbation theory. 
One writes $\tilde A$ as $A +E$, where $E$ is a random matrix with zero mean, and uses results from perturbation theory, such as the Davis--Kahan theorem;  see, for instance, \cite{WW1, JW1, JW2, nadler2008finite, zhu2022high}. This method makes use of modern estimates on the norm of $E$; see \cite{koltchinskii2017concentration, T1, srivastava2013covariance, tikhomirov2018sample}. For a representative result, see \cite[Section 8]{WW1}. Recently, there have been several attempts to refine this approach, either using the contour integral method from numerical analysis \cite{KX1, KL1, JW1, JW2,  mas2003perturbation} or exploiting the fact that the matrix $E$ is random, using tools from high-dimensional probability  \cite{OVK13, OVK22, CCF1, XZ1}.

The perturbation approach applies in very general settings and is more applicable in practice. 
To the best of our knowledge, all results obtained using this approach provide only {\it upper bounds}, and little has been known about their sharpness. 
In the general setting
($d,n,A$  arbitrary), it seems 
very hard to obtain
limiting theorems. Thus, determining the correct order of magnitude  seems to be  the next best thing 
that one may achieve.

{\bf \noindent Problem 1.} 
Determine the order of magnitude of $\|\tilde{u}_p \tilde{u}_p^T - u_p u_p^T\| $  and $\|\tilde{\Pi}_p - \Pi_p\|$.

\noindent We will also address the consistency problem:

{\bf \noindent Problem 2.}  Find necessary and sufficient conditions for 
the estimators $\tilde{u}_p$ and  $\tilde{\Pi}_p $ to be consistent.

\subsubsection{Our new results}  
Using the results and ideas developed in this paper, we are able to make progress on these problems in \cite{TranvuCovariance}. 
In particular, we can provide an answer to Problem 1, for a wide range of the parameters $d$ and $n$, under mild assumptions on the covariance matrix $A$. In this range, we determine the error $\|\tilde{u}_p \tilde{u}_p^T - u_p u_p^T\|$ ($\| \tilde \Pi_p - \Pi_p \|$) up to a constant factor. 

As corollaries,  we obtain answers for Problem 2, again in a wide range of $d, n$, and $A$. We would like to point out that in these new necessary and sufficient conditions,  the lower bound required on $n$, the number of samples, {\it does not} depend on $d$, the dimension. 
Instead, it depends on the effective rank and the signal-to-gap ratio of $A$, both of which can be significantly smaller than $d$.  

In what follows, we present a few representative results from \cite{TranvuCovariance}. For more discussion, see \cite{TranvuCovariance}.

{\bf \noindent Effective rank and signal-to-gap ratio.} Two quantities that play an important role in our analysis are the effective rank $r_{\mathrm{eff}}$ (also referred to as the intrinsic rank) and 
the signal-to-gap ratio  $\frac{\lambda_p}{\delta_{(p)}}$, where $\delta_{(p)} $ is the distance from $\lambda_p$ to the closest eigenvalue.

The notion of effective rank (also known as intrinsic rank) was introduced in \cite{Vershynin2012} and later studied systematically in the context of covariance estimation in \cite{BuneaXiao2015}. Since then, it has become a fundamental concept in high-dimensional statistics. Among its many applications, effective rank plays a key role in the formulation of matrix Bernstein inequalities; see, for example, \cite{RV2, koltchinskii2017concentration, rudelson2013hanson, martinsson2020randomized, T1, Ver1book}.
For matrix $A$ with singular values $\sigma_1 \ge \dots \ge \sigma_m \ge 0$, we set 
$$r_{\mathrm{eff}} (A) : = \frac{\sum_{i=1}^m \sigma_i }{\sigma_1}. $$
As  $A$ is positive semi-definite, its singular values and eigenvalues are the same, so 
\begin{equation} \label{defrstable} r_{\mathrm{eff}} (A) : = \frac{\sum_{i=1}^d \lambda_i }{\lambda_1}=\frac{\mathrm{Trace}\, A}{\lambda_1}. \end{equation}

Our study reveals that the relation between these two parameters and the magnitude of $n$ governs the order of magnitude of $\left\| \tilde{u}_p \tilde{u}_p^T - u_p u_p^T \right\|$ 
and $\|\tilde \Pi_p- \Pi_p\|$.  Our method is effective even if the signal-to-gap ratio $\frac{\lambda_p}{\delta_{(p)}}$ tends to infinity with $d$. In fact, this ratio can be as large as a polynomial in $d$. This is an important difference compared to many previous works, where it was assumed that this ratio is bounded (e.g., \cite{FW1, BYZ1, JP1, J1}). This means that we can obtain a good estimator for an eigenvector even when its corresponding eigenvalue is very close to the nearest one.

{\bf \noindent Models for the random vector $X$.} A standard model for $X$ is $X= A^{1/2} Y$, where $Y$ is a random vector such that 
$\E YY^T = I_d $. 
It is easy to check that $\E XX^T = A$. The generality of the model is determined by the assumptions on $A$ and $Y$.

For simplicity, we present our representative result for perturbation of the leading eigenvector when $A$ has a block form.  For results concerning more complicated  models and treatment of $\|\tilde u_p \tilde u_p^T -u_p u_p^T \|$ or $\|\tilde \Pi_p  - \Pi_p \|$, for a general $p$, see \cite{TranvuCovariance}.

\noindent\textbf{Treatment of the leading eigenvector when $A$ has a block form.} 
Our goal is to estimate
$\| \tilde u_1 \tilde u_1^T -u_1 u_1^T \| $. We assume that 
\(
X=A^{1/2}Y,
\)
where $A$ is a $d\times d$ positive semi-definite matrix of the form
\[
A=
\begin{pmatrix}
\Lambda_r & 0\\
0 & V_d
\end{pmatrix},
\] and $Y$ is a random vector with sub-Gaussian coordinates. 

This model is one of the most popular models for spiked covariance matrices, whose study was initiated in \cite{J1} and has generated considerable interest; see 
\cite{BY2, BYZ1, JW1, JW2, BS1,shen2013surprising, P1, J1} and the references therein. 
Here, $\Lambda_r$ is an $r\times r$ matrix whose eigenvalues are $\lambda_1, \dots, \lambda_r$, the largest $r$ eigenvalues  of $A$.
We call these eigenvalues the spikes of $A$. $V_d$ is a $(d-r)\times(d-r)$ matrix whose eigenvalues are $\lambda_{r+1}, \dots, \lambda_d$, which form the non-spiked part of the spectrum;  see, e.g., \cite[Chapter 11]{BYZ1} for more details. 
We do not make any assumption about the eigenvalues. They can grow with $d$, and can have different orders of magnitude. 

We write $Y= (y_1, \dots, y_d)$,  and assume that $y_i$ are independent (not necessarily iid) sub-Gaussian random variables with mean 0 and variance 1, satisfying $\max_{1 \leq i \leq d} \|y_i\|_{\psi_2}^2 \leq K$ for some $K \geq 1$.  We assume that the parameters $K$ and $r$ are bounded, i.e., \begin{equation} \label{cond1} K,r =O(1) . \end{equation}

To distinguish between the spiked and non-spiked eigenvalues, we assume that 
\begin{equation} \label{cond2} \lambda_{r+1} \le \frac{1}{2} \lambda_1 .  \end{equation}
The constant $\frac{1}{2}$ is chosen for convenience; 
it can be replaced by any constant $0 < c < 1$. We also assume that  
\begin{equation} \label{cond3} r_{\mathrm{eff}} + \frac{\lambda_1}{\delta_1} \ge C_0, \end{equation} where $C_0 $ is a sufficiently large  constant depending on $K$ and $r$, to be chosen.

We will consider three ranges of the effective rank, depending on the relation between $r_{\mathrm{eff}} $ and $\frac{\lambda_1}{\delta_1}$.  We show that under assumptions \eqref{cond1}, \eqref{cond2}, \eqref{cond3} (with a proper choice of $C_0$),  the following three theorems, corresponding to these ranges, hold. 

\begin{theorem} [Small effective rank] \label{thm:small} For any constants $\epsilon, K, r$, there are constants $C_0, C_1, C_2 $ so that the following holds. 
If $r_{\mathrm{eff}} \leq  \frac{r\lambda_1}{\delta_1}$,
then with probability at least $1-\epsilon$,  
    $$   \left\| \tilde{u}_1 \tilde{u}_1^T - u_1 u_1^T \right\| \leq C_1 \cdot \min \bigg\{ \frac{ \lambda_1 }{\sqrt{n} \delta_1}, 1 \bigg\},$$
\noindent and  with probability at least $\frac{1}{16}$,
$$ C_2 \cdot \min \bigg\{ \frac{ \lambda_1 }{\sqrt{n} \delta_1}, 1 \bigg\} \leq  \left\| \tilde{u}_1 \tilde{u}_1^T - u_1 u_1^T \right\| .$$
\end{theorem}
\begin{theorem} [Medium effective rank] \label{thm:medium} 
For any constants $\epsilon, K,r $, there are constants $C_0, C_1, C_2, C_3$ so that the following holds.  If $\frac{r\lambda_1}{\delta_1} < r_{\mathrm{eff}} \leq r \big(\frac{\lambda_1}{\delta_1}\big)^2 $ and $n \geq C_3 \cdot \frac{\lambda_1}{\delta_1} r_{\mathrm{eff}}$, then with probability at least $1-\epsilon$,
\[
\left\| \tilde{u}_1 \tilde{u}_1^T - u_1 u_1^T \right\| \leq C_1 \cdot \frac{ \lambda_1 }{\sqrt{n} \delta_1},
\]
and with probability at least $\frac{1}{16}$, 
   \[
  C_2\cdot \frac{\lambda_1}{\sqrt{n} \delta_1} \leq \left\| \tilde{u}_1 \tilde{u}_1^T - u_1 u_1^T \right\|.
   \]
\end{theorem}


\begin{theorem} [Large effective rank] \label{thm:large} 
 For any constants $\epsilon, K,r $, there are constants $C_0, C_1, C_2, C_3$ so that the following holds. If $ r_{\mathrm{eff}} > r \big(\frac{\lambda_1}{\delta_1}\big)^2 $ and  $n \geq C_3 \cdot  r_{\mathrm{eff}}  \frac{\lambda_1}{\delta_1}   $, then with probability at least $1-\epsilon$, 
\[
\left\| \tilde{u}_1 \tilde{u}_1^T - u_1 u_1^T \right\| \leq C_1\cdot \sqrt{\frac{r_{\mathrm{eff}}}{n}},
\]
and with probability at least $\frac{1}{16}$,
\[
C_2 \cdot \sqrt{\frac{r_{\mathrm{eff}}}{n}}  \leq \left\| \tilde{u}_1 \tilde{u}_1^T - u_1 u_1^T \right\| .
\]
\end{theorem}

 In all three theorems, the lower bound and upper bound match up to the constants involved.

{\bf \noindent Necessary and sufficient conditions for consistency.}  
As immediate corollaries of the above theorems, we obtain the following results, addressing  Problem 2.

\begin{corollary} \label{consistent1}
  Assume that  $r_{\mathrm{eff}} \leq  \frac{r\lambda_1}{\delta_1}$. Then 
    $  \left\| \tilde{u}_1 \tilde{u}_1^T - u_1 u_1^T \right\| =o(1) $
 if and only if 
    $\frac{ \lambda_1 }{\sqrt{n} \delta_1} =o(1). $
\end{corollary}


%
\begin{corollary} \label{consistent2}
    Assume that  $\frac{r\lambda_1}{\delta_1} < r_{\mathrm{eff}} \leq r \big(\frac{\lambda_1}{\delta_1}\big)^2 $ and $n \geq C_3 \cdot r_{\mathrm{eff}} \frac{\lambda_1}{\delta_1}$. Then
     $  \left\| \tilde{u}_1 \tilde{u}_1^T - u_1 u_1^T \right\| =o(1) $ 
     if and only if 
    $\frac{ \lambda_1 }{\sqrt{n} \delta_1} =o(1)$ . 
\end{corollary}

\begin{corollary} \label{consistent3}
    Assume that $ r_{\mathrm{eff}} > r \big(\frac{\lambda_1}{\delta_1}\big)^2 $ and $n \geq  C_3 \cdot r_{\mathrm{eff}} \frac{\lambda_1}{\delta_1}$. Then
     $  \left\| \tilde{u}_1 \tilde{u}_1^T - u_1 u_1^T \right\| =o(1) $
 if and only if
    $ \sqrt{\frac{r_{\mathrm{eff}}}{n}} =o(1)$ .
\end{corollary}

As mentioned earlier, the lower bounds on the number of samples $n$ in these results do not 
involve the dimension $d$. For example, in Corollary \ref{consistent1}, we need 
$$\textstyle \frac{n} {(\frac{\lambda_1}{\delta_1})^2 } \rightarrow \infty,  $$ 
\noindent while in Corollary \ref{consistent3},
we need $n \geq C_3 \cdot \frac{\lambda_1}{\delta_1} r_{\mathrm{eff}}$
and 
$$\frac{n}{r_{\mathrm{eff}} } \rightarrow \infty.$$
The latest result we found that provides a necessary and sufficient condition for consistency is \cite[Theorem 8]{JW2}. In this result,  $n$ is required to be super-quadratic in $d$, namely, $\frac{n}{d^2} \rightarrow \infty$.


\subsection {Matrix completion} 

Let $A$ be a large matrix of which only a small fraction of the entries are observed (say, one percent). The \emph{matrix completion problem} asks whether one can recover $A$ from these observations, under suitable structural assumptions. A famous example is the Netflix prize problem \cite{netflix}: the rows of $A$ represent viewers, the columns represent movies, and the entries are ratings. Since each viewer rates only a small fraction of the movies, the goal is to predict the missing ratings. Similar problems arise naturally in large recommendation systems.

Matrix completion is a fundamental problem in data science with many applications; see, for example, \cite{AFSU1, AEP1, BLWY1, AFKMS1, TK1, Si1, CRT1}. Recovery is possible only under appropriate assumptions on $A$. A standard assumption is that $A$ has low, or approximately low, rank. One also needs an incoherence condition; see the references above.

Let us explain how matrix completion fits into our framework. Following the standard model, suppose that each entry of $A$ is observed independently with probability $\rho$ (say, $\rho=.01$). Define $\widetilde A$ by replacing each unobserved entry by zero and scaling each observed entry by $1/\rho$. Then $\widetilde A$ has independent entries and
\(
\E \tilde A=A.
\)
Thus,
\(
\tilde A=A+E,
\)
where $A$ is (approximately) low rank and $E$ is a random matrix with mean zero. This is precisely the type of perturbation considered in this paper.

There are many approaches to matrix completion; see \cite{TranlinhVu} for a mini-survey. The first solution, due to Candes et al.\ \cite{CT1, CR1}, used convex optimization and required essentially only the basic low-rank and incoherence assumptions. However, convex optimization programs usually
run slowly in practice. The survey \cite{li2019survey} mentioned the interior point-based methods SDPT3 \cite{Sturm1999} and SeDuMi \cite{TohToddTutuncu1999}, which can take up to $O(|\Omega|^2 \cdot n^2)$ floating-point operations (FLOPs). Here, $|\Omega|=\rho n^2$ is the number of observed entries. Thus, the number of operations is $O(\rho n^6)$, which is too large even for moderate $n$. This motivates the search for faster methods.


A popular alternative is the spectral method.
One first computes a low-rank approximation of $\widetilde A$ and uses perturbation bounds to show that it approximates $A$ in operator norm or Frobenius norm. This initial approximation is then used as a base, from which one can make various refinements to obtain an entrywise ($l_{\infty}$ norm) approximation;  see \cite{CT1,KMO1,C1,CR1,AV}.

To the best of our knowledge, most results based on the spectral method depend on the condition number
\(
\kappa=\frac{\sigma_1(A)}{\sigma_r(A)},
\)
and deteriorate as $\kappa$ grows; see, for instance, the main results of \cite{CT1,KMO1}. In a recent paper \cite{AV}, A. Bhardwaj and the second author developed an analysis that removes the dependence on $\kappa$, at the cost of assuming that the singular values of $A$ are sufficiently well separated.

Using the methods introduced in this paper, together with some additional ideas, the second author and L. Tran \cite{TranlinhVu} obtained an entrywise recovery result that requires neither a condition-number bound nor a singular-value gap assumption. More precisely, they showed that a suitably chosen low-rank approximation of $\tilde A$ approximates $A$, with high probability, entrywise to any prescribed accuracy.  The interesting fact here is that the low-rank approximation achieves the goal directly, without further refinements. 
Informally speaking, what is going on here is that the error in the Frobenius norm spreads somewhat evenly on all entries (thanks to the incoherence property of $A$). Thus, the error on each entry is also small. 

To conclude this subsection, we present the symmetric versions of the algorithm and the main result in \cite{TranlinhVu}. Assume that $\|A\|_{\infty} \leq K$ and $\mathrm{rank}\,A \leq r_{\max}$. We know $K, \rho$ and $r_{\max}$ (but may not know the exact rank of $A$). In fact, it is sufficient to know $\rho$ up to a constant factor  (with high probability), which is easy to achieve using the standard large deviation bounds. 

\begin{algorithm}[H]
\caption{Recovering every entry of matrix $A$ within an absolute error $\varepsilon_0$}
\label{alg:thresholdandround}
\renewcommand{\thealgorithm}{}
\begin{algorithmic}[1] 

\STATE {\text{Compute the truncated SVD:} $\tilde{A}_{r_{\max}}:= \sum_{i=1}^{r_{\max}} \tilde{\sigma}_i \tilde{u}_i \tilde{v}_i^T.$}      
\STATE {\text{Find the largest index $s \leq r_{\max} -1$ such that} $\tilde{\sigma}_s -\tilde{\sigma}_{s+1} \geq 20 K \sqrt{\frac{2r_{\max} n}{\rho}}.$ } \\
\text{If there is no such $s$, choose $s=r_{\max}$.} 

\STATE {Output $\tilde{A}_s =\sum_{i=1}^s \tilde{\sigma}_i \tilde{u}_i \tilde{v}_i^T$.}
\end{algorithmic}
\end{algorithm}
The validity of the above algorithm is given by the symmetric version of \cite[Theorem~1.5]{TranlinhVu}.

\begin{theorem}
    There is a universal constant $C>0$ such that the following holds. Suppose $r_{\max} \leq \log ^2 n$. Assume that 
    $$\sigma_1 \geq 100 r K \sqrt{\frac{r_{\max}n}{\rho}}\,\,\text{and}\,\, \rho \geq \frac{C}{n} \max \left\{\log^4 n, \frac{r^3 K^2}{\varepsilon_0^2 } \left( 1+ \frac{\mu_0^2}{\log^2 n}\right) \right\} \log^6 n,$$
where $\mu_0:= \max_{i \in [n]} \frac{n}{r}\|e_i^T U\|^2$  is the coherence parameter \footnote{$e_1, e_2, \dots, e_n$ are the standard basis vectors; for instance,  $e_1 = (1,0,0, \dots, 0)$.}. Then with probability $1 -O(n^{-1})$, the  Algorithm \ref{alg:thresholdandround} recovers 
(simultaneously) every entry of $A$ within an error at most $\varepsilon_0$. 
\end{theorem}


The method we developed in this paper applies to 
other settings of the matrix completion problem. For instance, it can accommodate heterogeneous missingness, where each entry $(i,j)$ is observed with its own probability $\rho_{ij}$. We can also use it to handle recovery with respect to other norms.  For example, in a very recent work \cite{TranVuell1norm}, we apply this framework to matrix completion with heterogeneous missingness and obtain a recovery guarantee in the two-to-infinity norm, which is the natural norm to consider in various medical experiments.

\section{Proof of Theorem \ref{deterministicHRgeneral} }  
\label{section:proof1}


In this section, we prove our key result, Theorem \ref{deterministicHRgeneral}. We will show  that under (the extended) Assumption {\bf D0} (see \eqref{D0: full version}), we have 
\[
 \| \tilde  \Pi_{\tilde S} - {\Pi}_{S} \| \le  12 \sqrt{p}  \frac{\| E\|_{A,S, \bar \lambda} }{\delta_S} + 12 \sqrt p \frac{ \|E\|}{\bar \lambda}.
\]

From now on, we use \(N_{\bar{\lambda}}(S)^c\) and \(S^c\) to denote the complements of \(N_{\bar{\lambda}}(S)\) and \(S\), respectively, in \({1,2,\ldots,n}\). It is important that no eigenvalue indexed by \(S\) coincides with an eigenvalue indexed by \(S^c\); i.e.,
\(
\lambda_i\neq\lambda_j,\) for all \(i\in S,\ j\in S^c.
\)

Following the plan outlined in Section~\ref{section:method}, we first estimate
\(\|F_1\|\) and \(\|F_2\|\) in Subsections~\ref{F_1Es} and
\ref{subsec: F2}, respectively. These estimates contain all the main terms
appearing in our final bound, and their derivations also serve as a warm-up
for the more involved combinatorial analysis required to estimate
\(\|F_s\|\) for general \(s\geq 3\). We then establish the general estimate
for \(\|F_s\|\) in Subsection~\ref{F_2Es}, using three key lemmas:
Lemmas~\ref{LemmaCase1M(a,b)bound}, \ref{LemmaCase2M(a,b)bound}, and
\ref{LemmaCase3M(a,b)bound}. Finally, we combine these estimates to prove
Theorem~\ref{deterministicHRgeneral} in
Subsection~\ref{subsec: provevialemma}.

\subsection{Estimating the first term $ \Norm{F_1}$} \label{F_1Es} 
In this subsection, we estimate the first term $\| F_1 \|$. Before starting the estimation, we record a technical lemma, which is a simple consequence of the Cauchy integral formula \eqref{Cauchy0} and will be used throughout this subsection.
\begin{lemma} \label{integral} Let $a, b$ be two points not on the boundary of $\Gamma$. Then 
$$ \frac{1}{2 \pi \textbf{i}}   \int_{\Gamma} \frac{1}{(z-a) (z-b) } dz    =\frac{1 }{  b-a } $$ if exactly $b$ lies inside $\Gamma$ and $a$ lies outside $\Gamma$. In all other cases, the integral equals zero.  
\end{lemma} 
\noindent\textit{Proof of Lemma \ref{integral}.} Notice that 
\begin{equation*}  \frac{1}{2 \pi \textbf{i}}\int_{\Gamma} \frac{1}{(z-a) (z-b)} dz  = \frac{1}{a-b} \Big( \frac{1}{2 \pi \textbf{i}} \int_{\Gamma} \frac{1}{ z-a } dz -  \frac{1}{2 \pi \textbf{i}}\int_{\Gamma} \frac{1}{ z-b} dz \Big).  
\end{equation*} 
The claim follows from  (\ref{Cauchy0}).

Back to our estimation  of $\|F_1\|$, recall our strategy that we split the resolvent $(z-A)^{-1}$ into two parts:
$$(z-A)^{-1}= \sum_{i=1}^{n} \frac{1}{z-\lambda_i} u_i u_i^{T} = P+Q ,$$ where 
 \begin{equation} \label{P-HR} P :=\sum_{i \in N_{\bar{\lambda}}(S)} \frac{u_i u_i^T}{z-\lambda_i} \,\,\,\text{and}\,\,\, Q :=   \sum_{j \in N_{\bar{\lambda}}(S)^c} \frac{ u_j u_j^T}{z -\lambda_j}. 
 \end{equation} 
\noindent For all $s \ge 1$, we  can rewrite $F_s$ as 
\begin{equation} \label{F_sEFor-HR}
F_s= \int_{\Gamma} (P+Q)[ E(P+Q)]^s dz.
\end{equation}
In particular, $ F_1 =\int_{\Gamma} (P+Q) E(P+Q)dz $, which can be decomposed as
$$ \int_{\Gamma} \left( PEP + QEP+PEQ +QEQ  \right) dz = \int_{\Gamma} QEQ dz+ \int_{\Gamma} PEQ dz + \int_{\Gamma} QEP dz+\int_{\Gamma} PEP dz.$$

Since $\Gamma$ does not enclose any point from the set $\{ \lambda_j\}_{j \in N_{\bar{\lambda}}(S)^c}$,  $\frac{1}{2 \pi \textbf{i}}\int_{\Gamma} \frac{1}{(z-\lambda_{j})(z-\lambda_{j'})} dz =0 $ for any pair $ j,j' \in N_{\bar{\lambda}}(S)^c$. Therefore
\begin{equation} \label{estimate1-HR} \frac{1}{2 \pi \textbf{i}}\int_{\Gamma} QEQ  dz =  \frac{1}{2 \pi \textbf{i}}\sum_{j,j' \in N_{\bar{\lambda}}(S)^c} \int_{\Gamma}  \frac{u_j u_j^T E u_{j'} u_{j'}^T}{(z-\lambda_{j})(z-\lambda_{j'})} dz = 0. 
\end{equation}

Next, we have 
\begin{equation*}
\begin{split}
\frac{1}{2 \pi \textbf{i}} \int_{\Gamma} P E Q dz & = \frac{1}{2 \pi \textbf{i}} \sum_{ \substack{i \in N_{\bar{\lambda}}(S)\\ j \in N_{\bar{\lambda}}(S)^c}} \int_{\Gamma} \frac{1}{(z- \lambda_i)(z - \lambda_j)} u_i u_i^T E u_j u_j^T dz \\
& = \sum_{\substack{i \in S \\ j \in N_{\bar{\lambda}}(S)^c}} \frac{1}{\lambda_i - \lambda_j} u_i u_i^T E u_j u_j^T (\text{\,by Lemma \, \ref{integral}}) \\
& = \sum_{i \in S} u_i u_i^T E \bigg( \sum_{j \in N_{\bar{\lambda}}(S)^c} \frac{u_j u_j^T}{\lambda_i -\lambda_j} \bigg).
\end{split}
\end{equation*}
Therefore, by the definition of the spectral norm, 
$$\frac{1}{2 \pi}\Norm{ \int_{\Gamma} P EQ dz} = \max_{\|\textbf{v}\|=\|\textbf{w}\|=1} \sum_{i \in S}  \textbf{v}^T u_i   u_i^T E \bigg( \sum_{j \in N_{\bar{\lambda}}(S)^c} \frac{u_j u_j^T}{\lambda_i -\lambda_j} \bigg) \textbf{w} .$$
By the triangle inequality, the RHS is at most
\[
 \max_{\|\textbf{v}\|=\|\textbf{w}\|=1}\sum_{i \in S} |\textbf{v}^T u_i|\cdot \bigg|u_i^T E \bigg( \sum_{j \in N_{\bar{\lambda}}(S)^c} \frac{u_j u_j^T}{\lambda_i -\lambda_j} \bigg) \textbf{w}  \bigg| \leq  \max_{\|\textbf{v}\|=1} (\sum_{i \in S} |\textbf{v}^T u_i|)  \  \bigg\|E \bigg( \sum_{j \in N_{\bar{\lambda}}(S)^c} \frac{u_j u_j^T}{\lambda_i -\lambda_j} \bigg)\bigg\|.
\]
Furthermore, by definition of $\bar \lambda$, we have $ \big\|E \big( \sum_{j \in N_{\bar{\lambda}}(S)^c} \frac{u_j u_j^T}{\lambda_i -\lambda_j} \big)\big\| \leq \frac{\|E\|}{\bar \lambda}$. By the Cauchy--Schwarz inequality, for any $\|\textbf{v}\|=1$, $\sum_{i \in S} |\textbf{v}^T u_i| \leq \sqrt{|S| \cdot (\sum_{i \in S} |\textbf{v}^T u_i|^2) } \leq \sqrt{|S|}=\sqrt{p}.$ These estimates imply 
\begin{equation} \label{estimate2-HR}
\begin{split} 
 \frac{1}{2 \pi} \bigg\|\int_{\Gamma} PEQ dz \bigg\| &  \leq \frac{\sqrt{p} \|E\|}{\bar{\lambda}}.
\end{split}
\end{equation}
Similarly, 
\begin{equation} \label{estimate21-HR}
    \frac{1}{2 \pi} \bigg\|\int_{\Gamma} QEP dz \bigg\| \leq \frac{\sqrt{p} \|E\|}{\bar{\lambda}}.
\end{equation}

Now we consider the last term $\int_{\Gamma} PEP dz$, which can be expressed as 
$$ \int_{\Gamma} \sum_{i,j \in N_{\bar{\lambda}}(S)} \frac{1}{(z -\lambda_i)(z-\lambda_j)}  u_j u_j^T E u_i u_i^T dz =   \int_{\Gamma} \sum_{  i, j  \in N_{\bar{\lambda}}(S)} \frac{1}{(z -\lambda_i)(z-\lambda_j)}  u_j (u_j^T E u_i ) u_i^T dz .$$
By Lemma \ref{integral},   
$\int_{\Gamma } \frac{1}{(z -\lambda_i)(z-\lambda_j)} dz = \frac{2\pi \textbf{i} }{\lambda_j - \lambda_i} $  
if $\lambda _i$ lies inside $\Gamma$ and $\lambda_j$ lies outside $\Gamma$. In all other cases, the integral is zero. Therefore, 
$$\frac{1}{2\pi} \Norm{\int_{\Gamma} PEP dz} =  \max_{\|\textbf{v}\|=\|\textbf{w}\|=1} \bigg( \sum_{\substack{i \in S \\ j \in N_{\bar{\lambda}}(S) \setminus S}} \textbf{v}^T u_i \frac{(u_i^T E u_j)}{\lambda_i - \lambda_j} u_j^T \textbf{w} +  \sum_{\substack{i \in S \\ j \in N_{\bar{\lambda}}(S) \setminus S}} \textbf{v}^T u_j \frac{(u_j^T E u_i)}{\lambda_i - \lambda_j} u_i^T \textbf{w} \bigg).$$
By the triangle inequality, the RHS is at most 
\begin{equation*}
    \begin{split}
         & \max_{\|\textbf{v}\|=\|\textbf{w}\|=1} \bigg( \sum_{\substack{i \in S \\ j \in N_{\bar{\lambda}}(S) \setminus S}} \norm{\textbf{v}^T u_i} \norm{ \frac{(u_i^T E u_j)}{\lambda_i - \lambda_j}} \norm{ u_j^T \textbf{w}}+ \sum_{\substack{i \in S \\ j \in N_{\bar{\lambda}}(S) \setminus S}} \norm{\textbf{v}^T u_j} \norm{ \frac{(u_j^T E u_i)}{\lambda_i - \lambda_j}} \norm{ u_i^T \textbf{w}} \bigg) \\
  & \leq  2 \max_{\|\textbf{v}\|=\|\textbf{w}\|=1}\sum_{\substack{i \in S \\ j \in N_{\bar{\lambda}}(S) \setminus S}} \norm{\textbf{v}^T u_i} \norm{ \frac{(u_i^T E u_j)}{\lambda_i - \lambda_j}} \norm{ u_j^T \textbf{w}}. 
    \end{split}
\end{equation*}
By the definitions of $\delta_S$ and $x$ (see \eqref{xyz}), we have 
\[
|\lambda_i -\lambda_j| \geq \delta_S\,\,\text{and}\,\, |u_j^T E u_i| \leq x\,\,\text{for any $i\in S,  j \in N_{\bar{\lambda}}(S) \setminus S $}.
\]
Therefore, we have
\begin{equation}  \label{estimate3-HR}
\begin{split}
\frac{1}{2\pi} \Norm{\int_{\Gamma} PEP dz}
  & \leq \frac{2x}{\delta_S} \max_{\|\textbf{v}\|=\|\textbf{w}\|=1} \sum_{i \in S, j \in N_{\bar{\lambda}}(S) \setminus S} \norm{\textbf{v}^T u_i} \norm{ u_j^T \textbf{w}}  \\
  & = \frac{2x}{\delta_S} \max_{\|\textbf{v}\|=\|\textbf{w}\|=1} \bigg( \sum_{i \in S} \norm{\textbf{v}^T u_i} \bigg) \bigg( \sum_{j \in N_{\bar{\lambda}}(S) \setminus S} \norm{u_j^T \textbf{w}} \bigg) \\
  & \leq \frac{2 \sqrt{p(r-p)} x}{\delta_S} (\text{by the Cauchy--Schwarz inequality}).
  \end{split}
 \end{equation}

By \eqref{estimate1-HR}, \eqref{estimate2-HR},\eqref{estimate21-HR}, and \eqref{estimate3-HR}, we obtain  
  \begin{equation} \label{F1-HR} 
   \frac{1}{2\pi}\|F_1\| \le \frac{ 2 \sqrt{p}  \| E\|}{ \bar{\lambda}}   + \frac{  2 \sqrt{ p r} x} { \delta_S } . 
   \end{equation}

\subsection{Estimating the second term $\| F_2 \|$.} \label{subsec: F2} 



In this subsection, we bound \(\|F_2\|\). Together with the estimate obtained in the previous subsection, these bounds will produce all the terms appearing in Theorem~\ref{deterministicHRgeneral}. As before, we start from the expansion
$$
\begin{aligned}
F_2
&=\int_{\Gamma}(P+Q)E(P+Q)E(P+Q)\,dz\\
&=\int_{\Gamma}QEQEQ\,dz
+\int_{\Gamma}QEPEQ\,dz
+\int_{\Gamma}PEQEP\,dz
+\int_{\Gamma}PEPEP\,dz\\
&\quad
+\int_{\Gamma}PEQEQ\,dz
+\int_{\Gamma}PEPEQ\,dz
+\int_{\Gamma}QEQEP\,dz
+\int_{\Gamma}QEPEP\,dz.
\end{aligned}
$$
We divide these eight integrals into three types according to their two end factors: those with \(Q\) at both ends (e.g., \(QEQEQ\)), those with \(P\) at both ends (e.g., \(PEQEP\)), and those with \(P\) at one end and \(Q\) at the other (e.g., \(PEQEQ\)). Each type requires a different treatment. The arguments in this subsection can be viewed as a warm-up for the corresponding treatment of the three types of terms in \(F_s\) developed later; see Subsection~\ref{F_2Es} and Section~\ref{section: proof of main lemmas}.

Technically, we will expand each integral into a sum of terms of the form
\[ \textstyle
\int_{\Gamma}
\frac{dz}
{(z-\lambda_i)(z-\lambda_k)(z-\lambda_j)}
\,\;
u_i u_i^T E u_k u_k^T E u_j u_j^T.
\]
Each contour integral
\(
\int_{\Gamma}
\frac{1}
{(z-\lambda_i)(z-\lambda_k)(z-\lambda_j)}
\,dz
\)
is evaluated exactly using Lemma~\ref{contour1-HR} for \(s=2\) and
\(T\in\{1,2\}\), where \(T\) denotes the number of indices among
\(i,j,k\) that belong to \(S\), while \(s\), as throughout the paper,
denotes the number of factors \(E\). Except for $\int_{\Gamma} QEQEQ dz$ and $\int_{\Gamma}PEPEP dz$, we then group the resulting terms according to which of the indices \(i,j,k\) belong to \(N_{\bar\lambda}(S)^c\).

For the first term $\int_{\Gamma} Q E Q E Q dz $, since $\Gamma$ does not enclose any point from the set $\{ \lambda_j\}_{j \in N_{\bar{\lambda}}(S)^c}$,  $\int_{\Gamma} \frac{1}{(z-\lambda_{j_1})(z-\lambda_{j_2}) (z- \lambda_{j_3})} dz =0 $ for any $ j_1,j_2, j_3 \in N_{\bar{\lambda}}(S)^c$. Therefore,
\begin{equation} \label{estimate1-HRF2} \int_{\Gamma} QEQEQ  dz =  \sum_{j_1,j_2, j_3 \in N_{\bar{\lambda}}(S)^c} \int_{\Gamma}  \frac{u_{j_1} u_{j_1}^T E u_{j_2} u_{j_2}^T E u_{j_3} u_{j_3}^T}{(z-\lambda_{j_1})(z-\lambda_{j_2})(z-\lambda_{j_3})} dz = 0. 
\end{equation}

Next, we have 
\begin{equation*}
\begin{split}
\frac{1}{2 \pi \textbf{i}} \int_{\Gamma}QE P E Q dz & = \frac{1}{2 \pi \textbf{i}} \sum_{ \substack{i \in N_{\bar{\lambda}}(S)\\ j,j' \in N_{\bar{\lambda}}(S)^c}} \int_{\Gamma} \frac{1}{(z-\lambda_{j'})(z- \lambda_i)(z - \lambda_j)} u_{j'} u_{j'}^T E u_i u_i^T E u_j u_j^T dz \\
& = \sum_{\substack{i \in S \\ j,j' \in N_{\bar{\lambda}}(S)^c}} \frac{1}{(\lambda_i -\lambda_{j'})(\lambda_i - \lambda_j)} u_{j'} u_{j'}^T E u_i u_i^T E u_j u_j^T (\text{by Lemma \ref{contour1-HR} for $T=1, s=2$}) \\
& = \sum_{i \in S} \bigg( \sum_{j' \in N_{\bar{\lambda}}(S)^c} \frac{u_{j'} u_{j'}^T}{\lambda_i -\lambda_{j'}} \bigg) E u_i u_i^T E \bigg( \sum_{j \in N_{\bar{\lambda}}(S)^c} \frac{u_j u_j^T}{\lambda_i -\lambda_j} \bigg).
\end{split}
\end{equation*}
Therefore, by the triangle inequality,
\begin{equation*}
    \begin{split}
         \frac{1}{2 \pi} \Norm{\int_{\Gamma}QE P E Q dz} & \leq  \sum_{i\in S}  \bigg\| \sum_{j' \in N_{\bar{\lambda}}(S)^c} \frac{u_{j'} u_{j'}^T}{\lambda_i -\lambda_{j'}}\bigg\| \Norm{E u_i} \Norm{ u_i^T E}  \bigg\| \sum_{j \in N_{\bar{\lambda}}(S)^c} \frac{u_j u_j^T}{\lambda_i -\lambda_j}\bigg\| \\
        & \leq \sum_{i \in S} \frac{1}{\min_{j' \in N_{\bar{\lambda}}(S)^c}|\lambda_i -\lambda_{j'}|} \cdot \|E\| \cdot \|E\| \cdot \frac{1}{\min_{j \in N_{\bar{\lambda}}(S)^c}|\lambda_i -\lambda_{j}|}.
    \end{split}
\end{equation*}
By the definition of $\bar \lambda$ and by $ \frac{\sqrt{p} \|E\|}{\bar{\lambda}} \leq \frac{1}{12}$ (Assumption \textbf{D0}), we further have
\begin{equation} \label{QEPEQ}
 \frac{1}{2 \pi} \Norm{\int_{\Gamma}QE P E Q dz} \leq \sum_{i \in S} \frac{1}{\bar{\lambda}} \cdot \|E\| \cdot \|E\| \cdot \frac{1}{\bar \lambda}= \frac{p \|E\|^2}{\bar{\lambda}^2} \leq \frac{1}{12} \frac{\sqrt{p} \|E\|}{\bar{\lambda}}.  
    \end{equation}
The same estimates hold for $\frac{1}{2 \pi} \Norm{\int_{\Gamma}PE Q E Q dz}, \frac{1}{2 \pi} \Norm{\int_{\Gamma}QE Q E P dz}$. 

Now, we consider the term $\int_{\Gamma} PE  Q EP dz$, which can be expressed as
$$\frac{1}{2 \pi \textbf{i}} \int_{\Gamma}PE Q E P dz = \frac{1}{2 \pi \textbf{i}} \sum_{ \substack{i,i' \in N_{\bar{\lambda}}(S)\\ j \in N_{\bar{\lambda}}(S)^c}} \int_{\Gamma} \frac{1}{(z- \lambda_{i'})(z - \lambda_j)(z-\lambda_{i})} u_{i'} u_{i'}^T E u_j u_j^T E u_i u_i^T dz.$$
To compute the contour integral, we need to identify the $\lambda_i$ inside (outside) $\Gamma$. Thus, we split the set $N_{\bar{\lambda}}(S)$ into $S$ and $N_{\bar{\lambda}}(S) \setminus S$ and further decompose $\frac{1}{2 \pi \textbf{i}} \int_{\Gamma}PE Q E P dz$ into 4 sub-sums:
\begin{equation*}
    \begin{split}
        & \frac{1}{2 \pi \textbf{i}} \sum_{ \substack{i,i' \in N_{\bar{\lambda}}(S) \setminus S\\ j \in N_{\bar{\lambda}}(S)^c}} \int_{\Gamma} \frac{1}{(z- \lambda_{i'})(z - \lambda_j)(z-\lambda_{i})} u_{i'} u_{i'}^T E u_j u_j^T E u_i u_i^T dz \\
        & + \frac{1}{2 \pi \textbf{i}} \sum_{ \substack{i,i' \in S\\ j \in N_{\bar{\lambda}}(S)^c}} \int_{\Gamma} \frac{1}{(z- \lambda_{i'})(z - \lambda_j)(z-\lambda_{i})} u_{i'} u_{i'}^T E u_j u_j^T E u_i u_i^T dz 
    \end{split}
\end{equation*}
\begin{equation*}
    \begin{split}
& +\frac{1}{2 \pi \textbf{i}} \sum_{ \substack{i \in S \\ i' \in N_{\bar{\lambda}}(S) \setminus S\\ j \in N_{\bar{\lambda}}(S)^c}} \int_{\Gamma} \frac{1}{(z- \lambda_{i'})(z - \lambda_j)(z-\lambda_{i})} u_{i'} u_{i'}^T E u_j u_j^T E u_i u_i^T dz \\
&+ \frac{1}{2 \pi \textbf{i}} \sum_{ \substack{i \in N_{\bar{\lambda}}(S) \setminus S \\ i' \in S \\ j \in N_{\bar{\lambda}}(S)^c}} \int_{\Gamma} \frac{1}{(z- \lambda_{i'})(z - \lambda_j)(z-\lambda_{i})} u_{i'} u_{i'}^T E u_j u_j^T E u_i u_i^T dz.
    \end{split}
\end{equation*}
By Lemma~\ref{contour1-HR}, we can evaluate exactly the contour integral
\(
\int_{\Gamma}
\frac{dz}{(z-\lambda_{i'})(z-\lambda_j)(z-\lambda_i)}
\)
for each group of indices \((i,i',j)\). For
\(i,i'\in N_{\bar\lambda}(S)\setminus S\) and
\(j\in N_{\bar\lambda}(S)^c\), this integral vanishes, and hence the first
sub-sum is zero. Applying Lemma~\ref{contour1-HR} to the remaining three
sub-sums, with \(T=2\) for the second sub-sum and \(T=1\) for the third
and fourth sub-sums, we can rewrite them as
\begin{equation*} \label{splitPEQEP}
    \begin{split}
& \sum_{ \substack{i,i' \in S\\ j \in N_{\bar{\lambda}}(S)^c}} \frac{-1}{(\lambda_i -\lambda_j)(\lambda_{i'} -\lambda_j)} u_{i'} u_{i'}^T E u_j u_j^T E u_i u_i^T + \sum_{ \substack{i \in S \\ i' \in N_{\bar{\lambda}}(S) \setminus S\\ j \in N_{\bar{\lambda}}(S)^c}} \frac{1}{(\lambda_i -\lambda_{i'})(\lambda_i -\lambda_j)} u_{i'} u_{i'}^T E u_j u_j^T E u_i u_i^T  \\
&+  \sum_{ \substack{i \in N_{\bar{\lambda}}(S) \setminus S \\ i' \in S \\ j \in N_{\bar{\lambda}}(S)^c}} \frac{1}{(\lambda_{i'} -\lambda_i)(\lambda_{i'} -\lambda_j)} u_{i'} u_{i'}^T E u_j u_j^T E u_i u_i^T.         \end{split}
\end{equation*}

\noindent We bound the norm of the first remaining sub-sum as follows. First, by the definition of the spectral norm, 
\begin{equation*}
    \begin{split}
&\bigg\| \sum_{ \substack{i,i' \in S\\ j \in N_{\bar{\lambda}}(S)^c}} \frac{-1}{(\lambda_i -\lambda_j)(\lambda_{i'} -\lambda_j)} u_{i'} u_{i'}^T E u_j u_j^T E u_i u_i^T \bigg\| = \max_{\|\textbf{v}\|=\|\textbf{w}\|=1}  \sum_{ \substack{i,i' \in S\\ j \in N_{\bar{\lambda}}(S)^c}} \frac{  \textbf{v}^T u_{i'} u_{i'}^T E u_j u_j^T E u_i u_i^T \textbf{w}}{(\lambda_i -\lambda_j)(\lambda_{i'} -\lambda_j)}  \\
& =\max_{\|\textbf{v}\|=\|\textbf{w}\|=1}  \sum_{ \substack{i,i' \in S}}  \textbf{v}^T u_{i'} u_{i'}^T E \bigg(\sum_{ j \in N_{\bar{\lambda}}(S)^c} \frac{  u_j u_j^T }{(\lambda_i -\lambda_j)(\lambda_{i'} -\lambda_j)} \bigg) E u_i u_i^T \textbf{w}. 
\end{split}
\end{equation*}
By the triangle inequality, the RHS is at most 
\begin{equation*}
    \begin{split}
& \max_{\|\textbf{v}\|=\|\textbf{w}\|=1}  \sum_{ \substack{i,i' \in S}}  \norm{\textbf{v}^T u_{i'}} \bigg| u_{i'}^T E \bigg(\sum_{ j \in N_{\bar{\lambda}}(S)^c} \frac{  u_j u_j^T }{(\lambda_i -\lambda_j)(\lambda_{i'} -\lambda_j)} \bigg) E u_i \bigg| \norm{ u_i^T \textbf{w}} \\
& \leq \frac{p \|E\|^2}{\bar{\lambda}^2} (\text{by the definition of $\bar{\lambda}$ and the Cauchy--Schwarz inequality}).
  \end{split}
\end{equation*}

Next, we bound the norm of the second sub-sum. Similarly, we write
\begin{equation*}
    \begin{split}
& \bigg\|\sum_{ \substack{i \in S \\ i' \in N_{\bar{\lambda}}(S) \setminus S\\ j \in N_{\bar{\lambda}}(S)^c}} \frac{1}{(\lambda_i -\lambda_{i'})(\lambda_i -\lambda_j)} u_{i'} u_{i'}^T E u_j u_j^T E u_i u_i^T\bigg\| = \max_{\|\textbf{v}\|=\|\textbf{w}\|=1}  \sum_{ \substack{i \in S \\ i' \in N_{\bar{\lambda}}(S) \setminus S\\ j \in N_{\bar{\lambda}}(S)^c}} \frac{  \textbf{v}^T u_{i'} u_{i'}^T E u_j u_j^T E u_i u_i^T \textbf{w}}{(\lambda_i -\lambda_{i'})(\lambda_{i} -\lambda_j)}  \\
& =\max_{\|\textbf{v}\|=\|\textbf{w}\|=1}  \sum_{ \substack{i \in S \\ i' \in N_{\bar{\lambda}}(S) \setminus S} } \frac{\textbf{v}^T u_{i'} u_{i'}^T E}{\lambda_i -\lambda_{i'}} \bigg(\sum_{ j \in N_{\bar{\lambda}}(S)^c} \frac{  u_j u_j^T }{\lambda_i -\lambda_j} \bigg) E u_i u_i^T \textbf{w}, 
\end{split}
\end{equation*}
which, by the triangle inequality, is at most
\[
\max_{\|\textbf{v}\|=\|\textbf{w}\|=1}  \sum_{ \substack{i \in S \\ i' \in N_{\bar{\lambda}}(S) \setminus S} } \frac{|\textbf{v}^T u_{i'}|}{|\lambda_i -\lambda_{i'}|} \bigg| u_{i'}^T E \bigg(\sum_{ j \in N_{\bar{\lambda}}(S)^c} \frac{  u_j u_j^T }{\lambda_i -\lambda_j} \bigg) E u_i \bigg| \norm{ u_i^T \textbf{w}} 
\]
Since $i' \neq i$ and the definition of $y$ in \eqref{xyz}, we have $|u_{i'}^T E \big(\sum_{ j \in N_{\bar{\lambda}}(S)^c} \frac{  u_j u_j^T }{\lambda_i -\lambda_j} \big) E u_i| \leq y$. Together with $|\lambda_i -\lambda_{i'}| \geq \delta_S$, these estimates show that the norm of the second sub-sum is at most 
\[
\frac{y}{\delta_S} \cdot \max_{\|\textbf{v}\|=\|\textbf{w}\|=1}  \sum_{ \substack{i \in S \\ i' \in N_{\bar{\lambda}}(S) \setminus S} }\norm{\textbf{v}^T u_{i'}} \cdot \norm{ u_i^T \textbf{w}} \leq \frac{\sqrt{pr} y}{\delta_S} (\text{by the Cauchy--Schwarz inequality}). 
\]
A similar estimate holds for the norm of the third sub-sum. Putting the bounds on the norms of these sub-sums together, we obtain 
\begin{equation} \label{PEQEP}
    \frac{1}{2 \pi} \Norm{\int_{\Gamma} PEQEP dz} \leq \frac{p \|E\|^2}{\bar{\lambda}^2} + \frac{2 \sqrt{pr} y}{\delta_S}.
\end{equation}

Now we consider the fourth term $\int_{\Gamma} PEPEP dz$. By expansion, we have 
\begin{equation*}
 \int_{\Gamma} PEPEP dz = \sum_{i_1,i_2,i_3 \in N_{\bar{\lambda}}(S)} \int_{\Gamma} \frac{u_{i_1} u_{i_1}^T E u_{i_2} u_{i_2}^T E u_{i_3} u_{i_3}^T}{(z- \lambda_{i_1})(z-\lambda_{i_2})(z-\lambda_{i_3})} dz.  
\end{equation*}
Notice that the integral is non-zero when there is either one or two of  $\{i_1,i_2,i_3\}$ inside the contour $\Gamma$. Thus, there are  6 such scenarios, and we split $\int_{\Gamma} PEPEP dz$ into 6 corresponding subsums. We will bound the norm of each sub-sum by $\frac{\sqrt{pr} r x^2}{\delta_S^2}.$ We  present here the treatment for the case $i_1 \in S$ and $i_2, i_3 \in N_{\bar{\lambda}} (S) \setminus S$; the treatment of other cases are similar. We write 
\begin{equation*}
\begin{split}
  & \bigg\| \sum_{\substack{i_1 \in S\\i_2, i_3 \in N_{\bar{\lambda}} (S) \setminus S}} \int_{\Gamma} \frac{u_{i_1} u_{i_1}^T E u_{i_2} u_{i_2}^T E u_{i_3} u_{i_3}^T}{(z- \lambda_{i_1})(z-\lambda_{i_2})(z-\lambda_{i_3})} dz  \bigg\|  = \bigg\| \sum_{\substack{i_1 \in S\\i_2, i_3 \in N_{\bar{\lambda}} (S) \setminus S}} \frac{u_{i_1} u_{i_1}^T E u_{i_2} u_{i_2}^T E u_{i_3} u_{i_3}^T}{(\lambda_{i_1}-\lambda_{i_2})(\lambda_{i_1}-\lambda_{i_3})}  \bigg\| \\
  & = \max_{\|\textbf{v}\|=\|\textbf{w}\|=1} \sum_{\substack{i_1 \in S\\i_2, i_3 \in N_{\bar{\lambda}} (S) \setminus S}} \frac{ \textbf{v}^T u_{i_1} u_{i_1}^T E u_{i_2} u_{i_2}^T E u_{i_3} u_{i_3}^T \textbf{w}}{(\lambda_{i_1}-\lambda_{i_2})(\lambda_{i_1}-\lambda_{i_3})}. 
  \end{split}
\end{equation*}
By the triangle inequality, this is at most
\begin{equation} \label{eqF2PEPEP}
  \max_{\|\textbf{v}\|=\|\textbf{w}\|=1} \sum_{\substack{i_1 \in S\\i_2, i_3 \in N_{\bar{\lambda}} (S) \setminus S}} \frac{ \norm{\textbf{v}^T u_{i_1}} \norm{ u_{i_1}^T E u_{i_2}} \norm{ u_{i_2}^T E u_{i_3}} \norm{ u_{i_3}^T \textbf{w}}}{(\lambda_{i_1}-\lambda_{i_2})(\lambda_{i_1}-\lambda_{i_3})}.  
\end{equation}
By the definition of $x$ in \eqref{xyz}, $\norm{ u_{i_1}^T E u_{i_2}} \norm{ u_{i_2}^T E u_{i_3}} \leq x^2$. Thus, \eqref{eqF2PEPEP} is less than or equal to
\begin{equation*}
    \begin{split}
  \frac{r x^2}{\delta_S^2} \max_{\|\textbf{v}\|=\|\textbf{w}\|=1} \sum_{\substack{i_1 \in S\\ i_3 \in N_{\bar{\lambda}} (S) \setminus S}} \norm{\textbf{v}^T u_{i_1}} \norm{ u_{i_3}^T \textbf{w}}  \leq \frac{\sqrt{pr} r x^2}{\delta_S^2} (\text{by the Cauchy--Schwarz inequality}).      
    \end{split}
\end{equation*}
Therefore, 
\begin{equation} \label{PEPEP}
\frac{1}{2 \pi} \Norm{\int_{\Gamma} PEPEP dz} \leq \frac{6 \sqrt{pr} r x^2}{\delta_S^2}.
\end{equation}

Finally, we bound $\int_\Gamma PEPEQ dz$. Again, we expand the product and obtain 
$$\frac{1}{2 \pi \textbf{i}} \int_{\Gamma}PE P E Q dz = \frac{1}{2 \pi \textbf{i}} \sum_{ \substack{i,i' \in N_{\bar{\lambda}}(S)\\ j \in N_{\bar{\lambda}}(S)^c}} \int_{\Gamma} \frac{u_{i'} u_{i'}^T E u_i u_i^T E u_j u_j^T dz}{(z-\lambda_{i'})(z- \lambda_i)(z - \lambda_j)} ,$$
which can be split into
\begin{equation*} 
    \begin{split}
 & \frac{1}{2 \pi \textbf{i}} \sum_{ \substack{i,i' \in S\\ j \in N_{\bar{\lambda}}(S)^c}} \int_{\Gamma} \frac{ u_{i'} u_{i'}^T E u_j u_j^T E u_i u_i^T dz}{(z-\lambda_{i'})(z- \lambda_i)(z - \lambda_j)}  +\frac{1}{2 \pi \textbf{i}} \sum_{ \substack{i \in S \\ i' \in N_{\bar{\lambda}}(S) \setminus S\\ j \in N_{\bar{\lambda}}(S)^c}} \int_{\Gamma} \frac{u_{i'} u_{i'}^T E u_i u_i^T E u_j u_j^T dz }{(z-\lambda_{i'})(z- \lambda_i)(z - \lambda_j)} \\
&+ \frac{1}{2 \pi \textbf{i}} \sum_{ \substack{i \in N_{\bar{\lambda}}(S) \setminus S \\ i' \in S \\ j \in N_{\bar{\lambda}}(S)^c}} \int_{\Gamma} \frac{u_{i'} u_{i'}^T E u_i u_i^T E u_j u_j^T dz}{(z-\lambda_{i'})(z- \lambda_i)(z - \lambda_j)} .
    \end{split}
\end{equation*}

 \noindent Applying Lemma \ref{contour1-HR} to compute the integral, we further decompose $ \frac{1}{2 \pi \textbf{i}} \int_{\Gamma}PE P E Q dz  $ as
\begin{equation}\label{splitPEPEQ}
    \begin{split}    
& \sum_{ \substack{i,i' \in S\\ j \in N_{\bar{\lambda}}(S)^c}} \frac{-1}{(\lambda_i -\lambda_j)(\lambda_{i'} -\lambda_j)} u_{i'} u_{i'}^T E u_i u_i^T E u_j u_j^T  + \sum_{ \substack{i \in S \\ i' \in N_{\bar{\lambda}}(S) \setminus S\\ j \in N_{\bar{\lambda}}(S)^c}} \frac{1}{(\lambda_i -\lambda_{i'})(\lambda_i -\lambda_j)} u_{i'} u_{i'}^T E u_i u_i^T E u_j u_j^T  \\
&+  \sum_{ \substack{i \in N_{\bar{\lambda}}(S) \setminus S \\ i' \in S \\ j \in N_{\bar{\lambda}}(S)^c}} \frac{1}{(\lambda_{i'} -\lambda_i)(\lambda_{i'} -\lambda_j)} u_{i'} u_{i'}^T E u_i u_i^T E u_j u_j^T . 
    \end{split}
\end{equation}
\noindent We bound the norm of the first sum as follows. First, by the definition of the spectral norm, we write 
\begin{equation*}
    \begin{split}
& \bigg\|\sum_{ \substack{i,i' \in S\\ j \in N_{\bar{\lambda}}(S)^c}} \frac{-1}{(\lambda_i -\lambda_j)(\lambda_{i'} -\lambda_j)} u_{i'} u_{i'}^T E u_i u_i^T E u_j u_j^T \bigg\| = \max_{\|\textbf{v}\|=\|\textbf{w}\|=1} \sum_{ \substack{i,i' \in S\\ j \in N_{\bar{\lambda}}(S)^c}} \frac{\textbf{v}^T u_{i'} u_{i'}^T E u_i u_i^T E u_j u_j^T \textbf{w}}{(\lambda_i -\lambda_j)(\lambda_{i'} -\lambda_j)}  \\
& = \max_{\|\textbf{v}\|=\|\textbf{w}\|=1} \sum_{ \substack{i,i' \in S}} \textbf{v}^T u_{i'} u_{i'}^T E u_i u_i^T E \bigg(\sum_{j \in N_{\bar{\lambda}}^c} \frac{u_j u_j^T }{(\lambda_i -\lambda_j)(\lambda_{i'} -\lambda_j)}  \bigg) \textbf{w}.
\end{split}
\end{equation*}
By the triangle inequality, the absolute value of the RHS is at most
\begin{equation*}
    \begin{split}
& \max_{\|\textbf{v}\|=\|\textbf{w}\|=1} \sum_{ \substack{i,i' \in S}} \norm{\textbf{v}^T u_{i'}} \norm{ u_{i'}^T E u_i} \Norm{ u_i^T E}  \bigg\|\sum_{j \in N_{\bar{\lambda}}^c} \frac{u_j u_j^T }{(\lambda_i -\lambda_j)(\lambda_{i'} -\lambda_j)}  \textbf{w} \bigg\|  \\
& \leq \max_{\|\textbf{v}\|=\|\textbf{w}\|=1} \sum_{ \substack{i,i' \in S}} \norm{\textbf{v}^T u_{i'}} \norm{ u_{i'}^T E u_i} \|E\| \cdot  \frac{1}{\min_{j \in N_{\bar{\lambda}}^c} |\lambda_i -\lambda_j|\cdot|\lambda_{i'}-\lambda_j|}  \\
& \leq \frac{p x \|E\|}{\bar{\lambda}^2} \max_{\|\textbf{v}\|=1} \sum_{ \substack{i' \in S}} \norm{\textbf{v}^T u_{i'}} (\text{by definitions of $x$ and $\bar{\lambda}$)} \\
& \leq \frac{px}{\bar{\lambda}} \times \frac{\sqrt{p} \|E\|}{\bar{\lambda}} (\text{by the Cauchy--Schwarz inequality}).
    \end{split}
\end{equation*}
By a similar argument, we also obtain 
$$\bigg\|\sum_{ \substack{i \in S \\ i' \in N_{\bar{\lambda}}(S) \setminus S\\ j \in N_{\bar{\lambda}}(S)^c}} \frac{1}{(\lambda_i -\lambda_{i'})(\lambda_i -\lambda_j)} u_{i'} u_{i'}^T E u_i u_i^T E u_j u_j^T \bigg\| \leq \frac{\sqrt{pr} x}{\delta_S} \times \frac{\sqrt{p} \|E\|}{\bar{\lambda}},$$
and 
$$\bigg\|\sum_{ \substack{i \in N_{\bar{\lambda}}(S) \setminus S \\ i' \in S \\ j \in N_{\bar{\lambda}}(S)^c}} \frac{1}{(\lambda_{i'} -\lambda_i)(\lambda_{i'} -\lambda_j)} u_{i'} u_{i'}^T E u_i u_i^T E u_j u_j^T \bigg\| \leq \frac{\sqrt{pr}x}{\delta_S} \times \frac{\sqrt{p}\|E\|}{\bar{\lambda}}. $$
Therefore, by Assumption \textbf{D0}, we have 
\begin{equation} \label{PEPEQ}
\frac{1}{2 \pi} \Norm{\int_{\Gamma}PE P E Q dz } \leq 3 \times \frac{1}{12} \times \frac{\sqrt{p} \|E\|}{\bar{\lambda}} = \frac{1}{4}   \frac{\sqrt{p} \|E\|}{\bar{\lambda}}.
\end{equation}
The same estimate holds for $\frac{1}{2 \pi} \Norm{\int_{\Gamma} QEPEP dz}$.

Combining \eqref{estimate1-HRF2}, \eqref{QEPEQ} (and the line following that), \eqref{PEQEP}, \eqref{PEPEP}, and \eqref{PEPEQ} (and the line following that), we obtain  a bound for $\| F_2\| $ 
\begin{equation*}
\begin{split}
\frac{1}{2 \pi} \|F_2\| & \leq \frac{3}{12} \frac{\sqrt{p}\|E\|}{\bar{\lambda}} + \frac{p\|E\|^2}{\bar{\lambda}^2} + \frac{2 \sqrt{pr} y}{\delta_S} + \frac{6 \sqrt{pr} r x^2}{\delta_S^2}+ \frac{2}{4} \frac{\sqrt{p}\|E\|}{\bar{\lambda}} \\
& \leq \left( \frac{3}{12}+ \frac{\sqrt{p} \|E\|}{\bar{\lambda}} + \frac{1}{2} \right) \frac{\sqrt{p} \|E\|}{\bar{\lambda}} + \frac{6r x}{\delta_S} \times \frac{\sqrt{pr} x}{\delta_S} + \frac{2 \sqrt{pr} y}{\delta_S}\\
& \leq \frac{\sqrt{p} \|E\|}{\bar{\lambda}} + \frac{\sqrt{pr} x}{2\delta_S} + \frac{2 \sqrt{pr} y}{\delta_S} \,\,(\text{since}\,\,  \frac{\sqrt{p} \|E\|}{\bar{\lambda}}, \frac{r x}{\delta_S}  \leq \frac{1}{12} \,\, \text{by Assumption \textbf{D0}})   \\
& \leq 2 \sqrt{p} \left( \frac{\|E\|}{\bar{\lambda}} + \frac{\sqrt{r}x}{\delta_S} + \frac{\sqrt{r} y}{\delta_S}   \right).
\end{split}
\end{equation*}
\subsection{Estimating $ \Norm{ F_s} $  for a general $s$} \label{F_2Es} 
In this subsection, we present three technical lemmas for estimating $\Norm{ F_s} $, with a general $s$. These lemmas are the core of our proof.  Combining these lemmas with \eqref{bound0}, we obtain a proof of Theorem \ref{deterministicHRgeneral}, which we present in the next subsection. The proofs of the  lemmas 
are involved and will be presented in the next section.

Following \eqref{F_sEFor-HR}, we  expand $(P+Q)[ E(P+Q)]^s$ into the sum of $2^{s+1}$ operators, each of which is a product of  alternating $Q$-blocks and 
$P$-blocks 
$$(QEQE \dots QE)(PE PE \dots PE) \dots (PE PE \dots PE) ( QE  QE \dots Q), $$ where we allow the first and last blocks to be empty. 

We code each operator like this by recording the numbers of $Q$'s and the numbers of $P$'s in the respective blocks, from left to right. If there are $(k+1)$  $Q$-blocks and $k$ $P$-blocks, for some integer $k$,  we let 
$\alpha_1, \dots, \alpha_{k+1} $ and $\beta_1, \dots, \beta_k$ be these numbers. These numbers satisfy the following conditions 
$$ \alpha_1, \alpha_{k+1} \geq 0,$$
$$ \alpha_i, \beta_j \geq 1,\,\, \forall 1<i<k+1,\,\, \forall  1\le j \le k, $$
$$\alpha_1+...+\alpha_{k+1}+\beta_1 +...+\beta_k=s+1.$$
In what follows, we set $\alpha =(\alpha_1,...,\alpha_{k+1}) \in \mathbb{Z}^{k+1} $ and $\beta=(\beta_1,...,\beta_k) \in \mathbb{Z}^{k}$, and use $M(\alpha; \beta)$ to denote the corresponding operator.

\noindent \textit{Example:} given $s=10, k=2, \alpha_1=3, \alpha_2=2, \alpha_3=1, \beta_1=2, \beta_2=3$, we have 
 $$\int_{\Gamma} M(3,2,1;2,3) dz = \int_{\Gamma} \underset{\substack{\alpha_1=3\\ \text{numbers of $Q$}}}{\underbrace{(QEQEQE)}} \underset{\substack{\beta_1=2\\ \text{numbers of $P$}}} {\underbrace{(PEPE)}}\,\,\, \underset{\substack{\alpha_2=2 \\ \text{numbers of } \, Q}}{\underbrace{(QEQE)}} \underset{\substack{\beta_2=3 \\ \text{numbers of}\, P}}{\underbrace{(PEPEPE)}} \,\,\, \underset{\substack{\alpha_3=1 \\ \text{numbers of}\, Q}}{\underbrace{Q}} dz.$$

We are going to bound $\Norm{\int_{\Gamma}M(\alpha; \beta) dz} $ for each pair $(\alpha; \beta)$ separately, and use the triangle inequality to add up the bounds. 
For $s$ and $k$ fixed, we split the collection of pairs $(\alpha;\beta)$  into the following types, according to the values of $\alpha_1$ and $\alpha_{k+1}$:
\begin{itemize}
    \item Type I: $\alpha_1, \alpha_{k+1} > 0$.
    \item Type II: $\alpha_1, \alpha_{k+1} = 0$. 
    \item Type III: either $(\alpha_1 = 0, \alpha_{k+1} > 0)$ or $(\alpha_1 > 0, \alpha_{k+1} = 0)$. 
\end{itemize}
Let $s_1 :=  \sum_{i=1}^{k+1} \alpha_i$ and  $s_2 := \sum_{j=1}^{k} \beta_j $. We prove three lemmas bounding  $\Norm{\int_{\Gamma}M(\alpha; \beta) dz}$ with respect to the above 3 types.
We assume that Assumption {\bf D0}
 holds in all three lemmas. 
\begin{lemma} \label{LemmaCase1M(a,b)bound} For a pair $(\alpha; \beta)$ of Type I,   
\begin{equation} \label{Case1boundM(a,b)bound}
    \frac{1}{2 \pi} \Norm{\int_{\Gamma} M(\alpha;\beta) dz} \leq \left( \frac{ \sqrt{p}\|E\|}{\bar{\lambda}} + \frac{px}{\bar{\lambda}} \right) \frac{2^{s_2+s-1}}{12^{s-1}}.
\end{equation}
\end{lemma}

\begin{lemma}\label{LemmaCase2M(a,b)bound}
  For a pair $(\alpha; \beta)$ of Type II, 
    \begin{itemize}
        \item If $M(\alpha; \beta) \neq PEQEPE ... QEP$, i.e.,  $(s_1,s_2) \neq (k-1,k)$, then 
        \begin{equation} \label{Cas2M(a,b)boundgood}
     \frac{1}{2 \pi} \Norm{ \int_{\Gamma} M(\alpha; \beta) dz}  \le  \left(\frac{\sqrt{p} \|E\|}{\bar{\lambda}} + \frac{\sqrt{pr} x}{\delta_S} \right) \frac{2^{s_2+s-1}}{12^{s-1}}, 
 \end{equation}
 \item If $M(\alpha;\beta) = PEQEPEQE...EQEP $, i.e., $(s_1,s_2) = (k-1,k)$, then 
 \begin{equation} \label{Case2M(a,b)boundbad}
    \frac{1}{2 \pi} \Norm{ \int_{\Gamma} PEQEPEQE...EQEP dz}  \le \left( \frac{px}{\bar{\lambda}} + \frac{\sqrt{p} \|E\|}{\bar{\lambda}} \right) \frac{2^{3s/2} }{12^{s-1}}  +   \frac{(s+2) \sqrt{pr} y}{\delta_S} \left( \frac{2 \sqrt{r} w}{\sqrt{\bar{\lambda} \delta_S}} \right)^{s-2}. 
 \end{equation}
    \end{itemize}
     \end{lemma}

\begin{lemma} \label{LemmaCase3M(a,b)bound}
    For a pair $(\alpha; \beta)$ of Type III, 
    \begin{equation} \label{Case3M(a,b)bound}
\begin{split}
 \frac{1}{2 \pi} \Norm{\int_{\Gamma} M(\alpha;\beta) dz} &\leq \left( \frac { \sqrt{p}\|E\|}{\bar{\lambda}}  + \frac{p x}{\bar{\lambda}} \right) \frac{2^{s_2+s-1}}{12^{s-1}}.
\end{split}
\end{equation} 
\end{lemma}

\subsection {Proof of Theorem \ref{deterministicHRgeneral} using the lemmas.} \label{subsec: provevialemma} 
 For any integer $ 1 \le s_2 \leq s$,  the expansion of $F_s$ contains  $\binom{s+1}{s_2}$ operators $M(\alpha;\beta)$, each of which has exactly $s_2$ $P$ operators. A simple  consideration reveals that among these,  $\binom{s-1}{s_2}$ are of Type I,  $\binom{s-1}{s_2-2}$ are of Type II, and $2\binom{s-1}{s_2 -1}$ are of Type III. As shown in Section \ref{section:method}, the only term in the expansion with $s_2=0$ vanishes by the Cauchy integral formula.

 Putting \eqref{Case1boundM(a,b)bound}, \eqref{Cas2M(a,b)boundgood}, \eqref{Case2M(a,b)boundbad}, \eqref{Case3M(a,b)bound} together and summing over $s_2$,  we show that  $ \frac{1}{2\pi } \| F_s \|$ is at most 
 \begin{equation*}
 \begin{split}
  & \left[\sum_{s_2=1}^{s+1} \binom{s-1}{s_2}   \left( \frac{ \sqrt{p}\|E\|}{\bar{\lambda}} + \frac{px}{\bar{\lambda}} \right) \frac{2^{s_2+s-1}}{12^{s-1}}\right] +   \bigg[ \left( \frac{px}{\bar{\lambda}} + \frac{\sqrt{p} \|E\|}{\bar{\lambda}} \right) \frac{2^{3s/2} }{12^{s-1}}+  \frac{(s+2) \sqrt{pr} y}{\delta_S} \bigg( \frac{2 \sqrt{r} w}{\sqrt{\bar{\lambda} \delta_S}} \bigg)^{s-2} \bigg]+ \\
    & \left[ \sum_{s_2=1}^{s+1} \binom{s-1}{s_2-2}  \left(\frac{\sqrt{p} \|E\|}{\bar{\lambda}}  + \frac{\sqrt{pr} x}{\delta_S} \right) \frac{2^{s_2+s-1}}{12^{s-1}} \right] + \left[ \sum_{s_2=1}^{s+1} 2\binom{s-1}{s_2 -1} \left( \frac { \sqrt{p}\|E\|}{\bar{\lambda}}  + \frac{p x}{\bar{\lambda}} \right) \times \frac{2^{s_2+s-1}}{12^{s-1}}  \right],
  \end{split}     
 \end{equation*}  
 here $\binom{s-1}{l} = 0$ for either $l > s-1$ or $l < 0$.
 Consider the first term 
 \begin{equation*}
     \begin{split}
         \sum_{s_2=1}^{s+1} \binom{s-1}{s_2}   \left( \frac{ \sqrt{p}\|E\|}{\bar{\lambda}} + \frac{px}{\bar{\lambda}} \right) \frac{2^{s_2+s-1}}{12^{s-1}} & = \sum_{s_2=1}^{s+1} \binom{s-1}{s_2} 2^{s_2}  \left( \frac{ \sqrt{p}\|E\|}{\bar{\lambda}} + \frac{px}{\bar{\lambda}} \right) \frac{2^{s-1}}{12^{s-1}} \\
         & = \left( \sum_{s_2=1}^{s+1} \binom{s-1}{s_2} 2^{s_2} 1^{s-1-s_2}  \right) \times \left( \frac{ \sqrt{p}\|E\|}{\bar{\lambda}} + \frac{px}{\bar{\lambda}} \right) \frac{2^{s-1}}{12^{s-1}} \\
         & = 3^{s-1}  \times \left( \frac{ \sqrt{p}\|E\|}{\bar{\lambda}} + \frac{px}{\bar{\lambda}} \right) \frac{2^{s-1}}{12^{s-1}} \\
         & = \left( \frac{ \sqrt{p}\|E\|}{\bar{\lambda}} + \frac{px}{\bar{\lambda}} \right) \frac{1}{2^{s-1}}. 
     \end{split}
 \end{equation*}
Using a similar argument, we have 
$$\sum_{s_2=1}^{s+1} \binom{s-1}{s_2-2}  \left(\frac{\sqrt{p} \|E\|}{\bar{\lambda}} + \frac{\sqrt{pr} x}{\delta_S} \right) \frac{2^{s_2+s-1}}{12^{s-1}} = 4  \left(\frac{\sqrt{p} \|E\|}{\bar{\lambda}} + \frac{\sqrt{pr} x}{\delta_S} \right) \frac{1}{2^{s-1}},   $$
(here we pair up $2^{s_2-2}$ with $\binom{s-1}{s_2-2}$, the remaining factor is $4$),
and
$$\sum_{s_2=1}^{s+1} 2\binom{s-1}{s_2 -1} \left( \frac { \sqrt{p}\|E\|}{\bar{\lambda}}  + \frac{p x}{\bar{\lambda}} \right) \times \frac{2^{s_2+s-1}}{12^{s-1}} =  4 \left( \frac { \sqrt{p}\|E\|}{\bar{\lambda}}  + \frac{p x}{\bar{\lambda}} \right) \times \frac{1}{2^{s-1}},$$
(we pair up $2^{s_2-1}$ with $\binom{s-1}{s_2-1}$, the remaining factor is $2 \times 2 =4$). Finally, we have 
$$ \left[ \left( \frac{px}{\bar{\lambda}} + \frac{\sqrt{p} \|E\|}{\bar{\lambda}} \right) \frac{2^{3s/2} }{12^{s-1}}+  \frac{(s+2) \sqrt{pr} y}{\delta_S} \left( \frac{2 \sqrt{r} w}{\sqrt{\bar{\lambda} \delta_S}} \right)^{s-2} \right] \le \left[ \frac{2}{3} \left( \frac{px}{\bar{\lambda}} + \frac{\sqrt{p} \|E\|}{\bar{\lambda}} \right) \frac{1 }{4 ^{s-2}} + \frac{ \sqrt{pr} y}{\delta_S} \frac{s+2}{6^{s-2}}  \right],$$
where we use Assumption {\bf D0} to bound $\frac{\sqrt{r} w}{\sqrt{\bar{\lambda} \delta_S}} $ and the trivial estimate that $\frac{2^{3/2}}{12} \leq \frac{3}{12}=\frac{1}{4}$. It follows that 
\begin{equation} \label{boundFs} 
    \begin{split}
        \frac{1}{2 \pi } \Norm{F_s}  & \leq  \left( \frac{ \sqrt{p}\|E\|}{\bar{\lambda}} + \frac{px}{\bar{\lambda}} \right) \frac{1}{2^{s-1}} +  \left[ \frac{2}{3} \left( \frac{px}{\bar{\lambda}} + \frac{\sqrt{p} \|E\|}{\bar{\lambda}} \right) \frac{1 }{4 ^{s-2}} + \frac{ \sqrt{pr} y}{\delta_S} \frac{s+2}{6^{s-2}}  \right]  +  \\
    & + 4  \left(\frac{\sqrt{p} \|E\|}{\bar{\lambda}} + \frac{\sqrt{pr} x}{\delta_S} \right) \frac{1}{2^{s-1}} + 4 \left( \frac { \sqrt{p}\|E\|}{\bar{\lambda}}  + \frac{p x}{\bar{\lambda}} \right) \frac{1}{2^{s-1}}.
    \end{split}
\end{equation}

The RHS converges. By summing (over $s$), we can bound  $\frac{1}{2 \pi} \sum_{s\geq 2} \|F_s\|$ by 
 \begin{equation}  \label{F2-HR} 10 \frac{ \sqrt{p}\|E\|}{\bar{\lambda}} + 4 \frac{\sqrt{pr} x}{\delta_S}    + 6 \frac{px}{\bar{\lambda}}    + \frac{26}{5} \frac{\sqrt{pr} y}{\delta_S}.    
 \end{equation} 
(Here we omit the routine calculation and do not try to optimize the constants.) 
Combining \eqref{bound0}, \eqref{F1-HR}, \eqref{F2-HR}, we obtain 
 \begin{equation} \label{F3-HR} 
 \|\tilde{\Pi}_{\tilde{S}} - \Pi_S   \|\le \frac{ 2 \sqrt{p}  \| E\|}{ \bar{\lambda}}   + \frac{   2\sqrt{p r} x} { \delta_S } +  \left(10 \frac{ \sqrt{p}\|E\|}{\bar{\lambda}} +  6\frac{px}{\bar{\lambda}} + 4\frac{\sqrt{pr} x}{\delta_S}      + \frac{26}{5} \frac{\sqrt{pr} y}{\delta_S} \right) . 
 \end{equation} 
  Since  $\frac{p  x}{\bar{\lambda}} \leq \frac{ \sqrt{pr} x}{\delta_S}$, \eqref{F3-HR} implies 
  \begin{equation*} \label{F4-HR} 
\|\tilde{\Pi}_{\tilde{S}} - \Pi_S   \|
 \le \frac{12  \sqrt{p} \| E \| }{ \bar{\lambda}} + \frac{ 12 \sqrt{pr} x } {\delta_S} + 6 \frac{ \sqrt{p r} y }{ \delta_S } \leq 12 \sqrt{p} \left( \frac{\|E\|}{\bar{\lambda}} + \frac{\sqrt{r} x}{\delta_S} + \frac{\sqrt{r} y}{\delta_S}  \right),
 \end{equation*} completing the proof.

\section{Proofs of the main lemmas}  \label{section: proof of main lemmas}

In this section, we begin with a combinatorial lemma (Lemma~\ref{contour1-HR}), which plays a central role in our analysis. 
We then prove the two main lemmas, Lemma~\ref{LemmaCase1M(a,b)bound} and Lemma~\ref{LemmaCase2M(a,b)bound}, in the next two subsections. 
The proof of the third main lemma, Lemma~\ref{LemmaCase3M(a,b)bound}, follows a similar strategy and is presented separately in Appendix~\ref{Appx: proofLemmaCase3}. 
Finally, we conclude this section with the proof of Expansion~\eqref{TaylorEx} in Subsection~\ref{TaylorExpansion}.

\subsection{A combinatorial lemma} \label{subsec: combi lemma}
An important  part of our analysis is to 
compute the contour integral 
$$ \frac{1}{2\pi \textbf{i}}\int_{\Gamma} \prod_{j=1}^l \frac{1}{ z- \lambda_{i_j}}  , $$ for a sequence of indices 
$i_1, \dots, i_l$. There are many ways to do this, which could lead to different-looking (but equivalent) formulae.

The result in this section provides a combinatorial way to compute the integral. This will be essential in regrouping the terms after the expansion. 

We start with a few definitions. 
Fix a natural number $s$. Given a parameter $0 \le T \le s+1$, consider two sequences of (not necessarily different)  indices: $ X= [i_1, \dots, i_T]$ and $Y = [i_{T+1}, ..., i_{s+1}]$.

Let $\mathcal{L}(T,s)$ denote the set of all sequences of natural numbers $m_1,m_2,...,m_T$ such that $m_1+...+m_T=s$.  

For each sequence $L=[m_1,...,m_T] \in \mathcal{L}(T,s)$, we define a directed bipartite graph $G(X,Y|L)$,  with $X$ and $Y$ being the vertex sets,  as follows:

\begin{itemize}
\item The edges go from $X$ to $Y$: $\forall e \in \mathcal{E}(X,Y|L)$, $e = e^{+} \rightarrow e^{-}:  e^{+} \in X$ and $ e^{-} \in Y.$

\item For each $1 \leq j \leq T$, $i_j$ in $X$ is connected to $i_{j'}$ in $Y$, for all   $j'$ in the interval 
$$[  (s+1) - ((m_1-1)+...+(m_{j}-1) ),  (s+1) - ((m_1-1)+...+(m_{j-1}-1))].$$
\end{itemize}

In particular, $i_1$ is connected to all $i_{j'}$ with $j'$ in  $[ (s+1)- (m_1-1), s+1]$. Thus, if $m_1=1$, then $i_1$ is connected to 
only $i_{s+1}$.  It is clear that given $T $ and $L$, the graph is defined up to isomorphism. Denote the set of (directed) edges of the graph $G(X,Y|L)$ by $\mathcal{E}(X,Y| L)$.

\vskip2mm 
\noindent \textit{Example:} Given two sequences of indices: $X=[1,1,2]$ and $Y=[3,3,9]$. We have  $T=3, s=5$, $i_1=1, i_2=1, i_3=2, i_4=3, i_5=3, i_6=9$ and  a sequence $(m_1,m_2,m_3):=(1,1,3)$, we obtain: 

$$ G(X,Y| 1,1,3): \begin{tikzpicture}[->,shorten >=1pt,auto,node distance=3cm, thick,main node/.style={circle,draw}]

\draw  node[fill,circle,inner sep=0pt,minimum size=1pt] {};
  \node[main node] (1) {$i_1=1$};
  \node[main node] (2) [right of=1] {$i_2=1$};
  \node[main node] (3) [right of=2] {$i_3=2$};
  \node[main node] (4) [below of=1] {$i_4=3$};
  \node[main node] (5) [right of=4] {$i_5=3$};
  \node[main node] (6) [right of=5] {$i_6=9$};

  \path[every node/.style={font=\sffamily\small}]
  (1) edge [ right]  (6)
  (2) edge [ right] (6)
  (3) edge [ right]  (6)
 (3) edge [ right]  (5)
 (3) edge [ right]  (4)

            ;
            
\end{tikzpicture}.$$
\begin{remark} \label{G(X,Y)edgesproperty}
 $G(X,Y|L)$ has exactly $s$ edges and for every vertex $j$ in $Y$, there is at least one edge ending at $j$.  
\end{remark}



            

In what follows, we fix two sequences $X,Y$ of indices and the parameter $T$. Thus, for each sequence $L \in \mathcal{L}(T,s)$, we obtain a graph $G(X,Y| L)$.
Let $\lambda_{i_1}, \dots, \lambda_{i_{s+1}}$ be the sequence of eigenvalues corresponding to the sequence of indices $i_1,...,i_{s+1}$.  Define 
 $$w(X,Y| L) := \prod_{e \in \mathcal{E}(X,Y| L)} \frac{1}{(\lambda_{e^+}- \lambda_{e^-})}.$$ If either $X$ or $Y$ is the empty set, we simply set $w(X,Y| L)=0$ for any $L \in \mathcal{L}(T,s)$.

\noindent  \textit{Example:} For $X= [1,1,2]$ and $Y=[3,3,9]$ and $L=(m_1,m_2,m_3) =(1,1,3)$
 $$w(X,Y| L)= \frac{1}{(\lambda_1- \lambda_9)(\lambda_1- \lambda_9)(\lambda_2 -\lambda_9) (\lambda_2- \lambda_3) (\lambda_2 -\lambda_3)}.$$

\begin{lemma} \label{contour1-HR}  Consider a contour $\Gamma$, a parameter $0 \leq T \leq s+1$, and two sequences $X=[i_1,...,i_T], Y= [i_{T+1},...,i_{s+1}]$ of indices.
Assume that $\lambda_{i_1},...,\lambda_{i_T}$ are inside the contour $\Gamma$ and $\lambda_{i_{T+1}},...,\lambda_{i_{s+1}}$ are outside. 
Then, 
\begin{equation} \label{form1-HR}
\frac{1}{2 \pi \textbf{i}} \int_{\Gamma} \frac{1}{\prod_{j=1}^{s+1} (z- \lambda_{i_j})} dz = (-1)^{T+1}\sum_{L \in \mathcal{L}(T,s)}  w(X,Y| L).
\end{equation}
\end{lemma} 

\begin{remark} \label{remark: why residue is wrong}
From the complex analysis point of view, a natural way to compute 
the LHS of \eqref{form1-HR} is to use the residue theorem. However, in this approach, 
the resulting formula is very hard to analyze. The introduction of the combinatorial profile in \eqref{form1-HR} plays a critical role in our analysis, as it makes the number of terms controllable and many cancellations possible. 
 \end{remark}

\begin{proof}[Proof of Lemma \ref{contour1-HR}]
We use induction on $s$. In the base case  $s=1$,  $T$ is either $0, 1$ or $2$, and the claim follows from  Lemma \ref{integral}. 

 Now assume that for all $1\leq r \leq s$, the induction hypothesis holds. This means that given any two sequences $X_1=[i'_1, i'_2,..., i'_{T_1}], Y_1=[i'_{T_1+1},...,i'_{r}]$ of indices (their total length is at most $s$) such that $\lambda_{i'_1},...,\lambda_{i_{T_1}}$ are inside the contour $\Gamma$ and $\lambda_{i'_{T_1+1}},...,\lambda_{i'_{r}}$ are outside,   
$$\frac{1}{2 \pi \textbf{i}} \int_{\Gamma} \frac{1}{\prod_{j=1}^{r} (z- \lambda_{i'_j})} dz = (-1)^{T_1+1}\sum_{L \in \mathcal{L}(T_1,r-1)}  w(X_1,Y_1| L). $$
We will prove the claim for two sequences $X=[i_1,...,i_T], Y=[i_{T+1},...,i_{s+1}]$ (which have $s+1$ indices in total)
and a parameter $0\le T \le s+1$,
$$\frac{1}{2 \pi \textbf{i}} \int_{\Gamma} \frac{1}{\prod_{j=1}^{s+1} (z- \lambda_{i_j})} dz = (-1)^{T+1}\sum_{L \in \mathcal{L}(T,s)}  w(X,Y| L).$$

For each $l < k$, define $F(l,k):= \frac{1}{2 \pi \textbf{i}} \int_{\Gamma} \frac{1}{\prod_{j=l}^{k} (z-\lambda_{i_j})} dz.$ The LHS  is $F(1,s+1)$. To generate the edge \((i_1,i_{s+1})\), or equivalently the factor
\(
\frac{1}{\lambda_{i_1}-\lambda_{i_{s+1}}},
\)
we separate the two factors
\(
\frac{1}{(z-\lambda_{i_1})(z-\lambda_{i_{s+1}})}
\)
and use the identity
\(
\frac{1}{(z-\lambda_{i_1})(z-\lambda_{i_{s+1}})}
=
\frac{1}{\lambda_{i_1}-\lambda_{i_{s+1}}}
\left(
\frac{1}{z-\lambda_{i_1}}
-
\frac{1}{z-\lambda_{i_{s+1}}}
\right).
\)
Consequently,
\begin{equation} \label{IntAl0}
\begin{split}
F(1,s+1) & = \frac{1}{2 \pi \textbf{i}} \int_{\Gamma} \frac{1}{(z -\lambda_{i_1})(z-\lambda_{i_{s+1}})} \times \frac{1}{\prod_{j=2}^{s} (z-\lambda_{i_j})} dz \\
& = \frac{1}{2 \pi \textbf{i}} \times \frac{1}{(\lambda_{i_1} - \lambda_{i_{s+1}})} \int_{\Gamma} \left(  \frac{1}{z -\lambda_{i_1}} -\frac{1}{z-\lambda_{i_{s+1}}}\right) \frac{1}{\prod_{j=2}^{s} (z-\lambda_{i_j})} dz \\
& = \frac{1}{(\lambda_{i_1} -\lambda_{i_{s+1}})} \left(    \frac{1}{2 \pi \textbf{i}} \int_{\Gamma} \frac{1}{\prod_{j=1}^{s} (z-\lambda_{i_j})} dz - \frac{1}{2 \pi \textbf{i}} \int_{\Gamma} \frac{1}{\prod_{j=2}^{s+1} (z -\lambda_{i_j})} dz  \right) \\
& = \frac{1}{(\lambda_{i_1} -\lambda_{i_{s+1}})} \left[ F(1,s) - F(2,s+1)  \right].
\end{split}
\end{equation}

Since the lengths of $F(1,s)$ and $F_(2,s+1)$ are both $s$, we can apply the induction hypothesis here. Indeed, to handle $F(2,s+1)$, consider two sequences: $X':= \{ i_2,...,i_{T} \}$ and $ Y:= \{i_{T+1},...,i_{s+1}\}$, which have $s$ indices in total. Notice that for any sequence of natural numbers $m_2,...,m_T$ such that $m_2+...+m_T=s-1$, we have 
\begin{equation}
    \begin{split}
    & G(X,Y| 1,m_2,...,m_T) = G(X',Y| m_2,...,m_T) \cup \left(\vec{ e}_{i_1, i_{s+1}} \right), 
        \end{split}
\end{equation}
where $\vec{ e}_{i_1, i_{s+1}}$ is the directed edge from $i_1$ to $i_{s+1}$. Therefore, 
$$ w(X,Y|1,m_2,...,m_T) =\frac{1}{(\lambda_{i_1} -\lambda_{i_{s+1}})} \times w(X',Y| m_2,...,m_T).$$
Applying the induction hypothesis on  $\lambda_{i_2},...,\lambda_{i_{s+1}}$, we  obtain
\begin{equation} \label{IntAl1}
    \begin{split}
   \frac{1}{(\lambda_{i_1} -\lambda_{i_{s+1}})}  F(2,s+1) & = \sum_{\substack{
    m_2+...+m_T=s-1 } } \frac{ 1}{(\lambda_{i_1} - \lambda_{i_{s+1}}) } \times (-1)^{(T-1)+1} w(X',Y|m_2,...,m_T) \\
    & = \sum_{
    1+m_2+...+m_T=s } (-1)^{T} w(X,Y| 1,m_2,...,m_T).
    \end{split}
\end{equation}

Similarly, considering two other sequences $X:=\{ i_1,...,i_T \}$ and $Y':=\{i_{T+1},...,i_s\}$ (which have $s$ indices in total), and applying the induction hypothesis, we obtain 
\begin{equation} \label{IntAl2}
    \begin{split}
    \frac{1}{(\lambda_{i_1} -\lambda_{i_{s+1}})}  F(1,s) & = \sum_{\substack{m_1+...+m_T=s-1 }} \frac{ 1}{(\lambda_{i_1} - \lambda_{i_{s+1}}) } \times (-1)^{T+1}  w(X, Y'| m_1,...,m_T) \\
    & = \sum_{\substack{m_1+...+m_T=s \\ m_1 > 1}} (-1)^{T+1} w(X,Y|m_1,...,m_T).
    \end{split}
\end{equation}
Combining (\ref{IntAl0}), (\ref{IntAl1}) and (\ref{IntAl2}), we obtain 
\begin{equation*} 
\begin{split}
    F(1,s+1) & = \frac{1}{\lambda_{i_1} -\lambda_{i_{s+1}}} \left[ F(1,s) - F(2,s+1) \right] \\
    & = \sum_{\substack{m_1+...+m_T=s \\ m_1 > 1}} (-1)^{T+1} w(X,Y|m_1,...,m_T) - \sum_{ \substack{
    m_1+m_2+...+m_T=s \\
    m_1=1} } (-1)^{T} w(X,Y|m_1,m_2,...,m_T)\\
    & = \sum_{L \in \mathcal{L}(T,s)} (-1)^{T+1} w(X,Y|L).
    \end{split}
\end{equation*}
This verifies the claim for $s+1$, completing the proof. 
\end{proof} 
We will also use the following simple lemma several times in the upcoming subsections. 
\begin{lemma} \label{toytrick}
Let $a_1,a_2,a_3$ be positive numbers and $l_1,l_2,l_3$ be nonnegative integers. One has 
\begin{equation*}
    a_1^{l_1} a_2^{l_2} a_3^{l_3} \leq \max \{a_1,a_2,a_3\}^{l_1+l_2+l_3}.
\end{equation*}
 If $\max \{l_1,l_2\} \geq 1$, then 
 \begin{equation*}
    a_1^{l_1} a_2^{l_2} a_3^{l_3} \leq (a_1+a_2) \times \max\{a_1,a_2,a_3\}^{l_1+l_2+l_3-1}.  
 \end{equation*}
\end{lemma}

\subsection{Proof of Lemma \ref{LemmaCase1M(a,b)bound}}\label{section:proof2} Fix $(\alpha; \beta)$ of Type I, define
\begin{equation} \label{A_hB_hdefinition}
    \begin{split}
    & A_l := \alpha_1 + \alpha_2 + ... + \alpha_l \,\,\,\, \text{for each}\,\,\, 1 \leq l \leq k+1, \\
    & B_{\ell} := \beta_1 + \beta_2 + ... + \beta_{\ell} \,\,\,\, \text{for each} \,\,\, 1 \leq \ell \leq k, \\
    & A_0 = B_0 = 0.
    \end{split}
\end{equation}
Recall that $s_1= \sum_{i=1}^{k+1} \alpha_i = A_{k+1}$ and $s_2= \sum_{j=1}^k \beta_j = B_k.$ Note that for a pair $(\alpha; \beta)$ of Type I, since $\alpha_i \geq 1$ for all $1\leq i \leq k+1$, we have  $s_1 \geq k+1$.

Multiplying out all the terms in $P$ (there are $s_2$ of $P$'s) and all terms in $Q$ (there are $s_1$ of $Q$'s), we obtain the expansion 
$$ \int_{\Gamma} M(\alpha; \beta) dz= \sum_{ \substack{i_1,i_2,...,i_{s_2} \in N_{\bar{\lambda}}(S) \\ j_1,...,j_{s_1} \in N_{\bar{\lambda}}(S)^c}} \int_{\Gamma} \frac{dz} { \prod_{\ell=1}^{s_1} (z- \lambda_{j_{\ell}}) \prod_{l=1}^{s_2}(z-\lambda_{i_l})} \times $$
$$\left\lbrace \prod_{h=0}^{k-1} \left[ \prod_{l=1}^{\alpha_{h+1}} \left(u_{j_{l+A_{h}}} u_{j_{l+A_{h}}} ^T E  \right) \right] \times \left[ \prod_{l=1}^{\beta_{h+1}} \left(u_{i_{l+B_{h}}} u_{i_{\ell+B_{h}}}^T E \right) \right] \right\rbrace \times \left[ \prod_{l =1}^{\alpha_{k+1}-1}\left(  u_{j_{l + A_k}} u_{j_{l + A_k}}^T E \right)\right] u_{j_{s_1}} u_{j_{s_1}}^T , $$
where  the indices $j_{\ell}$'s represent the $Q$-operators, and the indices $i_l$'s represent the $P$-operators. To make the grouping of the indices clearer, let us present a small example. 

\noindent \textit{Example.} Consider $s=10, k=2$, and 
$\alpha= (3,2,1); \beta = (2,3)$. In this case 
$$s_1= 3+2+1=  6, s_2= 2+3 =5, $$ and we have 
 $$\int_{\Gamma} M(3,2,1;2,3) dz = \int_{\Gamma} \underset{\substack{\alpha_1=3\\ \text{$Q$-operators}}}{\underbrace{(QEQEQE)}} \underset{\substack{\beta_1=2\\ \text{$P$-operators}}} {\underbrace{(PEPE)}}\,\,\, \underset{\substack{\alpha_2=2 \\ \text{$Q$-operators } }}{\underbrace{(QEQE)}} \underset{\substack{\beta_2=3 \\ \text{$P$-operators}}}{\underbrace{(PEPEPE)}} \,\,\, \underset{\substack{\alpha_3=1 \\ \text{$Q$-operator}}}{\underbrace{Q}} dz$$
into 
 $$\sum_{ \substack{i_1,i_2,...,i_{5} \in N_{\bar{\lambda}}(S) \\ j_1,...,j_{6} \in N_{\bar{\lambda}}(S)^c}} \int_{\Gamma} \frac{dz} { \prod_{\ell=1}^{6} (z- \lambda_{j_{\ell}}) \prod_{l=1}^{5}(z-\lambda_{i_l})} \left[ (u_{j_1} u_{j_1}^T E)(u_{j_2} u_{j_2}^T E) (u_{j_3} u_{j_3}^T E)(u_{i_1} u_{i_1}^T E) (u_{i_2} u_{i_2}^T E) \right]$$
 $$ \times  \left[(u_{j_4} u_{j_4}^T E)(u_{j_5} u_{j_5}^T E) (u_{i_3} u_{i_3}^T E) (u_{i_4} u_{i_4}^T E) (u_{i_4} u_{i_4}^T E) \right] \times \left[ u_{j_6} u_{j_6}^T \right] . $$
From now on, to reduce unnecessary indices, we will write $\sum_{i_1,...,i_{s_2}}$ and $\sum_{j_1,...,j_{s_1}}$ instead of $\sum_{i_1,...,i_{s_2} \in N_{\bar{\lambda}}(S)}$ and $\sum_{j_1,...,j_{s_1} \in N_{\bar{\lambda}}(S)^c}$. The next idea is  to use Lemma \ref{contour1-HR} to treat $\int_{\Gamma} \frac{dz} { \prod_{\ell=1}^{s_1} (z- \lambda_{j_{\ell}}) \prod_{l=1}^{s_2}(z-\lambda_{i_l})}$. In order to do that, we need to know which eigenvalues are inside/outside the contour $\Gamma$. 

We have two sequences of indices: $[i_1,...,i_{s_2}]$, whose elements are in  $N_{\bar{\lambda}}(S)$ and $[j_1,...,j_{s_1}]$, whose elements are in $N_{\bar{\lambda}}(S)^c$. The order of the elements in the  sequence $[i_1,...,i_{s_2}]$ ($[j_1,...,j_{s_1}]$) is  inherited from the positions of $P$'s ($Q$'s resp.), counting from left to right. Now, we fix a sequence $i_1,...,i_{s_2} $ and split it exactly into the following two sub-sequences. The first one (denoted by $X$) is the subsequence of indices $i_{l}$ in $S$, i.e. $\lambda_{i_l}$ is inside contour $\Gamma$. The second one (denoted by $I_2= i_{l_1}, i_{l_2},..., i_{l_{|I_2|}}$) is obtained from $i_1,...,i_{s_2}$ by removing $X$. Define the sequence $Y$ as 
$$Y= [i_{l_1}, i_{l_2},..., i_{l_{|I_2|}}, j_1, j_2,..., j_{s_1}].$$
\noindent \textit{Example:} Continuing with the above example ($s=10, k=2$ and so on), we further set $N_{\lambda}(S)=\{1,2,3,4,5\}$ and $S=\{1,2,3\}$. Consider the sequence  $[i_1,i_2,i_3,i_4,i_5] = [1,4,2,2,3]$. In this case, 
$$X=[1,2,2,3] \,\,\, \text{and} \,\,\, Y =[4, j_1, j_2, j_3, j_4, j_5, j_6].$$
 Let $T =|X|$ and apply Lemma \ref{contour1-HR}, we have 
 \begin{equation}
     \begin{split}
  \frac{1}{2 \pi \textbf{i}}\int_{\Gamma} \frac{dz} { \prod_{\ell=1}^{s_1} (z- \lambda_{j_{\ell}}) \prod_{l=1}^{s_2}(z-\lambda_{i_l})} & = \sum_{L \in \mathcal{L}(T,s)} (-1)^{T+1} w(X,Y|L) \\
  & = \sum_{L \in \mathcal{L}(T,s)} (-1)^{T+1} \prod_{e \in \mathcal{E}(X,Y|L)} \frac{1}{ (\lambda_{e^+} - \lambda_{e^-})}.
     \end{split}
 \end{equation}
 Notice that once the sequence $[i_1, i_2, ...,i_{s_2}]$ and $L$ are  fixed, then the set $X$ is well-defined, and the shape of the graph $G(X, Y| L)$ is also determined, regardless of the particular values of 
 $j_l$.   In this case, we can write $\mathcal{E}_L$ instead of $\mathcal{E}(X,Y|L)$ (for the edge set of the graph). Let us illustrate with a concrete example.

\noindent \textit{Example:} Continuing with the previous example $(s=10, k=2$ and so on), we further set $N_{\lambda}(S)=\{1,2,3,4,5\}$ and $S=\{1,2,3\}$. Consider the sequence  $[i_1,i_2,i_3,i_4,i_5] = [1,4,2,2,3]$. In this case, 
$$X=[1,2,2,3] \,\,\, \text{and} \,\,\, Y =[4, j_1, j_2, j_3, j_4, j_5, j_6].$$
Thus, $T=4$. We set
$$L=(m_1, m_2, m_3, m_4)= (4,3,2,1).$$

The graph $G(1,2,2,3; 4,j_1,j_2,j_3,j_4,j_5,j_6|4,3,2,1)$ is 
$$ 
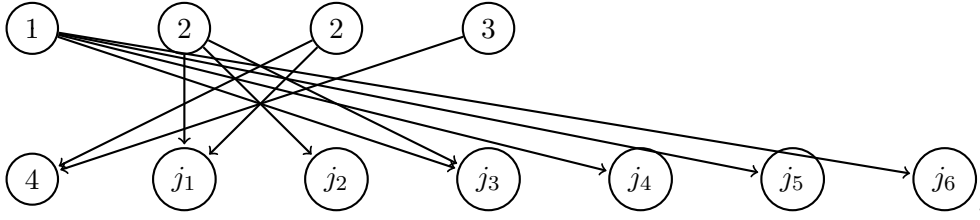
\begin{tikzpicture}[->,shorten >=1pt,auto,node distance=2cm, thick,main node/.style={circle,draw}]

\draw  node[fill,circle,inner sep=0pt,minimum size=1pt] {};
  \node[main node] (1) {$1$};
  \node[main node] (2) [right of=1] {$2$};
  \node[main node] (3) [right of=2] {$2$};
  \node[main node] (4) [right of=3] {$3$};
  \node[main node] (5) [below of=1] {$4$};
  \node[main node] (6) [right of=5] {$j_1$};
  \node[main node] (7) [right of=6] {$j_2$};
  \node[main node] (8) [right of=7] {$j_3$};
  \node[main node] (9) [right of=8] {$j_4$};
  \node[main node] (10) [right of=9] {$j_5$};
  \node[main node] (11) [right of=10] {$j_6$};

  \path[every node/.style={font=\sffamily\small}]
  (1) edge [ right]  (11)
  (1) edge [ right]  (10)
  (1) edge [ right]  (9)
  (1) edge [ right]  (8)
  (2) edge [ right] (8)
  (2) edge [ right] (7)
  (2) edge [ right] (6)
  (3) edge [ right]  (6)
 (3) edge [ right]  (5)
 (4) edge [ right]  (5)

            ;
            
\end{tikzpicture}.$$
The set of edges $\mathcal{E}(X,Y|L)$ is 
$$\left\lbrace (1 \rightarrow j_6), (1 \rightarrow j_5), (1 \rightarrow j_4), (1 \rightarrow j_3), (2 \rightarrow j_3), (2 \rightarrow j_2), (2 \rightarrow j_1), (2 \rightarrow j_1), (2 \rightarrow 4), (3 \rightarrow 4)    \right\rbrace ,$$
and $w(X,Y|L)$ equals
$$ \frac{1}{(\lambda_1 - \lambda_{j_6}) (\lambda_1 - \lambda_{j_5}) (\lambda_1 - \lambda_{j_4}) (\lambda_1 - \lambda_{j_3}) (\lambda_2 - \lambda_{j_3}) (\lambda_2 - \lambda_{j_2}) (\lambda_2 - \lambda_{j_1}) (\lambda_2 - \lambda_{j_1}) (\lambda_2 -\lambda_4)(\lambda_3 -\lambda_4) }.$$

Therefore, we can write 
\begin{equation}\label{M(a,b)Case1firstidentity}
\begin{split}
& \frac{1}{2\pi \textbf{i}} \int_{\Gamma} M(\alpha; \beta) dz  = \sum_{i_1,i_2,...,i_{s_2}} \sum_{L \in \mathcal{L}(T,s)}\sum_{j_1,...,j_{s_1}}  (-1)^{T+1} \prod_{e \in \mathcal{E}_L} \frac{1 }{ (\lambda_{e^+} - \lambda_{e^-})} \times    \\
& \bigg\{ \prod_{h=0}^{k-1} \bigg[ \prod_{l=1}^{\alpha_{h+1}} \left(u_{j_{l+A_{h}}} u_{j_{l+A_{h}}} ^T E  \right) \bigg] \times \bigg[ \prod_{l=1}^{\beta_{h+1}} \left(u_{i_{l+B_{h}}} u_{i_{l+B_{h}}}^T E \right) \bigg] \bigg\} \times \bigg[ \prod_{l =1}^{\alpha_{k+1}-1}\left(  u_{j_{l + A_k}} u_{j_{l + A_k}}^T E \right)\bigg] u_{j_{s_1}} u_{j_{s_1}}^T. 
\end{split}
\end{equation}

Now, we "redistribute" the factors in $ \prod_{e \in \mathcal{E}_L} \frac{1}{(\lambda_{e^+}-\lambda_{e^-})}$. 
By the definition of $Y$, $e^{-}$ either is in $I_2$ or equals $j_l$ for some $1\leq l \leq s_1$. Accordingly, we decompose the product as
$$\prod_{e \in \mathcal{E}_L} \frac{1}{(\lambda_{e^+}-\lambda_{e^-})}= \bigg[\prod_{\substack{e \in \mathcal{E}_L\\ e^{-} \in I_2}} \frac{1}{(\lambda_{e^+}-\lambda_{e^-})} \bigg] \times \prod_{l=1}^{s_1} \bigg[\prod_{\substack{e \in \mathcal{E}_L \\ e^{-} = j_l}} \frac{1}{(\lambda_{e^+}-\lambda_{e^-})} \bigg].$$
Intuitively, the first sub-product will give rise to factors of \(1/\delta_S\), while the remaining sub-products will give rise to factors of \(1/\bar{\lambda}\). 

\noindent\textit{Example:} Continuing with the above example, we can rewrite $w(X,Y|L)= \prod_{e \in \mathcal{E}_L} \frac{1}{(\lambda_{e^+}-\lambda_{e^-})}$ as 
\begin{equation*}
    \begin{split}
  & \frac{1}{(\lambda_2 -\lambda_4)(\lambda_3 -\lambda_4)} \times \frac{1}{(\lambda_2 - \lambda_{j_1}) (\lambda_2 - \lambda_{j_1})} \times \frac{1}{ (\lambda_2 - \lambda_{j_2})} \times \\
  &\frac{1}{ (\lambda_1 - \lambda_{j_3}) (\lambda_2 - \lambda_{j_3})} \times \frac{1}{(\lambda_1 - \lambda_{j_4})} \times \frac{1}{(\lambda_1 - \lambda_{j_5})} \times  \frac{1}{(\lambda_1 - \lambda_{j_6})}.       
    \end{split}
\end{equation*}

For each $1 \leq l \leq s_1$, we group the factor $\prod_{\substack{e \in \mathcal{E}_L \\ e^{-} = j_l}} \frac{1}{(\lambda_{e^+}-\lambda_{e^-})}$ with the rank-one projection $u_{j_l}  u_{j_l}^T$. We may view this factor as the ''\textit{perturbation weight}" associated with $u_{j_l}  u_{j_l}^T$. With this regrouping, we obtain
\begin{equation*}
\begin{split}
 \frac{1}{2\pi \textbf{i}} \int_{\Gamma} M(\alpha; \beta) dz 
 & = \sum_{i_1,i_2,...,i_{s_2}} \sum_{L \in \mathcal{L}(T,s)}\sum_{j_1,...,j_{s_1}}  (-1)^{T+1} \prod_{\substack{e \in \mathcal{E}_L \\ e^{-} \in I_2}} \frac{1 }{ (\lambda_{e^+} - \lambda_{e^-})} \times    \\
& \bigg\{ \prod_{h=0}^{k-1} \bigg[ \prod_{l=1}^{\alpha_{h+1}} \bigg( \frac{u_{j_{l+A_{h}}} u_{j_{l+A_{h}}}^T}{\displaystyle \prod_{\substack{e \in \mathcal{E}_L; \, e^{-} = j_{l+A_{h}}}} (\lambda_{e^{+}}- \lambda_{j_{l+A_{h}}})} E  \bigg) \bigg] \times \bigg[  \prod_{l=1}^{\beta_{h+1}} \left( u_{i_{l+B_{h}}} u_{i_{l+B_{h}}}^T E \right) \bigg] \bigg\} \times \\
& \times \bigg[  \prod_{l =1}^{\alpha_{k+1}-1} \bigg( \frac{u_{j_{l + A_k}} u_{j_{l + A_k}}^T}{\displaystyle \prod_{\substack{e \in \mathcal{E}_L;\, e^{-} = j_{l+A_{k}}}} (\lambda_{e^{+}}- \lambda_{j_{l+A_{k}}})} E \bigg) \bigg]\frac{u_{j_{s_1}} u_{j_{s_1}}^T}{\displaystyle\prod_{\substack{e \in \mathcal{E}_L;\, e^{-}= j_{s_1}}}(\lambda_{e^+}-\lambda_{j_{s_1}})}.
\end{split}
\end{equation*}
%
%
 \begin{remark} A key point is that
 the combinatorial Lemma~\ref{contour1-HR} allows us to treat 
$j_1,\dots,j_{s_1}$ purely as indices—i.e., with the same combinatorial profile—regardless of their specific values in $N_{\bar{\lambda}}(S)^c$. 
For each $1 \leq l \leq s_1$, rather than handling each interaction $E u_{j_l} u_{j_l}^T$ separately for $j_l \in N_{\bar{\lambda}}(S)^c$, we can evaluate the aggregated term
\(
    E \Big( \sum_{j_l \in N_{\bar{\lambda}}(S)^c} u_{j_l} u_{j_l}^T \Big).
\)
Consequently, instead of incurring a factor $(n-r)$ from separate interactions, we pay only a factor at most $\|E\|$ (and often less, since all interactions $E u_{i_k} u_{i_k}^T$ and $E u_{j_l} u_{j_l}^T$ are treated jointly in a single profile). 
This yields a substantial saving; for example, when $E$ is Wigner, $\|E\|=\Theta(\sqrt{n})$, and hence the total factor $(n-r)^{s_1}$ for all $j_1, j_2, \dots, j_{s_1}$ reduces to a factor at most $\|E\|^{s_1}=\Theta( n^{s_1/2})$. 
 \end{remark}
Indeed, we can distribute the sum over $j_1, j_2,..., j_{s_1}$ into each position of $u_{j_l} u_{j_l}^T$ in the product and obtain that $\frac{1}{2\pi \textbf{i}} \int_{\Gamma} M(\alpha; \beta) dz$ equals 
\begin{equation*}
    \begin{split}
     & \sum_{i_1,i_2,...,i_{s_2}} \sum_{L \in \mathcal{L}(T,s)} (-1)^{T+1} \prod_{\substack{e \in \mathcal{E}_L \\ e^{-} \in I_2}} \frac{1 }{ (\lambda_{e^+} - \lambda_{e^-})} \times    \\
& \bigg\{ \prod_{h=0}^{k-1} \bigg[ \prod_{l=1}^{\alpha_{h+1}}  \bigg( \sum_{j_{l+A_{h}} \in N_{\bar{\lambda}}(S)^c} \frac{u_{j_{l+A_{h}}} u_{j_{l+A_{h}}}^T}{  \displaystyle \prod_{\substack{e \in \mathcal{E}_L;\, e^{-} = j_{l+A_{h}}}} (\lambda_{e^{+}}- \lambda_{j_{l+A_{h}}})}  E  \bigg) \bigg] \times \bigg[ \bigg( \prod_{l=1}^{\beta_{h+1}} u_{i_{l+B_{h}}} u_{i_{l+B_{h}}}^T E \bigg) \bigg] \bigg\} \times
 \end{split}
\end{equation*}
\begin{equation*}
    \begin{split}
& \times \bigg[ \prod_{l =1}^{\alpha_{k+1}-1} \bigg( \sum_{j_{l+A_k} \in N_{\bar{\lambda}}(S)^c} \frac{u_{j_{l + A_k}} u_{j_{l + A_k}}^T}{\displaystyle \prod_{\substack{e \in \mathcal{E}_L \\ e^{-} = j_{l+A_{k}}}} (\lambda_{e^{+}}- \lambda_{j_{l+A_{k}}})} E \bigg) \bigg] \times \bigg(\sum_{j_{s_1} \in N_{\bar{\lambda}}(S)^c} \frac{u_{j_{s_1}} u_{j_{s_1}}^T}{\displaystyle \prod_{\substack{e \in \mathcal{E}_L \\ e^{-}= j_{s_1}}}(\lambda_{e^+}-\lambda_{j_{s_1}})} \bigg)\\
& =  \sum_{i_1,i_2,...,i_{s_2}} \sum_{L \in \mathcal{L}(T,s)} (-1)^{T+1} \prod_{\substack{e \in \mathcal{E}_L \\ e^{-} \in I_2}} \frac{1 }{ (\lambda_{e^+} - \lambda_{e^-})} \times   W \times  \\
& \times \bigg[ \prod_{l =1}^{\alpha_{k+1}-1} \bigg( \sum_{j_{l+A_k} \in N_{\bar{\lambda}}(S)^c} \frac{u_{j_{l + A_k}} u_{j_{l + A_k}}^T}{\displaystyle \prod_{\substack{e \in \mathcal{E}_L \\ e^{-} = j_{l+A_{k}}}} (\lambda_{e^{+}}- \lambda_{j_{l+ A_{k}}})} E \bigg) \bigg] \times \bigg(\sum_{j_{s_1} \in N_{\bar{\lambda}}(S)^c} \frac{u_{j_{s_1}} u_{j_{s_1}}^T}{\displaystyle \prod_{\substack{e \in \mathcal{E}_L \\ e^{-} = j_{s_1}}}(\lambda_{e^+}-\lambda_{j_{s_1}})} \bigg), 
    \end{split}
\end{equation*}
where

\begin{equation*}
    \begin{split}
    W & =  \bigg\{ \prod_{h=0}^{k-1} \bigg[ \prod_{l=1}^{\alpha_{h+1}}  \bigg( \sum_{j_{l+A_{h}} \in N_{\bar{\lambda}}(S)^c} \frac{u_{j_{l+A_{h}}} u_{j_{l+A_{h}}}^T}{ \displaystyle \prod_{\substack{e \in \mathcal{E}_L;\, e^{-} = j_{l+A_{h}}}} (\lambda_{e^{+}}- \lambda_{j_{l+A_{h}}})}  E  \bigg) \bigg] \times \bigg[ \bigg( \prod_{l=1}^{\beta_{h+1}} u_{i_{l+B_{h}}} u_{i_{l+B_{h}}}^T E \bigg) \bigg] \bigg\}.
    \end{split}
\end{equation*}

It is important to notice that the \textit{combinatorial profile} tells us exactly the perturbation weight for each  $u_k u_k^T$-projection for $k \in [n]$. With respect to each combinatorial profile, we then distribute the perturbation weights (p.w.) and group up the projections so that the interaction factors are either $u_i^T E u_{i'}\,\,\,\text{for} \,\,\,i, i' \in N_{\bar{\lambda}}(S) $ (interaction of $E$ with the important eigenvectors) or $E  \cdot \big[\sum_{j \in N_{\bar{\lambda}}(S)^c} (\textit{p.w. of} \,\, u_ju_j^T) \cdot u_j u_j^T \big]$  (interaction of $E$ with the \textit{unimportant} eigenspace).

Next, by the triangle inequality, we have
\begin{equation} \label{M(a,b)triangle}
\frac{1}{2 \pi} \Norm{\int_{\Gamma} M(\alpha;\beta) dz}  \leq \sum_{i_1,i_2,...,i_{s_2}} \sum_{L \in \mathcal{L}(T,s)}  \bigg|\prod_{\substack{e \in \mathcal{E}_L \\ e^{-} \in I_2}} \frac{1 }{ (\lambda_{e^+} - \lambda_{e^-})}\bigg| \times  \|W\| \times 
\end{equation}
\[
\times \bigg[ \prod_{l =1}^{\alpha_{k+1}-1} \bigg( \bigg\|\sum_{j_{l+A_k} \in N_{\bar{\lambda}}(S)^c} \frac{u_{j_{l + A_k}} u_{j_{l + A_k}}^T}{\displaystyle \prod_{\substack{e \in \mathcal{E}_L \\ e^{-} = j_{l+A_{k}}}} (\lambda_{e^{+}}- \lambda_{j_{l+A_{k}}})}\bigg\|  \|E\| \bigg) \bigg] \times \bigg\|\sum_{j_{s_1} \in N_{\bar{\lambda}}(S)^c} \frac{u_{j_{s_1}} u_{j_{s_1}}^T}{\displaystyle \prod_{\substack{e \in \mathcal{E}_L \\ e^{-} = j_{s_1}}}(\lambda_{e^+}-\lambda_{j_{s_1}})}\bigg\|.
\]

Now, we bound the terms (interactions) from left to right. First, 
\begin{equation} \label{Factor1}
  \bigg| \displaystyle\prod_{\substack{e \in \mathcal{E}_L \\ e^{-} \in I_2}} \frac{1}{ (\lambda_{e^{+}}-\lambda_{e^{-}})} \bigg| \leq \frac{1}{\delta_S^{|E(I_2)|}},
\end{equation} where $E(I_2)$ is the set of edges ending in $I_2$ in $G(X,Y|L)$. 

Next, we can split $\|W\|$ into (1) the interactions of $E$ with the important eigenvectors, (2) interactions of $E$ with the \textit{unimportant} eigenspaces, (3) the remaining factors of the form $\|E u_i\|$ with $i \in N_{\bar{\lambda}}(S)$. Indeed,   
\begin{equation*}
    \begin{split}
    W & = \prod_{h=0}^{k-1} \bigg[ \prod_{l=1}^{\alpha_{h+1}-1} \bigg(  \sum_{j_{l+A_{h}} \in N_{\bar{\lambda}}(S)^c} \frac{u_{j_{l+A_{h}}} u_{j_{l+A_{h}}}^T}{ \displaystyle \prod_{\substack{e \in \mathcal{E}_L;\, e^{-} = j_{l+A_{h}}}} (\lambda_{e^{+}}- \lambda_{j_{l+A_{h}}})}  E  \bigg) \bigg] \\
    & \times \bigg[ \bigg( \sum_{j_{\alpha_{h+1}+A_{h}} \in N_{\bar{\lambda}}(S)^c} \frac{u_{j_{\alpha_{h+1}+A_{h}}} u_{j_{\alpha_{h+1}+A_{h}}}^T}{\displaystyle \prod_{\substack{e \in \mathcal{E}_L; \, e^{-} = j_{\alpha_{h+1}+A_{h}}}} (\lambda_{e^{+}}- \lambda_{j_{\alpha_{h+1}+A_{h}}})} \bigg) E   \bigg( \prod_{l=1}^{\beta_{h+1}} u_{i_{l+B_{h}}} u_{i_{l+B_{h}}}^T E \bigg)\bigg] \\
    & =  \prod_{h=0}^{k-1} \bigg[ \prod_{l=1}^{\alpha_{h+1}-1} \bigg(  \sum_{j_{l+A_{h}} \in N_{\bar{\lambda}}(S)^c} \frac{u_{j_{l+A_{h}}} u_{j_{l+A_{h}}}^T}{ \displaystyle \prod_{\substack{e \in \mathcal{E}_L;\, e^{-} = j_{l+A_{h}}}} (\lambda_{e^{+}}- \lambda_{j_{l+A_{h}}})}  E  \bigg) \bigg] \\
       \end{split}
\end{equation*}
\begin{equation*}
    \begin{split}
        & \times \bigg[ \bigg( \sum_{j_{A_{h+1}} \in N_{\bar{\lambda}}(S)^c} \frac{u_{j_{A_{h+1}}} u_{j_{A_{h+1}}}^T}{ \displaystyle\prod_{\substack{e \in \mathcal{E}_L;\, e^{-} = j_{A_{h+1}}}} (\lambda_{e^{+}}- \lambda_{j_{A_{h+1}}})} \bigg) E   \bigg( \prod_{l=1}^{\beta_{h+1}} u_{i_{l+B_{h}}} u_{i_{l+B_{h}}}^T E \bigg)\bigg] \\
    & (\text{by the definitions in \eqref{A_hB_hdefinition}, we have $\alpha_{h+1}+A_h=A_{h+1}$}).
    \end{split}
\end{equation*}
Therefore, 
\begin{equation} \label{Wdefinitionbound}
    \begin{split}
    \|W\| & \leq   \prod_{h=0}^{k-1} \bigg[ \prod_{l=1}^{\alpha_{h+1}-1} \bigg( \bigg\| \sum_{j_{l+A_{h}} \in N_{\bar{\lambda}}(S)^c} \frac{u_{j_{l+A_{h}}} u_{j_{l+A_{h}}}^T}{ \displaystyle\prod_{\substack{e \in \mathcal{E}_L;\, e^{-} = j_{l+A_{h}}}} (\lambda_{e^{+}}- \lambda_{j_{l+A_{h}}})}\bigg\| \cdot \|E\|  \bigg) \bigg] \\
    & \times  \bigg\|\sum_{j_{A_{h+1}} \in N_{\bar{\lambda}}(S)^c} \frac{u_{j_{A_{h+1}}} u_{j_{A_{h+1}}}^T}{ \displaystyle \prod_{\substack{e \in \mathcal{E}_L;\, e^{-} = j_{A_{h+1}}}} (\lambda_{e^{+}}- \lambda_{j_{A_{h+1}}})} \bigg\| \cdot \bigg\| E \bigg( \prod_{l=1}^{\beta_{h+1}} u_{i_{l+B_{h}}} u_{i_{l+B_{h}}}^T E \bigg)\bigg\|.      \end{split}
\end{equation}
Notice that 
\begin{equation} \label{uEuB}
    \begin{split}
    \bigg\| E \bigg( \prod_{l=1}^{\beta_{h+1}} u_{i_{l+B_{h}}} u_{i_{l+B_{h}}}^T E \bigg) \bigg\| & = \bigg\| E u_{i_{1+B_{h}}} \bigg( \prod_{l=1}^{\beta_{h+1}-1}  (u_{i_{l+B_{h}}}^T E u_{i_{l+1+B_{h}}}) \bigg) u_{i_{B_{h+1}}}^T E  \bigg\| \\
    & \leq \Norm{E u_{i_{1+B_{h}}}} \cdot \Norm{u_{i_{B_{h+1}}}^T E } \times \prod_{l=1}^{\beta_{h+1}-1} \norm{ u_{i_{l+B_{h}}}^T E u_{i_{l+1+B_{h}}}} \\
    & \leq w^2 \cdot x^{\beta_{h+1}-1} (\text{by Definition \ref{xyz}}).
    \end{split}
\end{equation}
Moreover, for any fixed $1 \leq l \leq s_1$,
\begin{equation} \label{u_ju_j^Tsumbound}
    \bigg\| \sum_{j_l \in N_{\bar{\lambda}}(S)^c } \frac{u_{j_l} u_{j_l}^T}{ \displaystyle \prod_{\substack{e \in \mathcal{E}_L; \, e^{-} = j_l}} (\lambda_{e^{+}} -\lambda_{j_l})} \bigg\| \leq \frac{1}{\bar{\lambda}^{d(j_l)}}, 
\end{equation}
where $d(j_l)$ is the degree of $j_l$ in $G(X,Y|L)$ (this degree does not depend on the particular value of $j_l$). Applying \eqref{uEuB} and \eqref{u_ju_j^Tsumbound} on the RHS of \eqref{Wdefinitionbound}, we obtain
\begin{equation} \label{Wbound}
\begin{split}
  \|W\| & \leq \prod_{h =0}^{k-1} \left[ w^{2} x^{(\beta_{h+1}-1)} \|E\|^{(\alpha_{h+1}-1)} \left(\prod_{1 \leq l \leq \alpha_{h+1} } \frac{1}{\bar{\lambda}^{d(j_{l+A_h})}} \right) \right] \\
  & = \prod_{h =0}^{k-1} \left[ w^{2} x^{(\beta_{h+1}-1)} \|E\|^{(\alpha_{h+1}-1)} \left(\prod_{1+A_h \leq l \leq A_{h+1} } \frac{1}{\bar{\lambda}^{d(j_{l})}} \right) \right] \\
  & = w^{2k} x^{\sum_{h=0}^{k-1}(\beta_{h+1}-1)} \|E\|^{\sum_{h=0}^{k-1}(\alpha_{h+1}-1)} \prod_{h=0}^{k-1} \left(\prod_{1+A_h \leq l \leq A_{h+1} } \frac{1}{\bar{\lambda}^{d(j_{l})}} \right) \\
  & = w^{2k} x^{\sum_{l=1}^{k}(\beta_l -1)} \|E\|^{\sum_{l=1}^{k}(\alpha_l -1)} \prod_{1 \leq l \leq A_k} \frac{1}{\bar{\lambda}^{d(j_l)}} \\
  & =   \frac{w^{2k} x^{\sum_{l=1}^{k}(\beta_l -1)} \|E\|^{\sum_{l=1}^{k}(\alpha_l -1)}}{\bar{\lambda}^{\sum_{l=1}^{A_k} d(j_l)}}.
\end{split}
   \end{equation}

Arguing as in  \eqref{u_ju_j^Tsumbound}, we also have
\[
\bigg\|\sum_{j_{l+A_k} \in N_{\bar{\lambda}}(S)^c} \frac{u_{j_{l + A_k}} u_{j_{l + A_k}}^T}{\displaystyle \prod_{\substack{e \in \mathcal{E}_L; \, e^{-} = j_{l+A_{k}}}} (\lambda_{e^{+}}- \lambda_{j_{l+A_{k}}})} \bigg\| \leq \frac{1}{\bar{\lambda}^{d(j_{l+A_k})}} \,\,\text{for all}\,\, 1 \leq l \leq \alpha_{k+1}-1,
\]
and
\[
\bigg\|\sum_{j_{s_1} \in N_{\bar{\lambda}}(S)^c} \frac{u_{j_{s_1}} u_{j_{s_1}}^T}{\displaystyle\prod_{\substack{e \in \mathcal{E}_L; \, e^{-}= j_{s_1}}}(\lambda_{e^+}-\lambda_{j_{s_1}})} \bigg\| \leq \frac{1}{\bar{\lambda}^{d(j_{s_1})}}.
\]
And hence,
\begin{equation} \label{Factor34case1}
\begin{split}
     & \bigg[ \prod_{l =1}^{\alpha_{k+1}-1} \bigg( \bigg\|\sum_{j_{l+A_k} \in N_{\bar{\lambda}}(S)^c} \frac{u_{j_{l + A_k}} u_{j_{l + A_k}}^T}{\displaystyle \prod_{\substack{e \in \mathcal{E}_L \\ e^{-} = j_{l+A_{k}}}} (\lambda_{e^{+}}- \lambda_{j_{l+A_{k}}})} \bigg\|  \|E\| \bigg) \bigg] \times \bigg\| \sum_{j_{s_1} \in N_{\bar{\lambda}}(S)^c} \frac{u_{j_{s_1}} u_{j_{s_1}}^T}{\displaystyle \prod_{\substack{e \in \mathcal{E}_L \\ e^{-}= j_{s_1}}}(\lambda_{e^+}-\lambda_{j_{s_1}})}\bigg\| \\
& \leq \|E\|^{\alpha_{k+1}-1} \prod_{1 + A_k \leq l \leq s_1} \frac{1}{\bar{\lambda}^{d(j_l)}} = \frac{\|E\|^{\alpha_{k+1}-1}}{\bar{\lambda}^{\sum_{l=1+A_k}^{s_1} d(j_l)}}.
\end{split}
\end{equation}

Combining \eqref{M(a,b)triangle}, (\ref{Factor1}), (\ref{Wbound}), and (\ref{Factor34case1}), we finally obtain
\begin{equation} \label{beforesplitM(a,b)type1}
    \frac{1}{2 \pi } \Norm{\int_{\Gamma} M(\alpha;\beta) dz} \leq \sum_{i_1,...,i_{s_2}} \sum_{L \in \mathcal{L}(T,s)} \frac{w^{2k}  x^{\sum_{l=1}^k (\beta_l-1)} \|E\|^{\sum_{l=1}^{k+1} (\alpha_l-1)}}{\delta_S^{|E(I_2)|} \bar{\lambda}^{\sum_{l=1}^{s_1} d(j_l)}}.   
\end{equation}

 Set $r_c(T,L):= \sum_{l=1}^{s_1} (d(j_l)-1) \geq 0$. Then, by Remark \ref{G(X,Y)edgesproperty},
 \begin{equation} \label{GraphEdgeComputation}
   \sum_{l=1}^{s_1} d(j_l) = r_c(T,L)+s_1\,\,\, \text{and} \,\, |E(I_2)| = s - (r_c(T,L)+s_1) = s_2 -1-r_c(T,L).   
 \end{equation}
   We can rewrite the above inequality as
\begin{equation}
\begin{split}
\frac{1}{2 \pi}\Norm{\int_\Gamma M(\alpha;\beta) dz} & \leq \sum_{i_1,...,i_{s_2}} \sum_{L \in \mathcal{L}(T,s)} \frac{ \|E\|^{s_1-(k+1)} w^{2k} x^{s_2 - k}}{\delta_S^{s_2 -1-r_c(T,L)} \bar{\lambda}^{s_1+r_c(T,L)}}.
\end{split}
\end{equation}

Now, we need to estimate the number of \textit{combinatorial profiles}. Notice that if we fix $T$ and a sequence $L \in \mathcal{L}(T,s)$, then the shape of the graph $G(X,Y|L)$, and hence $r_c(T,L)$ are uniquely determined. Hence, if two different sequences $[i_1,...,i_{s_2}]$ and $[i'_1,...,i'_{s_2}]$ has the same number $T$ of elements in $S$, then their corresponding $r_c(T,L)$'s are the same for each fixed $L \in \mathcal{L}(T,s)$.  Moreover, when we fix $T$, the choices for $i_1,..., i_{s_2}$ are at most $ \binom{s_2}{T} p^{T} (r-p)^{s_2 -T} < \binom{s_2}{T} p^{T} r^{s_2 -T}$. Thus, \eqref{beforesplitM(a,b)type1} leads to 
\begin{equation} \label{splitM(a,b)type1}
\begin{split}
\frac{1}{2 \pi}\Norm{\int_\Gamma M(\alpha;\beta) dz} & \leq \sum_{T=1}^{s_2} \binom{s_2}{T}  p^{T} r^{s_2 -T} \sum_{L \in \mathcal{L}(T,s)} \frac{\|E\|^{s_1-(k+1)} w^{2k} x^{s_2 - k}}{\delta_S^{s_2 -1-r_c(T,L)} \bar{\lambda}^{s_1+r_c(T,L)}}  \\
& = M_1(\alpha;\beta) + M_2(\alpha;\beta), \,\, \text{where}\\
\end{split}
\end{equation}
\begin{equation*}
    \begin{split}
M_1(\alpha;\beta) &:= \sum_{T=1}^{s_2} \binom{s_2}{T} p^{T} r^{s_2 -T}  \sum_{\substack{L \in \mathcal{L}(T,s);\, r_c(T,L) \leq k-1}} \frac{ \|E\|^{s_1-(k+1)} w^{2k} x^{s_2 - k}}{\delta_S^{s_2 -1-r_c(T,L)} \bar{\lambda}^{s_1+r_c(T,L)}} , \\
M_2 (\alpha; \beta) &:= \sum_{T=1}^{s_2} \binom{s_2}{T} p^{T} r^{s_2 -T} \sum_{\substack{L \in \mathcal{L}(T,s);\, r_c(T,L) > k-1}} \frac{\|E\|^{s_1-(k+1)} w^{2k} x^{s_2 - k} }{\delta_S^{s_2 -1-r_c(T,L)} \bar{\lambda}^{s_1+r_c(T,L)}} . 
    \end{split}
\end{equation*}
In order to bound $M_1(\alpha;\beta)$ and $ M_2(\alpha;\beta)$ from above, we need the following property of $r_c(T,L)$. 
\begin{lemma} \label{r_cLemma} 
    $$r_c(T,L) \leq T-1.$$
\end{lemma}
\noindent\textit{Proof of Lemma \ref{r_cLemma}.} By Remark \ref{G(X,Y)edgesproperty}, we have
$\norm{E(I_2)} \geq |I_2|.$
Since $\norm{E(I_2)} = s_2 - (r_c(T,L)+1)$ and $|I_2| = s_2-T$, we obtain
$$r_c(T,L) \leq T-1.$$

Back to our estimation, we can bound $M_1(\alpha;\beta)$ from above as follows. Since $w \leq \|E\|$, $w^{2k} \leq w^{2(k-1 -r_c(T,L))} \|E\|^{2(r_c(T,L)+1)}$. Thus, 
\begin{equation} \label{Mab-bound-HR}
 \begin{split}
  M_1(\alpha;\beta) & \leq \sum_{T=1}^{s_2}\binom{s_2}{T}  p^{T} r^{s_2 -T} \sum_{\substack{L \in \mathcal{L}(T,s);\, r_c(T,L) \leq k-1}} \frac{ \|E\|^{s_1-(k+1)+2(r_c(T,L)+1)} w^{2(k-1-r_c(T,L))} x^{s_2 - k}}{\delta_S^{s_2 -1-r_c(T,L)} \bar{\lambda}^{s_1+r_c(T,L)}} \\
  & = \sum_{T=1}^{s_2} \binom{s_2}{T}   \sum_{\substack{L \in \mathcal{L}(T,s);\, r_c(T,L) \leq k-1}} p^T r^{s_2-T} \left( \frac{ \|E\|}{\bar{\lambda}} \right)^{T_1} \left( \frac{ x}{\delta_S} \right)^{T_2} \bigg( \frac{ w}{\sqrt {\bar{\lambda} \delta_S}  } \bigg)^{T_3},
 \end{split}
  \end{equation}
where 
$$T_1= s_1 - (k-1)+2r_c(T,L) \geq 2, \qquad T_2= s_2 - k,$$
$$T_3= 2(k-1 -r_c(T,L)).$$
Redistributing the factors $\sqrt{p}, r, \sqrt{r}$ to  $\frac{\|E\|}{\bar{\lambda}}, \frac{x}{\delta_S}, \frac{w}{\sqrt{\bar{\lambda} \delta_S}}$, respectively, we can rewrite the RHS of \eqref{Mab-bound-HR} as
$$ \sum_{T=1}^{s_2} \binom{s_2}{T}   \sum_{\substack{L \in \mathcal{L}(T,s);\, r_c(T,L) \leq k-1}} \left( \frac{p}{r} \right)^{T-(r_c(T,L)+1)} p^{\frac{k+1 -s_1}{2}}  \left( \frac{ \sqrt{p}\|E\|}{\bar{\lambda}} \right)^{T_1} \left( \frac{ rx}{\delta_S} \right)^{T_2} \bigg( \frac{ \sqrt r w}{\sqrt {\bar{\lambda} \delta_S}  } \bigg)^{T_3}.$$

Since $T \geq r_c(T,L)+1$ (Lemma  \ref{r_cLemma}) and $s_1-(k+1) \geq 0$ ($(\alpha,\beta)$ is Type I),
$$ 0 \leq \left( \frac{p}{r} \right)^{T-(r_c(T,L)+1)} \leq 1 \,\,\,\text{and}\,\,\, 0 \leq  p^{\frac{k+1 -s_1}{2}}  \leq 1.$$
We can omit these factors and obtain 
\begin{equation*}
    \begin{split}
  M_1(\alpha;\beta)&  \leq \sum_{T=1}^{s_2} \binom{s_2}{T}   \sum_{\substack{L \in \mathcal{L}(T,s);\, r_c(T,L) \leq k-1}}  \left( \frac{ \sqrt{p}\|E\|}{\bar{\lambda}} \right)^{T_1} \left( \frac{ rx}{\delta_S} \right)^{T_2} \bigg( \frac{ \sqrt r w}{\sqrt {\bar{\lambda} \delta_S}  } \bigg)^{T_3} \\
  & = \sum_{T=1}^{s_2} \binom{s_2}{T}   \sum_{\substack{L \in \mathcal{L}(T,s);\, r_c(T,L) \leq k-1}} \frac{\sqrt{p}\|E\|}{\bar{\lambda}} \times  \left( \frac{ \sqrt{p}\|E\|}{\bar{\lambda}} \right)^{T_1-1} \left( \frac{ rx}{\delta_S} \right)^{T_2} \bigg( \frac{ \sqrt r w}{\sqrt {\bar{\lambda} \delta_S}  } \bigg)^{T_3}.
    \end{split}
\end{equation*}
\noindent By the setting of $M_1(\alpha;\beta)$, $r_c(T,L) \leq k-1$, and then $T_3 \geq 0$. We also have $T_2 \geq 0$ and $T_1 -1 \geq 1$ by the definition of Type I $(\alpha,\beta)$. Therefore, by Lemma \ref{toytrick}, 
\begin{equation} \label{M_1a1ak+1>0}
    \begin{split}
 M_1(\alpha;\beta) 
 & \leq \sum_{T=1}^{s_2} \binom{s_2}{T}   \sum_{\substack{L \in \mathcal{L}(T,s): r_c(T,L) \leq k-1}}  \left( \frac{ \sqrt{p}\|E\|}{\bar{\lambda}} \right) \times   \max \left\lbrace  \frac{\sqrt{p}\|E\|}{\bar{\lambda}}, \frac{rx}{\delta_S} , \frac{ \sqrt r w}{\sqrt { \bar{\lambda} \delta_S }} \right\rbrace^{T_1+T_2+T_3-1} \\
 & = \bigg[\sum_{T=1}^{s_2}  \sum_{\substack{L \in \mathcal{L}(T,s): \\ r_c(T,L) \leq k-1}} \binom{s_2}{T} \bigg] \times \left( \frac{ \sqrt{p}\|E\|}{\bar{\lambda}} \right) \times   \max \left\lbrace  \frac{\sqrt{p}\|E\|}{\bar{\lambda}}, \frac{rx}{\delta_S} , \frac{ \sqrt r w}{\sqrt { \bar{\lambda} \delta_S }} \right\rbrace^{T_1+T_2+T_3-1} \\
 & \leq \bigg[\sum_{T=1}^{s_2}  \sum_{L \in \mathcal{L}(T,s)} \binom{s_2}{T} \bigg] \times \left( \frac{ \sqrt{p}\|E\|}{\bar{\lambda}} \right) \times   \max \left\lbrace  \frac{\sqrt{p}\|E\|}{\bar{\lambda}}, \frac{rx}{\delta_S} , \frac{ \sqrt r w}{\sqrt { \bar{\lambda} \delta_S }} \right\rbrace^{T_1+T_2+T_3-1}. 
 \end{split}
\end{equation}
Since $\norm{\mathcal{L}(T,s)} =\binom{s-1}{T-1} \leq 2^{s-1}$ for all $T \leq s$, we have 
\begin{equation} \label{s_2Lsumbound}
   \sum_{T=1}^{s_2}  \sum_{L \in \mathcal{L}(T,s)} \binom{s_2}{T} \leq 2^{s-1}  \sum_{T=1}^{s_2} \binom{s_2}{T} \leq 2^{s_2+s-1}. 
\end{equation}
Thus, \eqref{M_1a1ak+1>0} becomes 
\begin{equation}
    \begin{split}
    M_1(\alpha;\beta) & \leq 2^{s_2+s-1}  \left( \frac{ \sqrt{p}\|E\|}{\bar{\lambda}} \right) \times   \max \left\lbrace  \frac{\sqrt{p}\|E\|}{\bar{\lambda}}, \frac{rx}{\delta_S} , \frac{ \sqrt r w}{\sqrt { \bar{\lambda} \delta_S }} \right\rbrace^{T_1+T_2+T_3-1} \\
  & \leq \left( \frac{ \sqrt{p}\|E\|}{\bar{\lambda}} \right) \frac{2^{s_2+s-1}}{12^{T_1+T_2+T_3-1}} \,\,(\text{by Assumption}\,\,\textbf{D0}) \\
 & = \left( \frac{ \sqrt{p}\|E\|}{\bar{\lambda}} \right) \frac{2^{s_2+s-1}}{12^{s-1}} \,\, (\text{since}\,\, T_1+T_2+T_3 =s).      \end{split}
\end{equation}

Similarly, we can bound $M_2(\alpha;\beta)$ as follows. Since $w \leq \|E\|$, $w^{2k} \leq \|E\|^{2k}$ and hence 
\begin{equation}
    \begin{split} \label{M_2case1Rearrange}
   M_2(\alpha;\beta) & \leq \sum_{T=1}^{s_2} \binom{s_2}{T} p^{T} r^{s_2 -T} \sum_{\substack{L \in \mathcal{L}(T,s);\, r_c(T,L) > k-1}} \frac{\|E\|^{s_1+k-1} x^{s_2 - k} }{\delta_S^{s_2 -1-r_c(T,L)} \bar{\lambda}^{s_1+r_c(T,L)}} \\
   & = \sum_{T=1}^{s_2} \binom{s_2}{T} p^{T} r^{s_2 -T} \sum_{\substack{L \in \mathcal{L}(T,s);\, r_c(T,L) > k-1}} \left( \frac{ x}{\delta_S} \right)^{s_2-1-r_c(T,L)} \left(\frac{x}{\bar{\lambda}} \right)^{r_c(T,L)-(k-1)} \left( \frac{ \|E\|}{\bar{\lambda}} \right)^{k-1+s_1}.
    \end{split}
\end{equation}
\noindent Distributing the factors $r,p, \sqrt{p}$ to $\frac{x}{\delta_S}, \frac{x}{\bar{\lambda}}, \frac{\|E\|}{\bar{\lambda}}$ respectively, we can rewrite the RHS as 
$$ \sum_{T=1}^{s_2} \binom{s_2}{T} \sum_{\substack{L \in \mathcal{L}(T,s) \\ r_c(T,L) > k-1}} \left( \frac{p}{r}  \right)^{T-(r_c(T,L)+1)} p^{\frac{k+1 -s_1}{2}} \left( \frac{r x}{\delta_S} \right)^{s_2-1-r_c(T,L)} \left(\frac{px}{\bar{\lambda}} \right)^{r_c(T,L)-(k-1)} \left( \frac{ \sqrt{p} \|E\|}{\bar{\lambda}} \right)^{k-1+s_1}. $$
\noindent Since $T \geq r_c(T,L)+1$ (Lemma  \ref{r_cLemma}) and $s_1-(k+1) \geq 0$ ($(\alpha,\beta)$ is Type I), 
$$ 0 \leq \left( \frac{p}{r} \right)^{T-(r_c(T,L)+1)} \leq 1 \,\,\,\text{and}\,\,\, 0\leq p^{\frac{k+1 -s_1}{2}}  \leq 1.$$
We can omit these factors and obtain 
\begin{equation*}
    \begin{split}
 M_2(\alpha;\beta) & \leq \sum_{T=1}^{s_2} \binom{s_2}{T} \sum_{\substack{L \in \mathcal{L}(T,s) \\ r_c(T,L) > k-1}} \left( \frac{r x}{\delta_S} \right)^{s_2-1-r_c(T,L)} \left(\frac{px}{\bar{\lambda}} \right)^{r_c(T,L)-(k-1)} \left( \frac{ \sqrt{p} \|E\|}{\bar{\lambda}} \right)^{k-1+s_1} \\
 & = \sum_{T=1}^{s_2} \binom{s_2}{T} \sum_{\substack{L \in \mathcal{L}(T,s) \\ r_c(T,L) > k-1}} \frac{px}{\bar{\lambda}} \times  \left( \frac{r x}{\delta_S} \right)^{s_2-1-r_c(T,L)} \left(\frac{px}{\bar{\lambda}} \right)^{r_c(T,L)-(k-1)-1} \left( \frac{ \sqrt{p} \|E\|}{\bar{\lambda}} \right)^{k-1+s_1}.
    \end{split}
\end{equation*}
By the setting of $M_2(\alpha; \beta)$, $r_c(T,L) - (k-1) -1\geq 0$. By definitions of $s_1,s_2, r_c(T,L)$, we also have $s_2 -1-r_c(T,L),   k-1+s_1 \geq 0$. Applying Lemma \ref{toytrick} and the fact that 
\[
[s_2-1-r_c(T,L)]+[r_c(T,L)-(k-1)-1]+[k-1+s_1] = s_1+s_2 -2 = s-1,
\]
 we obtain
\begin{equation*}
\begin{split} 
M_2(\alpha;\beta)
& \leq \sum_{T=1}^{s_2} \binom{s_2}{T} \sum_{\substack{L \in \mathcal{L}(T,s) \\ r_c(T,L) > k-1}} \frac{px}{\bar{\lambda}} \times \max \left\lbrace \frac{rx}{\delta_S}, \frac{px}{\bar{\lambda}}, \frac{\sqrt{p}\|E\|}{\bar{\lambda}} \right\rbrace^{s-1} \\
& = \bigg[\sum_{T=1}^{s_2}  \sum_{\substack{L \in \mathcal{L}(T,s) \\ r_c(T,L) > k-1}} \binom{s_2}{T} \bigg] \times  \frac{px}{\bar{\lambda}} \times \max \left\lbrace \frac{rx}{\delta_S}, \frac{px}{\bar{\lambda}}, \frac{\sqrt{p}\|E\|}{\bar{\lambda}} \right\rbrace^{s-1} \\
& \leq \bigg[\sum_{T=1}^{s_2}  \sum_{\substack{L \in \mathcal{L}(T,s)}} \binom{s_2}{T} \bigg] \times  \frac{px}{\bar{\lambda}} \times \max \left\lbrace \frac{rx}{\delta_S}, \frac{px}{\bar{\lambda}}, \frac{\sqrt{p}\|E\|}{\bar{\lambda}} \right\rbrace^{s-1}.
\end{split}
\end{equation*}
Using \eqref{s_2Lsumbound} and the fact that $\frac{px}{\bar{\lambda}} \leq \frac{rx}{\delta_S}$, we finally have 
\begin{equation}\label{M_2a_1a_{k+1}>0}
 M_2(\alpha;\beta) \leq   2^{s_2+s-1}  \frac{p x}{\bar{\lambda}} \max \left\lbrace \frac{rx}{\delta_S}, \frac{\sqrt{p}\|E\|}{\bar{\lambda}} \right\rbrace^{s-1} \leq \frac{p x}{\bar{\lambda}} \times \frac{2^{s_2+s-1}}{12^{s-1}} \,\,(\text{by Assumption}\,\,\textbf{D0}).
\end{equation}
Combining \eqref{splitM(a,b)type1}, \eqref{M_1a1ak+1>0} and \eqref{M_2a_1a_{k+1}>0}, we obtain that  
\begin{equation*} 
    \frac{1}{2 \pi} \Norm{\int_{\Gamma} M(\alpha; \beta) dz} \leq \left( \frac{ \sqrt{p}\|E\|}{\bar{\lambda}} + \frac{px}{\bar{\lambda}} \right) \frac{2^{s_2+s-1}}{12^{s-1}}.
\end{equation*}
This proves Lemma \ref{LemmaCase1M(a,b)bound}.

\subsection{Proof of Lemma \ref{LemmaCase2M(a,b)bound}} \label{subsec: proofCase2}
Recall that in the setting of Lemma \ref{LemmaCase2M(a,b)bound}, $(\alpha;\beta)$ is of Type II, i.e., 
\[
\alpha_{1}=\alpha_{k+1} = 0, s_1 =\sum_{i=2}^k \alpha_{i} \geq k-1, \,\text{and}\,\, s_2=\sum_{j=1}^k \beta_j \geq k. 
\]
Similar to \eqref{M(a,b)Case1firstidentity}, we decompose $\frac{1}{2\pi \textbf{i}} \int_{\Gamma} M(\alpha; \beta) dz$ as
\[
\sum_{i_1,i_2,...,i_{s_2}} \sum_{L \in \mathcal{L}(T,s)} \sum_{j_1,...,j_{s_1}}  (-1)^{T+1} \prod_{\substack{e \in \mathcal{E}_L \\ e^{-} \in I_2}} \frac{1 }{ (\lambda_{e^+} - \lambda_{e^-})} \times \left[ \prod_{l=1}^{\beta_1} \left( u_{i_{l}} u_{i_{l}}^T E \right) \right] \times 
\]
\[
\bigg\{ \prod_{h=1}^{k-2} \bigg[ \prod_{l=1}^{\alpha_{h+1}} \bigg( \frac{u_{j_{l+A_{h}}} u_{j_{l+A_{h}}}^T}{\displaystyle \prod_{\substack{e \in \mathcal{E}_L;\, e^{-} = j_{l+A_{h}}}} (\lambda_{e^{+}}- \lambda_{j_{l+A_{h}}})} E  \bigg) \bigg] \times \bigg[  \prod_{l=1}^{\beta_{h+1}} \left( u_{i_{l+B_{h}}} u_{i_{l+B_{h}}}^T E \right) \bigg] \bigg\} \times
\]
\[
\times \bigg[  \prod_{l =1}^{\alpha_{k}} \bigg( \frac{u_{j_{l + A_{k-1}}} u_{j_{l + A_{k-1}}}^T}{\displaystyle\prod_{\substack{e \in \mathcal{E}_L;\, e^{-} = j_{l+A_{k-1}}}} (\lambda_{e^{+}}- \lambda_{j_{l+A_{k-1}}})} E \bigg) \bigg] \times \bigg[ \prod_{l =1}^{\beta_k-1} \left( u_{i_{l+B_{k-1}}} u_{i_{l+B_{k-1}}}^T E \right)  \bigg] u_{i_{s_2}} u_{i_{s_2}}^T.
\]
Similar to the previous subsection, we then distribute the sum over $j_1, j_2, ...,j_{s_1}$ according to  each position of $u_{j_l} u_{j_l}^T$ in the product, and obtain
\begin{equation*}
    \begin{split}
 & \frac{1}{2\pi \textbf{i}} \int_{\Gamma} M(\alpha; \beta) dz \\   
 & = \sum_{i_1,i_2,...,i_{s_2}} \sum_{L \in \mathcal{L}(T,s)} (-1)^{T+1} \prod_{\substack{e \in \mathcal{E}_L \\ e^{-} \in I_2}} \frac{1 }{ (\lambda_{e^+} - \lambda_{e^-})} \times \left[ \prod_{l=1}^{\beta_1} \left( u_{i_{l}} u_{i_{l}}^T E \right) \right] \times     \\
& \bigg\{ \prod_{h=1}^{k-2} \bigg[ \prod_{l=1}^{\alpha_{h+1}} \bigg( \sum_{j_{l+A_{h}} \in N_{\bar{\lambda}}(S)^c} \frac{u_{j_{l+A_{h}}} u_{j_{l+A_{h}}}^T}{ \displaystyle\prod_{\substack{e \in \mathcal{E}_L;\, e^{-} = j_{l+A_{h}}}} (\lambda_{e^{+}}- \lambda_{j_{l+A_{h}}})} E  \bigg) \bigg] \times \bigg[  \prod_{l=1}^{\beta_{h+1}} \left( u_{i_{l+B_{h}}} u_{i_{l+B_{h}}}^T E \right) \bigg] \bigg\} \times \\
& \times \bigg[  \prod_{l =1}^{\alpha_{k}} \bigg( \sum_{j_{l+A_{k-1}} \in N_{\bar{\lambda}}(S)^c}  \frac{u_{j_{l + A_{k-1}}} u_{j_{l + A_{k-1}}}^T}{\displaystyle\prod_{\substack{e \in \mathcal{E}_L;\, e^{-} = j_{l+A_{k-1}}}} (\lambda_{e^{+}}- \lambda_{j_{l+A_{k-1}}})} E \bigg) \bigg] \times \left[ \prod_{l =1}^{\beta_k-1} \left( u_{i_{l+B_{k-1}}} u_{i_{l+B_{k-1}}}^T E \right)  \right] u_{i_{s_2}} u_{i_{s_2}}^T. 
    \end{split}
\end{equation*}

Intuitively, estimating $\frac{1}{2\pi \textbf{i}} \int_{\Gamma} M(\alpha; \beta) dz$ for $(\alpha;\beta)$ of Type II can be viewed as a generalization of the estimation of $\frac{1}{2 \pi \textbf{i}} \int_{\Gamma} PEQEP dz$ arising from $F_2$. We split the sum above into two parts. The first gives rise to the term involving \(\frac{y}{\delta_S}\) and requires additional care, while the second can be handled by an argument similar to that of the previous subsection. 

Indeed, notice that once we fix the sequence $[i_1,...,i_{s_2}]$ and $L \in \mathcal{L}(T,s)$, the graph $G(X,Y|L)$ and $r_c(T,L)$ are determined. We can split the double sum $\frac{1}{2\pi \textbf{i}} \int_{\Gamma} M(\alpha;\beta) dz= \sum_{i_1,i_2,...,i_{s_2}} \sum_{L \in \mathcal{L}(T,s)} \left( \cdots \right) $ above into two sub-sums: 
\begin{equation} \label{Sigma1Sigma2split}
    \begin{split}
 & \Sigma_1 := \sum_{i_1,i_2,...,i_{s_2}} \sum_{\substack{L \in \mathcal{L}(T,s)\\ r_c(T,L)=0}}\left( \cdots \right),\, \text{and}\,\,\, \Sigma_2:= \sum_{i_1,i_2,...,i_{s_2}} \sum_{\substack{L \in \mathcal{L}(T,s) \\ r_c(T,L)> 0}} \left( \cdots \right).
    \end{split}
\end{equation}

We prove two lemmas bounding $\Norm{\Sigma_1}$ and $\Norm{\Sigma_2}$. We assume that Assumption \textbf{D0} holds in both lemmas. 
\begin{lemma} \label{lemma: Sigma1allcases} 
 If $(s_2,s_1) \neq (k,k-1)$, then 
 \begin{equation} \label{lemma: Sigma1goodcase}
     \Norm{\Sigma_1} \leq \left( \frac{ \sqrt{pr} x}{\delta_S} + \frac{\sqrt{p} \|E\|}{\bar{\lambda}} \right) \frac{\sum_{T = 1}^{s_2} \displaystyle \sum_{L \in \mathcal{L}(T,s);\, r_c(T,L)=0} \binom{s_2}{T} }{12^{s-1}}.
 \end{equation}
 If $(s_2,s_1) = (k,k-1)$, then 
 \begin{equation} \label{lemma: Sigma1badcase}
     \Norm{\Sigma_1} \leq \frac{(s+2) \sqrt{rp} y}{\delta_S} \left( \frac{2 \sqrt{r} w}{\sqrt{\bar{\lambda} \delta_S}} \right)^{s-2}.
 \end{equation}
\end{lemma}
\begin{lemma}\label{lemma: Sigma2allcases} 
    \begin{equation}
        \Norm{\Sigma_2} \leq \frac{p x}{\bar{\lambda}} \frac{ \sum_{T=2}^{s_2}  \sum_{\substack{L \in \mathcal{L}(T,s) \\r_c(T,L) > k-1}} \binom{s_2}{T} }{12^{s-1}} + \left( \frac{ \sqrt{p}\|E\|}{\bar{\lambda}} \right) \frac{\sum_{T=2}^{s_2}    \sum_{\substack{L \in \mathcal{L}(T,s) \\ 1 \leq r_c(T,L) \leq k-1}}\binom{s_2}{T}}{12^{s-1}}.
    \end{equation}
\end{lemma}

Together, Lemma \ref{lemma: Sigma2allcases} and \eqref{lemma: Sigma1goodcase} imply that if $(s_2,s_1) \neq (k,k-1)$, then 
\begin{equation*}
    \begin{split}
\Norm{\Sigma_1} + \Norm{\Sigma_2} & \leq \left( \frac{ \sqrt{pr} x}{\delta_S} + \frac{\sqrt{p} \|E\|}{\bar{\lambda}} \right) \frac{\sum_{T = 1}^{s_2}  \sum_{\substack{L \in \mathcal{L}(T,s)\\ r_c(T,L)=0}} \binom{s_2}{T} }{12^{s-1}} \\
& + 
\frac{p x}{\bar{\lambda}} \frac{ \sum_{T=2}^{s_2}  \sum_{\substack{L \in \mathcal{L}(T,s) \\r_c(T,L) > k-1}} \binom{s_2}{T} }{12^{s-1}} + \left( \frac{ \sqrt{p}\|E\|}{\bar{\lambda}} \right) \frac{\sum_{T=2}^{s_2}    \sum_{\substack{L \in \mathcal{L}(T,s) \\ 1 \leq r_c(T,L) \leq k-1}}\binom{s_2}{T}}{12^{s-1}} .
\end{split}
\end{equation*}
Replacing $\frac{px}{\bar{\lambda}}$ by $\frac{\sqrt{pr} x}{\delta_S}$ (the larger one) and  rearranging the RHS, we further obtain
\begin{equation}
    \begin{split}
\Norm{\Sigma_1} + \Norm{\Sigma_2} & \leq   \frac{\sqrt{pr} x}{\delta_S} \frac{ \sum_{T = 1}^{s_2}  \sum_{\substack{L \in \mathcal{L}(T,s)\\ r_c(T,L)=0}} \binom{s_2}{T} +\sum_{T=2}^{s_2}  \sum_{\substack{L \in \mathcal{L}(T,s); \\r_c(T,L) > k-1}} \binom{s_2}{T} }{12^{s-1}} + \\
& + \frac{\sqrt{p} \|E\|}{\bar{\lambda}} \frac{\sum_{T = 1}^{s_2}  \sum_{\substack{L \in \mathcal{L}(T,s) \\ r_c(T,L)=0}} \binom{s_2}{T} +\sum_{T=2}^{s_2}    \sum_{\substack{L \in \mathcal{L}(T,s); \\ 1 \leq r_c(T,L) \leq k-1}}\binom{s_2}{T} }{12^{s-1}} \\
& \leq  \frac{\sqrt{pr} x}{\delta_S} \frac{\sum_{T=1}^{s_2}  \sum_{\substack{L \in \mathcal{L}(T,s) }} \binom{s_2}{T} }{12^{s-1}} + \frac{\sqrt{p} \|E\|}{\bar{\lambda}} \frac{\sum_{T=1}^{s_2}  \sum_{\substack{L \in \mathcal{L}(T,s) }} \binom{s_2}{T} }{12^{s-1}} \\
& \leq \left( \frac{\sqrt{pr} x}{\delta_S} + \frac{\sqrt{p}\|E\|}{\bar{\lambda}} \right)  \frac{2^{s_2+s-1}}{12^{s-1}} \,\, (\text{by \eqref{s_2Lsumbound}}).
    \end{split}
\end{equation}
Therefore, if $(s_2,s_1) \neq (k,k-1)$, we have 
$$\frac{1}{2 \pi} \Norm{\int_{\Gamma} M(\alpha;\beta)} \leq \Norm{\Sigma_1} + \Norm{\Sigma_2} \leq \left( \frac{\sqrt{pr} x}{\delta_S} + \frac{\sqrt{p}\|E\|}{\bar{\lambda}}\right) \frac{2^{s_2+s-1}}{12^{s-1}}.$$
This proves the first part of Lemma \ref{LemmaCase2M(a,b)bound}. 
    
On the other hand, together, Lemma \ref{lemma: Sigma2allcases} and \eqref{lemma: Sigma1badcase} imply that if $(s_2,s_1) =(k,k-1)$, then
\begin{equation} \label{Simga_1+Sigma2badcase}
\begin{split}
  & \Norm{\Sigma_1} + \Norm{\Sigma_2} \leq \\
  & \frac{(s+2) \sqrt{rp} y}{\delta_S} \left( \frac{2 \sqrt{r} w}{\sqrt{\bar{\lambda} \delta_S}} \right)^{s-2} + \frac{p x}{\bar{\lambda}} \frac{ \sum_{T=2}^{s_2}  \sum_{\substack{L \in \mathcal{L}(T,s) \\r_c(T,L) > k-1}} \binom{s_2}{T} }{12^{s-1}} + \left( \frac{ \sqrt{p}\|E\|}{\bar{\lambda}} \right) \frac{\sum_{T=2}^{s_2}    \sum_{\substack{L \in \mathcal{L}(T,s) \\ 1 \leq r_c(T,L) \leq k-1}}\binom{s_2}{T}}{12^{s-1}}. 
\end{split} 
\end{equation}

Moreover, by \eqref{s_2Lsumbound},
$$ \max \left\lbrace \frac{ \sum_{T=2}^{s_2}  \sum_{\substack{L \in \mathcal{L}(T,s) \\r_c(T,L) > k-1}} \binom{s_2}{T} }{12^{s-1}}, \frac{\sum_{T=2}^{s_2}    \sum_{\substack{L \in \mathcal{L}(T,s) \\ 1 \leq r_c(T,L) \leq k-1}}\binom{s_2}{T}}{12^{s-1}} \right\rbrace \leq \frac{\sum_{T=1}^{s_2}  \sum_{\substack{L \in \mathcal{L}(T,s) }} \binom{s_2}{T}}{12^{s-1}} \leq \frac{2^{s_2+s-1}}{12^{s-1}}.  $$
We can rewrite \eqref{Simga_1+Sigma2badcase} into
\begin{equation*}
\Norm{\Sigma_1} + \Norm{\Sigma_2} \leq \frac{(s+2) \sqrt{rp} y}{\delta_S} \left( \frac{2 \sqrt{r} w}{\sqrt{\bar{\lambda} \delta_S}} \right)^{s-2} + \left(\frac{px}{\bar{\lambda}} + \frac{\sqrt{p} \|E\|}{\bar{\lambda}} \right) \frac{2^{s_2+s-1}}{12^{s-1}}. 
\end{equation*}
This implies that if $(s_2,s_1) = (k,k-1)$, 
$$\frac{1}{2 \pi} \Norm{\int_{\Gamma} M(\alpha;\beta)} \leq \Norm{\Sigma_1} + \Norm{\Sigma_2} \leq \left(\frac{px}{\bar{\lambda}} + \frac{\sqrt{p} \|E\|}{\bar{\lambda}} \right) \frac{2^{s_2+s-1}}{12^{s-1}}+ \frac{(s+2) \sqrt{rp} y}{\delta_S} \left( \frac{2 \sqrt{r} w}{\sqrt{\bar{\lambda} \delta_S}} \right)^{s-2}.  $$
In this case, $s_2=k, s_1 =k-1, s=s_1+s_2-1=2(k-1)$ and hence $s_2=s/2+1$, $s_2+s-1=3s/2$. The RHS can be rewritten as 
$$\left(\frac{px}{\bar{\lambda}} + \frac{\sqrt{p} \|E\|}{\bar{\lambda}} \right) \frac{2^{3s/2}}{12^{s-1}}+ \frac{(s+2) \sqrt{rp} y}{\delta_S} \left( \frac{2 \sqrt{r} w}{\sqrt{\bar{\lambda} \delta_S}} \right)^{s-2}.$$
This proves the second part of Lemma \ref{LemmaCase2M(a,b)bound}.
\subsubsection{Proof of Lemma \ref{lemma: Sigma2allcases}}
\noindent Recall that for  any square matrix $M$, 
\begin{equation} \label{normdefinition}
 \|M\|:= \max_{\|\textbf{v}\| = \|\textbf{w}\|=1} \textbf{v}^T M \textbf{w}.   
\end{equation}
\noindent Therefore, $\Norm{\Sigma_2}$ equals
\begin{equation*}
\begin{split}
 &\max_{\|\textbf{v}\|=\|\textbf{w}\|=1} \sum_{i_1,i_2,...,i_{s_2}} \sum_{L \in \mathcal{L}(T,s);\, r_c(T,L) >0} (-1)^{T+1} \prod_{\substack{e \in \mathcal{E}_L \\ e^{-} \in I_2}} \frac{1 }{ (\lambda_{e^+} - \lambda_{e^-})} \times \textbf{v}^T \left[ \prod_{l=1}^{\beta_1}  \left( u_{i_{l}} u_{i_{l}}^T E \right) \right] \times     \\
& \bigg\{ \prod_{h=1}^{k-2} \bigg[\prod_{l=1}^{\alpha_{h+1}} \bigg( \sum_{j_{l+A_{h}} \in N_{\bar{\lambda}}(S)^c} \frac{u_{j_{l+A_{h}}} u_{j_{l+A_{h}}}^T}{ \displaystyle\prod_{\substack{e \in \mathcal{E}_L;\, e^{-} = j_{l+A_{h}}}} (\lambda_{e^{+}}- \lambda_{j_{l+A_{h}}})} E  \bigg) \bigg] \times \bigg[  \prod_{l=1}^{\beta_{h+1}} \left( u_{i_{l+B_{h}}} u_{i_{l+B_{h}}}^T E \right) \bigg] \bigg\} \times \\
& \times \bigg[  \prod_{l =1}^{\alpha_{k}} \bigg( \sum_{j_{l+A_{k-1}} \in N_{\bar{\lambda}}(S)^c}  \frac{u_{j_{l + A_{k-1}}} u_{j_{l + A_{k-1}}}^T}{\displaystyle\prod_{\substack{e \in \mathcal{E}_L;\, e^{-} = j_{l+A_{k-1}}}} (\lambda_{e^{+}}- \lambda_{j_{l+A_{k-1}}})} E \bigg) \bigg] \times \bigg[ \prod_{l =1}^{\beta_k-1} \left( u_{i_{l+B_{k-1}}} u_{i_{l+B_{k-1}}}^T E \right)  \bigg] u_{i_{s_2}} u_{i_{s_2}}^T \textbf{w}.
\end{split}    
\end{equation*}
Thus, by the triangle inequality,
\begin{equation}\label{Sigma2splitstep0}
    \begin{split}
 & \Norm{\Sigma_2} \leq \max_{\|\textbf{v}\|=\|\textbf{w}\|=1} \sum_{i_1,i_2,...,i_{s_2}} \sum_{L \in \mathcal{L}(T,s);\, r_c(T,L) > 0} \norm{ \prod_{\substack{e \in \mathcal{E}_L \\ e^{-} \in I_2}} \frac{1 }{ (\lambda_{e^+} - \lambda_{e^-})}} \times \Norm{ \textbf{v}^T \prod_{l=1}^{\beta_1} \left( u_{i_{l}} u_{i_{l}}^T E \right) } \times     \\
& \bigg\| \prod_{h=1}^{k-2} \bigg[ \prod_{l=1}^{\alpha_{h+1}} \bigg( \sum_{j_{l+A_{h}} \in N_{\bar{\lambda}}(S)^c} \frac{u_{j_{l+A_{h}}} u_{j_{l+A_{h}}}^T}{ \displaystyle\prod_{\substack{e \in \mathcal{E}_L;\, e^{-} = j_{l+A_{h}}}} (\lambda_{e^{+}}- \lambda_{j_{l+A_{h}}})} E  \bigg) \bigg] \times \bigg[  \prod_{l=1}^{\beta_{h+1}} \left( u_{i_{l+B_{h}}} u_{i_{l+B_{h}}}^T E \right) \bigg] \bigg\| \times \\
& \times \bigg\| \bigg[  \prod_{l =1}^{\alpha_{k}} \bigg( \sum_{j_{l+A_{k-1}} \in N_{\bar{\lambda}}(S)^c}  \frac{u_{j_{l + A_{k-1}}} u_{j_{l + A_{k-1}}}^T}{\displaystyle\prod_{\substack{e \in \mathcal{E}_L;\, e^{-} = j_{l+A_{k-1}}}} (\lambda_{e^{+}}- \lambda_{j_{l+A_{k-1}}})} E \bigg) \bigg] \times \bigg[ \prod_{l =1}^{\beta_k-1} \bigg( u_{i_{l+B_{k-1}}} u_{i_{l+B_{k-1}}}^T E \bigg)  \bigg] u_{i_{s_2}} u_{i_{s_2}}^T \textbf{w}\bigg\|.
    \end{split}
\end{equation}
We proceed in a fashion similar to the previous subsection. First,
\begin{equation} \label{Factor1,case2}
  \bigg|\frac{1}{\prod_{\substack{e \in \mathcal{E}_L \\ e^{-} \in I_2}} (\lambda_{e^{+}}-\lambda_{e^{-}})} \bigg| \leq \frac{1}{\delta_S^{|E(I_2)|}}, 
\end{equation}
in which $E(I_2)$ is the set of edges ending in $I_2$ in $G(X,Y|L)$.  
Second, 
\begin{equation} \label{Factor2,case2}
    \begin{split}
     \Norm{ \textbf{v}^T \prod_{l=1}^{\beta_1} \left( u_{i_{l}} u_{i_{l}}^T E \right) } &  \leq  \norm{\textbf{v}^T u_{i_1}} \times \left[\prod_{l =1}^{\beta_1-1} \norm{u_{i_{l}}^T E u_{i_{l+1}} } \right] \times \Norm{u_{i_{\beta_1}}^T E} \\
     & \leq \norm{\textbf{v}^T u_{i_1}} x^{\beta_1 -1} w \,(\text{by Definition \ref{xyz}}).
    \end{split}
\end{equation}
Third, by repeating the procedure of bounding $\|W\|$ from the previous subsection, we obtain 

\begin{equation} \label{Factor3,case2}
\begin{split}
  & \bigg\| \prod_{h=1}^{k-2} \bigg[ \prod_{l=1}^{\alpha_{h+1}} \bigg( \sum_{j_{l+A_{h}} \in N_{\bar{\lambda}}(S)^c} \frac{u_{j_{l+A_{h}}} u_{j_{l+A_{h}}}^T}{ \displaystyle\prod_{\substack{e \in \mathcal{E}_L;\, e^{-} = j_{l+A_{h}}}} (\lambda_{e^{+}}- \lambda_{j_{l+A_{h}}})} E  \bigg) \bigg] \times \bigg[  \prod_{l=1}^{\beta_{h+1}} \left( u_{i_{l+B_{h}}} u_{i_{l+B_{h}}}^T E \right) \bigg] \bigg\|  \\
  & \leq \prod_{h =1}^{k-2} \bigg[ w^{2} x^{(\beta_{h+1}-1)} \|E\|^{(\alpha_{h+1}-1)} \bigg(\prod_{1+ A_{h} \leq l \leq A_{h+1} } \frac{1}{\bar{\lambda}^{d(j_l)}} \bigg) \bigg] \\
  & = \frac{w^{2(k-2)} x^{\sum_{l=2}^{k-1}(\beta_l -1)} \|E\|^{\sum_{l=2}^{k-1}(\alpha_l-1)}}{\bar{\lambda}^{\sum_{l=1}^{A_{k-1}} d(j_l)}} \,\, (\text{since}\,\, \alpha_1 =0\,\,\,\text{and then}\,\,1+A_1=1+0=1).
\end{split}    
\end{equation}
Finally, since $\alpha_k +A_{k-1}=A_k$,  we have 
\begin{equation*} 
    \begin{split}
        & \bigg\| \bigg[  \prod_{l =1}^{\alpha_{k}} \bigg( \sum_{j_{l+A_{k-1}} \in N_{\bar{\lambda}}(S)^c}  \frac{u_{j_{l + A_{k-1}}} u_{j_{l + A_{k-1}}}^T}{\displaystyle\prod_{\substack{e \in \mathcal{E}_L;\, e^{-} = j_{l+A_{k-1}}}} (\lambda_{e^{+}}- \lambda_{j_{l+A_{k-1}}})} E \bigg) \bigg] \times \left[ \prod_{l =1}^{\beta_k-1} \left( u_{i_{l+B_{k-1}}} u_{i_{l+B_{k-1}}}^T E \right)  \right] u_{i_{s_2}} u_{i_{s_2}}^T \textbf{w} \bigg\| \\
        & = \bigg\| \bigg[  \prod_{l =1}^{\alpha_{k}-1} \bigg( \sum_{j_{l+A_{k-1}} \in N_{\bar{\lambda}}(S)^c}  \frac{u_{j_{l + A_{k-1}}} u_{j_{l + A_{k-1}}}^T}{\displaystyle\prod_{\substack{e \in \mathcal{E}_L \\ e^{-} = j_{l+A_{k-1}}}} (\lambda_{e^{+}}- \lambda_{j_{l+A_{k-1}}})} E \bigg) \bigg] \times  \bigg( \sum_{j_{A_{k}} \in N_{\bar{\lambda}}(S)^c}  \frac{u_{j_{ A_{k}}} u_{j_{ A_{k}}}^T}{\displaystyle\prod_{\substack{e \in \mathcal{E}_L \\ e^{-} = j_{A_{k}}}} (\lambda_{e^{+}}- \lambda_{j_{A_{k}}})} E \bigg) \\
        & \times  \left[ \prod_{l =1}^{\beta_k-1} \left( u_{i_{l+B_{k-1}}} u_{i_{l+B_{k-1}}}^T E \right)  \right] \times u_{i_{s_2}}  u_{i_{s_2}}^T \textbf{w} \bigg\| . 
         \end{split}
\end{equation*}
We can rearrange the RHS into
\begin{equation*}
    \begin{split}
  &  \bigg\| \bigg[  \prod_{l =1}^{\alpha_{k}-1} \bigg( \sum_{j_{l+A_{k-1}} \in N_{\bar{\lambda}}(S)^c}  \frac{u_{j_{l + A_{k-1}}} u_{j_{l + A_{k-1}}}^T}{\displaystyle\prod_{\substack{e \in \mathcal{E}_L \\ e^{-} = j_{l+A_{k-1}}}} (\lambda_{e^{+}}- \lambda_{j_{l+A_{k-1}}})} E \bigg) \bigg] \times  \bigg( \sum_{j_{A_{k}} \in N_{\bar{\lambda}}(S)^c}  \frac{u_{j_{ A_{k}}} u_{j_{ A_{k}}}^T}{\displaystyle\prod_{\substack{e \in \mathcal{E}_L \\ e^{-} = j_{A_{k}}}} (\lambda_{e^{+}}- \lambda_{j_{A_{k}}})} \bigg) \\
        & \times \bigg\{ E \bigg[ \prod_{l =1}^{\beta_k-1} \left( u_{i_{l+B_{k-1}}} u_{i_{l+B_{k-1}}}^T E \right)  \bigg]  u_{i_{s_2}} \bigg\} \times  (u_{i_{s_2}}^T \textbf{w}) \bigg\|,\,\,\text{which is}
     \end{split}
\end{equation*}
\begin{equation*}
    \begin{split}
        & \bigg\| \bigg[  \prod_{l =1}^{\alpha_{k}-1} \bigg( \sum_{j_{l+A_{k-1}} \in N_{\bar{\lambda}}(S)^c}  \frac{u_{j_{l + A_{k-1}}} u_{j_{l + A_{k-1}}}^T}{\displaystyle\prod_{\substack{e \in \mathcal{E}_L \\ e^{-} = j_{l+A_{k-1}}}} (\lambda_{e^{+}}- \lambda_{j_{l+A_{k-1}}})} E \bigg) \bigg] \times  \bigg( \sum_{j_{A_{k}} \in N_{\bar{\lambda}}(S)^c}  \frac{u_{j_{ A_{k}}} u_{j_{ A_{k}}}^T}{\displaystyle\prod_{\substack{e \in \mathcal{E}_L \\ e^{-} = j_{A_{k}}}} (\lambda_{e^{+}}- \lambda_{j_{A_{k}}})} \bigg) \\
        & \times  (E u_{i_{1+B_{k-1}}}) \times \bigg[ \prod_{l =1}^{\beta_k-1} \bigg( u_{i_{l+B_{k-1}}}^T E u_{i_{l+1+B_{k-1}}}  \bigg)  \bigg]   \times  (u_{i_{s_2}}^T \textbf{w}) \bigg\|.
    \end{split}
\end{equation*}
By the sub-multiplicative property of the spectral norm, this is at most
\begin{equation*}
    \begin{split}
        & \bigg[ \prod_{l =1}^{\alpha_{k}-1} \bigg\| \sum_{j_{l+A_{k-1}} \in N_{\bar{\lambda}}(S)^c}  \frac{u_{j_{l + A_{k-1}}} u_{j_{l + A_{k-1}}}^T}{\displaystyle\prod_{\substack{e \in \mathcal{E}_L \\ e^{-} = j_{l+A_{k-1}}}} (\lambda_{e^{+}}- \lambda_{j_{l+A_{k-1}}})}\bigg\| \| E\| \bigg]  \times  \bigg\| \sum_{j_{A_{k}} \in N_{\bar{\lambda}}(S)^c}  \frac{u_{j_{ A_{k}}} u_{j_{ A_{k}}}^T}{\displaystyle\prod_{\substack{e \in \mathcal{E}_L \\ e^{-} = j_{A_{k}}}} (\lambda_{e^{+}}- \lambda_{j_{A_{k}}})} \bigg\|  \\
        & \times  \Norm{E u_{i_{1+B_{k-1}}}} \times \bigg[ \prod_{l =1}^{\beta_k-1} \norm{ u_{i_{l+B_{k-1}}}^T E u_{i_{l+1+B_{k-1}}} }  \bigg]   \times  \norm{u_{i_{s_2}}^T \textbf{w}}.
        \end{split}
\end{equation*}
By \eqref{u_ju_j^Tsumbound} and Definition \ref{xyz}, the long product above is at most
$$ \left[ \|E\|^{(\alpha_{k}-1)} \left(\prod_{1+ A_{k-1} \leq l \leq A_{k} } \frac{1}{\bar{\lambda}^{d(j_l)}} \right)  w x^{(\beta_{k}-1)} \right] \times \norm{u_{i_{s_2}}^T \textbf{w}}  = \frac{w x^{\beta_k-1} \|E\|^{\alpha_k-1} \times \norm{u_{i_{s_2}}^T \textbf{w}}  }{\bar{\lambda}^{\sum_{l=1+A_{k-1}}^{A_k} d(j_l)}}.$$
  It implies 
 \begin{equation} \label{Factor4,case2}
 \begin{split}
& \bigg\| \bigg[  \prod_{l =1}^{\alpha_{k}} \bigg( \sum_{j_{l+A_{k-1}} \in N_{\bar{\lambda}}(S)^c}  \frac{u_{j_{l + A_{k-1}}} u_{j_{l + A_{k-1}}}^T}{\displaystyle\prod_{\substack{e \in \mathcal{E}_L;\, e^{-} = j_{l+A_{k-1}}}} (\lambda_{e^{+}}- \lambda_{j_{l+A_{k-1}}})} E \bigg) \bigg] \times \left[ \prod_{l =1}^{\beta_k-1} \left( u_{i_{l+B_{k-1}}} u_{i_{l+B_{k-1}}}^T E \right)  \right] u_{i_{s_2}} u_{i_{s_2}}^T \textbf{w}\bigg\|\\
& \leq \frac{w x^{\beta_k-1} \|E\|^{\alpha_k-1} \times \norm{u_{i_{s_2}}^T \textbf{w}}  }{\bar{\lambda}^{\sum_{l=1+A_{k-1}}^{A_k} d(j_l)}}.
 \end{split}     
 \end{equation}

Combining \eqref{Sigma2splitstep0}, \eqref{Factor1,case2}, \eqref{Factor2,case2}, \eqref{Factor3,case2}, and \eqref{Factor4,case2}, we obtain 
\begin{equation} \label{Mstep1bound}
\begin{split}
 \Norm{\Sigma_2} & \leq \max_{\|\textbf{v}\|=\|\textbf{w}\|=1} \sum_{i_1,...,i_{s_2}} \sum_{\substack{L \in \mathcal{L}(T,s)\\ r_c(T,L) > 0}}  \norm{\textbf{v}^T u_{i_1}}\norm{u_{i_{s_2}}^T \textbf{w}} \frac{w^{2(k-1)}  x^{\sum_{l=1}^k (\beta_l-1)} \|E\|^{\sum_{l=2}^{k} (\alpha_l-1)}}{\delta_S^{|E(I_2)|} \bar{\lambda}^{\sum_{l=1}^{s_1} d(j_l)}}  \\
  & = \max_{\|\textbf{v}\|=\|\textbf{w}\|=1}  \sum_{i_1,...,i_{s_2}} \sum_{\substack{L \in \mathcal{L}(T,s) \\ r_c(T,L) > 0}}  \norm{\textbf{v}^T u_{i_1}}\norm{u_{i_{s_2}}^T \textbf{w}} \frac{\|E\|^{s_1-(k-1)} w^{2(k-1)} x^{s_2 - k}}{\delta_S^{s_2 -1-r_c(T,L)} \bar{\lambda}^{s_1+r_c(T,L)}} \, (\text{by \eqref{GraphEdgeComputation}}) \\
  & = \max_{\|\textbf{v}\|=\|\textbf{w}\|=1} \sum_{T=2}^{s_2} \sum_{\substack{L \in \mathcal{L}(T,s) \\ r_c(T,L) > 0}}  \bigg[ \sum_{\substack{i_1,...,i_{s_2} \\ |X(i_1,...,i_{s_2})|=T}} \norm{\textbf{v}^T u_{i_1}}\norm{u_{i_{s_2}}^T \textbf{w}} \bigg]  \frac{\|E\|^{s_1-(k-1)} w^{2(k-1)} x^{s_2 - k}}{\delta_S^{s_2 -1-r_c(T,L)} \bar{\lambda}^{s_1+r_c(T,L)}},
  \end{split}
     \end{equation}
here recall that $X(i_1,...,i_{s_2})$ is the sub-sequence of $[i_1,...,i_{s_2}]$ whose elements are in $S$. Note that in this case $r_c(T, L) \geq 1$, by Lemma \ref{r_cLemma}, $T \geq 2$ (that is why our sum runs from $T=2$ to $T=s_2$).

\noindent Moreover, 
\begin{equation*}
    \begin{split}
 \bigg[ \sum_{\substack{i_1,...,i_{s_2} \\ |X(i_1,...,i_{s_2})|=T}} \norm{\textbf{v}^T u_{i_1}}\norm{u_{i_{s_2}}^T \textbf{w}}\bigg] & =   \bigg[ \sum_{\substack{i_1,...,i_{s_2} \\ i_1, i_{s_2} \in S \\ |X(i_1,...,i_{s_2})|=T}} \norm{\textbf{v}^T u_{i_1}}\norm{u_{i_{s_2}}^T \textbf{w}} \bigg] + \bigg[ \sum_{\substack{i_1,...,i_{s_2} \\ i_1 \in S, i_{s_2} \notin S\\ |X(i_1,...,i_{s_2})|=T}} \norm{\textbf{v}^T u_{i_1}}\norm{u_{i_{s_2}}^T \textbf{w}} \bigg] \\
 & + \bigg[ \sum_{\substack{i_1,...,i_{s_2} \\ i_1 \notin S, i_{s_2} \in S \\ |X(i_1,...,i_{s_2})|=T}} \norm{\textbf{v}^T u_{i_1}}\norm{u_{i_{s_2}}^T \textbf{w}} \bigg] + \bigg[ \sum_{\substack{i_1,...,i_{s_2} \\ i_1, i_{s_2}\notin S \\ |X(i_1,...,i_{s_2})|=T}} \norm{\textbf{v}^T u_{i_1}}\norm{u_{i_{s_2}}^T \textbf{w}}\bigg]. 
   \end{split}
\end{equation*}
When $i_1, i_{s_2} \in S$ and $|X(i_1,...,i_{s_2})|=T \geq 2$, the choices for $[i_2,...,i_{s_2-1}]$ is at most $ \binom{s_2 -2}{T-2} p^{T-2} r^{s_2 -T}$.  Thus, by the Cauchy--Schwarz inequality, 
\begin{equation} \label{vUUwboundArgument}
    \begin{split}
  \sum_{\substack{i_1,...,i_{s_2} \\ i_1, i_{s_2} \in S \\ |X(i_1,...,i_{s_2})|=T}} \norm{\textbf{v}^T u_{i_1}}\norm{u_{i_{s_2}}^T \textbf{w}}  & \leq  \binom{s_2-2}{T-2} p^{T-2} r^{s_2 -T} \times \bigg[\sum_{i_1 \in S} \norm{\textbf{v}^T u_{i_1}} \bigg] \bigg[ \sum_{i_{s_2} \in S} \norm{u_{i_{s_2}}^T \textbf{w}} \bigg] \\
   & \leq \binom{s_2-2}{T-2}    p^{T-2} r^{s_2 -T} \times \sqrt{p} \sqrt{p} = \binom{s_2-2}{T-2}  p^{T-1} r^{s_2 -T}.     
    \end{split}
\end{equation}
By a similar argument, we also have 
\begin{equation} \label{vWSquareArgumentCase2}
 \sum_{\substack{i_1,...,i_{s_2} \\ i_1 \in S, i_{s_2} \notin S\\ |X(i_1,...,i_{s_2})|=T}} \norm{\textbf{v}^T u_{i_1}}\norm{u_{i_{s_2}}^T \textbf{w}}  \leq \binom{s_2-2}{T-1} p^{T-1} r^{s_2-T-1} \sqrt{pr}.
\end{equation}
\begin{equation} \label{vWSquareArgumentCase3}
 \sum_{\substack{i_1,...,i_{s_2} \\ i_1 \notin S, i_{s_2} \in S \\ |X(i_1,...,i_{s_2})|=T}} \norm{\textbf{v}^T u_{i_1}}\norm{u_{i_{s_2}}^T \textbf{w}}  \leq \binom{s_2-2}{T-1} p^{T-1} r^{s_2-T-1} \sqrt{pr}.
\end{equation}
\begin{equation}  \label{vWSquareArgumentCase4}
 \sum_{\substack{i_1,...,i_{s_2} \\ i_1, i_{s_2}\notin S \\ |X(i_1,...,i_{s_2})|=T}} \norm{\textbf{v}^T u_{i_1}}\norm{u_{i_{s_2}}^T \textbf{w}} \leq \binom{s-2}{T} p^T r^{s-2-T} \sqrt{r} \sqrt{r}= \binom{s_2-2}{T} p^T r^{s_2-1-T}. 
\end{equation}
\noindent Combining \eqref{vUUwboundArgument}, \eqref{vWSquareArgumentCase2}, \eqref{vWSquareArgumentCase3}, and \eqref{vWSquareArgumentCase4}, we obtain
\begin{equation} \label{vWSquareArgumentT>1}
\begin{split}
 \sum_{\substack{i_1,...,i_{s_2} \\ |X(i_1,...,i_{s_2})|=T}} \norm{\textbf{v}^T u_{i_1}}\norm{u_{i_{s_2}}^T \textbf{w}} &  \leq \binom{s_2-2}{T-2}  p^{T-1} r^{s_2 -T} + 2 \binom{s_2-2}{T-1} \sqrt{pr} p^{T-1} r^{s_2-T-1} + \\
 & +\binom{s_2-2}{T}  p^T r^{s_2 - T-1} \\
 & \leq \binom{s_2}{T} p^{T-1} r^{s_2 -T}.    
\end{split}   
\end{equation}
 Combining \eqref{Mstep1bound} and \eqref{vWSquareArgumentT>1}, we have 
\begin{equation} \label{splitSigma2M(a,b)}
\begin{split}
\Norm{\Sigma_2} & \leq \sum_{T=2}^{s_2} \binom{s_2}{T} p^{T-1} r^{s_2 -T} \sum_{L \in \mathcal{L}(T,s);\, r_c(T,L) > 0} \frac{ \|E\|^{s_1-(k-1)} w^{2(k-1)} x^{s_2 - k}}{\delta_S^{s_2 -1-r_c(T,L)} \bar{\lambda}^{s_1+r_c(T,L)}}  \\
& = M_1(\alpha; \beta) + M_2(\alpha;\beta), \, \text{where}
\end{split}
\end{equation}
$$M_1(\alpha;\beta) := \sum_{T=2}^{s_2} \binom{s_2}{T} p^{T-1} r^{s_2 -T} \sum_{\substack{L \in \mathcal{L}(T,s);\, 1\leq r_c(T,L) \leq k-1}} \frac{\|E\|^{s_1-(k-1)} w^{2(k-1)} x^{s_2 - k} }{\delta_S^{s_2 -1-r_c(T,L)} \bar{\lambda}^{s_1+r_c(T,L)}}, $$
$$M_2(\alpha;\beta) := \sum_{T=2}^{s_2} \binom{s_2}{T} p^{T-1} r^{s_2 -T}\sum_{\substack{L \in \mathcal{L}(T,s);\, r_c(T,L) > k-1}} \frac{\|E\|^{s_1-(k-1)} w^{2(k-1)} x^{s_2 - k} }{\delta_S^{s_2 -1-r_c(T,L)} \bar{\lambda}^{s_1+r_c(T,L)}}.$$

Similar to the arguments of bounding $M_1(\alpha;\beta), M_2(\alpha;\beta)$ for $(\alpha;\beta)$ of Type I from the previous subsection, we obtain the following lemmas.
\begin{lemma} \label{M1(a,b)case2part1} For any pair $(\alpha;\beta)$ of Type II,
    \begin{equation} 
  M_1(\alpha; \beta) \leq \left( \frac{ \sqrt{p}\|E\|}{\bar{\lambda}} \right) \frac{\sum_{T=2}^{s_2}   \displaystyle \sum_{\substack{L \in \mathcal{L}(T,s);\, 1 \leq r_c(T,L) \leq k-1}}\binom{s_2}{T}}{12^{s-1}}. 
\end{equation}
\end{lemma}
\begin{lemma}\label{M2(a,b)case2part1} For any pair $(\alpha;\beta)$ of Type II,
    \begin{equation}
\begin{split} 
M_2(\alpha;\beta) \leq \frac{p x}{\bar{\lambda}} \frac{ \sum_{T=2}^{s_2}  \displaystyle\sum_{\substack{L \in \mathcal{L}(T,s);\,r_c > k-1}} \binom{s_2}{T} }{12^{s-1}}. 
\end{split}
\end{equation}
\end{lemma}
\noindent The proofs of these lemmas are deferred to Appendix \ref{section: sublemmas}. Combining \eqref{splitSigma2M(a,b)} with Lemma \ref{M1(a,b)case2part1} and Lemma \ref{M2(a,b)case2part1}, we finally obtain 
\begin{equation*}
    \begin{split}
\Norm{\Sigma_2} & \leq M_1(\alpha;\beta) + M_2(\alpha;\beta) \\
& \leq \frac{p x}{\bar{\lambda}} \frac{ \sum_{T=2}^{s_2}  \displaystyle\sum_{\substack{L \in \mathcal{L}(T,s);\, r_c > k-1}} \binom{s_2}{T} }{12^{s-1}} + \left( \frac{ \sqrt{p}\|E\|}{\bar{\lambda}} \right) \frac{\sum_{T=2}^{s_2}    \displaystyle \sum_{\substack{L \in \mathcal{L}(T,s);\, 1 \leq r_c(T,L) \leq k-1}}\binom{s_2}{T}}{12^{s-1}}.
    \end{split}
\end{equation*}
This proves Lemma \ref{lemma: Sigma2allcases}.

\subsubsection{Proof of Lemma \ref{lemma: Sigma1allcases} - first part} Recall that in this case, $(s_2,s_1) \neq (k,k-1)$. Together with definition of $(s_2,s_1)$, we have 
\begin{equation}\label{s_2s_1condition}
    \max \{s_2 -k, s_1-(k-1)\} \geq 1.
\end{equation}
Back to our estimation of $\Norm{\Sigma_1}$. First, repeating the argument on bounding $\Norm{\Sigma_2}$ up to \eqref{Mstep1bound}, we also obtain
\begin{equation} \label{Sigma_1Step1}
    \begin{split}
 \Norm{\Sigma_1} \leq   \max_{\|\textbf{v}\|=\|\textbf{w}\|=1} \sum_{T=1}^{s_2} \sum_{\substack{L \in \mathcal{L}(T,s) \\ r_c(T,L) = 0}}  \bigg[ \sum_{\substack{i_1,...,i_{s_2} \\ |X(i_1,...,i_{s_2})|=T}} \norm{\textbf{v}^T u_{i_1}}\norm{u_{i_{s_2}}^T \textbf{w}} \bigg]  \frac{\|E\|^{s_1-(k-1)} w^{2(k-1)} x^{s_2 - k}}{\delta_S^{s_2 -1-r_c(T,L)} \bar{\lambda}^{s_1+r_c(T,L)}}.      
    \end{split}
\end{equation}
Similar to \eqref{vWSquareArgumentT>1}, if $T \geq 2$, 
\begin{equation*}
 \sum_{\substack{i_1,...,i_{s_2} \\ |X(i_1,...,i_{s_2})|=T}} \norm{\textbf{v}^T u_{i_1}}\norm{u_{i_{s_2}}^T \textbf{w}} \leq \binom{s_2}{T} p^{T-1} r^{s_2-T}.   
\end{equation*}
If $T=1$, the first case $i_1,i_{s_2} \in S$ cannot occur, combining \eqref{vWSquareArgumentCase2}, \eqref{vWSquareArgumentCase3}, and \eqref{vWSquareArgumentCase4}, we obtain
\begin{equation} \label{vWSquareArgumentT=1}
\begin{split}
 \sum_{\substack{i_1,...,i_{s_2} \\ |X(i_1,...,i_{s_2})|=1}} \norm{\textbf{v}^T u_{i_1}}\norm{u_{i_{s_2}}^T \textbf{w}} & = \bigg[ \sum_{\substack{i_1,...,i_{s_2} \\ i_1 \in S, i_{s_2} \notin S\\ |X(i_1,...,i_{s_2})|=1}} \norm{\textbf{v}^T u_{i_1}}\norm{u_{i_{s_2}}^T \textbf{w}} \bigg] \\
 & + \bigg[ \sum_{\substack{i_1,...,i_{s_2} \\ i_1 \notin S, i_{s_2} \in S \\ |X(i_1,...,i_{s_2})|=1}} \norm{\textbf{v}^T u_{i_1}}\norm{u_{i_{s_2}}^T \textbf{w}} \bigg] + \bigg[ \sum_{\substack{i_1,...,i_{s_2} \\ i_1, i_{s_2}\notin S \\ |X(i_1,...,i_{s_2})|=1}} \norm{\textbf{v}^T u_{i_1}}\norm{u_{i_{s_2}}^T \textbf{w}}\bigg] \\
 & \leq 2 \binom{s_2-2}{0} \sqrt{pr} r^{s_2-2}+ \binom{s_2-2}{1} p r^{s_2-2}  \\
 & \leq s_2 \sqrt{pr} r^{s_2-2}. 
\end{split}      
\end{equation}
Combining \eqref{Sigma_1Step1}, \eqref{vWSquareArgumentT>1}, and \eqref{vWSquareArgumentT=1}, we obtain 
\begin{equation} \label{Sigma1split}
    \begin{split}
\Norm{\Sigma_1} & \leq M_3(\alpha;\beta) + M_4(\alpha;\beta), \, \text{where}\\
& M_3(\alpha;\beta):= \sum_{L \in \mathcal{L}(1,s)} s_2 \sqrt{pr} r^{s_2-2}  \frac{ \|E\|^{s_1-(k-1)} w^{2(k-1)} x^{s_2 - k}}{\delta_S^{s_2 -1} \bar{\lambda}^{s_1}} \,\,\,(T=1, r_c(T,L)=0), \\
& M_4(\alpha;\beta):= \sum_{T=2}^{s_2} \sum_{\substack{L \in \mathcal{L}(T,s) \\ r_c(T,L) = 0}} \binom{s_2}{T} p^{T-1} r^{s_2-T}    \frac{\|E\|^{s_1-(k-1)} w^{2(k-1)} x^{s_2 - k}}{\delta_S^{s_2 -1} \bar{\lambda}^{s_1}}.
    \end{split}
\end{equation}

The remaining task is to bound $M_3(\alpha;\beta)$ and $M_4(\alpha;\beta)$. First, we can rewrite $M_3(\alpha;\beta)$ as 
$$  M_3(\alpha;\beta) =  s_2 \sum_{L \in \mathcal{L}(1,s)} \sqrt{\frac{p}{r}} \left(\frac{rx}{\delta_S} \right)^{s_2 -k} \left( \frac{\|E\|}{\bar{\lambda}} \right)^{s_1-(k-1)} \left( \frac{ \sqrt{r} w}{\sqrt{\bar{\lambda} \delta_S}} \right)^{2(k-1)}.$$
Notice that by definition of $(s_1,s_2)$,  $s_2 \geq k, s_1 \geq k-1$, all the exponents are nonnegative. Applying Lemma \ref{toytrick} under condition \eqref{s_2s_1condition}, we have
\begin{equation}\label{Case2M_3}
\begin{split}
    M_3(\alpha;\beta) & \leq  s_2 \sum_{L \in \mathcal{L}(1,s)} \sqrt{\frac{p}{r}} \left( \frac{rx}{\delta_S} + \frac{\|E\|}{\bar{\lambda}} \right)  \max \left\lbrace \frac{\|E\|}{\bar{\lambda}}, \frac{rx}{\delta_S}, \frac{\sqrt{r}w}{ \sqrt{\bar{\lambda} \delta_S}} \right\rbrace^{s_1+s_2-2}  \\
    & = s_2 \sum_{L \in \mathcal{L}(1,s)} \sqrt{\frac{p}{r}} \left( \frac{rx}{\delta_S} + \frac{\|E\|}{\bar{\lambda}} \right)  \max \left\lbrace \frac{\|E\|}{\bar{\lambda}}, \frac{rx}{\delta_S}, \frac{\sqrt{r}w}{ \sqrt{\bar{\lambda} \delta_S}} \right\rbrace^{s-1} \,\,(\text{since}\,\, s_1+s_2=s+1) \\
    & = s_2 \sum_{L \in \mathcal{L}(1,s)}  \left(  \frac{\sqrt{pr} x}{\delta_S} + \frac{\sqrt{p} \|E\|}{\sqrt{r} \bar{\lambda}} \right)  \max \left\lbrace \frac{\|E\|}{\bar{\lambda}}, \frac{rx}{\delta_S}, \frac{\sqrt{r}w}{ \sqrt{\bar{\lambda} \delta_S}} \right\rbrace^{s-1}. 
    \end{split}    
\end{equation}
By Assumption \textbf{D0}, the RHS is at most
$$ s_2 \sum_{L \in \mathcal{L}(1,s)} \left( \frac{\sqrt{pr} x}{\delta_S} + \frac{\sqrt{p} \|E\|}{\sqrt{r} \bar{\lambda}} \right) \frac{1}{12^{s-1}}  \leq s_2 \sum_{L \in \mathcal{L}(1,s)} \left( \frac{\sqrt{pr} x}{\delta_S} + \frac{\sqrt{p} \|E\|}{ \bar{\lambda}} \right) \frac{1}{12^{s-1}}.$$

Next, $M_4(\alpha;\beta)$ can be rewritten as 
$$\sum_{T = 2}^{s_2} \binom{s_2}{T} \sum_{L \in \mathcal{L}(T,s);\, r_c(T,L)=0}  \left( \frac{p}{r} \right)^{T-1} \left(\frac{rx}{\delta_S} \right)^{s_2 -k} \left( \frac{\|E\|}{\bar{\lambda}} \right)^{s_1-(k-1)} \left( \frac{ \sqrt{r} w}{\sqrt{\bar{\lambda} \delta_S}} \right)^{2(k-1)}.$$
By the definition of $M_4(\alpha;\beta)$, $T \geq 2$ and hence $\left( \frac{p}{r} \right)^{T-1} \leq \frac{p}{r}$. Hence,
\begin{equation} \label{M4removed(p/r)}
 M_4(\alpha;\beta) \leq \sum_{T = 2}^{s_2} \binom{s_2}{T} \sum_{L \in \mathcal{L}(T,s);\, r_c(T,L)=0} \frac{p}{r} \left(\frac{rx}{\delta_S} \right)^{s_2 -k} \left( \frac{\|E\|}{\bar{\lambda}} \right)^{s_1-(k-1)} \left( \frac{ \sqrt{r} w}{\sqrt{\bar{\lambda} \delta_S}} \right)^{2(k-1)}. 
\end{equation}
Note that by the definition of $(s_1,s_2)$, all exponents are non-negative. Applying Lemma \ref{toytrick} under condition \eqref{s_2s_1condition}, we have
\begin{equation} \label{toytrickApplicaiton}
\begin{split}
   \left(\frac{rx}{\delta_S} \right)^{s_2 -k} \left( \frac{\|E\|}{\bar{\lambda}} \right)^{s_1-(k-1)} \left( \frac{ \sqrt{r} w}{\sqrt{\bar{\lambda} \delta_S}} \right)^{2(k-1)} & \leq \left( \frac{\|E\|}{\bar{\lambda}} +  \frac{r x}{\delta_S} \right) \max \left\lbrace \frac{\|E\|}{\bar{\lambda}}, \frac{rx}{\delta_S}, \frac{\sqrt{r}w}{ \sqrt{\bar{\lambda} \delta_S}} \right\rbrace^{s_1+s_2-2} \\
   &= \left( \frac{\|E\|}{\bar{\lambda}} +  \frac{r x}{\delta_S} \right) \max \left\lbrace \frac{\|E\|}{\bar{\lambda}}, \frac{rx}{\delta_S}, \frac{\sqrt{r}w}{ \sqrt{\bar{\lambda} \delta_S}} \right\rbrace^{s-1}. 
\end{split}
\end{equation}
Combining \eqref{M4removed(p/r)} and \eqref{toytrickApplicaiton}, we have
\begin{equation} \label{M4(a,b)}
\begin{split}
 M_4(\alpha;\beta)   & \leq   \sum_{T = 2}^{s_2} \binom{s_2}{T} \sum_{L \in \mathcal{L}(T,s);\, r_c(T,L)=0}  \frac{p}{r} \left( \frac{\|E\|}{\bar{\lambda}} +  \frac{r x}{\delta_S} \right) \max \left\lbrace \frac{\|E\|}{\bar{\lambda}}, \frac{rx}{\delta_S}, \frac{\sqrt{r}w}{ \sqrt{\bar{\lambda} \delta_S}} \right\rbrace^{s-1}    \\
 & =  \sum_{T = 2}^{s_2} \binom{s_2}{T} \sum_{L \in \mathcal{L}(T,s);\, r_c(T,L)=0}   \left( \frac{p\|E\|}{r \bar{\lambda}} +  \frac{p x}{\delta_S} \right) \max \left\lbrace \frac{\|E\|}{\bar{\lambda}}, \frac{rx}{\delta_S}, \frac{\sqrt{r}w}{ \sqrt{\bar{\lambda} \delta_S}} \right\rbrace^{s-1} \\
  &\leq   \left( \frac{px}{\delta_S} + \frac{p\|E\|}{r\bar{\lambda}} \right) \frac{\sum_{T = 2}^{s_2}  \displaystyle\sum_{L \in \mathcal{L}(T,s);\, r_c(T,L)=0} \binom{s_2}{T} }{12^{s-1}} (\text{by Assumption \textbf{D0}}).
\end{split}    
\end{equation}

 Together \eqref{Sigma1split}, \eqref{Case2M_3} and \eqref{M4(a,b)} imply that if $(s_2,s_1) \neq (k,k-1)$, then 
\begin{equation} \label{goodcasebound-HR}
\begin{split}
 \Norm{\Sigma_1} &\leq M_3(\alpha;\beta)+M_4(\alpha;\beta) \\ 
& \leq s_2 \sum_{L \in \mathcal{L}(1,s)} \left( \frac{\sqrt{pr} x}{\delta_S} + \frac{\sqrt{p} \|E\|}{ \bar{\lambda}} \right) \frac{1}{12^{s-1}} + \left( \frac{px}{\delta_S} + \frac{p\|E\|}{r\bar{\lambda}} \right) \frac{\sum_{T = 2}^{s_2}  \displaystyle\sum_{L \in \mathcal{L}(T,s);\, r_c(T,L)=0} \binom{s_2}{T} }{12^{s-1}} 
\end{split}
\end{equation} 
\begin{equation*}
    \begin{split}
    & \leq s_2 \sum_{L \in \mathcal{L}(1,s)} \left( \frac{\sqrt{pr} x}{\delta_S} + \frac{\sqrt{p} \|E\|}{ \bar{\lambda}} \right) \frac{1}{12^{s-1}} + \left( \frac{ \sqrt{pr} x}{\delta_S} + \frac{\sqrt{p} \|E\|}{\bar{\lambda}} \right) \frac{\sum_{T = 2}^{s_2}  \displaystyle\sum_{L \in \mathcal{L}(T,s);\, r_c(T,L)=0} \binom{s_2}{T} }{12^{s-1}} \\
& = \left( \frac{\sqrt{pr} x}{\delta_S} + \frac{\sqrt{p} \|E\|}{ \bar{\lambda}} \right) \frac{\sum_{T = 1}^{s_2}  \displaystyle\sum_{L \in \mathcal{L}(T,s);\, r_c(T,L)=0} \binom{s_2}{T}}{12^{s-1}}.    
    \end{split}
\end{equation*}
This proves the first part of Lemma \ref{lemma: Sigma1allcases}.

\subsubsection{Proof of Lemma \ref{lemma: Sigma1allcases} - second part}
Now we consider the very special case when $(s_2,s_1) = (k,k-1)$. In this case,
\begin{equation}\label{alpha-beta-particularcase}
\begin{split}
     & M(\alpha,\beta)= PEQEPEQEPEQ \dots PEQEP, \\
     & \alpha_1=\alpha_{k+1}=0, \alpha_2=...=\alpha_k=\beta_1=...=\beta_k=1, \\
     & s= s_1+s_2 -1 = 2k-2.
\end{split}
   \end{equation}
 Again, write $P =\sum _{i \in N_{\bar{\lambda}}(S)} \frac{1}{ z -\lambda_i } u_i u_i^T $ and $Q = \sum_{j \in N_{\bar{\lambda}}(S)^c} \frac{u_j u_j^T}{z - \lambda_j}$. Multiplying out all the terms, we can express 
$PEQEPEQ \cdots PEQEP$ as 
$$  \sum_{i_1,i_2,...,i_k \in N_{\bar{\lambda}}(S)} \sum_{j_1,j_2,...,j_{k-1} \in N_{\bar{\lambda}}(S)^c} \frac{1}{ (z-\lambda_{i_k}) \prod_{l=1}^{k-1} [ (z-\lambda_{i_l})(z-\lambda_{j_l})]} u_{i_1} u_{i_k}^T \prod_{l=1}^{k-1} (u_{i_l}^T E u_{j_{l}} u_{j_{l}}^T E u_{i_{l+1}}). $$
Following the procedure above, we can expand $\frac{1}{2 \pi \textbf{i}}\int_{\Gamma} PEQEPEQ...PEQEP dz$ as 
\begin{equation} \label{PEQEPEQfirstequation}
 \begin{split}
   \sum_{i_1,i_2,...,i_k } \sum_{j_1,j_2,...,j_{k-1}} \sum_{L \in \mathcal{L}(T,2k-2)} (-1)^{T+1} w(X,Y|L) u_{i_1} \left[ \prod_{l=1}^{k-1} (u_{i_l}^T E u_{j_{l}} u_{j_{l}}^T E u_{i_{l+1}}) \right] u_{i_k}^T.
\end{split}   
\end{equation}
Splitting $ w(X,Y|L)=\prod_{e \in \mathcal{E}_L} \frac{1}{(\lambda_{e^+}-\lambda_{e^-})}= \left[\prod_{\substack{e \in \mathcal{E}_L\\ e^{-} \in I_2}} \frac{1}{(\lambda_{e^+}-\lambda_{e^-})} \right] \times \prod_{l=1}^{s_1} \left[\prod_{\substack{e \in \mathcal{E}_L \\ e^{-} = j_l}} \frac{1}{(\lambda_{e^+}-\lambda_{e^-})} \right]$ and distributing the sum over $j_1,...,j_{k-1}$ according to each position of $u_{j_l} u_{j_l}^T$ in the product,  \eqref{PEQEPEQfirstequation} becomes 
\begin{equation*}
    \sum_{i_1,i_2,...,i_{k}} \sum_{\substack{L \in \mathcal{L}(T,s)}}   \frac{(-1)^{T+1} }{ \displaystyle\prod_{\substack{e \in \mathcal{E}_L \\ e^{-} \in I_2}} (\lambda_{e^+} - \lambda_{e^-})} u_{i_1} \prod_{l=1}^{k-1} \bigg[ u_{i_l}^T E \bigg( \sum_{j_l \in N_{\bar{\lambda}}(S)^c} \frac{u_{j_l} u_{j_l}^T}{\displaystyle\prod_{ \substack{e \in \mathcal{E}_L \\ e^{-}=j_l}} (\lambda_{e^+} -\lambda_{j_l})} \bigg) E u_{i_{l+1}} \bigg] u_{i_k}^T.
\end{equation*}
Therefore, in this particular case, $\Sigma_1$ (defined in \eqref{Sigma1Sigma2split}) becomes
\begin{equation}
    \begin{split}
    \sum_{i_1,i_2,...,i_{k}} \sum_{\substack{L \in \mathcal{L}(T,s)\\ r_c(T,L)=0}}  \frac{(-1)^{T+1} }{ \displaystyle \prod_{\substack{e \in \mathcal{E}_L \\ e^{-} \in I_2}}(\lambda_{e^+} - \lambda_{e^-})} u_{i_1} \prod_{l=1}^{k-1} \bigg[ u_{i_l}^T E \bigg( \sum_{j_l \in N_{\bar{\lambda}}(S)^c} \frac{u_{j_l} u_{j_l}^T}{\displaystyle\prod_{ \substack{e \in \mathcal{E}_L \\ e^{-}=j_l}} (\lambda_{e^+} -\lambda_{j_l})} \bigg) E u_{i_{l+1}} \bigg] u_{i_k}^T.
    \end{split}
\end{equation}

The trivial bound for each factor $\bigg|u_{i_l}^T E \bigg( \sum_{j_l \in N_{\bar{\lambda}}(S)^c} \frac{u_{j_l} u_{j_l}^T}{\prod_{ \substack{e \in \mathcal{E}_L \\ e^{-}=j_l}} (\lambda_{e^+} -\lambda_{j_l})} \bigg) E u_{i_{l+1}} \bigg|$ is $\frac{w^2}{\bar{\lambda}}$. However, given a combinatorial profile with $r_{c}(T,L)=0$, we can bound at least one factor by $y$ (see Definition \ref{xyz}).

 Recall that $r_c(T,L):= \sum_{l=1}^{s_1} (d(j_l)-1)$. So the condition $r_c(T,L)=0$ implies that $d(j_l)=1$ for all $1 \leq l \leq k-1$. This means that $j_l$ is connected to a unique vertex $e^{+}_l$ in the graph $G(X,Y|L)$. Therefore, $\Sigma_1$ equals 
 $$\sum_{i_1,i_2,...,i_{k}} \sum_{\substack{L \in \mathcal{L}(T,s)\\ r_c(T,L)=0}}  \frac{(-1)^{T+1} }{\displaystyle\prod_{\substack{e \in \mathcal{E}_L \\ e^{-} \in I_2}}  (\lambda_{e^+} - \lambda_{e^-})} u_{i_1} \prod_{l=1}^{k-1} \left[ u_{i_l}^T E \left( \sum_{j_l \in N_{\bar{\lambda}}(S)^c} \frac{u_{j_l} u_{j_l}^T}{ (\lambda_{e^+_l} -\lambda_{j_l})} \right) E u_{i_{l+1}} \right] u_{i_k}^T. $$
 Similar to the argument of bounding $\Norm{\Sigma_1}$ from the previous case, we have
\begin{equation} \label{Sigma1bound1}
    \begin{split}
             \Norm{\Sigma_1}   &\leq \sum_{ \substack{i_1,i_2,...,i_k \\ L \in \mathcal{L}(T,2k-2)\\ r_c(T,L)=0 }} \max_{\|\textbf{v}\|=\|\textbf{w}\|=1} \norm{\textbf{v}^T u_{i_1}}\norm{u_{i_{k}}^T \textbf{w}} \frac{1}{\delta_S^{|E(I_2)|}} \prod_{l=1}^{k-1} \norm{ u_{i_l}^T E \left( \sum_{j_l \in N_{\bar{\lambda}}(S)^c}\frac{u_{j_l}u_{j_l}^T}{\lambda_{e^{+}_l} - \lambda_{j_l}}\right) E u_{i_{l+1}}} \\
        & \leq  \sum_{ \substack{i_1,i_2,...,i_k \\ L \in \mathcal{L}(T,2k-2)\\ r_c(T, L)=0 }} \max_{\|\textbf{v}\|=\|\textbf{w}\|=1} \norm{\textbf{v}^T u_{i_1}}\norm{u_{i_{k}}^T \textbf{w}} \frac{1}{\delta_S^{k-1}} \prod_{l=1}^{k-1} \norm{ u_{i_l}^T E \left( \sum_{j_l \in N_{\bar{\lambda}}(S)^c} \frac{u_{j_l}u_{j_l}^T}{\lambda_{e^{+}_l} - \lambda_{j_l}}\right) E u_{i_{l+1}}},
    \end{split}
\end{equation}
since  $|E(I_2)|=s_2-1-r_c(T,L)=k-1$.

Consider the special case when all the elements of $X$ are the same, namely,  $X=[i,i,...,i]$, with $k$ repetitions,  $Y=[j_1,...,j_{k-1}]$, and $L \in \mathcal{L}(T,2k-2)$. 
In this case, the graph $G(X,Y|L)$  is nontrivial if $ i \in S$ and  $T=k$. Since $|Y|=k-1 = \frac{|\mathcal{E}(X,Y|L)|}{2}$, there exists at least one vertex $j_l$ such that $d(j_l) \geq 2$, which implies that  $r_c(X,Y|L) > 0$.  Therefore, the RHS of \eqref{Sigma1bound1} does not contain any sequence $[i_1, i_2, ...,i_k] = [i,i,..,i]$. We call a sequence $[i_1, i_2, ...,i_k] \neq [i,i,..,i]$ "\textit{non-uniform}", and rewrite \eqref{Sigma1bound1} as
$$\Norm{\Sigma_1} \leq   \sum_{ \substack{ \textit{non-uniform} \\i_1,i_2,...,i_k \\ L \in \mathcal{L}(T,2k-1);\, r_c(k, L)=0 }}  \max_{\|\textbf{v}\|=\|\textbf{w}\|=1} \norm{\textbf{v}^T u_{i_1}}\norm{u_{i_{k}}^T \textbf{w}}  \frac{1}{\delta_S^{k-1}} \prod_{l=1}^{k-1} \norm{ u_{i_l}^T E \left( \sum_{j_l \in N_{\bar{\lambda}}(S)^c} \frac{u_{j_l}u_{j_l}^T}{\lambda_{e^{+}_l} - \lambda_{j_l}}\right) E u_{i_{l+1}}}.$$

Observe that for any $1 \leq l \leq k-1$, by Definition \ref{xyz},
$$ \norm{u_{i_l}^T E \left( \sum_{j_{l} \in N_{\bar{\lambda}}(S)^c} \frac{u_{j_{l}} u_{j_{l}}^T}{\lambda_{e^{+}_{l}}-\lambda_{j_{l}}} \right) E u_{i_{l+1}}} \leq \Norm{u_{i_l}^T E } \Norm{ \sum_{j_{l} \in N_{\bar{\lambda}}(S)^c} \frac{u_{j_{l}} u_{j_{l}}^T}{\lambda_{e^{+}_{l}}-\lambda_{j_{l}}}} \Norm{E u_{i_{l+1}}} \leq \frac{w^2}{\bar{\lambda}}.$$
Furthermore, since $[i_1,...,i_k] \neq [i,...,i]$, there exists at least one pair $i_{l'} \neq i_{l'+1}$, and hence 
$$ \norm{u_{i_l'}^T E \left( \sum_{j_{l'} \in N_{\bar{\lambda}}(S)^c} \frac{u_{j_{l'}} u_{j_{l'}}^T}{\lambda_{e^{+}_{l'}}-\lambda_{j_{l'}}} \right) E u_{i_{l'+1}}} \leq y \,\,\, (\text{by Definition \ref{xyz}}).$$
Together, these observations imply that  
\begin{equation}
    \prod_{l=1}^{k-1} \norm{ u_{i_l}^T E \left( \sum_{j_l \in N_{\bar{\lambda}}(S)^c} \frac{u_{j_l}u_{j_l}^T}{\lambda_{e^{+}_l} - \lambda_{j_l}}\right) E u_{i_{l+1}}} \leq \frac{y w^{2(k-2)}}{\bar{\lambda}^{k-2}}.
\end{equation}
\noindent Therefore,  \eqref{Sigma1bound1} implies 
$$\Norm{\Sigma_1} \leq   \sum_{ \substack{ \textit{non-uniform} \\ i_1,i_2,...,i_k   }} \sum_{ \substack{ L \in \mathcal{L}(T,2k-2)\\ r_c(k, L)=0}}  \max_{\|\textbf{v}\|=\|\textbf{w}\|=1} \norm{\textbf{v}^T u_{i_1}}\norm{u_{i_{k}}^T \textbf{w}}  \frac{1}{\delta_S^{k-1}} \frac{y w^{2(k-2)}}{\bar{\lambda}^{k-2}}. $$

\noindent By the combinatorial fact that 
$ | \mathcal{L}(T,2k-2)| = \binom{(2k-2)-1}{T-1} = \binom{2k-3}{T-1} \le 2^{2k-2}$ for $T \le 2k-1$,  the last bound on $\Norm {\Sigma_1} $ implies  
\begin{equation} \label{Sigma1badcaseFinalstep}
    \begin{split}
     \Norm{\Sigma_1} & \leq  \sum_{ \substack{\textit{non-uniform} \\i_1,i_2,...,i_k  \\  }} 2^{2k-2} \max_{\|\textbf{v}\|=\|\textbf{w}\|=1} \norm{\textbf{v}^T u_{i_1}}\norm{u_{i_{k}}^T \textbf{w}}  \frac{y w^{2(k-2)}}{\bar{\lambda}^{k-2} \delta_S^{k-1}} \\
     & =  \bigg(\sum_{ \substack{\textit{non-uniform} \\ i_1,i_2,...,i_k  \\  }}  \max_{\|\textbf{v}\|=\|\textbf{w}\|=1} \norm{\textbf{v}^T u_{i_1}}\norm{u_{i_{k}}^T \textbf{w}} \bigg) 2^{2k-2} \frac{y w^{2(k-2)}}{\bar{\lambda}^{k-2} \delta_S^{k-1}}. 
      \end{split}
\end{equation}
As there is at least one index $i_l$ such that $i_l \in S$, we have 
\begin{equation} \label{Sig1vWSquaresplit}
    \begin{split}
   \sum_{ \substack{\textit{non-uniform} \\ i_1,i_2,...,i_k  \\  }}  \max_{\|\textbf{v}\|=\|\textbf{w}\|=1} \norm{\textbf{v}^T u_{i_1}}\norm{u_{i_{k}}^T \textbf{w}}  & \leq \sum_{l=1}^{k} \sum_{ \substack{i_1,i_2,...,i_k  \\ i_l \in S  }}  \max_{\|\textbf{v}\|=\|\textbf{w}\|=1} \norm{\textbf{v}^T u_{i_1}}\norm{u_{i_{k}}^T \textbf{w}} \\  
   & = \sum_{l=2}^{k-1} \sum_{ \substack{i_1,i_2,...,i_k  \\ i_l \in S  }}  \max_{\|\textbf{v}\|=\|\textbf{w}\|=1} \norm{\textbf{v}^T u_{i_1}}\norm{u_{i_{k}}^T \textbf{w}} \\
   & +  \sum_{ \substack{i_1,i_2,...,i_k  \\ i_1 \in S  }}  \max_{\|\textbf{v}\|=\|\textbf{w}\|=1} \norm{\textbf{v}^T u_{i_1}}\norm{u_{i_{k}}^T \textbf{w}} \\
   & +  \sum_{ \substack{i_1,i_2,...,i_k  \\ i_k \in S  }}  \max_{\|\textbf{v}\|=\|\textbf{w}\|=1} \norm{\textbf{v}^T u_{i_1}}\norm{u_{i_{k}}^T \textbf{w}}. 
    \end{split}
\end{equation}
Similar to the argument in \eqref{vUUwboundArgument} and \eqref{vWSquareArgumentT=1}, we have 
\begin{equation} \label{Sig1vWSquare1}
\sum_{ \substack{i_1,i_2,...,i_k  \\ i_l \in S  }}  \max_{\|\textbf{v}\|=\|\textbf{w}\|=1} \norm{\textbf{v}^T u_{i_1}}\norm{u_{i_{k}}^T \textbf{w}} \leq  \sqrt{r} \sqrt{r} p r^{k-3} \leq \sqrt{pr} r^{k-2}, \,\,\text{for any}\,\,2 \leq l \leq k-1,  
\end{equation}
\begin{equation} \label{Sig1vWSquare2}
  \sum_{ \substack{i_1,i_2,...,i_k \\ i_1 \in S  }}  \max_{\|\textbf{v}\|=\|\textbf{w}\|=1} \norm{\textbf{v}^T u_{i_1}}\norm{u_{i_{k}}^T \textbf{w}} \leq \sqrt{p} \sqrt{r} r^{k-2} = \sqrt{pr} r^{k-2}, 
\end{equation}
\begin{equation} \label{Sig1vWSquare3}
   \sum_{ \substack{i_1,i_2,...,i_k \\ i_k \in S  }}  \max_{\|\textbf{v}\|=\|\textbf{w}\|=1} \norm{\textbf{v}^T u_{i_1}}\norm{u_{i_{k}}^T \textbf{w}} \leq  \sqrt{r} \sqrt{p} r^{k-2} = \sqrt{pr} r^{k-2}.
\end{equation}
Together \eqref{Sig1vWSquaresplit}, \eqref{Sig1vWSquare1}, \eqref{Sig1vWSquare2}, and \eqref{Sig1vWSquare3} imply that 
\begin{equation} \label{Sig1vWSquarebound}
  \sum_{ \substack{\textit{non-uniform} \\ i_1,i_2,...,i_k  \\  }}  \max_{\|\textbf{v}\|=\|\textbf{w}\|=1} \norm{\textbf{v}^T u_{i_1}}\norm{u_{i_{k}}^T \textbf{w}} \leq k \sqrt{pr} r^{k-2} = k \sqrt{pr} r^{k-2}.  
\end{equation}

Combining \eqref{Sigma1badcaseFinalstep} and \eqref{Sig1vWSquarebound}, we obtain 
 $$\Norm{\Sigma_1} \leq  k \sqrt{rp} r^{k-2} \times  2^{2k-2} \frac{y w^{2(k-2)}}{\bar{\lambda}^{k-2} \delta_S^{k-1}}. $$

As $k=s/2+1$, we can write this as 
\begin{equation*}
    \Norm{\Sigma_1} \leq \frac{(s+2) \sqrt{pr} y}{\delta_S} \left( \frac{2 \sqrt{r} w}{\sqrt{\bar{\lambda} \delta_S}} \right)^{s-2}.
\end{equation*}

  This proves the second part of Lemma \ref{lemma: Sigma1allcases}.

\subsection{Proof of Expansion \eqref{TaylorEx} } \label{TaylorExpansion} To conclude the proof of Theorem \ref{deterministicHRgeneral}, we still need to prove that (under Assumption \textbf{D0})  the series $\sum_{s=1}^{\infty} \|(z-A)^{-1} [E(z-A)^{-1}]^{s}\|$ converges. 

 Define $G_s:= \|(z-A)^{-1} [E(z-A)^{-1}]^{s}\|.$ By the triangle inequality, we have 
 \begin{equation} \label{F^0_s}
    G_s \leq \sum_{(\alpha;\beta)} \|M(\alpha;\beta)\|,
 \end{equation}
 in which $\alpha, \beta, M(\alpha;\beta)$ are as defined in Subsection \ref{F_2Es}.
\begin{lemma} \label{ConvergenceM(a,b)} For any pair $(\alpha,\beta)$,
\begin{equation}
   \Norm{M(\alpha,\beta)} \leq \frac{2r}{ 6^s \delta_S}. 
\end{equation}
    \end{lemma} 
Combining \eqref{F^0_s}, Lemma \ref{ConvergenceM(a,b)}, and the fact that there are exactly $2^{s+1}$ pairs $(\alpha,\beta)$ in the expansion of $G_s$, we obtain 
$$G_s \leq 2^{s+1} \times \frac{2r}{  6^s \delta_S} = \frac{4 r}{3^s \delta_S }.$$
Therefore, 
$$\sum_{s=1}^{\infty} \|(z-A)^{-1} [E(z-A)^{-1}]^{s}\| \leq \sum_{s=1}^{\infty} \frac{4 r}{ 3^s \delta_S} = \frac{2 r}{\delta_S }.$$ 
So, our remaining duty is to prove Lemma \ref{ConvergenceM(a,b)}. We first bound $\Norm{M(\alpha;\beta)}$ for $(\alpha;\beta)$ of Type I; the proofs for the other types are similar.  Multiplying out all the terms in the operators $P$ that appear in $M(\alpha;\beta)$ (there are $s_2$ such operators),  we decompose $M(\alpha;\beta)$ as

\begin{equation*}
\begin{split}
& \sum_{i_1,i_2,...,i_{s_2}} \frac{1}{\prod_{l=1}^{s_2}(z-\lambda_{i_l})} \left\lbrace \prod_{h=0}^{k-1} \left[ \prod_{l=1}^{\alpha_{h+1}} \left( Q E  \right) \right] \times \left[  \prod_{l=1}^{\beta_{h+1}} \left( u_{i_{l+B_{h}}} u_{i_{l+B_{h}}}^T E \right) \right] \right\rbrace \times \left[  \prod_{l =1}^{\alpha_{k+1}-1} \left( Q E \right) \right] Q \\
 & = \sum_{i_1,i_2,...,i_{s_2}} \frac{1}{\prod_{l=1}^{s_2}(z-\lambda_{i_l})} \left\lbrace \prod_{h=0}^{k-1} \left[ \prod_{l=1}^{\alpha_{h+1}-1} \left( Q E  \right) \right] QE \left[ \prod_{l=1}^{\beta_{h+1}} \left( u_{i_{l+B_{h}}} u_{i_{l+B_{h}}}^T E \right) \right] \right\rbrace \times \left[  \prod_{l =1}^{\alpha_{k+1}-1} \left( Q E \right) \right] Q \\
 & = \sum_{i_1,i_2,...,i_{s_2}} \frac{1}{\prod_{l=1}^{s_2}(z-\lambda_{i_l})} \left\lbrace \prod_{h=0}^{k-1} \left[ \prod_{l=1}^{\alpha_{h+1}-1} \left( Q E  \right) \right] Q \left[ E \prod_{l=1}^{\beta_{h+1}} \left( u_{i_{l+B_{h}}} u_{i_{l+B_{h}}}^T E \right) \right] \right\rbrace \times \left[  \prod_{l =1}^{\alpha_{k+1}-1} \left( Q E \right) \right] Q.
\end{split}
\end{equation*}
Thus, by the triangle inequality,  
 \begin{equation} \label{ConveM(a,b)form}
\begin{split}
&  \Norm{M(\alpha; \beta)}  \leq  \\
 & \sum_{i_1,i_2,...,i_{s_2}} \norm{\frac{1}{\prod_{l=1}^{s_2}(z-\lambda_{i_l})}} \times  \prod_{h=0}^{k-1} \left[ \prod_{l=1}^{\alpha_{h+1}-1} \| Q E\| \right] \|Q\| \Norm{ E \prod_{l=1}^{\beta_{h+1}} \left( u_{i_{l+B_{h}}} u_{i_{l+B_{h}}}^T E \right) }  \times \left[  \prod_{l =1}^{\alpha_{k+1}-1} \| Q E\| \right] \|Q\|.
\end{split}
\end{equation}
On the other hand,  Assumption \textbf{D0} implies  that $|z -\lambda_{j_l}| \ge \bar{\lambda}/2$ for all $z$ on the contour and $j_l \in N_{\bar{\lambda}}(S)^c$ and that 
 $|z - \lambda_i | \ge \delta_S/2$ for all $z$ on the contour and  eigenvalue $\lambda_i, i \in N_{\bar{\lambda}}(S) $. Thus, we have 
 \begin{equation} \label{z-lambda-delta}
 \norm{\frac{1}{\prod_{l=1}^{s_2}(z-\lambda_{i_l})}}  \leq \frac{2^{s_2}}{\delta_S^{s_2}}    
 \end{equation}
and 
\begin{equation} \label{Qbound}
    \|Q\| \leq \frac{2}{\bar{\lambda}}.
\end{equation}
Moreover, the estimate in \eqref{uEuB}  shows that
\begin{equation} \label{ConvuEu}
   \norm{ E \prod_{l=1}^{\beta_{h+1}} \left( u_{i_{l+B_{h}}} u_{i_{l+B_{h}}}^T E \right) } \leq x^{\beta_{h+1}-1} w^2. 
\end{equation}
Combining \eqref{ConveM(a,b)form}, \eqref{z-lambda-delta}, \eqref{ConvuEu}, and the fact that there are exactly $s_1-(k+1)$ factors $\|QE\|,$ we obtain 
\begin{equation*} 
    \begin{split}
    \Norm{M(\alpha;\beta)} & \leq \sum_{i_1,...,i_{s_2}} \frac{x^{s_2-k}  w^{2k} \|QE\|^{s_1-1-k} \|Q\|^{k+1} 2^{s_2}}{\delta_S^{s_2}} \\
    & \leq \sum_{i_1,...,i_{s_2}} \frac{x^{s_2-k}  (\|Q\| \cdot \|E\|)^{s_1-1-k} w^{2k} \|Q\|^{k+1} 2^{s_2}}{\delta_S^{s_2}} \\
    & = \sum_{i_1,...,i_{s_2}} \frac{x^{s_2-k}   \|E\|^{s_1-1-k} w^{2k} \|Q\|^{s_1} 2^{s_2}}{\delta_S^{s_2}}. 
 \end{split}
\end{equation*}
Since $w \leq \|E\|$, we further obtain
\begin{equation} \label{MabCon0}
    \begin{split}
  \Norm{M(\alpha;\beta)}   &  \leq \sum_{i_1,...,i_{s_2}} \frac{x^{s_2-k} \|E\|^{s_1+1-k} w^{2(k-1)} \|Q\|^{s_1} 2^{s_2}}{\delta_S^{s_2}}  \\
    & \leq \sum_{i_1,...,i_{s_2}} \frac{x^{s_2-k} \|E\|^{s_1+1-k} w^{2(k-1)} 2^{s_2+s_1}}{ \bar{\lambda}^{s_1} \delta_S^{s_2}} \,(\text{by \eqref{Qbound}}). 
 \end{split}
\end{equation}
Note that  there are $r^{s_2}$ choices for $i_1,...,i_{s_2}$ and $s_1+s_2=s+1$. The last formula in \eqref{MabCon0} becomes
\begin{equation*} 
    \begin{split}
    & 2^{s+1} r^{s_2} \frac{x^{s_2-k} \|E\|^{s_1+1-k} w^{2(k-1)}}{\bar{\lambda}^{s_1} \delta_S^{s_2}}  = 2^{s+1} \left( \frac{r x}{\delta_S} \right)^{s_2-k} \left( \frac{\sqrt{r} w}{\sqrt{\bar{\lambda} \delta_S}} \right)^{2(k-1)} \left( \frac{\|E\|}{\bar\lambda} \right)^{s_1+1-k}  \frac{r}{ \delta_S}.
\end{split}
\end{equation*}
Applying Lemma \ref{toytrick}, we finally have
\begin{equation*} \label{MabCon}
    \begin{split}
  \Norm{M(\alpha;\beta)}    & \leq \frac{ 2^{s+1} r}{\delta} \max \left\lbrace \frac{rx}{\delta_S}, \frac{\|E\|}{\bar{\lambda}}, \frac{\sqrt{r} w}{\sqrt{\bar{\lambda} \delta_S}} \right\rbrace^{s_2-k+2(k-1)+s_1-(k-1)} \\
    & = \frac{2^{s+1} r}{\delta} \max \left\lbrace \frac{rx}{\delta_S}, \frac{\|E\|}{\bar{\lambda}}, \frac{\sqrt{r} w}{\sqrt{\bar{\lambda} \delta_S}} \right\rbrace^{s_2+s_1-1} \, \\
    & \leq \frac{2^{s+1} r}{\delta_S} \left(\frac{1}{12} \right)^{s_2+s_1-1} \,\, ( \text{by Assumption \textbf{D0}}) \\
    & = \frac{ 2^{s+1} r}{ 12^{s}\delta_S} \,\, (\text{since $s_1+s_2-1=s$}) \\
    & = \frac{2r}{ 6^{s} \delta_S}.
 \end{split}
\end{equation*}   
 
This completes the proof of Lemma \ref{ConvergenceM(a,b)} and then the proof of Theorem \ref{deterministicHRgeneral}. 
 
\section{Extensions to singular subspaces} \label{sec: extension}
In this section, we discuss two extensions of our eigenspace perturbation bounds that are
naturally related to singular subspaces. We first extend our results to rectangular matrices,
using symmetrization to obtain perturbation bounds for the corresponding left and right
singular subspaces. We then return to the symmetric setting and consider eigenspaces
associated with the leading singular values of $A$, or equivalently, with the eigenvalues of
largest absolute value.
\subsection{Extension to rectangular matrices} \label{section:recMa}


 In this subsection, we extend our theorems to the non-symmetric case. Our object now is an   $m \times n$ matrix $A$
 of rank $r_A \le \min \{m, n\}$. Consider the singular value decomposition of $A$
 $$ A:= U \Sigma V^T, $$ where $\Sigma = \text{diag}(\sigma_1,...,\sigma_{r_A})>0$ is a diagonal matrix containing the non-zero singular values $\sigma_1 \geq \sigma_2 \geq ... \geq \sigma_{r_A}$. The columns of the matrices $U=(u_1,..., u_{r_A})$ and $V=(v_1,...,v_{r_A})$ are the  left and right singular vectors of $A$, respectively.  By definition 
$$U^T U =V^T V=I_{r_A}.$$

Let $S$ be a subset of $\{1,2,..., r_A\}$. Using our notation, we set  
\begin{equation}
\Pi_{S}^{left}=\sum_{i \in S} u_i u_i^T, \qquad \Pi_{S}^{right} = \sum_{i \in S} v_i v_i^T.
\end{equation}
In particular, for $S=\{1,...,p\}$, we simply write $\Pi_{p}^{left}\, (\Pi_{p}^{right})$.
Let $E$ be an $m \times n$ ``noise" matrix and set $\tilde{A}=A+E$. Define $\tilde{\Pi}_{S}^{left}, \tilde{\Pi}_{S}^{right}$ accordingly. Our goal is to bound the perturbations of the left and right singular spaces, i.e. 
$$\Norm{\tilde{\Pi}_{S}^{left} -\Pi_{S}^{left}}, \Norm{\tilde{\Pi}_{S}^{right}-\Pi_{S}^{right}}. $$

As in the case of eigenspaces, a classical upper bound on the perturbation of the left and right singular subspaces is provided by Wedin's theorem, the rectangular analogue of the Davis--Kahan theorem; see \cite[Theorem V.4.4]{SS1}. 
 \begin{theorem}[Wedin \cite{W1} - modified version \footnote{In the original version, $\delta_S:= \min_{i\in S, j \notin S}|\sigma_i -\tilde{\sigma}_j|$ .}] Let $\delta_S:= \min_{i \in S, j \notin S}|\sigma_i -\sigma_j|$. Then, 
 \begin{equation}
\max \left\lbrace \Norm{\tilde{\Pi}_{S}^{left} -\Pi_{S}^{left}}, \Norm{\tilde{\Pi}_{S}^{right}-\Pi_{S}^{right}} \right\rbrace \leq \frac{2\|E\|}{\delta_S}.     
 \end{equation}
 \end{theorem}

In this paper, to obtain a refinement of Wedin's theorem, we symmetrize $A$ and $\tilde A$ and apply our theorems in Section \ref{section: main}. There are two natural symmetrization approaches, each with its own advantages.

{\it Additive symmetrization.}
Set  $\mathcal{A}:= \begin{pmatrix}
0 & A \\
A^T & 0
\end{pmatrix}, \mathcal{E}:= \begin{pmatrix}
0 & E \\
E^T & 0
\end{pmatrix}$. 
It is well known that the nonzero eigenvalues of 
$\mathcal{A}$
are $\sigma_1 , \dots, \sigma_{r_A}, -\sigma_1 , \dots, -\sigma_{r_A}$ and with corresponding eigenvectors 
$\frac{1}{\sqrt{2}}(u_1, v_1), \dots , \frac{1}{\sqrt{2}}(u_{r_A}, v_{r_A}), \frac{1}{\sqrt{2}}(u_1, -v_1), \dots,\frac{1}{\sqrt{2}}(u_{r_A}, - v_{r_A} )$ (where $\frac{1}{\sqrt{2}}(u,v)$ is the concatenation of $u$ and $v$). Theorem~\ref{deterministicHRgeneral} can then be applied directly to
\(\mathcal{A}\) and \(\mathcal{A}+\mathcal{E}\). The drawback of this
approach is that the left and right singular vectors are coupled through
the concatenation. Consequently, the resulting argument controls the two
singular subspaces simultaneously, leading to a bound on
\[
\max\left\{
\Norm{\tilde{\Pi}_{S}^{left}
      -\Pi_S^{left}},
\Norm{\tilde{\Pi}_{ S}^{right}
      -\Pi_S^{right}}
\right\}.
\]

{\it Multiplicative symmetrization.}
Alternatively, consider  $M= AA^T$ and $\tilde{M}:= (A+E)(A^T+E^T)$. The nonzero eigenvalues of \(M\) are
\(\sigma_1^2,\ldots,\sigma_{r_A}^2\), and the corresponding eigenvectors
are precisely the left singular vectors of \(A\). Thus, this approach
allows us to study the perturbation of the left singular subspace
directly. Similarly, using \(M=A^T A\) yields a corresponding treatment
of the right singular subspace.
The price for this separation is that the perturbation becomes
\[
\bar E
:=
\tilde M-M
=
AE^T+EA^T+EE^T,
\]
which has a more complicated structure than the original perturbation
\(E\). Consequently, applying our results requires additional estimates
for \(\bar E\), and the corresponding version of Assumption~{\bf D0}
becomes more involved.

In this paper, we restrict our attention to additive symmetrization,
which allows for a direct application of our general perturbation
theorems. The multiplicative approach requires a more technical analysis
and will be investigated separately in future work.

\subsubsection{Results using additive symmetrization} \label{subsec: RecSym}
Following Theorem \ref{deterministicHRgeneral}, we choose a quantity $\bar{\sigma}$ and let $N_{\bar{\sigma}}(S)$ be the $\bar{\sigma}$-neighborhood of $S$:
\begin{equation}
\begin{split}
 & N_{\bar{\sigma}}(S):=\{i  |\sigma_i - \sigma_j| \leq \bar{\sigma} \,\, \text{for some}\,\, j \in S\},\\
 & p:=|S|, r=|N_{\bar{\sigma}}(S)|,\,\,\text{and}\,\,\delta_S= \min_{\sigma_i \in S, \sigma_j \notin S} |\sigma_i -\sigma_j|.
\end{split}
\end{equation}

Now we have the analogue of Definition \ref{xyz}. 
\begin{definition}
\end{definition}
\begin{itemize}
\item $x:= \max_{i,j} |u_i^T E v_j|$ where $i,j \in N_{\bar{\sigma}}(S)$.
\item $y:= \max_{i,j,k} \left\lbrace \norm{v_{i}^T E^T \left(\sum_{l \notin N_{\bar{\sigma}(S)}} \frac{u_l u_{l}^T}{\sigma_k - \sigma_l} \right) E v_{j}},\norm{ u_{i}^T E \left( \sum_{l \notin N_{\bar{\sigma}(S)}} \frac{v_l v_{l}^T}{\sigma_k - \sigma_l} \right) E^T u_{j}} \right\rbrace $ where $i< j \in N_{\bar{\sigma}(S)}$ and $k \in S$.
\item $w := \max_{i,j} \left\lbrace \|E v_i\|, \|E^T u_j\| \right\rbrace$ where $i,j \in  N_{\bar{\sigma}}(S)$.
\end{itemize}
Set 
$$\|E\|_{A,S, \bar{\sigma}} = \sqrt{r} x + \sqrt{r} y.$$
Next are the analogues of Assumptions {\bf C0} and  {\bf D0}: 
\begin{equation}
(\textbf{C1}):\,\,\, \delta_p \geq \max\{24 rx, \frac{576 rw^2}{\bar{\sigma}}\}.
\end{equation}
\begin{equation}
(\textbf{D1}):\,\,\, \delta_S \geq \max\{24 rx, \frac{576 rw^2}{\bar{\sigma}}\}.
\end{equation}
Similar to Subsection \ref{subsection:generalS}, we draw a contour $\Gamma$ which encircles those singular values with indices in $S$ (viewed as positive eigenvalues of $\mathcal A $) such that the distance from any eigenvalue of $\mathcal A$ to $\Gamma$ is at least $\delta_S/2$, and the distance from $\sigma_j, j \notin N_{\bar{\sigma}}(S)$ to $\Gamma$ is at least $\bar{\sigma}/2$. Let $\tilde{S}$ be the set of eigenvalues of $\tilde{A}$ inside $\Gamma$.  Denote 
$$\tilde{\Pi}_{\tilde{S}}^{left}=\sum_{i \in \tilde{S}} \tilde{u}_i \tilde{u}_i^T \,\,\, \text{and} \,\,\, \tilde{\Pi}_{\tilde{S}}^{right}= \sum_{i \in \tilde{S}} \tilde{v}_i \tilde{v}_i^T.$$
\begin{theorem} \label{deterministicRecHR} Under Assumption \textbf{D1}, 
\begin{equation}
\max \left\lbrace \Norm{\tilde{\Pi}_{\tilde{S}}^{left} -\Pi_{S}^{left}}, \Norm{\tilde{\Pi}_{\tilde{S}}^{right}-\Pi_{S}^{right}} \right\rbrace \leq 24 \sqrt{p}  \left(\frac{\|E\|_{A,S, \bar{\sigma}}}{\delta_S} + \frac{ \|E\|}{\bar{\sigma}} \right).
\end{equation}
\end{theorem}

We also obtain the following corollaries for $S=\{1,2,...,p\}$ and $S=\{ \pi(1),...,\pi(p)\}.$
\begin{corollary} \label{cor: recPip}  Under Assumption \textbf{C1}, 
\begin{equation}
\max \left\lbrace \Norm{\tilde{\Pi}_{p}^{left} -\Pi_{p}^{left}}, \Norm{\tilde{\Pi}_{p}^{right}-\Pi_{p}^{right}} \right\rbrace \leq 24 \sqrt{p}  \left(\frac{\|E\|_{A,S, \bar{\sigma}}}{\delta_p} + \frac{ \|E\|}{\bar{\sigma}}  \right).
\end{equation}
\end{corollary}
Using the trivial bound $y \leq \frac{\|E\|^2}{\bar{\sigma}}$ and hence $\|E\|_{A,S, \bar{\sigma}} \leq \sqrt{r} x + \frac{\sqrt{r} \|E\|^2}{\bar{\sigma}} $, we obtain 
\begin{corollary} Under Assumption \textbf{C1}, 
\begin{equation}
\max \left\lbrace \Norm{\tilde{\Pi}_{p}^{left} -\Pi_{p}^{left}}, \Norm{\tilde{\Pi}_{p}^{right}-\Pi_{p}^{right}} \right\rbrace \leq 24\sqrt{p}  \left( \frac{\sqrt{r}x}{\delta_p} + \frac{\sqrt{r} \|E\|^2}{\delta_p \bar{\sigma}}+ \frac{ \|E\|}{\bar{\sigma}}  \right).
\end{equation}
\end{corollary} 

The constants $12, 144$ (from the results for symmetric matrices) get replaced by $24, 576$ due to the symmetrization; see Subsection \ref{proof: reccase}.

\subsection{Eigenspaces associated with the leading singular values}

The preceding subsection treats singular subspaces of rectangular matrices. We now consider a closely related problem in the symmetric setting. Since the singular values of a symmetric matrix $A$ are the absolute values of its eigenvalues, the subspace associated with the leading singular values may contain eigenvectors from both ends of the spectrum. This setting provides an important special case of Theorem~\ref{deterministicHRgeneral}.

Let $S$ correspond to the leading $p$ singular values of $A$; that is,
\[
\{|\lambda_i|:i\in S\}
=
\{\sigma_1,\sigma_2,\ldots,\sigma_p\}
\qquad\text{for some }1\le p\le n.
\]
More precisely, we can write $$
S=S_1\cup S_2,
$$
where
$$
S_1:=\{1,\ldots,k-1\},
\qquad
S_2:=\{n,n-1,\ldots,n-(p-k)\},
$$
for some \(1\le k\le p\). Consequently, $N_{\bar{\lambda}}(S)=\{1,2,...,r_{+} , n-r_{-}+1, n- r_{-}+2,...,n\}$  for properly chosen $r_{+}$ and $r_{-}$.  Set  $r:= r_{+} + r_{-}$ and denote $\bar{\delta}_p:= \min \{ \delta_{k-1}, \delta_{n-(p-k)-1} \}$.
We obtain the following results. 
\begin{theorem} \label{deterministicHRsingular} Under Assumption {\bf D0}, 
\begin{equation} \label{finalboundHRsingular}  \|  \tilde{\Pi}_{(p)} - \Pi_{(p)}  \| \le  12 \sqrt{p} \left(\frac{\|E\|_{A,S, \bar{\lambda}}}{\bar{\delta}_p} + \frac{ \|E\|}{\bar{\lambda}}  \right), 
\end{equation} given that the RHS is less than 1. 
\end{theorem}

Similarly, by the trivial bound $y \leq \frac{\|E\|^2}{\bar{\lambda}}$, we obtain the following corollaries. 
\begin{corollary} \label{leadingprojection1-cor1} 
 Under  Assumption {\bf D0}, 
\begin{equation} \label{leadingprojectionbound1}  \|  \tilde{\Pi}_{(p)} - \Pi_{(p)}  \| \le  12 \sqrt{p} \left( \frac{\sqrt{r}x}{\bar{\delta}_p} + \frac{ \sqrt{r} \| E\|^2 }{\bar{\delta}_p \bar \lambda}+ \frac{ \|E\|}{\bar{\lambda}} \right). \end{equation} 
 \end{corollary} 
 
 \begin{corollary} \label{leadingprojection1-cor2}  Let $r= \mathrm{rank} \, A$. Under Assumption {\bf D0}, 
\begin{equation} \label{leadinglowrank}   \|  \tilde{\Pi}_{(p)} -\Pi_{(p)} \| \le  12 \sqrt{p} \left(  \frac{\sqrt{r}x}{\bar{\delta}_p} + \frac{ \sqrt{r} \bar y}{\bar{\delta}_p \sigma_p} + \frac{r^{3/2} x^2}{ \bar{\delta}_p \sigma_p} + \frac{ \|E\|}{\sigma_p} \right). \end{equation} 
\end{corollary}

\section{Proofs of the main theorems in Section \ref{section:randomnoise} and Section \ref{sec: extension} } \label{section:proof3}

In this section, we prove the main results stated in Sections~\ref{section:randomnoise} and
\ref{sec: extension}. We begin with the
results for random perturbations (Subsection~\ref{subsec: proof random}), where the main task is to obtain suitable probabilistic bounds for the quantities appearing in our deterministic perturbation theorem. Several technical
moment computations used in these arguments are deferred to Appendix \ref{AppendixA}. In Subsection~\ref{proof: reccase}, we then turn to the results concerning the leading singular subspace from Section~\ref{sec: extension}. 
\subsection{Proof of Theorems \ref{mainTh01} and \ref{mainTh-002}: Perturbation of eigenspaces with random noise} \label{subsec: proof random}
To prove Theorem~\ref{mainTh-002}, we first bound the quantities \(x\) and \(y\) for the random perturbation \(E\), and then apply Theorem~\ref{deterministicHR}. We begin with the following lemma, which provides an upper bound for \(x\).

\begin{lemma} \label{uEv} Let $E=(\xi_{kl})_{k,l=1}^{n}$ be a $n \times n$ real symmetric random matrix where 
$$ \{ \xi_{kl}: 1\leq k \leq l \leq n\}$$
is a collection of independent random variables each with mean zero. Further assume 
$$\sup_{1 \leq k \leq l \leq n} | \xi_{kl}| \leq K, \sup_{1 \leq k \leq l \leq n} \Big( \E (\xi_{kl}^2)= \sigma_{kl}^2 \Big) \leq \sigma^2.$$
with probability $1$. Then for any $1\leq i \leq j \leq r$ and every $t >0$
$$\P( |u_i^T E u_j| \geq t) \leq \exp \left( \frac{-t^2/2}{2\sigma^2+ Kt} \right).$$
\end{lemma}

\begin{proof}
Note that for any $1 \leq i,j \leq r$, 
$$u_i^T E u_j = \sum_{k,l} u_{i k} \xi_{kl} u_{jl} = \sum_{k<l} (u_{ik} u_{jl} +u_{il} u_{jk}) \xi_{kl} + \sum_{k} u_{ik} u_{jk} \xi_{kk}.$$
Therefore, $\E(u_i^T E u_j) = 0.$ Furthermore,
$$ \Var (u_i^T E u_j) = \sum_{k<l} (u_{ik} u_{jl} +u_{il} u_{jk})^2 \sigma_{kl}^2 + \sum_{k} u_{ik}^2 u_{jk}^2 \sigma_{kk}^2$$
$$ \leq \sum_{k<l} 2(u_{ik}^2 u_{jl}^2 +u_{il}^2 u_{jk}^2) \sigma_{kl}^2 + \sum_{k} u_{ik}^2 u_{jk}^2 \sigma_{kk}^2 \leq 2 \sum_{k,l} u_{ik}^2 u_{jl}^2 \sigma_{kl}^2.$$
By the assumed upper bound on $\sigma_{kl}^2$, we obtain
$$ \Var (u_i^T E u_j) \leq 2 \sigma^2 \sum_{k,l} u_{ik}^2 u_{jl}^2 = 2 \sigma^2 \sum_{k} u_{ik}^2 \sum_{l} u_{jl}^2 = 2 \sigma^2.$$
Moreover, by the Cauchy-Schwarz inequality and the assumption that $\max_{1 \leq k\leq l}|\xi_{kl}| \leq K$, we have $|(u_{ik} u_{jl}+ u_{il}u_{jk}) \xi_{kl}| \leq K,$ for all $k<l$ and $i,j$. 

These upper bounds on the variances and the off-diagonal terms allow us to apply the Bernstein inequality and get
\begin{equation} 
\P(|u^TEv| \geq t) \leq 2\exp \left( \frac{-t^2/2}{2\sigma^2 + Kt} \right), \, \forall t >0.
\end{equation}
\end{proof}

To obtain the $y$-upper bound, we  will use the moment method to bound $G_{ij,\bar p}:= u_i^T E \left( \sum_{l > r} \frac{u_l u_l^T}{\lambda_{\bar p} -\lambda_l} \right) E u_j $ for a fixed triple $(\bar p,i,j)$: $1 \leq i < j \leq r, 1 \leq \bar p \leq p$ and then apply the union bound on all feasible triples $(i,j,\bar{p})$.

Recall that 
$$s_{jj}^{(\bar p)} := (\lambda_p -\lambda_{r+1} ) \cdot \sum_{l >r} \frac{u_{lj}^2}{\lambda_{\bar p}-\lambda_l},\,\,\, \text{for each}\,\,\, 1 \leq j \leq n, 1 \leq \bar p \leq p,$$
and
$$ \mu_{A,E}:=  \max_{\substack{1 \leq h < l  \leq r \\ 1\leq \bar p \leq p}} \bigg\{\bigg| \sum_{i=1}^n u_{hi} u_{li} \sum_{j=1}^n \sigma_{ji}^2 s_{jj}^{(\bar p)} \bigg|\bigg\}.$$
We obtain the following lemma bounding the first and second moments of $G_{12,\bar p}$. 
\begin{lemma} \label{MomentofG12}
Let $E=(\xi_{ij})_{1 \leq i,j \leq n}$ be a symmetric random $n \times n$ matrix with independent, zero-mean upper-triangular entries. We further assume that with probability $1$, 
$$\norm{\xi_{ij}} \leq K, \,\, \sigma_{ij}^2 \leq \sigma^2, \,\,\, \forall 1 \leq i,j \leq n.$$
One has 
\begin{equation} \label{second moment G12p}
   \E  G_{12,\bar p}^2 \leq (\E G_{12,\bar p})^2 + \frac{8n \sigma^2 (\sigma^2+K^2)}{(\lambda_p -\lambda_{r+1})^2},
\end{equation}
and 
\begin{equation} \label{first  moment G12p}
\begin{split}
     |\E G_{12,\bar{p}}| \leq \frac{\mu_{A,E} +\sqrt{n} \sigma^2}{\lambda_p- \lambda_{r+1}}.
\end{split}
    \end{equation}
    In particular, if the upper triangular entries of $E$ are iid, then $\E G_{12,\bar p}=0$.
\end{lemma}
Thus, by \eqref{second moment G12p} and Chebyshev's inequality,  we obtain: 
\begin{corollary} \label{yboundCorHR} \label{u_iE()Eu_j} For any $t >0$ 
$$\P \bigg(| G_{12,\bar p} -\E G_{12,\bar p} | \ge \frac{t}{\lambda_p -\lambda_{r+1}}\bigg) \le \frac{ 8n \sigma^2 (\sigma^2+K^2) } {t^2} . $$ 
In particular, if the upper triangular entries of $E$ are iid, then $\E G_{12,\bar p}=0$ and 
$$\P \bigg(| G_{12,\bar p} | \ge \frac{t}{\lambda_p -\lambda_{r+1}}\bigg) \le \frac{ 8n \sigma^2 (\sigma^2 + K^2)} {t^2} . $$ 
\end{corollary} 
The computations of $G_{12,\bar p}$'s moments are deferred to Appendix \ref{AppendixA}. Corollary \ref{yboundCorHR} implies that
\begin{equation} \label{G12K}
 \P \bigg(|G_{12,\bar p}| \leq \frac{t}{\lambda_p -\lambda_{r+1}} +|\E G_{12,\bar p}|\bigg) \geq   \P \bigg(| G_{12,\bar p} -\E G_{12,\bar p} | \leq \frac{t}{\lambda_p -\lambda_{r+1}}\bigg) \geq 1 - \frac{8 n \sigma^2 (\sigma^2+K^2)}{t^2}.
\end{equation}

Moreover, by \eqref{first  moment G12p}, we further have  
\begin{equation} \label{PG12kmu to G12k}
    \P\bigg(|G_{12,\bar p}| \leq \frac{t+ \sigma^2 \sqrt{n}+\mu_{A,E}}{\lambda_p -\lambda_{r+1}}\bigg) \geq \P \bigg(|G_{12,\bar p}| \leq \frac{t}{\lambda_p -\lambda_{r+1}} +|\E G_{12,\bar p}|\bigg).
\end{equation}
Together, \eqref{G12K} and \eqref{PG12kmu to G12k} imply that for any $t > 0$,
\[
\P\bigg(|G_{12,\bar p}| \leq \frac{t+\sqrt{n}\sigma^2 +\mu_{A,E}}{\lambda_p -\lambda_{r+1}}\bigg) \geq 1 - \frac{8n \sigma^2 (\sigma^2+K^2)}{t^2}.
\]
Thus, for any $t > 2 \sqrt{n} \sigma^2$, we have
\begin{equation} \label{G12Kmu}
    \P\bigg(|G_{12,\bar p}| \leq \frac{t +\mu_{A,E}}{\lambda_p -\lambda_{r+1}}\bigg) \geq 1 - \frac{8n \sigma^2 (\sigma^2+K^2)}{(t/2)^2} = 1 - \frac{32n \sigma^2 (\sigma^2+K^2)}{t^2}.
\end{equation}

Now, we are ready to prove Theorem \ref{mainTh-002}. Let $\Omega$ be the compound event defined by
\begin{equation*}
    \begin{split}
   \Omega & = \bigg\{ \text{Assumption}\,\textbf{D2}\,\,\text{holds} \bigg\} \cap      \left\lbrace |u_i^T E u_j| \leq t_1, \forall 1 \leq i,j \leq r \right\rbrace \\
   & \cap  \bigg\{\bigg|u_i^T E \bigg( \sum_{l >r} \frac{u_l u_l^T}{\lambda_{\bar p} -\lambda_l} \bigg) E u_j \bigg| \leq \frac{t_2+\mu_{A,E}}{\lambda_p -\lambda_{r+1}}, \forall 1 \leq i < j \leq r, 1 \leq \bar p \leq p \bigg\}.
    \end{split}
\end{equation*}
By Lemma \ref{uEv}, \eqref{G12Kmu}\footnote{If $t_2 < 2\sqrt{n} \sigma^2$, then the claimed probability lower bound is non-positive, and hence the assertion is trivial.}, and that Assumption \textbf{D2} holds with probability $1-\alpha_1$, $\Omega$ happens with probability at least 
$$1 - \alpha_1 - 2r^2 \exp\left( \frac{-t_1^2/2}{2 \sigma^2 + Kt_1} \right) - \frac{32r(r-1)p  n \sigma^2 (\sigma^2+ K^2)}{t_2^2}.$$ 
Now, applying Theorem \ref{deterministicHR} on the event $\Omega$, we have
$$\Norm{\tilde{\Pi}_p - \Pi_p  } \leq 12 \sqrt{p} \left( \epsilon_1 + \sqrt{r} t_1 \epsilon_2 + \frac{ \sqrt{r} (t_2+\mu_{A,E})}{\|E\|} \epsilon_1 \epsilon_2 \right)\,\,\text{as desired}.$$

The proof of Theorem~\ref{mainTh01} follows the same procedure as above, except for the
following modification. Instead of bounding \(y\) by directly estimating the first and
second moments of \(G_{12,\bar p}\), we first use
\[
y
\leq
\max_{1\leq i<j\leq r}
\frac{|u_i^T E^2u_j|}{\lambda_p}
+
\frac{rx^2}{\lambda_p}.
\]
We then bound
\[
\max_{1\leq i<j\leq r}|u_i^T E^2u_j|
\]
using the moment method.  Compared with the proof of
Theorem~\ref{mainTh-002}, this approach introduces the additional term
\(rx^2/\lambda_p\). On the other hand, the moment computations for
\(u_i^T E^2u_j\) are considerably simpler than those for
\(G_{ij,\bar p}\); see Lemma~\ref{MomentofG12lowrank} and
Corollary~\ref{u_iEEu_j}.

\begin{lemma} \label{MomentofG12lowrank}
    Let $E=(\xi_{ij})_{1 \leq i,j \leq n}$ be a random symmetric $n \times n$ matrix with independent, zero-mean upper-triangular entries. One has 
    \begin{equation}
        \E \left(E u_1 \cdot E u_2 \right)^2 \leq 5nm_4 + \mu^2,
    \end{equation}
    here, $\mu:= \E E u_1 \cdot E u_2.$

    In particular, if $E$ is regular, then $\mu = 0$.
\end{lemma}

\begin{corollary} \label{u_iEEu_j} For any $t >0$, 
$$\P (| E u_1 \cdot E u_2 -\mu  | \ge t) \le \frac{ 5nm_4 } {t^2} . $$ 
In particular, if $E$ is regular, then
$$\P (| E u_1 \cdot E u_2 | \ge t) \le \frac{ 5nm_4 } {t^2} . $$ 
\end{corollary} 
The proofs of Lemma \ref{MomentofG12lowrank} and Corollary \ref{u_iEEu_j} are also deferred to Appendix \ref{AppendixA}.
\subsection{Proof of Theorem \ref{deterministicRecHR}: Perturbation of singular spaces} \label{proof: reccase}
In this subsection, given Theorem \ref{deterministicHRgeneral}, we are going to prove Theorem \ref{deterministicRecHR}. 

First, we apply Additive symmetrization, considering the $(m+n) \times (m+n)$ matrices
   $$\mathcal{A}:= \begin{pmatrix}
   0 & A \\
   A^T & 0
   \end{pmatrix}, \,\,  \,\, \mathcal{E}:= \begin{pmatrix}
   0 & E \\
   E^T & 0
   \end{pmatrix}, \,\, \text{and} \,\, \tilde{\mathcal{A}}:=\mathcal{A}+\mathcal{E}. $$
     The non-zero eigenvalues of $\mathcal{A}$ are $\pm \sigma_1,...,\pm \sigma_{r_A}$, with corresponding eigenvectors 
   $$\textbf{u}_j:=\frac{1}{\sqrt{2}} (u_j^T, v_j^T)^T, \,\,\, \text{and} \,\, \textbf{u}_{j+r}:= \frac{1}{\sqrt{2}} (u_j^T,-v_j^T)^T.$$
   Thus $\mathcal{A}$ admits the spectral decomposition 
   $$\mathcal{A}= \mathcal{U} \mathcal{D} \mathcal{U}^T,$$
   where $\mathcal{U}:=(\textbf{u}_1,...,\textbf{u}_{2{r_A}})$ is the matrix with columns $\textbf{u}_1,...,\textbf{u}_{2{r_A}}$. We expand $\Gamma$ into $\Gamma'$ which contains exactly $S'= S \cup -S$ and define  $\mathcal{U}_{S'}:=(\textbf{u}_i)_{i \in S'}$. By the decomposition of $\tilde{\mathcal{A}} = \tilde{\mathcal{U}} \tilde{\mathcal{D}} \tilde{\mathcal{U}}^T$, we also define $\tilde{\mathcal{U}}$, $ \tilde{S}'= \tilde{S} \cup - \tilde{S}, $ and  $\tilde{\mathcal{U}}_{\tilde{S'}}$ respectively.

 Note that for all $i,j \in N_{\bar{\sigma}}(S)$, 
 $$ 2\textbf{u}_i \mathcal{E} \textbf{u}_j= u^T_i E v_j + v_i^T E^T u_j,$$
 $$ 2\textbf{u}_i \mathcal{E} \textbf{u}_{j+r_A}= - u^T_i E v_j + v_i^T E^T u_j,$$
  $$ 2\textbf{u}_{i+r_A} \mathcal{E} \textbf{u}_{j+r_A}= - (u^T_i E v_j + v_i^T E^T u_j).$$
  For each $ 1 \leq i \leq 2r_A$, we define 
  \begin{equation*}
      \bar{i}:=\begin{cases}
      i \,\,\,& \text{if} \,\, i \leq r_A \\
 i-r_A \,\, &\text{if} \,\, i > r_A
  \end{cases}.
  \end{equation*}
     Therefore, for any $1 \leq i,j \leq 2r_A$ and $\bar{i}, \bar{j} \in N_{\bar{\sigma}}(S)$, 
   \begin{equation} \label{RecC(2)x}
   |\textbf{u}_i^T \mathcal{E} \textbf{u}_j| \leq \frac{1}{2} \left( |u_{\bar{i}}^T E v_{\bar{j}}| + |v_{\bar{i}}^T E^T u_{\bar{j}}| \right) \leq x.
   \end{equation}
   
   Similarly,  we also obtain that for all $1 \leq i < j \leq 2r_A$ and $\bar{i}, \bar{j} \in N_{\bar{\sigma}}(S)$, 
 \begin{equation*}
 \begin{split}
  4 \textbf{u}_i^T \mathcal{E} \left( \sum_{ \substack{l:\,\bar{l} \notin N_{\bar{\sigma}}(S)} }\frac{\textbf{u}_l \textbf{u}_l^T}{\sigma_k - \sigma_l}  \right)  \mathcal{E} \textbf{u}_j &  = \pm v_{\bar{i}}^T E^T \sum_{\substack{l \\\bar{l} \notin N_{\bar{\sigma}}(S)}} \frac{u_{\bar{l}} u_{\bar{l}}^T}{\sigma_k - \sigma_l} E v_{\bar{j}} \pm u_{\bar{i}}^T E \sum_{\substack{l \\\bar{l} \notin N_{\bar{\sigma}}(S)}} \frac{v_{\bar{l}} v_{\bar{l}}^T}{\sigma_k - \sigma_l} E^T u_{\bar{j}} \\
  & = 2 \left( \pm v_{\bar{i}}^T E^T \sum_{\substack{l \notin N_{\bar{\sigma}}(S)}} \frac{u_l u_{l}^T}{\sigma_k - \sigma_l} E v_{\bar{j}} \pm u_{\bar{i}}^T E \sum_{\substack{l \notin N_{\bar{\sigma}}(S)}} \frac{v_l v_{l}^T}{\sigma_k - \sigma_l} E^T u_{\bar{j}}\right).
 \end{split}
 \end{equation*}
  
   Thus, for all $ 1 \leq \bar{i} < \bar{j} \leq r$ and $\sigma_k \in S$, by the triangle inequality,  
   \begin{equation} \label{RecC(3)y}
   \bigg|\textbf{u}_i^T \mathcal{E} \bigg( \sum_{ \substack{l:\,\bar{l} \notin N_{\bar{\sigma}}(S)} }\frac{\textbf{u}_l \textbf{u}_l^T}{\sigma_k - \sigma_l}  \bigg)  \mathcal{E} \textbf{u}_j\bigg| \leq \frac{1}{2} \bigg( \bigg|v_{\bar{i}}^T E^T \sum_{\substack{l \notin N_{\bar{\sigma}}(S)}} \frac{u_l u_{l}^T}{\sigma_k - \sigma_l} E v_{\bar{j}}\bigg|+\bigg| u_{\bar{i}}^T E \sum_{\substack{l \notin N_{\bar{\sigma}}(S)}} \frac{v_l v_{l}^T}{\sigma_k - \sigma_l} E^T u_{\bar{j}} \bigg|\bigg)  \leq y.
   \end{equation}
  In fact, we still obtain Lemma \ref{lemma: Sigma1allcases} - second part (Section \ref{section:proof2}) if we consider the particular case $(\bar{i_1},...,\bar{i_k}) = (i,...,i), 1 \leq i \leq r$ instead of the case $(i_1,...,i_k)=(i,...,i)$ from the original proof. Hence, (\ref{RecC(2)x}), (\ref{RecC(3)y}), and the fact that $\|\mathcal{E}\|=\|E\|$ are enough for us to apply Theorem \ref{deterministicHRgeneral} on $\mathcal{A}$ and $\mathcal{E}$. It means that we do not need to require $\norm{\textbf{u}_i^T \mathcal{E} \left( \sum_{ \substack{l \\\bar{l} \notin N_{\bar{\sigma}}(S)} }\frac{\textbf{u}_l \textbf{u}_l^T}{\sigma_k - \sigma_l}  \right)  \mathcal{E} \textbf{u}_{i+r_A}} \leq y.$ \\ 
Since $|\mathcal{S'}|=2p$ and $|N_{\bar{\sigma}}(\mathcal{S}')|=2r$, we obtain
\begin{equation} \label{recEQ1}
 \Norm{\mathcal{U}_{S'}\mathcal{U}_{S'}^T - \tilde{\mathcal{U}}_{\tilde{S'}} \tilde{\mathcal{U}}_{\tilde{S'}}^T} \leq 12 \sqrt{2p} \left( \frac{ \|E\|}{\bar{\sigma}} + \frac{\sqrt{2r} x}{\delta_S} + \frac{\sqrt{2r} y}{ \delta_S} \right) \leq 24 \sqrt{p} \left( \frac{ \|E\|}{\bar{\sigma}} + \frac{\sqrt{r} x}{\delta_S} + \frac{\sqrt{r} y}{ \delta_S} \right).
\end{equation}

 Furthermore, by \cite[Lemma 26]{OVK13},  
 \begin{equation} \label{recEQ2}
 \max \left\lbrace \Norm{\tilde{\Pi}_{\tilde{S}}^{left} -\Pi_{S}^{left}}, \Norm{\tilde{\Pi}_{\tilde{S}}^{right}-\Pi_{S}^{right}} \right\rbrace \leq  \Norm{\mathcal{U}_{S'}\mathcal{U}_{S'}^T - \tilde{\mathcal{U}}_{\tilde{S'}} \tilde{\mathcal{U}}_{\tilde{S'}}^T} .
 \end{equation}
 Combining (\ref{recEQ1}) and (\ref{recEQ2}),  we complete the proof. 
 
\bibliographystyle{amsrefs}
\bibliography{RefO}
\appendix
\section{Moment Computations} \label{AppendixA}

In this appendix, we present the moment computations used in Section~\ref{section:proof3} to derive the probabilistic bounds in Theorems~\ref{mainTh01} and~\ref{mainTh-002}. Throughout, $E$ is a symmetric random matrix with independent, mean-zero upper-diagonal entries. In Subsection~\ref{subsec: moment easy}, we estimate the first and second moments of
\[
Eu_i\cdot Eu_j=u_i^T E^2u_j,
\]
which are used in the proof of Theorem~\ref{mainTh01}. In Subsection~\ref{subsec: moment hard}, we carry out the analogous (more complicated) computation for
\[
u_i^T E
\left(
\sum_{l>r}\frac{u_lu_l^T}{\lambda_{\bar p}-\lambda_l}
\right)
Eu_j,
\]
which yields the moment estimates needed for Theorem~\ref{mainTh-002}. In both cases, the main task is to exploit the independence and symmetry of the entries of $E$ to identify the non-vanishing terms in the corresponding moment expansions.
\subsection{Moment computation of $ \norm{E u_i \cdot E u_j} $} \label{subsec: moment easy}
To obtain a more friendly version (for applications), we first bound $ \norm{ Eu_i \cdot Eu_j}  $.
Without loss of generality, take $i=1$ and $j=2$. We have 
$$\bar G_{12} :=  Eu_1 \cdot Eu_2 = u_1^T E^2 u_2 = \sum_{i,j} u_{1i} u_{2j} \cdot (\sum_l \xi_{li} \xi_{lj}). $$
Since the upper-triangular entries are independent and centered, the expectation comes from the diagonal term. This implies
$$ \E \bar G_{12} = \sum_{i,j} u_{1i} u_{2j} \cdot (\sum_l \E \xi_{li} \xi_{lj})  = \sum_i u_{1i} u_{2i} \sum_l \E \xi_{li}^2. $$

Let $S_i = \sum_l \E \xi_{li}^2 $. If the random matrix $E$ is regular, then $S_i =S_j$ for all $1\leq i,j \leq n$, and hence the expectation of $\bar G_{12} $ is zero, as $u_1$ and $u_2$ are orthogonal. 
If $E$ is not a regular random matrix, then 
$$\mu:=  \E \bar G_{12} = \sum_i  S_i u_{1i} u_{2i} , $$ with $S_i$ defined as above. In geometrical term, $S_i$ is the expectation of the square of the length of the $i$th row of $E$. 

We now turn to the second moment 
$$\bar G_{12}^2 = \sum_{\substack{i,j,l \\ i',j',l'}} \xi_{il} \xi_{lj} \xi_{i' l'} \xi_{l' j'} u_{1i} u_{2j} u_{1i'} u_{2j'},\,\,\text{and hence}\,\,\,\E \bar G_{12}^2 = \sum_{\substack{i,j,l \\ i',j',l'}}u_{1i} u_{2j} u_{1i'} u_{2j'} \E  \xi_{il} \xi_{lj} \xi_{i' l'} \xi_{l' j'} . $$
Since $\E(\xi_{ij}) =0$, most of the expectations on the RHS are zero. The nontrivial ones can be divided into four cases.

\textit{Case 1. } (Fourth moment case)  $\{i,l \}=\{l ,j \}=\{i', l' \}=\{l',j'\}.$  This can only happen if
$i=j=i'=j', l=l'$ or $i=j=l', i'=j'=l.$

In the first scenario, we receive: $   \sum_i  (u_{1i } u_{2i})^2  \sum_l \E \xi_{il}^4  \le h_{12} n m_4 . $

In the second scenario, we receive: $ \sum_{i,l} u_{1i} u_{2i} u_{1l} u_{2l}  \E \xi_{il}^4 .$
As $u_{1i} u_{2i} u_{1l} u_{2l} \le \frac{1}{2} ( (u_{1i} u_{2i } )^2 + (u_{1l } u_{2l} )^2 $, the last quantity can also be bounded by $h_{12} n m_4 $. Thus, the contribution in Case 1 is at most 
$2 h_{12} nm _4 $.

\textit{Case 2.} $\{i,l \}=\{l ,j \} \neq \{i', l' \}=\{l',j'\}.$ This could only happen if  $i=j, i'=j'$ and $\{i,l \} \neq \{i',l'\}$.  The contribution from this case is 
$$\sum_{i,i', l, l'} u_{1i}u_{2i} u_{1i'} u_{2i'} \E \xi_{il}^2 \E \xi_{i' l'}^2 -D , $$ where $D$ collect all diagonal terms with $\{i,l \} = \{ i' l' \} $. Without any restriction on the indices, the sum $\sum_{i,i', l, l'} u_{1i}u_{2i} u_{1i'} u_{2i'} \E \xi_{ij}^2 \E \xi_{i' l'}^2 $  equals 
$$ \Big ( \sum_{i,l} u_{1i}u_{2i} \E \xi_{il}^2 \Big) ^2  = \Big ( \sum_i  u_{1i}u_{2i} \sum_l \E \xi_{il}^2 \Big)^2 = \mu^2 .$$  
The term $D$ can be bounded (in absolute value) by  $h_{12} n m_4$ similarly to Case 1. Thus, the contribution of this case is at most 
$$h_{12} n m_4 + \mu^2. $$

\textit{Case 3.}  $\{i,l \}=\{i' ,l' \} \neq \{l, j \}=\{l',j'\}.$ It happens if $(i=i', l=l', j=j', i \neq j)$. The contribution from this case is 
$$\sum_{i,j,l}  \E \xi_{il}^2 \xi_{jl}^2 u_{1i}^2 u_{2j}^2 =\sum _{i,j} u_{1i}^2 u_{2j}^2 \sum_l \E \xi_{il}^2 \xi_{jl}^2. $$
By the Cauchy--Schwarz inequality,
$$  \sum_l \E \xi_{il}^2 \xi_{jl}^2 \le \frac{1}{2} \sum _l \E ( \xi_{il}^4+ \xi_{jl}^4)  \le nm_4 . $$ Furthermore, as $u_1, u_2$ are unit vectors
$\sum _{i,j} u_{1i}^2 u_{2j}^2 = 1$. Thus, the contribution in this case is bounded from above by $n m_4$.

\textit{Case 4.} $\{i,l \}=\{l' ,j' \} \neq \{l, j \}=\{i',l'\}.$ It happens if  $(i=j', i'=j, l=l', i \neq j)$. The contribution in this case is 
$$\bigg|\sum_{i,j,l} u_{1i}u_{2i} u_{1j} u_{2j} \E \xi_{il}^2 \xi_{jl}^2 \bigg| \leq \sum _{i,j} | u_{1i}u_{2i} u_{1j} u_{2j}| \sum_l   
\E \xi_{il}^2 \xi_{jl}^2 . $$
By the argument in the previous case, 
$$\sum_l   
\E \xi_{il}^2 \xi_{jl}^2  \le n m_4. $$ Furthermore, by the Cauchy-Schwarz inequality, 
$$ \sum _{i,j} | u_{1i}u_{2i} u_{1j} u_{2j}|  =\sum_i | u_{1i} u_{2i} | \sum _j |u_{1j} u_{2j }| \le 1 . $$ Thus, the contribution from this case is bounded by $n m_4$.

Putting all the cases together,  
$$\E \bar G_{12}^2 \le 3 h_{12} n m_4 + 2n m_4 + \mu^2  \le 5 nm_4 + \mu^2. $$ 
Recall that $\E \bar G_{12} =\mu$, so the variance of $\bar G_{12} $ is at most $5nm_4$. By Chebyshev's inequality, we obtain 
\begin{corollary} For any $t >0$ 
$$\P (| \bar G_{12} -\mu  | \ge t) \le \frac{ 5nm_4 } {t^2} . $$ 

In particular, if $E$ is regular, then $\mu=0$ and 

$$\P (| \bar G_{12} | \ge t) \le \frac{ 5nm_4 } {t^2} . $$ 

\end{corollary} 

\subsection{Moment computation of $ \big| u_i^T E \left( \sum_{l > r} \frac{u_l u_l^T}{\lambda_{\bar p} -\lambda_l} \right) E u_j \big| $} \label{subsec: moment hard}
  In this subsection, via moment method, we are going to bound $ \norm{u_i^T E \left( \sum_{l > r} \frac{u_l u_l^T}{\lambda_{\bar p} -\lambda_l} \right) E u_j} $ for a fixed triple $(\bar p,i,j)$: $1 \leq i < j \leq r, 1 \leq \bar p \leq p$, i.e. $S:=\{1,2,...,p\}, N_{\bar{\lambda}}(S)=\{1,2,...,r\}.$ Here, we will work with only the assumption that $E$ is a symmetric $(K,\sigma)$-bounded random matrix with independent upper-triangular entries.  A similar computation also applies for a general subset $S$ and its corresponding $N_{\bar{\lambda}}(S)$.
 
 Without loss of generality, we take $i=1$, $j=2$ and fix $1 \leq  \bar p \leq p$. For convenience, for each $l > r$ and $1\leq k,m \leq n$, we define
 \[
 \bar{\sigma}_l := \lambda_{\bar p} - \lambda_l\,\text{and}\,\,  s_{km}:= \sum_{l > r}\frac{u_{lk} u_{lm}}{\bar{\sigma}_l}.
 \]
 Thus, we can rewrite $G_{12,\bar{p}}$ as
 \(
 G_{12,\bar p}= u_{1}^T E \left( \sum_{l > r} \frac{u_l u_l^T}{\bar{\sigma}_l} \right) E u_2,
 \) and expand it into
\begin{equation} \label{G12p expansion}
\begin{split}
G_{12,\bar p} &= \sum_{l > r} \sum_{k,i,j,m} u_{1k} \xi_{ki} u_{li} u_{lj} \xi_{jm} u_{2m} \frac{1}{\bar{\sigma}_l} = \sum_{l > r} \sum_{k,i,j,m} u_{1k}u_{2m} \frac{u_{li} u_{lj}}{\bar{\sigma}_l} \xi_{ki} \xi_{jm}.
\end{split}
\end{equation}

Before going to our estimation, we present here an identity that we will use later in this section. 
\begin{lemma} \label{lemma: skm}
    \[
    \sum_{k,m} s_{km}^2 = \sum_{l >r} \frac{1}{\bar{\sigma}_l^2}.
    \]
\end{lemma}
\textit{Proof of Lemma \ref{lemma: skm}.} Consider matrix $S:= (s_{km})_{1\leq k, m \leq n}$. Since $u_{r+1}, \cdots, u_n$ are orthonormal,  $S$ admits the spectral decomposition: 
\[
S = [u_{r+1} \,u_{r+2}\,\cdots u_n] \mathrm{Diag}\bigg[\frac{1}{\bar{\sigma}_{r+1}}, \frac{1}{\bar{\sigma}_{r+2}}, \dots, \frac{1}{\bar{\sigma}_n} \bigg]  [u_{r+1} \,u_{r+2}\,\cdots u_n]^T.
\]
Therefore, 
\[
\sum_{k,m} s_{km}^2 = \|S\|_F^2 = \sum_{l >r} \frac{1}{\bar{\sigma}_l^2}\,\,\text{as desired}.
\]
\subsubsection{Bounding $|\E G_{12,\bar{p}}|$: Proof of \eqref{first  moment G12p}} By \eqref{G12p expansion}, we have
\[
\E G_{12,\bar p}  = \sum_{l > r} \sum_{k,i,j,m} u_{1k}u_{2m} \frac{u_{li} u_{lj}}{\bar{\sigma}_l} \E \xi_{ki} \xi_{jm}.
\]
Since $\E \xi_{ki} \xi_{jm}$ is nontrivial when $\{k,i\} = \{j,m\}$, we have
\begin{equation}
\begin{split}
\E G_{12,\bar p} & = \sum_{l > r} \sum_{\{k,i\} = \{j,m\}} u_{1k}u_{2m} \frac{u_{li} u_{lj}}{\bar{\sigma}_l} \E \xi_{ki} \xi_{jm} \\
& = \sum_{l > r} \sum_{k \neq i} u_{1k}u_{2i} \frac{u_{li} u_{lk}}{\bar{\sigma}_l} \E \xi_{ki}^2 + \sum_{l > r} \sum_{k,i} u_{1k}u_{2k} \frac{u_{li} u_{li}}{\bar{\sigma}_l} \E \xi_{ki}^2. 
\end{split}
\end{equation}
By the triangle inequality, we have 
\[
|\E G_{12,\bar p}| \leq \bigg| \sum_{l > r} \sum_{k \neq i} u_{1k}u_{2i} \frac{u_{li} u_{lk}}{\bar{\sigma}_l} \E \xi_{ki}^2 \bigg| + \bigg| \sum_{l > r} \sum_{k,i} u_{1k}u_{2k} \frac{u_{li} u_{li}}{\bar{\sigma}_l} \E \xi_{ki}^2\bigg|.
\]
By the definition of $\mu_{A,E}$, the second term is bounded by $\frac{\mu_{A,E}}{\lambda_p -\lambda_{r+1}}$. Thus, to prove \eqref{first  moment G12p}, we need to show that 
\[
\bigg| \sum_{l > r} \sum_{k \neq i} u_{1k}u_{2i} \frac{u_{li} u_{lk}}{\bar{\sigma}_l} \E \xi_{ki}^2 \bigg|= \bigg| \sum_{l > r} \sum_{k \neq i} u_{1k}u_{2i} \frac{\sigma_{ki}^2 u_{li} u_{lk}}{\bar{\sigma}_l} \bigg|   \leq \frac{\sqrt{n} \sigma^2}{\lambda_p -\lambda_{r+1}}.
\]
For $1 \leq k,i \leq n$, define 
\[
M_{ki}:= \begin{cases} \sigma_{ki}^2 \cdot s_{ki}
      \,\,\,& \text{if $k \neq i$} \\
    0 &\text{if $k=i$}
\end{cases}.
\]
We can rewrite $\big| \sum_{l > r} \sum_{k \neq i} u_{1k}u_{2i} \frac{\sigma_{ki}^2 u_{li} u_{lk}}{\bar{\sigma}_l} \big| $ as 
\[
|\sum_{k,i} u_{1k} u_{2i} M_{ki}| = | u_1^T M u_2| \leq \|M\| \leq \|M\|_F,\,\,\text{where}\,\, M:=(M_{ki})_{1 \leq k,i \leq n}.
\]
Moreover, by Lemma \ref{lemma: skm} and the definition of $\sigma$ that $\max_{1 \leq k,i \leq n}\sigma_{ki}^2 \leq \sigma^2$, we have
\begin{equation}
    \begin{split}
        \|M\|_F^2= \sum_{k \neq i} \sigma_{ki}^4 \cdot s_{ki}^2  \leq \sigma^4 \sum_{k\neq i} s_{ki}^2  \leq \sigma^4 \sum_{k, i} s_{ki}^2  = \sigma^4 \cdot \bigg(\sum_{l >r} \frac{1}{\bar\sigma_l^2}  \bigg).
    \end{split}
\end{equation}
Since $\bar{\sigma}_l \geq (\lambda_p -\lambda_{r+1}) > 0$ for all $l > r$, we have $\|M\|_F^2 \leq \frac{(n-r) \sigma^4}{(\lambda_p -\lambda_{r+1})^2}$, and hence 
\[
\bigg| \sum_{l > r} \sum_{k \neq i} u_{1k}u_{2i} \frac{u_{li} u_{lk}}{\bar{\sigma}_l} \E \xi_{ki}^2 \bigg| \leq \|M\|_F \leq \frac{\sqrt{n} \sigma^2}{\lambda_p -\lambda_{r+1}} \,\,\text{as desired}. 
\]
\begin{remark} 
If the upper triangular entries of $E$ are iid, then, by \eqref{G12p expansion} and the fact that $u_l, u_i, u_k$ are orthogonal, $\E G_{12} =0$.
\end{remark}

\subsubsection{Bounding $\E (G_{12,\bar{p}})^2$: Proof of \eqref{second moment G12p}}
For the second moment, we have 
\begin{equation}
\begin{split}
\E G_{12,\bar p}^2& = \sum_{l,l' >r} \sum_{\substack{k,i,j,m \\ k',i',j',m'}} u_{1k}u_{2m} \frac{u_{li} u_{lj}}{\bar{\sigma}_l}  u_{1k'}u_{2m'} \frac{u_{l'i'} u_{l'j'}}{\bar{\sigma}_{l'}} \E \xi_{ki} \xi_{jm} \xi_{k'i'} \xi_{j'm'} \\
&= M_1+M_2+M_3+M_4,
\end{split}
\end{equation}
in which 
\[
M_1:= \sum_{l,l'>r} \sum_{\substack{ \{k,i\}=\{j,m\} \neq \\
\{k',i'\}=\{j',m'\}}} (u_{1k}u_{2m} \frac{u_{li} u_{lj}}{\bar{\sigma}_l}  u_{1k'}u_{2m'} \frac{u_{l'i'} u_{l'j'}}{\bar{\sigma}_{l'}} \E \xi_{ki}^2 \E \xi_{k'i'}^2 ),  
\]
\[
M_2: = \sum_{l,l'>r} \sum_{\substack{ \{k,i\}=\{k',i'\} \neq\\\{j,m\}=\{j',m'\}}} (u_{1k}u_{2m} \frac{u_{li} u_{lj}}{\bar{\sigma}_l}  u_{1k'}u_{2m'} \frac{u_{l'i'} u_{l'j'}}{\bar{\sigma}_{l'}} \E \xi_{ki}^2 \E \xi_{jm}^2), 
\]
\[
M_3:= \sum_{l,l'>r} \sum_{\substack{ \{k,i\}=\{j',m'\} \neq \\\{k',i'\}=\{j,m\}}} (u_{1k}u_{2m} \frac{u_{li} u_{lj}}{\bar{\sigma}_l}  u_{1k'}u_{2m'} \frac{u_{l'i'} u_{l'j'}}{\bar{\sigma}_{l'}} \E \xi_{ki} ^2 \E \xi_{jm}^2),
\]
and
\[
M_4:= \sum_{l,l'>r} \sum_{\substack{ \{k,i\}=\{j,m\}=\\
\{k',i'\}=\{j',m'\}}} (u_{1k}u_{2m} \frac{u_{li} u_{lj}}{\bar{\sigma}_l}  u_{1k'}u_{2m'} \frac{u_{l'i'} u_{l'j'}}{\bar{\sigma}_{l'}} \E \xi_{ki}^4 ).
\]
In fact, if we define
$$M_5:= M_1 + M_4 - (\E G_{12,\bar p})^2 = \sum_{l,l'} \sum_{\substack{ \{k,i\}=\{j,m\}=\\
\{k',i'\}=\{j',m'\}}} (u_{1k}u_{2m} \frac{u_{li} u_{lj}}{\bar{\sigma}_l}  u_{1k'}u_{2m'} \frac{u_{l'i'} u_{l'j'}}{\bar{\sigma}_{l'}})(\E \xi_{ki}^4 - (\E \xi_{ki}^2)^2),$$
then
\begin{equation} \label{G12power2split}
    \E G_{12,\bar p}^2 = (\E G_{12,\bar p})^2 + M_5 + M_2+M_3.
\end{equation}

Thus, to estimate $\E G_{12,\bar p}^2$, we only need to bound $M_2$ and $M_5$, since $M_3$ is $M_2$ with indices permutation. 
\subsubsection{ Bounding $M_2$}
We have
\begin{equation}
\begin{split}
M_2  &= \sum_{l,l'>r} \sum_{\substack{ \{k,i\} \neq\\\{j,m\}}} \sigma_{ki}^2 \sigma_{jm}^2 \left( \frac{(u_{1k})^2 (u_{2m})^2 u_{li} u_{lj} u_{l'i} u_{l'j}}{\bar{\sigma}_l \bar{\sigma}_{l'}} + \frac{(u_{1k})^2 u_{2m} u_{2j} u_{li} u_{lj} u_{l'i} u_{l'm} }{\bar{\sigma}_l \bar{\sigma}_{l'}} \right) \\
& + \sum_{l,l'>r} \sum_{\substack{ \{k,i\}\neq\\\{j,m\}}} \sigma_{ki}^2 \sigma_{jm}^2 \left( \frac{u_{1k}u_{1i} (u_{2m})^2 u_{li} u_{lj} u_{l'k} u_{l' j}}{\bar{\sigma}_l \bar{\sigma}_{l'}} + \frac{u_{1k} u_{1i} u_{2m} u_{2j} u_{li} u_{lj} u_{l'k} u_{l'm} }{\bar{\sigma}_l \bar{\sigma}_{l'}} \right).
\end{split}
\end{equation}
The first sub-sum can be bounded as follows. First, by the notation that 
$s_{ij} = \sum_{l >r} \frac{u_{li} u_{lj}}{\bar{\sigma}_l},$ we rewrite the first sub-sum into
\begin{equation*} 
\begin{split}
\sum_{l,l'>r} \sum_{\substack{ \{k,i\} \neq\\\{j,m\}}} \sigma_{ki}^2 \sigma_{jm}^2 \frac{(u_{1k})^2 (u_{2m})^2 u_{li} u_{lj} u_{l'i} u_{l'j}}{\bar{\sigma}_l \bar{\sigma}_{l'}} & \leq \sum_{k,i,j,m} \sigma_{ki}^2 \sigma_{jm}^2 u_{1k}^2 u_{2m}^2 s_{ij}^2 \leq \sigma^4 \sum_{k,i,j,m}  u_{1k}^2 u_{2m}^2 s_{ij}^2.
\end{split}
\end{equation*}
The first inequality follows from the fact that $s_{ij}^2 \geq 0$, while the second is true by the definition that $\max_{1\leq i,j \leq n} \sigma_{ij}^2 \leq \sigma^2.$ By Lemma \ref{lemma: skm}, we rewrite the RHS as
\begin{equation} \label{M_21}
    \begin{split}
    \sigma^4 \cdot \sum_{k,m} u_{1k}^2 u_{2m}^2 \sum_{i,j} s_{ij}^2 = \sigma^4\cdot \sum_{k,m} u_{1k}^2 u_{2m}^2  \cdot \sum_{l > r} \frac{1}{\bar{\sigma}_l^2} =\sigma^4 \sum_{l>r} \frac{1}{\bar{\sigma}_l^2} \leq \frac{\sigma^4 (n-r)}{(\lambda_{p} -\lambda_{r+1})^2}.     
    \end{split}
\end{equation}
The second equation follows from the fact that $u_1, u_2$ are unit vectors.

Next, we bound the second sub-sum as follows. Similarly, using $s_{ij}$, we rewrite the second sub-sum into
\begin{equation*}
\begin{split}
\sum_{l,l'>r} \sum_{\substack{ \{k,i\} \neq\\\{j,m\}}} \sigma_{ki}^2 \sigma_{jm}^2 \frac{(u_{1k})^2 u_{2m} u_{2j} u_{li} u_{lj} u_{l'i} u_{l'm} }{\bar{\sigma}_l \bar{\sigma}_{l'}} &=  \sum_{k} u_{1k}^2 \sum_{i} \sigma_{ki}^2 \sum_{j,m} (u_{2j} s_{ij}) \sigma_{jm}^2 (u_{2m} s_{im}) \\
& = \sum_{k} u_{1k}^2 \sum_{i} \sum_{j,m} ( \sigma_{ki} u_{2j} s_{ij}) \sigma_{jm}^2 ( \sigma_{ki}  u_{2m} s_{im}).
\end{split}
\end{equation*}
For each $k$, consider the $n \times 1$ vectors $U_{i,k} = [\sigma_{ki} u_{2j} s_{ij}]^T_{1 \leq j \leq n}$ and the $n \times n$ matrix $\Sigma=(\sigma_{jm}^2)_{1 \leq j,m \leq n}.$ The absolute value of the RHS above can be rewritten as 
\begin{equation*} 
\begin{split}
|\sum_{k} u_{1k}^2 \sum_{i} \sum_{j,m} ( \sigma_{ki} u_{2j} s_{ij}) \sigma_{jm}^2 ( \sigma_{ki}  u_{2m} s_{im})| &= \sum_{k} u_{1k}^2 \sum_{i} |U_{i,k}^T \Sigma U_{i,k}| \leq \sum_{k} u_{1k}^2 \| \Sigma \|  \sum_{i} \|U_{i,k}\|^2.
\end{split}
\end{equation*}
Since $\|\Sigma\| \leq n \sigma^2$, the RHS is at most $ n\sigma^2 \sum_{k} u_{1k}^2 \sum_{i} \|U_{i,k}\|^2$, which can be expanded using the definitions of $s_{ij}$ and $U_{i,k}$ as
\begin{equation} \label{M_22}
    \begin{split}
 n \sigma^2 \sum_{k} u_{1k}^2 \sum_{i} \sum_{j} \sigma_{ki}^2 u_{2j}^2 s_{ij}^2 
& \leq n \sigma^4 \sum_{k} u_{1k}^2 \sum_{i} \sum_{j} \sum_{l,l' >r} u_{2j}^2 \frac{u_{li} u_{lj} u_{l'i} u_{l'j}}{\bar{\sigma}_l \bar{\sigma}_{l'}} \\
& = n \sigma^4 \sum_{j} u_{2j}^2 \sum_{l>r} \frac{ u_{lj}^2}{\bar{\sigma}_l^2} \,\,\,(\text{since $\sum_{i} u_{li} u_{l'i} =\textbf{1}_{l=l'}$}) \\
& \leq n \sigma^4 \sum_{j} u_{2j}^2 \frac{\sum_{l>r} u_{lj}^2}{(\lambda_p - \lambda_{r+1})^2} \leq \frac{n \sigma^4}{(\lambda_p - \lambda_{r+1})^2}.      
    \end{split}
\end{equation}

The third sub-sum is estimated by a similar argument. Finally, the fourth sub-sum can be bounded as follows. First, we rewrite its absolute value as 
\begin{equation*}
\begin{split}
\bigg|\sum_{l,l'>r} \sum_{\substack{ \{k,i\}\neq\\\{j,m\}}} \sigma_{ki}^2 \sigma_{jm}^2 \frac{u_{1k} u_{1i} u_{2m} u_{2j} u_{li} u_{lj} u_{l'k} u_{l'm} }{\bar{\sigma}_l \bar{\sigma}_{l'}} \bigg| & = \bigg|\sum_{k,i,j,m} u_{1k} \sigma_{ki}^2 u_{1i} u_{2m} \sigma_{jm}^2 u_{2j} s_{ij} s_{km}\bigg| \\
&=\bigg|\sum_{k,i} (u_{1k} \sigma_{ki}^2 u_{1i}) \sum_{j,m} (u_{2m} s_{km} \sigma_{jm}^2 (u_{2j} s_{ij})\bigg|.
\end{split}
\end{equation*}
Define 
$$\Sigma':=(\sigma_{ki}^2 r_{ki})_{1 \leq k,i\leq n}\,\,\,\text{and}\,\,\, r_{ki}:= \sum_{j,m} (u_{2m} s_{km} \sigma_{jm}^2 (u_{2j} s_{ij}).$$
The RHS is rewritten as
$$\bigg|\sum_{k,i} u_{1k} ( \sigma_{ki}^2 r_{ki})u_{1i}\bigg|= |u_1^T \Sigma' u_1| \leq  \| \Sigma'\| \leq  \|\Sigma'\|_{F}.$$
Moreover, since $\sigma_{ki} \leq \sigma$ for all $k,i$, we have
\begin{equation*}
\begin{split}
\|\Sigma'\|_{F}^2 & \leq \sigma^4 \sum_{1 \leq k,i \leq n} r_{ki}^2 = \sigma^4 \sum_{1 \leq k,i \leq n} |V_{2,k}^T \Sigma V_{2,i}|^2 \leq \sigma^4 \|\Sigma\|^2 \sum_{k,i}  \|V_{2,k}\|^2 \cdot \|V_{2,i}\|^2,
\end{split}
\end{equation*}
where 
$$V_{2,k}:=[ u_{2m} s_{km}]^T_{1 \leq m \leq n}, V_{2,i} :=[ u_{2m} s_{im}]^T_{1 \leq m \leq n}, \text{and}\,\,\Sigma=(\sigma_{jm}^2)_{1 \leq j,m \leq n}.$$
On the other hand, 
$$\sum_{k,i}  \|V_{2,k}\|^2 \cdot \|V_{2,i}\|^2 = (\sum_k \sum_{m} u_{2m}^2 s_{km}^2)^2.$$ 
Similarly to the above estimation, we further rewrite $\sum_k \sum_{m} u_{2m}^2 s_{km}^2$ into
\begin{equation*}
\begin{split}
 \sum_{m} u_{2m}^2 \sum_{k} s_{km}^2 =\sum_{m} u_{2m}^2 \sum_{l,l'>r} \frac{ (\sum_{k} u_{lk} u_{l'k}) u_{lm} u_{l'm}}{\bar{\sigma}_l \bar{\sigma}_{l'}} = \sum_{m} u_{2m}^2  \sum_{l >r} \frac{u_{lm}^2}{\bar{\sigma}_l^2} \leq \frac{\sum_{m} u_{2m}^2}{(\lambda_p -\lambda_{r+1})^2}.
\end{split}
\end{equation*}
Since $\|\Sigma\|^2 \leq n^2 \sigma^4$, we finally obtain that the absolute value of the fourth sub-sum is at most
\begin{equation} \label{M_23}
\|\Sigma'\| \leq \frac{n \sigma^4}{(\lambda_p -\lambda_{r+1})^2}.
\end{equation}
Combining (\ref{M_21}), (\ref{M_22}), the line after \eqref{M_22}, and (\ref{M_23}), we obtain
\begin{equation} \label{M_2}
|M_2| \leq \frac{4 n\sigma^4}{(\lambda_p -\lambda_{r+1})^2}. 
\end{equation}
By a similar argument, we also have 
\begin{equation} \label{M_3}
    |M_3| \leq \frac{4 n\sigma^4}{(\lambda_p -\lambda_{r+1})^2}. 
\end{equation}
\subsubsection{ Bounding $M_5$}
Denote $\sigma_{4,ki}=\E \xi_{ki}^4 - (\E \xi_{ki}^2)^2$ and $\bar{\sigma}_4 =\max_{1\leq k,i \leq n} \sigma_{4,ki}$. \\
We have
$M_5= \sum_{l,l'>r} \sum_{\substack{ \{k,i\}=\{j,m\}=\\
\{k',i'\}=\{j',m'\}}} (...)(\E \xi_{ki}^4 - (\E \xi_{ki}^2)^2)$, which equals
\begin{equation*}
\begin{split}
 \sum_{l,l' > r} \frac{1}{\bar{\sigma}_l \bar{\sigma}_{l'}} \sum_{k,i} \sigma_{4,ki} \times \big[& u_{1k}u_{2i} u_{li} u_{lk} u_{1k}u_{2i} u_{l'i} u_{l'k}+ \\  & u_{1k}u_{2i} u_{li} u_{lk} u_{1i}u_{2i} u_{l'k} u_{l'k}+ \\
& u_{1k}u_{2k} u_{li} u_{li} u_{1k}u_{2i} u_{l'i} u_{l'k}+\\
& u_{1k}u_{2i} u_{li} u_{lk} u_{1k}u_{2k} u_{l'i} u_{l'i}+\\
& u_{1k}u_{2i} u_{li} u_{lk} u_{1i}u_{2k} u_{l'k} u_{l'i}+\\
& u_{1k}u_{2k} u_{li} u_{li} u_{1k}u_{2k} u_{l'i} u_{l'i}+\\
& u_{1k}u_{2k} u_{li} u_{li} u_{1i}u_{2i} u_{l'k} u_{l'k}+ \\
& u_{1k}u_{2k} u_{li} u_{li} u_{1i}u_{2k} u_{l'k} u_{l'i}      \big]. 
\end{split}
\end{equation*}
We rearrange this expression into
\begin{equation*}
\begin{split}
\sum_{l,l' > r} \frac{1}{\bar{\sigma}_l \bar{\sigma}_{l'}} \sum_{k,i} \sigma_{4,ki} \times  
 & \big[(u_{1k})^2 u_{lk}   u_{l'k} (u_{2i})^2 u_{li} u_{l'i} + 
\\  & u_{1k}u_{lk} (u_{l'k})^2 (u_{2i})^2 u_{li}  u_{1i}  + \\
& (u_{1k})^2 u_{2k} u_{l'k} (u_{li})^2 u_{2i} u_{l'i} +\\
& (u_{1k})^2 u_{lk} u_{2k}  u_{2i} u_{li}  (u_{l'i})^2+\\
& (u_{1k}u_{lk}u_{2k}u_{l'k}) (u_{2i} u_{li}  u_{1i}  u_{l'i})+\\
& (u_{1k})^2 (u_{2k})^2 (u_{li})^2  (u_{l'i})^2 +\\
& u_{1k}u_{2k} (u_{l'k})^2 (u_{li})^2 u_{1i}u_{2i} + \\
& u_{1k} (u_{2k})^2 u_{l'k} (u_{li})^2 u_{1i}  u_{l'i}      \big]. 
\end{split}
\end{equation*}
Here, we have $8$ sub-sums. We show that the absolute value of each sub-sum is at most $\frac{n \bar{\sigma}_4}{(\lambda_p-\lambda_{r+1})^2}$. Combining these estimates with the fact that $\bar{\sigma}_4 \leq \max_{1 \leq k, i \leq n} \E \xi_{ki}^4 \leq  \sigma^2 K^2$, we obtain
\begin{equation} \label{M_5}
|M_5| \leq \frac{8n \bar{\sigma}_4}{(\lambda_p -\lambda_{r+1})^2} \leq \frac{8n \sigma^2 K^2}{(\lambda_p -\lambda_{r+1})^2}.
\end{equation}

Without loss of generality, we only present how to work with the first two sub-sums; the others can be estimated by a similar argument. We bound the first one as follows.  By the triangle inequality, we have
$$\bigg|\sum_{l,l' > r} \frac{1}{\bar{\sigma}_l \bar{\sigma}_{l'}} \sum_{k,i} \sigma_{4,ki} 
(u_{1k})^2 u_{lk}   u_{l'k} (u_{2i})^2 u_{li} u_{l'i}\bigg| \leq \sum_{l,l' > r} \frac{1}{\bar{\sigma}_l \bar{\sigma}_{l'}} \sum_{k,i} \sigma_{4,ki} 
  |u_{1k}|^2 |u_{lk}|   |u_{l'k}| |u_{2i}|^2 |u_{li}| |u_{l'i}|.$$
  The RHS can be rewritten as 
  $$\sum_{k,i} \sigma_{4,ki} |u_{1k}|^2 |u_{2i}|^2 \left( \sum_{l>r} \frac{|u_{lk}| |u_{li}|}{\bar{\sigma}_l} \right)^2  \leq \sum_{k,i} \sigma_{4,ki} |u_{1k}|^2 |u_{2i}|^2 \left( \sum_{l=1}^n \frac{|u_{lk}| |u_{li}|}{\bar{\sigma}_l} \right)^2.$$
By the Cauchy--Schwarz inequality, it is at most 
$$\sum_{k,i}\sigma_{4,ki} |u_{1k}|^2 |u_{2i}|^2 \frac{1}{(\lambda_p -\lambda_{r+1})^2} \left( \sum_{l} u_{lk}^2 \right) \left( \sum_{l} u_{li}^2 \right) = \sum_{k,i}\sigma_{4,ki} |u_{1k}|^2 |u_{2i}|^2 \frac{1}{(\lambda_p -\lambda_{r+1})^2}.$$
Since $\bar{\sigma}_4 =\max_{1\leq k,i \leq n} \sigma_{4,ki}$, the RHS is at most
$$ \bar{\sigma}_4 \sum_{k,i} |u_{1k}|^2 |u_{2i}|^2 \frac{1}{(\lambda_p -\lambda_{r+1})^2} =\frac{\bar{\sigma}_4 (\sum_{k} u_{1k}^2)(\sum_{i} u_{2i}^2)}{(\lambda_p -\lambda_{r+1})^2} = \frac{\bar{\sigma}_4}{(\lambda_p -\lambda_{r+1})^2}.$$

The second sub-sum can be bounded as follows. By the triangle inequality, we have
$$\bigg|\sum_{l,l' > r} \frac{1}{\bar{\sigma}_l \bar{\sigma}_{l'}} \sum_{k,i} \sigma_{4,ki}  u_{1k}u_{lk} (u_{l'k})^2 (u_{2i})^2 u_{li}  u_{1i} \bigg| \leq \sum_{l,l' > r}  \sum_{k,i} \sigma_{4,ki}  |u_{1k}| |u_{1i}| (u_{2i})^2 \frac{|u_{lk} | |u_{li}|}{\bar{\sigma}_l} \frac{(u_{l'k})^2}{\bar{\sigma}_{l'}}.$$
Since 
$$\sum_{l>r} \frac{|u_{lk} | |u_{li}|}{\bar{\sigma}_l} \leq \sum_{l=1}^n \frac{|u_{lk} | |u_{li}|}{\bar{\sigma}_l} \leq \frac{1}{(\lambda_p -\lambda_{r+1})},$$
$$ \sum_{l'>r}\frac{(u_{l'k})^2}{\bar{\sigma}_{l'}} \leq \frac{1}{(\lambda_p -\lambda_{r+1})},\,\,\text{and}\,\, \sigma_{4,ki} \leq \bar{\sigma}_4,$$
the RHS is at most
$$ \frac{\bar{\sigma}_4}{(\lambda_p -\lambda_{r+1})^2} \left( \sum_{k} |u_{1k}| \right) \left( \sum_{i} |u_{1i}| u_{2i}^2 \right) \leq \frac{\bar{\sigma}_4}{(\lambda_p -\lambda_{r+1})^2} \left( \sum_{k} |u_{1k}| \right) \left( \sum_{i} |u_{1i}| \right) \leq \frac{n \bar{\sigma}_4}{(\lambda_p -\lambda_{r+1})^2}.$$
The last inequality follows from the Cauchy--Schwarz inequality that
$$\sum_{k} |u_{1k}| \leq \sqrt{n}, \sum_{i} |u_{1i}| \leq \sqrt{n}.$$
The remaining sub-sums can be handled similarly and also bounded by $\frac{n \bar{\sigma}_4}{(\lambda_p -\lambda_{r+1})^2}$. We complete the estimation of $|M_5|$.

Combining \eqref{G12power2split}, \eqref{M_2}, \eqref{M_3}, and \eqref{M_5} , we finally obtain 
$$ \E G_{12,\bar p}^2 \leq (\E G_{12,\bar{p}})^2 + M_2+M_3+M_5  \leq (\E G_{12,\bar{p}})^2 + \frac{n \sigma^2 (8\sigma^2+8K^2)}{(\lambda_p -\lambda_{r+1})^2}.$$
This proves \eqref{second moment G12p}.

\section{ Proof of the third main lemma - Lemma \ref{LemmaCase3M(a,b)bound}} \label{Appx: proofLemmaCase3}
In the setting of this lemma,  $s_1 \geq k$. Without loss of generality, it is sufficient to handle the case that $\alpha_1 =0, \alpha_{k+1} > 0$. Similar to Subsection \ref{section:proof2} and Subsection \ref{subsec: proofCase2}, we have 
\begin{equation*}
  \begin{split}
\frac{1}{2\pi \textbf{i}} \int_{\Gamma} M(\alpha; \beta) dz & = \sum_{i_1,i_2,...,i_{s_2}} \sum_{L \in \mathcal{L}(T,s)} \sum_{j_1,...,j_{s_1}}  (-1)^{T+1} \prod_{\substack{e \in \mathcal{E}_L \\ e^{-} \in I_2}} \frac{1 }{ (\lambda_{e^+} - \lambda_{e^-})} \times \left[ \prod_{l=1}^{\beta_1} \left( u_{i_{l}} u_{i_{l}}^T E \right) \right] \times   \\
& \bigg\{ \prod_{h=1}^{k-1} \bigg[ \prod_{l=1}^{\alpha_{h+1}} \bigg( \frac{u_{j_{l+A_{h}}} u_{j_{l+A_{h}}}^T}{ \displaystyle\prod_{\substack{e \in \mathcal{E}_L;\, e^{-} = j_{l+A_{h}}}} (\lambda_{e^{+}}- \lambda_{j_{l+A_{h}}})} E  \bigg) \bigg] \times \bigg[  \prod_{l=1}^{\beta_{h+1}} \left( u_{i_{l+B_{h}}} u_{i_{l+B_{h}}}^T E \right) \bigg] \bigg\} \times \\
& \times \bigg[  \prod_{l =1}^{\alpha_{k+1}-1} \bigg( \frac{u_{j_{l + A_k}} u_{j_{l + A_k}}^T}{\displaystyle\prod_{\substack{e \in \mathcal{E}_L\\ e^{-} = j_{l+A_{k}}}} (\lambda_{e^{+}}- \lambda_{j_{l+A_{k}}})} E \bigg) \bigg]\frac{u_{j_{s_1}} u_{j_{s_1}}^T}{\displaystyle\prod_{\substack{e \in \mathcal{E}_L \\ e^{-}= j_{s_1}}}(\lambda_{e^+}-\lambda_{j_{s_1}})}.
\end{split}  
\end{equation*}
Once again, we distribute the sum over $j_1, j_2, ...,j_{s_2}$ according to  each position of $u_{j_l} u_{j_l}^T$ in the product, and obtain
\[
\frac{1}{2\pi \textbf{i}} \int_{\Gamma} M(\alpha; \beta) dz  = \sum_{i_1,i_2,...,i_{s_2}} \sum_{L \in \mathcal{L}(T,s)} (-1)^{T+1} \prod_{\substack{e \in \mathcal{E}_L \\ e^{-} \in I_2}} \frac{1 }{ (\lambda_{e^+} - \lambda_{e^-})} \times \left[ \prod_{l=1}^{\beta_1} \left( u_{i_{l}} u_{i_{l}}^T E \right) \right] \times 
\]
\[
\bigg\{ \prod_{h=1}^{k-1} \bigg[ \prod_{l=1}^{\alpha_{h+1}} \bigg( \sum_{j_{l+A_{h}} \in N_{\bar{\lambda}}(S)^c} \frac{u_{j_{l+A_{h}}} u_{j_{l+A_{h}}}^T}{ \displaystyle\prod_{\substack{e \in \mathcal{E}_L;\, e^{-} = j_{l+A_{h}}}} (\lambda_{e^{+}}- \lambda_{j_{l+A_{h}}})} E  \bigg) \bigg] \times \bigg[  \prod_{l=1}^{\beta_{h+1}} \left( u_{i_{l+B_{h}}} u_{i_{l+B_{h}}}^T E \right) \bigg] \bigg\} \times 
\]
\[
\times \bigg[ \prod_{l =1}^{\alpha_{k+1}-1} \bigg( \sum_{j_{l+A_k} \in N_{\bar{\lambda}}(S)^c} \frac{u_{j_{l + A_k}} u_{j_{l + A_k}}^T}{\displaystyle\prod_{\substack{e \in \mathcal{E}_L \\ e^{-} = j_{l+A_{k}}}} (\lambda_{e^{+}}- \lambda_{j_{l+A_{k}}})} E \bigg) \bigg] \times \bigg(\sum_{j_{s_1} \in N_{\bar{\lambda}}(S)^c} \frac{u_{j_{s_1}} u_{j_{s_1}}^T}{\displaystyle\prod_{\substack{e \in \mathcal{E}_L \\ e^{-}= j_{s_1}}}(\lambda_{e^+}-\lambda_{j_{s_1}})} \bigg).
\]
Next, by \eqref{normdefinition}, we rewrite $\frac{1}{2\pi } \Norm{\int_{\Gamma} M(\alpha; \beta) dz} $ as
\[
\max_{\|\textbf{v}\|=\|\textbf{w}\|=1}  \sum_{i_1,i_2,...,i_{s_2}} \sum_{L \in \mathcal{L}(T,s)} (-1)^{T+1} \prod_{\substack{e \in \mathcal{E}_L \\ e^{-} \in I_2}} \frac{1 }{ (\lambda_{e^+} - \lambda_{e^-})} \times \left[ \textbf{v}^T \prod_{l=1}^{\beta_1} \left( u_{i_{l}} u_{i_{l}}^T E \right) \right] \times  
\]
\[
\bigg\{ \prod_{h=1}^{k-1} \bigg[ \prod_{l=1}^{\alpha_{h+1}} \bigg( \sum_{j_{l+A_{h}} \in N_{\bar{\lambda}}(S)^c} \frac{u_{j_{l+A_{h}}} u_{j_{l+A_{h}}}^T}{ \displaystyle\prod_{\substack{e \in \mathcal{E}_L;\, e^{-} = j_{l+A_{h}}}} (\lambda_{e^{+}}- \lambda_{j_{l+A_{h}}})} E  \bigg) \bigg] \times \bigg[  \prod_{l=1}^{\beta_{h+1}} \bigg( u_{i_{l+B_{h}}} u_{i_{l+B_{h}}}^T E \bigg) \bigg] \bigg\} \times
\]
\[
\times \bigg[ \prod_{l =1}^{\alpha_{k+1}-1} \bigg( \sum_{j_{l+A_k} \in N_{\bar{\lambda}}(S)^c} \frac{u_{j_{l + A_k}} u_{j_{l + A_k}}^T}{\displaystyle\prod_{\substack{e \in \mathcal{E}_L \\ e^{-} = j_{l+A_{k}}}} (\lambda_{e^{+}}- \lambda_{j_{l+A_{k}}})} E \bigg) \bigg] \times \bigg(\sum_{j_{s_1} \in N_{\bar{\lambda}}(S)^c} \frac{u_{j_{s_1}} u_{j_{s_1}}^T}{\displaystyle\prod_{\substack{e \in \mathcal{E}_L \\ e^{-}= j_{s_1}}}(\lambda_{e^+}-\lambda_{j_{s_1}})} \bigg) \textbf{w}.
\]
\noindent By the triangle inequality, the RHS is at most 
\begin{equation*}
    \begin{split}
 & \max_{\|\textbf{v}\|=\|\textbf{w}\|=1} \sum_{i_1,i_2,...,i_{s_2}} \sum_{L \in \mathcal{L}(T,s)} \norm{ \prod_{\substack{e \in \mathcal{E}_L \\ e^{-} \in I_2}} \frac{1 }{ (\lambda_{e^+} - \lambda_{e^-})}} \times \bigg\| \textbf{v}^T \prod_{l=1}^{\beta_1} \left( u_{i_{l}} u_{i_{l}}^T E \right) \bigg\| \times     \\
& \bigg\| \prod_{h=1}^{k-1} \bigg[ \prod_{l=1}^{\alpha_{h+1}} \bigg( \sum_{j_{l+A_{h}} \in N_{\bar{\lambda}}(S)^c} \frac{u_{j_{l+A_{h}}} u_{j_{l+A_{h}}}^T}{ \displaystyle\prod_{\substack{e \in \mathcal{E}_L;\, e^{-} = j_{l+A_{h}}}} (\lambda_{e^{+}}- \lambda_{j_{l+A_{h}}})} E  \bigg) \bigg] \times \bigg[  \prod_{l=1}^{\beta_{h+1}} \bigg( u_{i_{l+B_{h}}} u_{i_{l+B_{h}}}^T E \bigg) \bigg] \bigg\| \times \\
& \times \bigg\| \prod_{l =1}^{\alpha_{k+1}-1} \bigg( \sum_{j_{l+A_k} \in N_{\bar{\lambda}}(S)^c} \frac{u_{j_{l + A_k}} u_{j_{l + A_k}}^T}{\displaystyle\prod_{\substack{e \in \mathcal{E}_L \\ e^{-} = j_{l+A_{k}}}} (\lambda_{e^{+}}- \lambda_{j_{l+A_{k}}})} E \bigg)\bigg\| \times \bigg\|\sum_{j_{s_1} \in N_{\bar{\lambda}}(S)^c} \frac{u_{j_{s_1}} u_{j_{s_1}}^T}{\displaystyle\prod_{\substack{e \in \mathcal{E}_L \\ e^{-}= j_{s_1}}}(\lambda_{e^+}-\lambda_{j_{s_1}})} \textbf{w}\bigg\|.
 \end{split}
\end{equation*}
Since $\|w\|=1$, this is at most 
\begin{equation*}
    \begin{split}
&  \max_{\|\textbf{v}\|=1} \sum_{i_1,i_2,...,i_{s_2}} \sum_{L \in \mathcal{L}(T,s)} \norm{ \prod_{\substack{e \in \mathcal{E}_L \\ e^{-} \in I_2}} \frac{1 }{ (\lambda_{e^+} - \lambda_{e^-})}} \times \bigg\| \textbf{v}^T \prod_{l=1}^{\beta_1} \left( u_{i_{l}} u_{i_{l}}^T E \right) \bigg\| \times     \\
& \bigg\|\prod_{h=1}^{k-1} \bigg[ \prod_{l=1}^{\alpha_{h+1}} \bigg( \sum_{j_{l+A_{h}} \in N_{\bar{\lambda}}(S)^c} \frac{u_{j_{l+A_{h}}} u_{j_{l+A_{h}}}^T}{ \displaystyle\prod_{\substack{e \in \mathcal{E}_L;\, e^{-} = j_{l+A_{h}}}} (\lambda_{e^{+}}- \lambda_{j_{l+A_{h}}})} E  \bigg) \bigg] \times \bigg[  \prod_{l=1}^{\beta_{h+1}} \left( u_{i_{l+B_{h}}} u_{i_{l+B_{h}}}^T E \right) \bigg] \bigg\| \times \\
& \times \bigg\| \prod_{l =1}^{\alpha_{k+1}-1} \bigg( \sum_{j_{l+A_k} \in N_{\bar{\lambda}}(S)^c} \frac{u_{j_{l + A_k}} u_{j_{l + A_k}}^T}{\displaystyle\prod_{\substack{e \in \mathcal{E}_L \\ e^{-} = j_{l+A_{k}}}} (\lambda_{e^{+}}- \lambda_{j_{l+A_{k}}})} E \bigg)\bigg\| \times \bigg\|\sum_{j_{s_1} \in N_{\bar{\lambda}}(S)^c} \frac{u_{j_{s_1}} u_{j_{s_1}}^T}{\displaystyle\prod_{\substack{e \in \mathcal{E}_L \\ e^{-}= j_{s_1}}}(\lambda_{e^+}-\lambda_{j_{s_1}})}\bigg\|. 
    \end{split}
\end{equation*}
Using  the estimates in \eqref{Factor1,case2}, \eqref{Factor2,case2}, \eqref{Factor3,case2}, and \eqref{Factor34case1}, we obtain 
\begin{equation} \label{Mstep1boundcase3}
\begin{split} 
  \frac{1}{2 \pi } \Norm{\int_{\Gamma} M(\alpha;\beta) dz} & \leq \max_{\|\textbf{v}\|=1} \sum_{i_1,...,i_{s_2}} \sum_{L \in \mathcal{L}(T,s)}  \norm{\textbf{v}^T u_{i_1}}\frac{w^{2k-1}  x^{\sum_{l=1}^k (\beta_l-1)} \|E\|^{\sum_{l=2}^{k+1} (\alpha_l-1)}}{\delta_S^{|E(I_2)|} \bar{\lambda}^{\sum_{l=1}^{s_1} d(j_l)}}  \\
  & = \max_{\|\textbf{v}\|=1}  \sum_{i_1,...,i_{s_2}} \sum_{L \in \mathcal{L}(T,s)}  \norm{\textbf{v}^T u_{i_1}}\frac{\|E\|^{s_1-k} w^{2k-1} x^{s_2 - k}}{\delta_S^{s_2 -1-r_c(T,L)} \bar{\lambda}^{s_1+r_c(T,L)}} \\
  & = \max_{\|\textbf{v}\|=1} \sum_{T=1}^{s_2} \sum_{L \in \mathcal{L}(T,s)}  \left[ \sum_{\substack{i_1,...,i_{s_2} \\ |X(i_1,...,i_{s_2})|=T}} \norm{\textbf{v}^T u_{i_1}} \right]  \frac{\|E\|^{s_1-k} w^{2k-1} x^{s_2 - k}}{\delta_S^{s_2 -1-r_c(T,L)} \bar{\lambda}^{s_1+r_c(T,L)}}.
\end{split}
     \end{equation}
Moreover, similar to the argument in \eqref{vUUwboundArgument}, by the Cauchy--Schwarz inequality, we have
\begin{equation*} 
    \begin{split} 
\sum_{\substack{i_1,...,i_{s_2} \\ |X(i_1,...,i_{s_2})|=T}} \norm{\textbf{v}^T u_{i_1}} & = \sum_{\substack{i_1,...,i_{s_2} \\ |X(i_1,...,i_{s_2})|=T \\ i_1 \in S}} \norm{\textbf{v}^T u_{i_1}} + \sum_{\substack{i_1,...,i_{s_2} \\ |X(i_1,...,i_{s_2})|=T \\ i_1 \in N_{\bar{\lambda}}(S) \setminus S}} \norm{\textbf{v}^T u_{i_1}} \\
& \leq \binom{s_2-1}{T-1} \sqrt{p} p^{T-1} r^{s_2-T} + \binom{s_2-1}{T} \sqrt{r} p^{T} r^{s_2-T-1}.
\end{split}
\end{equation*}
Since $\sqrt{r} p^{T} r^{s_2-T-1} = p^T r^{s_2-T-1/2} \leq p^{T-1/2} r^{s_2-T}$ and $\binom{s_2 -1}{T-1}+\binom{s_2-1}{T} = \binom{s_2}{T}$, we have
\begin{equation} \label{vUboundcase3}
 \sum_{\substack{i_1,...,i_{s_2} \\ |X(i_1,...,i_{s_2})|=T}} \norm{\textbf{v}^T u_{i_1}} \leq \binom{s_2}{T} p^{T-1/2} r^{s_2-T}.   
\end{equation}
\noindent Combining \eqref{Mstep1boundcase3} and \eqref{vUboundcase3}, we obtain 
\begin{equation} \label{splitM(a,b)lastcase}
\begin{split}
\frac{1}{2 \pi}\Norm{\int_\Gamma M(\alpha;\beta) dz} & \leq  \sum_{T=1}^{s_2} \binom{s_2}{T} p^{T-1/2} r^{s_2 -T} \sum_{L \in \mathcal{L}(T,s)} \frac{\|E\|^{s_1-k} w^{2k} x^{s_2 - k} }{\delta_S^{s_2 -1-r_c} \bar{\lambda}^{s_1+r_c}}  \\
& = M_1(\alpha;\beta) + M_2(\alpha;\beta), {\rm where} \\
M_1(\alpha;\beta) & :=  \sum_{T=1}^{s_2} \binom{s_2}{T} p^{T-1/2} r^{s_2 -T} \sum_{\substack{L \in \mathcal{L}(T,s);\, r_c(T,L) \leq k-1}} \frac{\|E\|^{s_1-k} w^{2k-1} x^{s_2 - k} }{\delta_S^{s_2 -1-r_c(T,L)} \bar{\lambda}^{s_1+r_c(T,L)}} , \\
M_2 (\alpha;\beta) & :=  \sum_{T=1}^{s_2} \binom{s_2}{T} p^{T-1/2} r^{s_2 -T}\sum_{\substack{L \in \mathcal{L}(T,s);\, r_c(T,L) > k-1}} \frac{ \|E\|^{s_1-k} w^{2k-1} x^{s_2 - k}}{\delta_S^{s_2 -1-r_c(T,L)} \bar{\lambda}^{s_1+r_c(T,L)}} .
\end{split}
\end{equation}
We obtain the following lemmas, whose proofs are deferred to the next appendix (Appendix \ref{section: sublemmas}). 
\begin{lemma} \label{M_1Case3boundfinal} For any pair $(\alpha;\beta)$ of Type III,
    \begin{equation}
        M_1(\alpha;\beta) \leq \frac { \sqrt{p}\|E\|}{\bar{\lambda}} \times \frac{2^{s_2+s-1}}{12^{s-1}}.
    \end{equation}
\end{lemma}
\begin{lemma} \label{M_2Case3boundfinal} For any pair $(\alpha;\beta)$ of Type III,
  \begin{equation}
      M_2(\alpha;\beta) \leq \frac{p x}{\bar{\lambda}} \frac{2^{s_2+s-1}}{12^{s-1}}.
  \end{equation}  
\end{lemma}
Combining \eqref{splitM(a,b)lastcase}, Lemma \ref{M_2Case3boundfinal} and Lemma \ref{M_1Case3boundfinal}, we finally have 
\begin{equation}
\begin{split}
 \frac{1}{2 \pi} \Norm{\int_{\Gamma} M(\alpha;\beta) dz} &\leq M_1(\alpha;\beta)+ M_2(\alpha;\beta) \\ 
 & \leq \left( \frac { \sqrt{p}\|E\|}{\bar{\lambda}}  + \frac{p x}{\bar{\lambda}} \right) \times \frac{2^{s_2+s-1}}{12^{s-1}}.
\end{split}
\end{equation} 
This proves Lemma \ref{LemmaCase3M(a,b)bound}. 

\section{Proofs of some technical claims} \label{section: sublemmas}
In this section, we prove the technical results mentioned in Section \ref{section:randomnoise}, Section \ref{section: proof of main lemmas}, and Appendix \ref{Appx: proofLemmaCase3}: Lemmas \ref{M1(a,b)case2part1}, \ref{M2(a,b)case2part1}, \ref{M_1Case3boundfinal}, \ref{M_2Case3boundfinal}, and Corollary \ref{lowrank-random1}.

\subsection{Proof of Lemma \ref{M1(a,b)case2part1}}
We can bound $M_1(\alpha;\beta)$ in this case as follows. We replace $w^{2(k-1)}$ by $w^{2(k-1-r_c(T,L)} \|E\|^{2r_c(T,L)}$ and obtain 
$$ M_1(\alpha;\beta) \leq \sum_{T=2}^{s_2} \binom{s_2}{T} p^{T-1} r^{s_2 -T} \sum_{\substack{L \in \mathcal{L}(T,s);\, 1 \leq r_c(T,L) \leq k-1}} \frac{\|E\|^{s_1-(k-1)+2r_c(T,L)} w^{2(k-1-r_{c}(T,L))} x^{s_2 - k} }{\delta_S^{s_2 -1-r_c(T,L)} \bar{\lambda}^{s_1+r_c(T,L)}}.   $$
Similar to the argument in \eqref{Mab-bound-HR}, we can arrange the factors in the RHS and get
\begin{equation}\label{M_3Case2splitsum}
 \begin{split}
 M_1(\alpha;\beta) & \leq  \sum_{T=2}^{s_2} \binom{s_2}{T}   \sum_{\substack{L \in \mathcal{L}(T,s);\, 1\leq r_c(T,L) \leq k-1}} p^{T-1} r^{s_2-T} \left( \frac{ \|E\|}{\bar{\lambda}} \right)^{T_1} \left( \frac{ x}{\delta_S} \right)^{T_2} \left( \frac{ w}{\sqrt {\bar{\lambda} \delta_S}  } \right)^{T_3}  \\
 & =  \sum_{T=2}^{s_2} \binom{s_2}{T}   \sum_{\substack{L \in \mathcal{L}(T,s);\, 1 \leq r_c(T,L) \leq k-1}} \left( \frac{p}{r} \right)^{T-(r_c(T,L)+1)} p^{\frac{k-1 -s_1}{2}}  \left( \frac{ \sqrt{p}\|E\|}{\bar{\lambda}} \right)^{T_1} \left( \frac{ rx}{\delta_S} \right)^{T_2} \left( \frac{ \sqrt r w}{\sqrt {\bar{\lambda} \delta_S}  } \right)^{T_3},
 \end{split}
  \end{equation}
  where 
$$T_1= s_1 - (k-1)+2r_c(T,L),$$
$$T_2= s_2 - k,$$
$$T_3= 2(k-1 -r_c(T,L)).$$
Since $(\alpha,\beta)$ is type II, $s_1 \geq k-1$ and hence $0 \leq p^{\frac{k-1-s_1}{2}} \leq 1$. By Lemma \ref{r_cLemma}, we also have $0 \leq \left( \frac{p}{r} \right)^{T-(r_c(T,L)+1)} \leq 1$. We can omit these factors and obtain 
\begin{equation}
    \begin{split}
M_1(\alpha;\beta)&  \leq \sum_{T=2}^{s_2} \binom{s_2}{T}   \sum_{\substack{L \in \mathcal{L}(T,s);\, 1 \leq r_c(T,L) \leq k-1}}  \left( \frac{ \sqrt{p}\|E\|}{\bar{\lambda}} \right)^{T_1} \left( \frac{ rx}{\delta_S} \right)^{T_2} \left( \frac{ \sqrt r w}{\sqrt {\bar{\lambda} \delta_S}  } \right)^{T_3} \\
& = \sum_{T=2}^{s_2} \binom{s_2}{T}   \sum_{\substack{L \in \mathcal{L}(T,s);\, 1 \leq r_c(T,L) \leq k-1}} \frac{\sqrt{p}\|E\|}{\bar{\lambda}} \times  \left( \frac{ \sqrt{p}\|E\|}{\bar{\lambda}} \right)^{T_1-1} \left( \frac{ rx}{\delta_S} \right)^{T_2} \left( \frac{ \sqrt r w}{\sqrt {\bar{\lambda} \delta_S}  } \right)^{T_3}.
    \end{split}
\end{equation}
 By the setting of $M_1(\alpha;\beta)$ for $(\alpha,\beta)$ of type II, we have $T_2,T_3 \geq 0$. Particularly, $T_1 \geq 2 r_c(T,L) \geq 2$. Therefore, by Lemma \ref{toytrick}, $M_1(\alpha;\beta)$ is less than or equal to 
\begin{equation*}
\begin{split}
& \sum_{T=2}^{s_2} \binom{s_2}{T}   \sum_{\substack{L \in \mathcal{L}(T,s) \\ 1 \leq r_c(T,L) \leq k-1}}   \left( \frac{ \sqrt{p}\|E\|}{\bar{\lambda}} \right)  \max \left\lbrace \frac{\|E\|}{\bar{\lambda}}, \frac{rx}{\delta_S}, \frac{\sqrt{r}w}{ \sqrt{\bar{\lambda} \delta_S}} \right\rbrace^{T_1+T_2+T_3-1}\\
 & = \sum_{T=2}^{s_2} \binom{s_2}{T}   \sum_{\substack{L \in \mathcal{L}(T,s) \\ 1 \leq r_c(T,L) \leq k-1}}   \left( \frac{ \sqrt{p}\|E\|}{\bar{\lambda}} \right)  \max \left\lbrace \frac{\|E\|}{\bar{\lambda}}, \frac{rx}{\delta_S}, \frac{\sqrt{r}w}{ \sqrt{\bar{\lambda} \delta_S}} \right\rbrace^{s-1},   
\end{split}
\end{equation*}
since $T_1+T_2+T_3=s$. Moreover, by Assumption \textbf{D0}, this RHS is at most 
\begin{equation*} 
  \left( \frac{ \sqrt{p}\|E\|}{\bar{\lambda}} \right) \frac{\sum_{T=2}^{s_2}   \displaystyle \sum_{\substack{L \in \mathcal{L}(T,s);\, 1 \leq r_c(T,L) \leq k-1}}\binom{s_2}{T}}{12^{s-1}}.  
\end{equation*} 
This proves Lemma \ref{M1(a,b)case2part1}.

\subsection{Proof of Lemma \ref{M2(a,b)case2part1}}
We bound $M_2(\alpha;\beta)$ in this case as follows. Replacing $w^{2(k-1)}$ by $\|E\|^{2(k-1)}$, we obtain
 $$M_2(\alpha;\beta) \leq \sum_{T=2}^{s_2} \binom{s_2}{T} p^{T-1} r^{s_2 -T}\sum_{\substack{L \in \mathcal{L}(T,s);\, r_c(T,L) > k-1}} \frac{\|E\|^{s_1+(k-1)} x^{s_2 - k} }{\delta_S^{s_2 -1-r_c(T,L)} \bar{\lambda}^{s_1+r_c(T,L)}}, $$
 whose RHS can be rewritten as 
 $$\sum_{T=2}^{s_2} \binom{s_2}{T} \sum_{\substack{L \in \mathcal{L}(T,s);\, r_c(T,L) > k-1}} p^{T-1} r^{s_2 -T}  \left( \frac{ x}{\delta_S} \right)^{s_2-1-r_c(T,L)} \left(\frac{x}{\bar{\lambda}} \right)^{r_c(T,L)-(k-1)} \left( \frac{  \|E\|}{\bar{\lambda}} \right)^{k-1+s_1}. $$
\noindent Distributing each factors $r,p, \sqrt{p}$ to $\frac{x}{\delta_S}, \frac{x}{\bar{\lambda}}, \frac{\|E\|}{\bar{\lambda}}$, respectively, we can rewrite the RHS as 
$$ \sum_{T=2}^{s_2} \binom{s_2}{T} \sum_{\substack{L \in \mathcal{L}(T,s) \\ r_c(T,L) > k-1}} \left( \frac{p}{r}  \right)^{T-(r_c(T,L)+1)} p^{\frac{k-1 -s_1}{2}} \left( \frac{r x}{\delta_S} \right)^{s_2-1-r_c(T,L)} \left(\frac{px}{\bar{\lambda}} \right)^{r_c(T,L)-(k-1)} \left( \frac{ \sqrt{p} \|E\|}{\bar{\lambda}} \right)^{k-1+s_1}. $$
\noindent Since $T \geq r_c(T,L)+1$ (Lemma  \ref{r_cLemma}) and $s_1-(k-1) \geq 0$ ($(\alpha,\beta)$ is Type II), 
$$ 0 \leq \left( \frac{p}{r} \right)^{T-(r_c(T,L)+1)} \leq 1 \,\,\,\text{and}\,\,\, 0 \leq  p^{\frac{k-1 -s_1}{2}}  \leq 1.$$
We can omit these factors and obtain 
\begin{equation*}
    \begin{split}
 M_2(\alpha;\beta) & \leq \sum_{T=2}^{s_2} \binom{s_2}{T} \sum_{\substack{L \in \mathcal{L}(T,s) \\ r_c(T,L) > k-1}} \left( \frac{r x}{\delta_S} \right)^{s_2-1-r_c(T,L)} \left(\frac{px}{\bar{\lambda}} \right)^{r_c(T,L)-(k-1)} \left( \frac{ \sqrt{p} \|E\|}{\bar{\lambda}} \right)^{k-1+s_1} \\
 & = \sum_{T=2}^{s_2} \binom{s_2}{T} \sum_{\substack{L \in \mathcal{L}(T,s) \\ r_c(T,L) > k-1}} \frac{px}{\bar{\lambda}} \times  \left( \frac{r x}{\delta_S} \right)^{s_2-1-r_c(T,L)} \left(\frac{px}{\bar{\lambda}} \right)^{r_c(T,L)-(k-1)-1} \left( \frac{ \sqrt{p} \|E\|}{\bar{\lambda}} \right)^{k-1+s_1}.
    \end{split}
\end{equation*}
\noindent By the setting of $M_2(\alpha; \beta)$, $r_c(T,L) - (k-1) -1\geq 0$. By definitions of $s_1,s_2, r_c(T,L)$, we also have $s_2 -1-r_c(T,L), k-1+s_1 \geq 0$. Applying Lemma \ref{toytrick} and using 
\[
s_2-1-r_c(T,L)+r_c(T,L)-(k-1)-1+k-1+s_1 = s_1+s_2-2= s -1,
\]
we finally obtain
\begin{equation}
\begin{split} 
M_2(\alpha;\beta)
& = \sum_{T=2}^{s_2} \binom{s_2}{T} \sum_{\substack{L \in \mathcal{L}(T,s) \\ r_c(T,L) > k-1}} \frac{px}{\bar{\lambda}} \times \max \left\lbrace \frac{rx}{\delta_S}, \frac{px}{\bar{\lambda}}, \frac{\sqrt{p}\|E\|}{\bar{\lambda}} \right\rbrace^{s-1} \\ 
 & \leq  \sum_{T=2}^{s_2} \binom{s_2}{T} \sum_{\substack{L \in \mathcal{L}(T,s);\, r_c(T,L) > k-1}}   \frac{p x}{\bar{\lambda}} \max \left\lbrace \frac{rx}{\delta_S}, \frac{\sqrt{p}\|E\|}{\bar{\lambda}} \right\rbrace^{s-1}  (\text{since}\,\, \frac{px}{\bar{\lambda}} \leq \frac{rx}{\delta_S}) \\
& \leq \frac{p x}{\bar{\lambda}} \frac{ \sum_{T=2}^{s_2}  \displaystyle\sum_{\substack{L \in \mathcal{L}(T,s);\,r_c > k-1}} \binom{s_2}{T} }{12^{s-1}} \,\,(\text{by Assumption}\,\,\textbf{D0}).
\end{split}
\end{equation}
This proves Lemma \ref{M2(a,b)case2part1}.
\subsection{Proof of Lemma \ref{M_1Case3boundfinal}}
\noindent  We bound  $M_1(\alpha;\beta)$ as follows. Similarly, replacing $w^{2k-1}$ by $w^{2k-2-2r_c(T,L)} \|E\|^{2r_c(T,L)+1}$, we have 
$$M_1(\alpha;\beta) \leq \sum_{T=1}^{s_2} \binom{s_2}{T} p^{T-1/2} r^{s_2 -T} \sum_{\substack{L \in \mathcal{L}(T,s) \\ r_c(T,L) \leq k-1}} \frac{\|E\|^{s_1-(k-1)+2r_c(T,L)} w^{2(k-1-r_c(T,L))} x^{s_2 - k} }{\delta_S^{s_2 -1-r_c(T,L)} \bar{\lambda}^{s_1+r_c(T,L)}}.$$
Rearranging the factors in the RHS, we further obtain
\begin{equation*} 
 \begin{split}
 M_1(\alpha;\beta) & \leq  \sum_{T=1}^{s_2} \binom{s_2}{T}   \sum_{\substack{L \in \mathcal{L}(T,s);\, r_c(T,L) \leq k-1}} p^{T-1/2} r^{s_2-T} \left( \frac{ \|E\|}{\bar{\lambda}} \right)^{T_1} \left( \frac{ x}{\delta_S} \right)^{T_2} \left( \frac{ w}{\sqrt {\bar{\lambda} \delta_S}  } \right)^{T_3}  \\
 & =   \sum_{T=1}^{s_2} \binom{s_2}{T}   \sum_{\substack{L \in \mathcal{L}(T,s);\, r_c(T,L) \leq k-1}} \left( \frac{p}{r} \right)^{T-(r_c(T,L)+1)} p^{\frac{k -s_1}{2}}  \left( \frac{ \sqrt{p}\|E\|}{\bar{\lambda}} \right)^{T_1} \left( \frac{ rx}{\delta_S} \right)^{T_2} \left( \frac{ \sqrt r w}{\sqrt {\bar{\lambda} \delta_S}  } \right)^{T_3},
 \end{split}
  \end{equation*}
where 
$$T_1= s_1 - (k-1)+2r_c(T,L),$$
$$T_2= s_2 - k,$$
$$T_3= 2(k-1 -r_c(T,L)).$$
Since $T \geq r_c(T,L)+1, s_1 \geq k$, and then $0 \leq \left( \frac{p}{r} \right)^{T-(r_c(T,L)+1)} p^{\frac{k -s_1}{2}} \leq 1$,  we have
\begin{equation*}
    \begin{split}
  M_1(\alpha;\beta) & \leq \sum_{T=1}^{s_2} \binom{s_2}{T}   \sum_{\substack{L \in \mathcal{L}(T,s);\, r_c(T,L) \leq k-1}}   \left( \frac{ \sqrt{p}\|E\|}{\bar{\lambda}} \right)^{T_1} \left( \frac{ rx}{\delta_S} \right)^{T_2} \left( \frac{ \sqrt r w}{\sqrt {\bar{\lambda} \delta_S}  } \right)^{T_3} \\
  & = \sum_{T=1}^{s_2} \binom{s_2}{T}   \sum_{\substack{L \in \mathcal{L}(T,s);\, r_c(T,L) \leq k-1}} \frac{ \sqrt{p}\|E\|}{\bar{\lambda}} \times  \left( \frac{ \sqrt{p}\|E\|}{\bar{\lambda}} \right)^{T_1-1} \left( \frac{ rx}{\delta_S} \right)^{T_2} \left( \frac{ \sqrt r w}{\sqrt {\bar{\lambda} \delta_S}  } \right)^{T_3}.
    \end{split}
\end{equation*}
By the setting of $M_1(\alpha;\beta)$,  $r_c(T,L) \leq k-1$ and hence $T_3 \geq 0$. We also have $T_2 \geq 0$ by definition of $s_2$. Since $(\alpha,\beta)$ is type III, $s_1 \geq k$ and hence $T_1 -1 \geq 0$. Therefore, by Lemma \ref{toytrick}, we finally obtain  
\begin{equation*} 
\begin{split}
  M_1(\alpha;\beta) & \leq  \sum_{T=1}^{s_2} \binom{s_2}{T}   \sum_{\substack{L \in \mathcal{L}(T,s);\, r_c(T,L) \leq k-1}}  \frac { \sqrt{p}\|E\|}{\bar{\lambda}} \times \left\lbrace \frac{\sqrt{p}\|E\|}{\bar{\lambda}}, \frac{rx}{\delta_S}, \frac{\sqrt{r} w}{\sqrt{\bar{\lambda} \delta_S}} \right\rbrace^{T_1+T_2+T_3-1}\\
  & = \sum_{T=1}^{s_2} \binom{s_2}{T}   \sum_{\substack{L \in \mathcal{L}(T,s);\, r_c(T,L) \leq k-1}}  \frac { \sqrt{p}\|E\|}{\bar{\lambda}} \times \left\lbrace \frac{\sqrt{p}\|E\|}{\bar{\lambda}}, \frac{rx}{\delta_S}, \frac{\sqrt{r} w}{\sqrt{\bar{\lambda} \delta_S}} \right\rbrace^{s-1} \,(\text{since}\, T_1+T_2+T_3=s). 
  \end{split}    
\end{equation*}
By Assumption \textbf{D0}, the RHS is at most
\begin{equation*} 
\begin{split}
  \sum_{T=1}^{s_2} \binom{s_2}{T}   \sum_{\substack{L \in \mathcal{L}(T,s);\, r_c(T,L) \leq k-1}}  \frac { \sqrt{p}\|E\|}{\bar{\lambda}} \times \frac{1}{12^{s-1}}   &= \left(\sum_{T=1}^{s_2}   \sum_{\substack{L \in \mathcal{L}(T,s) \\ r_c(T,L) \leq k-1}}   \binom{s_2}{T}\right) \frac { \sqrt{p}\|E\|}{\bar{\lambda}} \times \frac{1}{12^{s-1}}\\
  & \leq \left(\sum_{T=1}^{s_2}   \sum_{\substack{L \in \mathcal{L}(T,s)}}   \binom{s_2}{T}\right) \frac { \sqrt{p}\|E\|}{\bar{\lambda}} \times \frac{1}{12^{s-1}}\\
  & \leq \frac { \sqrt{p}\|E\|}{\bar{\lambda}} \times \frac{2^{s_2+s-1}}{12^{s-1}} \,\,(\text{by \eqref{s_2Lsumbound}}).
\end{split}    
\end{equation*}
This proves Lemma \ref{M_1Case3boundfinal}.

\subsection{Proof of Lemma \ref{M_2Case3boundfinal}}
We bound $M_2(\alpha;\beta)$ in this case as follows. Replacing  $w^{2k-1}$ by $\|E\|^{2k-1}$, we get 
\begin{equation*}
    \begin{split}
   M_2 (\alpha;\beta) & \leq  \sum_{T=1}^{s_2} \binom{s_2}{T} p^{T-1/2} r^{s_2 -T}\sum_{\substack{L \in \mathcal{L}(T,s);\, r_c(T,L) > k-1}} \frac{ \|E\|^{s_1+k-1} x^{s_2 - k}}{\delta_S^{s_2 -1-r_c(T,L)} \bar{\lambda}^{s_1+r_c(T,L)}}.
    \end{split}
\end{equation*}
Similar to the argument in \eqref{M_2case1Rearrange}, the RHS can be rewritten as
\begin{equation}
\begin{split} \label{M_2a_1=1a_{k+1}=0}
 \sum_{T=1}^{s_2} \binom{s_2}{T} \sum_{\substack{L \in \mathcal{L}(T,s) \\ r_c(T,L) > k-1}} \left( \frac{p}{r}  \right)^{T-(r_c(T,L)+1)} p^{\frac{k -s_1}{2}} \left( \frac{r x}{\delta_S} \right)^{s_2-1-r_c(T,L)} \left(\frac{px}{\bar{\lambda}} \right)^{r_c(T,L)-(k-1)} \left( \frac{ \sqrt{p} \|E\|}{\bar{\lambda}} \right)^{k-1+s_1}.
\end{split}
\end{equation}
By Lemma \ref{r_cLemma}, we have $T -(r_c(T,L)+1) \geq 0$ and hence $ 0\leq  \left( \frac{p}{r}  \right)^{T-(r_c(T,L)+1)} \leq 1$. Since $(\alpha,\beta)$ is of Type III, $s_1 \geq k$ and then $0 \leq p^{\frac{k-s_1}{2}} \leq 1$. Omitting these factors, we obtain 
\begin{equation*}
\begin{split} 
M_2(\alpha;\beta) & \leq  \sum_{T=1}^{s_2} \binom{s_2}{T} \sum_{\substack{L \in \mathcal{L}(T,s) \\ r_c(T,L) > k-1}} \left( \frac{r x}{\delta_S} \right)^{s_2-1-r_c(T,L)} \left(\frac{px}{\bar{\lambda}} \right)^{r_c(T,L)-(k-1)} \left( \frac{ \sqrt{p} \|E\|}{\bar{\lambda}} \right)^{k-1+s_1} \\
& = \sum_{T=1}^{s_2} \binom{s_2}{T} \sum_{\substack{L \in \mathcal{L}(T,s) \\ r_c(T,L) > k-1}} \frac{px}{\bar{\lambda}} \times \left( \frac{r x}{\delta_S} \right)^{s_2-1-r_c(T,L)} \left(\frac{px}{\bar{\lambda}} \right)^{r_c(T,L)-(k-1)-1} \left( \frac{ \sqrt{p} \|E\|}{\bar{\lambda}} \right)^{k-1+s_1}.
\end{split}
\end{equation*}
By the setting of $M_2(\alpha;\beta)$, $r_c(T,L) -(k-1)-1 \geq 0$. We also have $s_2-1-r_c(T,L) \geq s_2 -T \geq 0$ and $ k-1+s_1 \geq 0$. By Lemma \ref{toytrick}, we further obtain
\begin{equation*}
    \begin{split}
M_2(\alpha;\beta) & \leq \sum_{T=1}^{s_2} \binom{s_2}{T} \sum_{\substack{L \in \mathcal{L}(T,s) \\ r_c(T,L) > k-1}}    \frac{p x}{\bar{\lambda}} \max \left\lbrace \frac{rx}{\delta_S}, \frac{px}{\bar{\lambda}},
\frac{\sqrt{p}\|E\|}{\bar{\lambda}} \right\rbrace^{s_1+s_2-2}\\
& =\sum_{T=1}^{s_2} \binom{s_2}{T} \sum_{\substack{L \in \mathcal{L}(T,s) \\ r_c(T,L) > k-1}}    \frac{p x}{\bar{\lambda}} \max \left\lbrace \frac{rx}{\delta_S}, \frac{px}{\bar{\lambda}},
\frac{\sqrt{p}\|E\|}{\bar{\lambda}} \right\rbrace^{s-1} 
\end{split}
\end{equation*}
Since $\frac{px}{\bar{\lambda}} \leq \frac{rx}{\delta_S}$, the RHS simplifies to
\begin{equation*}
    \begin{split}
  \sum_{T=1}^{s_2} \binom{s_2}{T} \sum_{\substack{L \in \mathcal{L}(T,s) \\ r_c(T,L) > k-1}}    \frac{p x}{\bar{\lambda}} \max \left\lbrace \frac{rx}{\delta_S},
\frac{\sqrt{p}\|E\|}{\bar{\lambda}} \right\rbrace^{s-1} 
& = \left(  \sum_{T=1}^{s_2}  \sum_{\substack{L \in \mathcal{L}(T,s) \\ r_c(T,L) > k-1}} \binom{s_2}{T} \right) \times \frac{p x}{\bar{\lambda}} \max \left\lbrace \frac{rx}{\delta_S},
\frac{\sqrt{p}\|E\|}{\bar{\lambda}} \right\rbrace^{s-1} \\
& \leq \left(  \sum_{T=1}^{s_2}  \sum_{\substack{L \in \mathcal{L}(T,s)}} \binom{s_2}{T} \right) \times \frac{p x}{\bar{\lambda}} \max \left\lbrace \frac{rx}{\delta_S},
\frac{\sqrt{p}\|E\|}{\bar{\lambda}} \right\rbrace^{s-1} \\
& \leq   2^{s+s-1}  \frac{p x}{\bar{\lambda}} \max \left\lbrace \frac{rx}{\delta_S}, \frac{\sqrt{p}\|E\|}{\bar{\lambda}} \right\rbrace^{s-1} (\text{by \eqref{s_2Lsumbound} }). 
   \end{split}
\end{equation*}
By Assumption \textbf{D0}, we finally obtain
$M_2(\alpha;\beta) \leq \frac{p x}{\bar{\lambda}} \frac{2^{s_2+s-1}}{12^{s-1}}, $   
 proving Lemma \ref{M_2Case3boundfinal}.
 \subsection{Proof of Corollary \ref{lowrank-random1}} \label{app: 5.3 to concretebound}
In this subsection, we prove Corollary~\ref{lowrank-random1} using
Theorem~\ref{mainTh01}. Recall that the entries
\(\{\xi_{ij}\}_{1\leq i\leq j\leq n}\) of the random perturbation \(E\)
are sub-Gaussian and satisfy
\begin{equation}\label{subGauss condition}
    \max_{1\leq i\leq j\leq n}\|\xi_{ij}\|_{\psi_2}\leq C.
\end{equation}
By the truncation argument in Remark \ref{remark: truncation}, with probability at least \(1-n^{-1}\),
\(E\) is \(K\)-bounded with
\[
K=3C\sqrt{\log n}.
\]
Moreover, \eqref{subGauss condition} implies
\begin{equation}\label{sigma m4}
    \sigma^2
    :=\max_{1\leq i,j\leq n}\E\xi_{ij}^2
    \leq C^2,
    \qquad
    m_4
    :=\frac{1}{n}\max_{1\leq i\leq n}
    \sum_{j=1}^n\E\xi_{ij}^4
    \leq 4C^4.
\end{equation}
Furthermore, with probability \(1-o(1)\),
\begin{equation}\label{E bound subGaussian}
    \|E\|\leq c_0C\sqrt n
\end{equation}
for some absolute constant \(c_0>0\). Thus, with probability at least
\(1-n^{-1}\), \(E\) is
\((3C\sqrt{\log n},C)\)-bounded.

We now apply Theorem~\ref{mainTh01} with
\[
t_1=12C\sqrt{\log n}\log(r\log n),
\qquad
t_2=10C^2r\sqrt{n\log n}.
\]
If
\begin{equation}\label{delta condition 0}
 \textstyle   \delta_p
    \geq
    \max\left\{
        144Cr\sqrt{\log n}\log(r\log n),
        \frac{144r\|E\|^2}{|\lambda_p|}
    \right\}
\end{equation}
holds with probability \(1-\alpha_1\), then Theorem~\ref{mainTh01},
with
\(
\epsilon_1=\frac{\|E\|}{|\lambda_p|},
\epsilon_2=\frac{1}{\delta_p},
\)
gives
\[ \textstyle
\Norm{\tilde{\Pi}_p-\Pi_p}
\leq
12\sqrt p
\left(
\frac{\|E\|}{|\lambda_p|}
+
\frac{12C\sqrt r\sqrt{\log n}\log(r\log n)}{\delta_p}
+
\frac{
Cr^{3/2}
\left[
10\sqrt{n\log n}
+
144\log n\,\log^2(r\log n)
\right]
}{
|\lambda_p|\delta_p
}
\right)
\]
with probability at least
\[ \textstyle
1-\alpha_1
-r^2\exp\left(
-\frac{144\log n\,\log^2(r\log n)}
{4+72\log n\,\log(r\log n)}
\right)
-\frac{m_4}{40C^4\log n}
\geq
1-o(1)-\frac{11}{\log n}.
\]

By \eqref{E bound subGaussian}, if
\[
C_1>2\max\{144C,144c_0^2C^2\},
\]
then the assumption of Corollary~\ref{lowrank-random1} implies
\eqref{delta condition 0} with probability \(1-o(1)-n^{-1}\).
Consequently,
\[ \textstyle
\Norm{\tilde{\Pi}_p-\Pi_p}
\leq
12\sqrt p
\left(
\frac{c_0C\sqrt n}{|\lambda_p|}
+
\frac{12C\sqrt r\sqrt{\log n}\log(r\log n)}{\delta_p}
+
\frac{
Cr^{3/2}
\left[
10\sqrt{n\log n}
+
144\log n\,\log^2(r\log n)
\right]
}{
|\lambda_p|\delta_p
}
\right)
\]
with probability at least
\[ \textstyle
1-o(1)-n^{-1}-\frac{11}{\log n}
=1-o(1).
\]

It remains to simplify the last two terms. Since
\(
\log n\,\log^2(r\log n)
\leq
2\log^3 n
\leq
c_1\sqrt{n\log n}
\)
for some absolute constant \(c_1>1\), we have
\[ \textstyle
\frac{
Cr^{3/2}
\left[
10\sqrt{n\log n}
+
144\log n\,\log^2(r\log n)
\right]
}{
|\lambda_p|\delta_p
}
\leq
\frac{160c_1Cr^{3/2}\sqrt{n\log n}}
{|\lambda_p|\delta_p}.
\]
Moreover, the assumption
\(
\delta_p\geq C_1r^{3/2}\log n
\)
implies
\(
\frac{r^{3/2}\sqrt{\log n}}{\delta_p}
\leq
\frac{1}{C_1\sqrt{\log n}}
\leq \frac{1}{C_1}.
\)
Hence, for some absolute constant \(c_2>0\),
\[ \textstyle
\frac{
Cr^{3/2}
\left[
10\sqrt{n\log n}
+
144\log n\,\log^2(r\log n)
\right]
}{
|\lambda_p|\delta_p
}
\leq
\frac{c_2C\sqrt n}{|\lambda_p|}.
\]
Finally, since
\(
\sqrt{\log n}\log(r\log n)
\leq 2\log^{3/2}n,
\)
we obtain
\[
\Norm{\tilde{\Pi}_p-\Pi_p}
\leq
C_2\sqrt p
\left(
\frac{\sqrt n}{|\lambda_p|}
+
\frac{\sqrt r\,\log^{3/2}n}{\delta_p}
\right),
\]
where
\[
C_2=12C\max\{c_0+c_2,24\}.
\]
This completes the proof of Corollary~\ref{lowrank-random1}.

\section{Sharpness of Theorem \ref{deterministicHRgeneral}} \label{section:sharpness} 
In this section, we fix the $n \times n $ matrix 
$$A:=  \begin{pmatrix}
\lambda & 0 & 0 &  \cdots 0 \\
0& \lambda -\delta & 0 & \cdots 0 \\
0 & 0 & 0 & \cdots 0 \\
... \\
0 & 0 & 0 & \cdots 0 
\end{pmatrix}, $$
for some $\lambda > \delta > 0$. Clearly, $\lambda_1= \lambda, \lambda_2 = \lambda -\delta$ are the first  and second eigenvalues of $A$ with the corresponding eigenvectors $e_1, e_2$, respectively. We will show that our bound for  $\Norm{\tilde{u}_1 \tilde{u}_1^T - u_1^T u_1^T}= \Norm{\tilde{\Pi}_1 -\Pi_1}$ is sharp by setting up three different noise matrices $E$, in which each term from Theorem \ref{deterministicHRgeneral} dominates and represents the true value of the perturbation. We set $\bar \lambda = \lambda$. Hence, $r=2$ and $N_{\bar{\lambda}}=\{1,2\}.$

In the following subsections, we will compute $\Norm{\tilde{u}_1 \tilde{u}_1^T - u_1 u_1^T}= \Norm{\tilde{\Pi}_1 -\Pi_1}$ by the sine formula:
\begin{equation}
    \Norm{\tilde{u}_1 \tilde{u}_1^T - u_1^T u_1^T}= \norm{ \sin \angle (\tilde{u}_1, u_1)} = \sqrt{1 - \cos^2 \angle (\tilde{u}_1, u_1) } = \sqrt{1 - \left\langle \tilde{u}_1, u_1 \right\rangle^2}.
\end{equation}
Therefore, given our setting $u_1:= e_1 =[1 \, 0 \, 0 \, \cdots \, 0]^T$ and that $\tilde{u}_1= [\tilde{u_{11}} \,\,\, \tilde{u}_{12} \,\,\, \tilde{u_{13}} \,\,\, \cdots \,\,\,\tilde{u}_{1n}]^T,$ we have
 \begin{equation} \label{eigenspaceformula}
    \Norm{\tilde{u}_1 \tilde{u}_1^T - u_1^T u_1^T}= \sqrt{\sum_{j >1} \tilde{u}_{1j}^2 }.  
 \end{equation}
\subsection{Necessity of the first term} \label{subsec: need adterm} Consider 
\[
E=  \begin{pmatrix}
 0 & 0 &  \mu &0 & \cdots  0 \\
0 & 0 & 0 & 0 &  \cdots 0 \\
\mu & 0 & 0 & 0& \cdots 0 \\
0 & 0 & 0 &0 & \cdots 0 \\
... \\
0 & 0 & 0 & 0 & \cdots 0 \\
\end{pmatrix},\,\text{for some $0 <\mu <\lambda$, and hence,}\,\,\tilde{A} =  \begin{pmatrix}
\lambda & 0 & \mu &  \cdots 0 \\
0& \lambda -\delta & 0 & \cdots 0 \\
\mu & 0 & 0 & \cdots 0 \\
... \\
0 & 0 & 0 & \cdots 0 
\end{pmatrix}.
\]
 By a routine computation, we can show that 
the first eigenvalue $\tilde{\lambda}_1= \frac{\lambda + \sqrt{\lambda^2 + 4 \mu^2}}{2}$,  with  corresponding eigenvector
 $$\tilde{u}_1= [\tilde{u}_{11} \,\,\, 0 \,\,\, \tilde{u}_{13} \,\,\, 0 \,\,\,0 ... 0]^T,$$ where 
 $$ \tilde{u}_{11}:= \frac{\lambda + \sqrt{\lambda^2 + 4 \mu^2}}{2 \sqrt{ \tilde{\lambda}_1^2+ \mu^2}}  \,\,\, \text{and}\,\,\, \tilde{u}_{13}:= \frac{\mu}{\sqrt{\tilde{\lambda}_1^2+\mu^2}}.$$
 
  \noindent By \eqref{eigenspaceformula}, we have 
 \begin{equation}
 \begin{split}
 \Norm{\tilde{u}_1 \tilde{u}_1^T - u_1^T u_1^T}& = |\tilde{u_{13}}| =\frac{\mu}{\sqrt{\tilde{\lambda}_1^2+\mu^2}} = \Theta \left( \frac{\mu}{\lambda} \right)= \Theta \left( \frac{\|E\|}{\lambda} \right).
 \end{split}
\end{equation}   
Moreover, in this setting $x = y = 0$ (see Definition \ref{xyz}). The first term $\frac{\| E \| } {\bar \lambda} $ obviously dominates and gives the correct order of magnitude of   $\Norm{\tilde{u}_1 \tilde{u}_1^T - u_1 u_1^T}$.

\subsection{Necessity of the second term} Consider
\[
E=  \begin{pmatrix}
 0& \mu &  0 &0 & \cdots  0 \\
\mu & 0 & 0 & 0 &  \cdots 0 \\
0 & 0 & \sqrt{\lambda} & 0& \cdots 0 \\
0 & 0 & 0 &0 & \cdots 0 \\
... \\
0 & 0 & 0 & 0 & \cdots 0 \\
\end{pmatrix},\,\,\text{for some $0 < \mu \ll \delta < \lambda,$ and hence,}\,\tilde{A} =   \begin{pmatrix}
\lambda & \mu & 0 &  \cdots 0 \\
 \mu & \lambda -\delta & 0 & \cdots 0 \\
0 & 0 & \sqrt{\lambda} & \cdots 0 \\
... \\
0 & 0 & 0 & \cdots 0 
\end{pmatrix}.
\]
By a routine computation, one can show that the first eigenvalue $\tilde{\lambda}_1= \left( \lambda -\frac{\delta}{2} \right)+ \frac{\sqrt{\delta^2 +4 \mu^2}}{2}$ with the corresponding eigenvector
 $$\tilde{u}_1= [\tilde{u}_{11} \,\,\, \tilde{u}_{12} \,\,\, 0  \,\,\, 0 \,\,\,0 ... 0]^T,$$ where 
 $$ \tilde{u}_{11}:= \frac{\delta + \sqrt{\delta^2 +4\mu^2}}{2\sqrt{(\tilde{\lambda}_1-\lambda+\delta)^2+\mu^2}}  \,\,\, \text{and}\,\,\, \tilde{u}_{12}:= \frac{\mu}{\sqrt{(\tilde{\lambda}_1-\lambda+\delta)^2+\mu^2}}.$$
\noindent By  \eqref{eigenspaceformula},
 \begin{equation}
 \begin{split}
\Norm{\tilde{u}_1 \tilde{u}_1^T - u_1^T u_1^T}& = |\tilde{u_{12}}| =\frac{\mu}{\sqrt{(\tilde{\lambda}_1-\lambda+\delta)^2+\mu^2}} = \Theta \left( \frac{\mu}{\delta} \right)= \Theta \left( \frac{x}{\delta} \right).
 \end{split}
\end{equation}

On the other hand, in this case, it is easy to check that $y =0$ and $\|E\|=\Theta(\sqrt{\lambda})$. Therefore, $$\max \left\lbrace \frac{\|E\|}{\lambda}, \frac{y}{ \delta} \right\rbrace = \Theta \left( \frac{\|E\|}{\lambda} \right) = \Theta \left( \frac{1}{\sqrt{\lambda}}\right).$$
We  choose $\lambda = n^{2}, \delta=n^{1/2}, \mu=n^{1/4}$. In this setting,  the second term  $\frac{x}{\delta} $ obviously dominates and gives the correct order of magnitude of  $\Norm{\tilde{u}_1 \tilde{u}_1^T - u_1 u_1^T}$.

\subsection{Necessity of the third term} Consider
\[
E= \begin{pmatrix}
 0 & 0 &  0  & \cdots   0 & 1 \\
0 & 0 & 0  &  \cdots  0 & 1 \\
0 & 0 & 0 & \cdots  0 & 0 \\
... \\
0 & 0 & 0  & \cdots  0 & 0 \\
1 & 1 & 0  & \cdots  0 & 0 \\
\end{pmatrix},\,\text{and hence}\,\tilde{A}:=  \begin{pmatrix}
 \lambda & 0 &  0 & \cdots   0 & 1 \\
0 & \lambda -\delta & 0 &  \cdots  0 & 1 \\
0 & 0 & 0 &  \cdots  0 & 0 \\
... \\
0 & 0 & 0 &  \cdots  0 & 0 \\
1 & 1 & 0 &  \cdots  0 & 0 \\
\end{pmatrix}.
\]
 In this case, $x = 0$, $y=\Theta \left(\frac{1}{\lambda} \right),$ and $\|E\|=\Theta(1)$. Choose $\lambda= n^{\epsilon_1}, \delta= n^{-\epsilon}$ for some $\epsilon_1 > \epsilon > 0$. Thus, 
$$\frac{\|E\|}{\lambda}=\Theta \left( \frac{1}{\lambda} \right), \frac{x}{\delta}=0, \frac{y}{ \delta} = \Theta \left( \frac{1}{\lambda \delta} \right). $$
 
As $\delta  \ll 1$, it is clear that the term $\frac{y}{ \delta}$ is the dominant term. 
 
 Let $\tilde{\lambda}_1$ be the first eigenvalue of $\tilde{A}$. The corresponding first eigenvector has the form 
$$\tilde{u}_1 =[ \tilde{u}_{11} \,\,\, \tilde{u}_{12} \,\, 0\,\, \cdots \,\, 0 \,\, \tilde{u}_{1n}],$$
 in which $\tilde{u}_{11}, \tilde{u}_{12}, \tilde{u}_{1n}$ satisfy the following system of equations: 
 \begin{equation}
     \begin{split}
         & \lambda \tilde{u}_{11} + \tilde{u}_{1n} = \tilde{\lambda}_1 \tilde{u}_{11}, \\
         & (\lambda -\delta) \tilde{u}_{12} + \tilde{u}_{1n} = \tilde{\lambda}_1 \tilde{u}_{12}, \\
         & \tilde{u}_{11} + \tilde{u}_{12} = \tilde{\lambda}_1 \tilde{u}_{1n}, \\
         & \tilde{u}_{11}^2 + \tilde{u}_{12}^2 + \tilde{u}_{1n}^2 = 1.
     \end{split}
 \end{equation}
 These equations imply that 
 $$(\tilde{\lambda}_1 - \lambda) \tilde{u}_{11} = \tilde{u}_{1n}, (\tilde{\lambda}_1 -\lambda+\delta) \tilde{u}_{12} = \tilde{u}_{1n},$$
 and then, 
 $$1 = \tilde{u}_{11}^2 + \tilde{u}_{12}^2 + \tilde{u}_{1n}^2 = \tilde{u}_{1n}^2 \left( \frac{1}{(\tilde{\lambda}_1 -\lambda)^2} + \frac{1}{(\tilde{\lambda}_1 -\lambda+\delta)^2} +1\right). $$
 Thus, 
 \begin{equation} \label{eigenspaceformula2}
    \tilde{u}_{1n}^2 = \frac{ (\tilde{\lambda}_1 -\lambda)^2 }{1 + (\tilde{\lambda}_1 -\lambda)^2 + \frac{(\tilde{\lambda}_1 -\lambda)^2}{(\tilde{\lambda}_1 -\lambda+\delta)^2}} \,\,\, \text{and} \,\,\, \tilde{u}_{12}^2 = \frac{ \frac{(\tilde{\lambda} -\lambda)^2}{(\tilde{\lambda} -\lambda+\delta)^2} }{1 + (\tilde{\lambda}_1 -\lambda)^2 + \frac{(\tilde{\lambda}_1 -\lambda)^2}{(\tilde{\lambda}_1 -\lambda+\delta)^2}}.  
 \end{equation}
Together \eqref{eigenspaceformula} and \eqref{eigenspaceformula2} imply that
  \begin{equation} \label{contradictionassumption}
     \begin{split}
& \Norm{\tilde{u}_1 \tilde{u}_1^T - u_1^T u_1^T}^2 = \frac{(\tilde{\lambda}_1 -\lambda)^2 + \frac{(\tilde{\lambda}_1 -\lambda)^2}{(\tilde{\lambda}_1 -\lambda+\delta)^2} }{1+(\tilde{\lambda}_1 -\lambda)^2 + \frac{(\tilde{\lambda}_1 -\lambda)^2}{(\tilde{\lambda}_1 -\lambda+\delta)^2}}.
     \end{split}
 \end{equation}
 On the other hand, since $\tilde{\lambda}_1$ is the largest root of the equation $\det (x I_n - \tilde{A})$, by Schur's formula for determinant, $\tilde{\lambda}_1 > \lambda$ satisfies the following equation:
\begin{equation} \label{tildelambda_1}
    \tilde{\lambda}_1 = \frac{1}{\tilde{\lambda}_1 -\lambda} + \frac{1}{\tilde{\lambda}_1 -\lambda+\delta}.  
\end{equation}
 Since $\tilde{\lambda}_1 > \lambda = n^{\epsilon_1}$, Equation \eqref{tildelambda_1} requires that $\tilde{\lambda}_1=\lambda+o(1),$ and that 
 $$ \frac{1}{\tilde{\lambda}_1-\lambda} < \tilde{\lambda}_1 < \frac{2}{\tilde{\lambda}_1 -\lambda}.$$
 These requirements imply that 
 $$\tilde{\lambda}_1= \Theta(n^{\epsilon_1}) \,\,\, \text{and} \,\,\, \tilde{\lambda}_1-\lambda=\Theta(n^{-\epsilon_1}).$$
 Therefore, we can simplify \eqref{contradictionassumption} as 
$$\Norm{\tilde{u}_1 \tilde{u}_1^T - u_1^T u_1^T}^2 = \Theta \bigg( \frac{1}{\left( \frac{\tilde{\lambda}_1 -\lambda+\delta}{\tilde{\lambda}_1-\lambda} \right)^2 +1  } \bigg) = \Theta \left(\frac{1}{(1+ n^{\epsilon_1 -\epsilon})^2+1}  \right) = \Theta \left( \frac{1}{ (\lambda \delta)^2}\right). $$
This means that 
$$\Norm{\tilde{u}_1 \tilde{u}_1^T - u_1^T u_1^T} = \Theta \left( \frac{y}{\delta}  \right).  $$
This proves the claim.

\subsection{Necessity of the assumptions} In this section, we show that our assumptions (such as {\bf C0, D0}, etc.) are also sharp (up to a constant factor). 

We focus on {\bf C0}. In this assumption, the term $ \frac{\sqrt{p}\|E\|}{\bar{\lambda}}$ and $ \frac{ \sqrt{pr} x}{\delta}$  appears in the bound. Thus, even to make the bound meaningful, we need to assume that these two terms are bounded by a small constant.

 For the condition that $\delta \gg \frac{r w^2}{\bar{\lambda}}$, let us consider  matrix $A:= \tiny\begin{pmatrix}
    \lambda_1 & 0 & \cdots & 0 \\
    0 & \lambda_2 & \cdots &  0 \\
     0 & 0 & \cdots & 0 \\
    ... \\
    0 & 0 & \cdots & 0
\end{pmatrix},$  with $\lambda_1 = c_1 \sqrt{n}, \lambda_2= c_2 \sqrt{n}$ for some positive constants $c_1 > c_2 >0$.

Let $E:= \begin{pmatrix}
    0 & 0 \\
    0 & W_{n-1}
\end{pmatrix}$, in which $W_{n-1}$ is $(n-1) \times (n-1)$ Wigner matrix. We also investigate the perturbation bound of $\sin (u_1, \tilde u_1 )=\|\tilde{u}_1 \tilde{u}_1^T - u_1 u_1^T\|$. The condition $\delta \gg \frac{r w^2}{\bar{\lambda}}$ becomes $\delta \gg \frac{\|E\|^2}{\lambda_2}$.

Results from  random matrix theory \cite{B-GN1, KY1} show that  $\tilde{\lambda}_2 \approx \lambda_2 + \frac{\|E\|^2}{4 \lambda_2}$ when $n$ tends to infinity; see \cite{KY1, B-GN1, B-GGM1}.

If $\frac{\|E\|^2}{4 \lambda_2} > \lambda_1 - \lambda_2 = \delta$, then $\tilde{\lambda}_2$ becomes  the leading eigenvalue. By 
our construction,  $u_1$ and $\tilde{u}_1$ are perpendicular.  This means that in order to obtain any meaningful bound on 
$\sin (u_1, \tilde u_1 )$, we must assume 
$$ \frac{\|E\|^2}{4 \lambda_2} < \delta \,\,\,\text{as desired.}$$

\end{document}